\documentclass[11pt]{amsart}         
\usepackage{amscd}                   
\usepackage{amssymb} 

\newtheorem{thm}{Theorem}[section]
\newtheorem{lem}[thm]{Lemma}

\newtheorem{cor}[thm]{Corollary}
\newtheorem{claim}[thm]{Claim}
\newtheorem{prop}[thm]{Proposition}

\theoremstyle{definition}

\renewcommand{\thecase}{}
\newtheorem{conj}[thm]{Conjecture}

\newtheorem{defn}[thm]{Definition}

\newtheorem{rmk}[thm]{Remark} 
\renewcommand{\thestep}{}

\theoremstyle{remark}

\makeatletter
\def\alphenumi{
  \def\theenumi{\alph{enumi}}
  \def\p@enumi{\theenumi}
  \def\labelenumi{(\@alph\c@enumi)}}
\makeatother

\makeatletter
\def\thecase{\@arabic\c@case}
\makeatother

\numberwithin{equation}{section}

\makeatletter
\def\thestep{\@arabic\c@step}
\makeatother

\newenvironment{pf}{\begin{proof}[\proofname]}{\end{proof}}
\newenvironment{pf*}[1]{\begin{proof}[#1]}{\end{proof}}

\newcommand\barA{{\bar A}}

\newcommand\barB{{\bar{B}}}
\newcommand\barsB{{\bar{\mathcal{B}}}}

\newcommand\barGl{{\bar{\textrm{Gl}}}}

\newcommand\barM{{\bar{M}}}

\newcommand\barOm{{\bar\Omega}}

\newcommand\barP{{\overline{P}}}

\newcommand\ubarRR{{\underline{\mathbb{R}}}}

\newcommand\AAA{\mathbb{A}}

\newcommand\CC{\mathbb{C}}

\newcommand\GG{\mathbb{G}}
\newcommand\HH{\mathbb{H}}

\newcommand\NN{\mathbb{N}}

\newcommand\PP{\mathbb{P}}
\newcommand\QQ{\mathbb{Q}}
\newcommand\RR{\mathbb{R}}

\newcommand\UU{\mathbb{U}}

\newcommand\ZZ{\mathbb{Z}}

\newcommand\bdelta{{\boldsymbol{\delta}}}

\newcommand\bga{{\boldsymbol{\gamma}}}

\newcommand\bchi{{\boldsymbol{\chi}}}

\newcommand\blambda{{\boldsymbol{\lambda}}}

\newcommand\bvarphi{{\boldsymbol{\varphi}}}

\newcommand\bFr{{\mathbf{Fr}}}

\newcommand\bG{{\mathbf{G}}}

\newcommand\bL{{\mathbf{L}}}

\newcommand\bp{{\mathbf{p}}}

\newcommand\bx{{\mathbf{x}}}

\newcommand\by{{\mathbf{y}}}

\newcommand{\cov}{\nabla}

\newcommand{\rd}{\partial}

\newcommand\half{{\textstyle{\frac{1}{2}}}}
\newcommand\third{{{\frac{1}{3}}}}

\newcommand\quarter{{\textstyle{\frac{1}{4}}}}

\newcommand\threequarter{{\textstyle{\frac{3}{4}}}}

\newcommand\fB{{\mathfrak{B}}}

\newcommand\fg{{\mathfrak{g}}}

\newcommand\fH{{\mathfrak{H}}}

\newcommand\fS{{\mathfrak{S}}}

\newcommand\fV{{\mathfrak{V}}}

\newcommand\al{\alpha}
\newcommand\be{\beta}

\newcommand\de{\delta}
\newcommand\eps{\varepsilon}

\newcommand\Ga{\Gamma}
\newcommand\la{\lambda}
\newcommand\La{\Lambda}
\newcommand\ka{\kappa}

\newcommand\Om{\Omega}
\newcommand\si{\sigma}
\newcommand\Si{\Sigma}

\newcommand\so{{\mathfrak{s}\mathfrak{o}}}

\newcommand\su{{\mathfrak{s}\mathfrak{u}}}

\newcommand\GL{\operatorname{GL}}

\newcommand\PU{\operatorname{PU}}

\newcommand\SO{\operatorname{SO}}

\newcommand\SU{\operatorname{SU}}
\newcommand\U{\operatorname{U}}

\newcommand\less{\setminus}

\newcommand{\8}{\infty}

\newcommand\cp{\mathbb{C}\mathbb{P}}

\newcommand\CS{\operatorname{CS}}

\newcommand\dist{\operatorname{dist}}

\newcommand\End{\operatorname{End}}

\newcommand\Fr{\operatorname{Fr}}

\newcommand\Gl{\operatorname{Gl}}

\newcommand\Ker{\operatorname{Ker}}

\newcommand\loc{\operatorname{loc}}

\newcommand\Met{\operatorname{Met}}

\newcommand\PD{\operatorname{PD}}

\newcommand\Ric{\operatorname{Ric}}

\newcommand\Ran{\operatorname{Ran}}

\newcommand\Spec{\operatorname{Spec}}

\newcommand\Stab{\operatorname{Stab}}

\newcommand\supp{\operatorname{supp}}

\newcommand\Sym{\operatorname{Sym}}

\newcommand\tr{\operatorname{tr}}

\newcommand\Vol{\operatorname{Vol}}

\newcommand\ad{{\mathrm{ad}\,}}

\newcommand\id{{\mathrm{id}}}

\newcommand\sA{{\mathcal{A}}}
\newcommand\sB{{\mathcal{B}}}
\newcommand\sC{{\mathcal{C}}}

\newcommand\sG{{\mathcal{G}}}
\newcommand\sH{{\mathcal{H}}}

\newcommand\sO{{\mathcal{O}}}

\newcommand\sU{{\mathcal{U}}}

\newcommand\sZ{{\mathcal{Z}}}

\newcommand\tg{{\tilde g}}

\newcommand\tM{{\tilde M}}

\begin{document}
\title[Donaldson Invariants and Wall-Crossing Formulas. I]
{Donaldson Invariants and Wall-Crossing Formulas, I:
Continuity of Gluing Maps}
\author[Paul M. N. Feehan]{Paul M. N. Feehan}
\address{Department of Mathematics\\
Ohio State University\\
Columbus, OH 43210}
\curraddr{School of Mathematics\\
Institute for Advanced Study\\
Olden Lane\\
Princeton, NJ 08540}
\email{feehan@math.ohio-state.edu and feehan@math.ias.edu}
\author[Thomas G. Leness]{Thomas G. Leness}
\address{Department of Mathematics\\
Florida International University\\
Miami, FL 33199}
\email{lenesst@fiu.edu}
\dedicatory{}
\subjclass{}
\thanks{The first author was supported in part by an NSF Mathematical 
Sciences Postdoctoral Fellowship under grant DMS 9306061 and by NSF grants
DMS 9704174 and, through the Institute for Advanced Study, by DMS 9729992}
\date{First version: December 10, 1998. math.DG/9812060. This version:
February 1, 1999.}  
\keywords{}
\begin{abstract}
The present article is the first in a series whose ultimate goal is to
prove the Kotschick-Morgan conjecture concerning the wall-crossing formula
for the Donaldson invariants of a four-manifold with $b^+=1$. The
conjecture asserts that, in essence, the wall-crossing terms due to changes
in the metric depend at most on the homotopy type of the four-manifold and
the degree of the invariant. We prove the first half of a general gluing
theorem endowed with several important enhancements not present in previous
treatments due to Taubes, Donaldson, and Mrowka. In particular, we show
that the gluing maps, with obstructions to deformation permitted, extend to
maps on the closure of the gluing data which are continuous with respect to
the Uhlenbeck topology. The second half of the proof of our general gluing
theorem is proved in the sequel \cite{FLKM2}.
\end{abstract}

\maketitle
\setcounter{tocdepth}{1}
\tableofcontents



\section{Introduction} 
\label{subsec:Introduction}
The underlying problem which motivates the present work and its sequels 
\cite{FLKM2}, \cite{FLKM3}, \cite{FLKM4} is
to prove the strong form of the Kotschick-Morgan conjecture for the
Donaldson invariants of a four-manifold with $b^+=1$ and $b_1=0$, namely
Conjecture 6.2.2 in \cite{KotschickMorgan}.  {}From work of G\"ottsche we
know that this (in fact, even the weaker Conjecture 6.2.1 of
\cite{KotschickMorgan}) allows one to completely determine the
wall-crossing formula for such four-manifolds.  Though the two problems are
not directly related, the computational and analytical problems arising in
attempts to prove the wall-crossing formula are analogous, but still
significantly simpler, than those appearing in our concurrent program to
prove Witten's conjecture \cite{FLGeorgia}, \cite{FL2}.

Our program to relate the Donaldson and Seiberg-Witten invariants
\cite{FL1}, \cite{FLGeorgia}, \cite{FL2}, \cite{FL3} requires a gluing 
theory, for $\PU(2)$ monopoles, possessing most of the properties of
the gluing theory established here and in \cite{FLKM2} for the easier
case of anti-self-dual connections.  Thus, we shall present the
results we need in the considerably simpler case of anti-self-dual
connections separately, so we can address elsewhere the difficulties
peculiar to $\PU(2)$ monopoles \cite{FL3}.  The other purpose of the
present article and its sequels is, of course, to prove the
Kotschick-Morgan conjecture.  That too bears on our program to prove
the relation between Donaldson and Seiberg-Witten invariants, as it
allows us to dispose of the technical difficulties in
the general method underlying the expected proof of
equivalence in an easier case. 

The principal difficulties arising in attempts to prove the
Kotschick-Morgan conjecture are due to `bubbling', a
phenomenon which also accounts for the problems arising when one
considers either $\PU(2)$ monopoles or pseudo-holomorphic curves in
symplectic manifolds (of any dimension). 

\subsection{Statement of results} 
The main purpose of the present article and its companion \cite{FLKM2}
is to prove a general gluing theorem
for anti-self-dual connections.  Naturally, many statements of
gluing theorems have appeared before in the literature
\cite{DonConn}, \cite{DK}, \cite{DS}, \cite{MrowkaThesis}, 
\cite{TauSelfDual}, \cite{TauIndef}, \cite{TauFrame}, \cite{TauStable}, so
the reader may reasonably ask what is new and different about the
present version.  The gluing-theorem statements that are closest to those
needed for a proof of the Kotschick-Morgan conjecture are given by
Theorems 3.4.10 and 3.4.17 in \cite{FrM}: both assertions are based on
work of Taubes, principally \cite{TauSelfDual}, \cite{TauIndef},
\cite{TauFrame}, \cite{TauStable}, though the assertions 
of \cite{FrM} go beyond what has been already proved in or is easily
adapted from the existing gluing-theory literature. In any event, the
statement of the gluing theorem given in
\cite{FrM} is not sufficient to prove the Kotschick-Morgan conjecture; some 
of the additional requirements are briefly explained below in \S
\ref{subsec:Sequels}.

The present article is intended to bridge the gap between the
gluing-theory literature cited above and an assertion of
\cite[Theorems 3.4.10 \& 3.4.17]{FrM}, as we describe below.

Taubes' gluing maps \cite{TauStable} are defined on the level of
connections, rather than ideal connections, and it is not clear from
\cite{TauSelfDual}, \cite{TauIndef}, \cite{TauFrame}, \cite{TauStable}
that they necessarily extend continuously to the closure
$\bar{\Gl}_\kappa^{w,+}(X,\bar\sU_0,\bar\Sigma,\lambda_0)$ of the gluing
data bundle $\Gl_\kappa^{w,+}(X,\sU_0,\Sigma,\lambda_0)$. See \S
\ref{subsec:GluingData} for the definitions of the gluing data and its
closure defined by a choice of open neighborhood $\sU_0$ in
$M_{\kappa_0}^w(X,C)$ of the `background' connection $[A_0]$, together with
a smooth stratum $\Sigma\subset\Sym^\ell(X)$ and small positive constant
$\lambda_0$. Let
$$
\bar{\sB}_\kappa^w(X) 
:=
\bigcup_{\ell=0}^{\lfloor\kappa\rfloor}
\sB_{\kappa-\ell}^w(X)\times\Sym^\ell(X),
$$ 
where $\lfloor\kappa\rfloor$ denotes the greatest integer less than or
equal to $\kappa$. Our main result is a proof of the continuity assertion
in \cite{FrM}:

\begin{thm} 
\label{thm:GluingTheorem} 
Let $X$ be a closed, oriented $C^\8$ manifold with generic $C^\8$
Riemannian metric $g$. Let $w\in H^2(X;\ZZ)$ be a class such that no
$\SO(3)$ bundle over $X$ with second Stiefel-Whitney class equal to
$w\pmod{2}$ admits a flat connection. Let $E$ be a rank-two,
Hermitian vector bundle over $X$ with $c_1(E)=w$ and $c_2(E)-c_1(E)^2/4 =
\kappa\geq 1$.  Let $0<\ell<\lfloor\kappa\rfloor$ be an integer, let
$\Sigma\subset\Sym^\ell(X)$ be a smooth stratum, let
$C\subset\GL(T^*X)$ be a compact, connected, smooth
submanifold-with-boundary containing the identity, let $g\in \Met(X)$,
and let $[A_0]\in M_{\kappa-\ell}^w(X,g)$. Then there are
\begin{itemize}
\item 
A smooth {\em gluing map\/}
$\bga:\Gl_\kappa^{w,+}(X,\sU_0,\Sigma,\lambda_0)\to\sB_\kappa^{w,*}(X)$ which
extends to a continuous map 
$\bga:\bar{\Gl}_\kappa^{w,+}(X,\bar\sU_0,\bar\Sigma,\lambda_0)\to
\bar{\sB}_\kappa^w(X)$. Here, $\lambda_0$ is a sufficiently small positive 
constant, depending at most on $X,g,\kappa$.
\item 
A smooth {\em Kuranishi obstruction section\/} $\bchi$ of a smooth {\em
Kuranishi obstruction vector bundle\/} $\bga^*\fV$ over
$\Gl_\kappa^{w,+}(X,\sU_0,\Sigma,\lambda)$, where $\fV$ is a continuous vector
bundle over an open neighborhood $\sU\subset\bar{\sB}_\kappa^w(X)\times C$
of $[A_0]\times\Sym^\ell(X)$ which restricts to a smooth bundle over
$\sU\cap(\sB_\kappa^{w,*}(X)\times C)$. The bundle $\fV$ has real rank
equal to the real dimension of $H^2_{A_0}$. The section $\bchi$ is the
restriction of a continuous section on $\sU$.
\end{itemize}
Together, the gluing map and obstruction section have the property that
$$
\bga(\bchi^{-1}(0)) \subset \barM_\kappa^w(X,C)\cap\sU,
$$
where $\barM_\kappa^w(X,C)$ is the Uhlenbeck
compactification of the universal moduli space $M_\kappa^w(X,C)$ of
connections on $\su(E)$, anti-self-dual with respect to the family of
metrics $\{f(g):f\in C\}\subset \Met(X)$.
\end{thm}

The gluing maps and obstruction sections are defined in \S
\ref{sec:Splicing} and \S \ref{sec:Existence}. As explained in \S
\ref{subsec:KuranishiReducibleASD}, 
when $b_1(X)=0$, we need only consider the cases where the reducible
background point $[A_0,g_0]$ in $M_{\kappa_0}^w(X,C)$ either (i) varies in
a family for which $H^2_{A_0,g_0}=0$ or (ii) has $H^2_{A_0,g_0}\neq 0$, but
$H^1_{A_0,g_0}=0$ and so is an isolated point of $M_{\kappa_0}^w(X,C)$. For
this reason, we can assume without loss of generality that $\sU_0 \subset
M_{\kappa_0}^w(X,C)$, rather than taking, say, $\sU_0 \subset \{[A,g]\in
\sB_{\kappa_0}^w(X)\times C: \|F_A^{+,g}\|_{L^{\sharp,2}(X,g)} < \eps\}$.

\begin{rmk}
In the statements of the gluing theorem given in \cite{FrM} (see Theorems
3.4.10, 3.4.17, and 3.4.18), it is assumed that the gluing data bundle
$\Gl_\kappa^{w,+}(X,\sU_0,\Sigma,\lambda_0)\to\Sigma$ is replaced by its
restriction $\Gl_\kappa^w(X,\sU_0,K,\lambda_0)\to K$, where
$K\Subset\Sigma$ is a precompact open subset. We make no such restriction
here, so the points $\bx\in\Sigma$ are allowed to approach the diagonals of
the symmetric product $\Sym^\ell(X)$ containing the non-compact smooth
stratum $\Sigma$.
\end{rmk}

There are several important properties which are required of a useful
gluing theory --- in any context, but particularly for applications to
the Kotschick-Morgan conjecture --- which are {\em not\/} asserted by
Theorem \ref{thm:GluingTheorem}. These extensions are discussed in \S
\ref{subsec:Sequels} and are derived in the companion article
\cite{FLKM2}.

We shall need the case where $C$ is the image of a path $\gamma$
of metrics on $T^*X$ and the connection $A_0$ is a non-flat reducible,
with stabilizer $S^1$. In this instance, an open neighborhood $\sU\cap
M_\kappa^w(X,\gamma)$ of $[A_0]\times\Sym^\ell(X)$ will be covered by
a collection of images of gluing maps, one gluing map per smooth
stratum $\Sigma$ of $\Sym^\ell(X)$.

To satisfy the constraint on the second Stiefel-Whitney class, we use
the now standard trick of Morgan and Mrowka of requiring the existence of
a spherical class $\beta\in H_2(X;\ZZ)$ such that $\langle
c_1(E),\beta\rangle\neq 0 \pmod{2}$. This condition can always be
satisfied by blowing up if necessary and so, as far as the computation
of Donaldson invariants is concerned, there is no loss in generality,
but the structure of the moduli space is considerably simpler as flat
connections are excluded from the Uhlenbeck compactification
$\barM_\kappa^w(X,g)$. For generic paths of metrics, $\gamma$, most, but
not all, non-flat reducibles occurring in $\barM_\kappa^w(X,\gamma)$
can be assumed to have vanishing obstruction spaces, so this explains
the appearance of the obstruction map in the statement of Theorem
\ref{thm:GluingTheorem}. After work on the continuity assertions of
Theorem \ref{thm:GluingTheorem} had been written up, we learned that Graham
Taylor \cite{GTaylor} had considered the issue of continuity, in a more
restricted form, of the gluing maps defined by Taubes in \cite{TauIndef}.

\subsection{Outline}
We outline the contents of the remainder of the article.  The Kotschick-Morgan
conjecture and how its solution bears on the expected solution to
Witten's conjecture are explained in \S \ref{sec:WallCrossing}.

In \S \ref{sec:Splicing} we define the `splicing map' for connections
via the usual cut and paste techniques. Of course, this part is
entirely standard: the main new issue is that for later purposes it is
necessary to keep track of the constants defining the splicing map
itself. The splicing maps yield continuous maps (smooth away from
singularities) from bundles of gluing data, one bundle per smooth
stratum $\Sigma\subset\Sym^\ell(X)$, into the quotient space
$\sB_\kappa^{w,*}(X)$ whose image lies suitably close to the stratum
$M_{\kappa-\ell}^w(g)\times\Sigma$ of $\barM_\kappa^w(g)$. The
question of exactly how `close' is addressed in \S
\ref{sec:Regularity} and \S \ref{sec:Decay}.

In \S \ref{sec:Regularity} we record the regularity theory for the
anti-self-dual equation. The exposition is modeled on the
corresponding theory for the $\PU(2)$ monopole equations and, indeed,
applies to any quasi-linear first-order elliptic equation with a
quadratic non-linearity on a four-manifold. The results are used
repeatedly throughout the present article and its sequels, so precise
statements are given. The general pattern is that one typically has
uniform, global $L^2_1$ estimates for solutions to the anti-self-dual
equation (obtained, for example, by deforming an approximate, spliced
solution), but when we need estimates with respect to stronger norms
--- uniform with respect to a non-compact parameter --- we must get by
with only local estimates.

We shall almost always want estimates --- for spliced connections or their
curvatures --- whose constants depend at most on the $L^2$ norms of the
curvatures, or the $L^\8$ norms of the self-dual curvatures, of the
component connections restricted to suitable open sets. There are different
ways these estimates can be achieved, depending on the geometry of the
region in $X$ where estimates are required. One option is to use decay
estimates, due in varying degrees of generality to Donaldson,
Groisser-Parker, Morgan-Mrowka-Ruberman, or R\aa de. However, since our
real interest is in the case of $\PU(2)$ monopoles, we provide in
\S \ref{sec:Decay} a slightly simpler approach here which avoids the 
use of such decay estimates or their generalizations. This section
concludes with an estimate for the self-dual component of the
curvature of a spliced connection; as this is a typical application of
decay-type estimates for anti-self-dual connection we present the
details in this case and then leave subsequent calculations of this
type for the reader.

Splicing anti-self-dual connections produces a family, $A$, of
approximately anti-self-dual connections on $X$ and thus a family of
Laplacians, $d_A^+d_A^{+,*}$, on $L^2(X,\Lambda^+\otimes\su(E))$. We
shall need upper and lower bounds for the small eigenvalues of this
Laplacian, along with estimates for its eigenfunctions and variations
with respect to $A$. These and related estimates are derived in \S
\ref{sec:Eigenvalue2}, the methods and results being modeled on those
of Taubes \cite{TauSelfDual}, \cite{TauStable}.

In \S \ref{sec:Existence} we prove the existence and regularity of
solutions to the anti-self-dual equation, the essential point being to
provide estimates for the solutions using a slightly stronger system
of norms than those of \cite{TauSelfDual}, \cite{TauStable}: the
availability of slightly stronger estimates is important in
\cite{FLKM2}. Our basic requirement on the Hermitian bundle $E$ is
that no $SO(3)$ bundle with second Stiefel-Whitney class equal to
$c_1(E)\pmod{2}$ admits a flat connection; we do not require that $X$
be simply-connected or even that $b_1(X)=0$.  Taking $w=c_1(E)$, this
condition ensures that the Uhlenbeck compactification
$\barM_\kappa^w(X,g)$ contains no flat connections and, in particular,
no connections $A_0$ with non-vanishing $H^2_{A_0}$ except reducibles.

In \S \ref{sec:Continuity} we are finally
in a position to prove the main result of the present article, namely
the continuity of the gluing map with respect to Uhlenbeck limits: in
general, we need to know that if either a proper subset of the scales
(used to attach the connections from $S^4$) tends to zero or if the
$\SU(2)$ connections on $S^4$ converge to an Uhlenbeck limit, then the
sequence of families of glued-up solutions converge to a family of
ideal glued-up solutions. For splicing maps, this property is an
immediate consequence of the definition, which is purely local, but
for gluing maps it is not obvious because of the global dependence of
the solutions on the gluing data via Green's operators. The only
continuity statement that is a direct consequence of Taubes' original
construction is that if {\em all\/} the scales tend to zero then the
glued-up solutions converge to the obvious ideal points in the
Uhlenbeck compactification.

\subsection{Sequels and related articles}
\label{subsec:Sequels}
Some remarks are in order concerning the contents of the companion
articles \cite{FLKM2}, \cite{FLKM3}, \cite{FLKM4}, together with the
articles \cite{FL3} and \cite{FL4} concerning gluing theory for $\PU(2)$
monopoles. 

\subsubsection{Surjectivity of gluing maps}
The reader will note that the main result proved here is less than the
`first half' of a desired `gluing theorem'. To be of use for parameterizing
neighborhoods of ideal reducibles, we would also need to know the
following:
\begin{enumerate}
\item
The gluing map $\bga$ is a smooth embedding of
$\Gl_\kappa^{w,+}(X,\sU_0,\Sigma,\lambda_0)$ and a topological embedding of
the Uhlenbeck closure of the gluing data,
$\bar{\Gl}_\kappa^{w,+}(X,\sU_0,\Sigma,\lambda_0)$.
\item
The image of
$\bar{\Gl}_\kappa^{w,+}(X,\sU_0,\Sigma,\lambda_0)\cap\bchi^{-1}(0)$ under
$\bga$ is an open subset of $\barM^w_\kappa(X,C)$;
\item
The space $\barM^w_\kappa(X,C)$ has a finite covering by such open subsets.
\end{enumerate}
These remaining properties (1)--(3) are proved in \cite{FLKM2},
together with certain $C^1$ differentiability properties of $\bga$,
and comprise the `second half' of the gluing theorem. 
It is important to note that preceding three
gluing-map properties are not simple consequences
of the gluing method itself and their justification constitutes the more
difficult half of the proof of the full gluing theorem. In particular,
despite an extensive literature on gluing theory for anti-self-dual
connections, the existing accounts (see, for example,
\cite{DonConn}, \cite{DK}, \cite{MrowkaThesis}, 
\cite{TauSelfDual}, \cite{TauIndef}, \cite{TauFrame}, \cite{TauStable})
only address very special cases which do not capture all of the
difficulties one encounters when attempting to solve the general
problem. 

\subsubsection{Bubble-tree compactifications and manifolds-with-corners
structures} The proof of injectivity and surjectivity of the gluing maps
relies on estimates for the differentials of the gluing maps: these
estimates in turn yield further information concerning the $C^1$
differentiable structure of the bubble-tree compactification of
$M_\kappa^w(X,g)$. Indeed, related work of Taubes in \cite{TauFloer}
suggests that the bubble-tree compactification should be a $C^1$ manifold
with corners (as we would anticipate from known results concerning the
moduli space of stable maps on a symplectic manifold --- an orbifold).  In
our case, the $C^1$ differentiability properties of the gluing maps
essentially imply that the bubble-tree compactification of
$M_\kappa^w(X,C)$ can be given the structure of a $C^1$ manifold with
corners \cite{FLKM3}.

\subsubsection{Intersection theory and completion of the proof of the
Kotschick-Morgan conjecture} In work in preparation, \cite{FLKM4}, we apply
the preceding gluing results to compute the wall-crossing contribution of
ideal reducibles in a one-parameter family of anti-self-dual moduli spaces,
$\barM_\kappa^w(X,g_t)$, $t\in[0,1]$, to the change in the Donaldson
invariants of $X$. Specifically, we show that the changes only depend on
the homotopy type of $X$.  The importance of the continuity assertion for
gluing maps and that some condition on the differentials of the gluing maps
would be needed in order to compare overlaps of different gluing-map images
--- with applications to a proof of the Kotschick-Morgan conjecture ---
first arose in discussions with John Morgan and Peter Ozsv\'ath on
stratified spaces \cite{MorganPrivate}.

\subsubsection{Gluing $\PU(2)$ monopoles}
It remains to point out that techniques for gluing anti-self-dual
conections do not readily carry over to the case of $\PU(2)$ monopoles. 
For example, when using the splicing construction to produce families of
almost anti-self-dual connections, $A$, we find that the norms
$\|F_A^+\|_{L^{\sharp,2}(X)}$ are uniformly less than a small constant: this
alone implies that one can obtain key estimates such as that of Lemma
\ref{lem:LinftyL22CovLapEstv} by
rearrangement arguments which are typical of those used by Taubes in
\cite[\S 5]{TauStable} (see also \cite[\S 5]{FeehanSlice}). On the other
hand $\|F_A^+\|_{L^{\sharp,2}(X)}$ will not generally be small in the case of
almost $\PU(2)$ monopoles and more sophisticated approaches to deriving
necessary estimates (such as Lemma \ref{lem:LinftyL22CovLapEstv}) are
required. This and many other difficulties are discussed in 
\cite{FL3}, \cite{FL4}. 

An earlier version of this article did not fully take advantage of
rearrangement arguments which are possible when gluing anti-self-dual
connections rather than $\PU(2)$ monopoles. For example, the previous
statement of Lemma \ref{lem:LinftyL22CovLapEstv} implied that one needed a
uniform bound on $\|F_A^{+,g}\|_{L^\8(X)}$, though the current, sharper
statement of Lemma \ref{lem:LinftyL22CovLapEstv} (which does not require
this bound) had been known to the first author a year prior to the
discovery of $\PU(2)$ monopoles. Hence, our introduction to the previous
version of this article incorrectly stated that one needed both (i) some
uniform bound on $\|F_A^{+,g}\|_{L^p}$, $p>2$, and (ii)
$\|F_A^{+,g}\|_{L^{\sharp,2}(X)}$ to be uniformly small when deforming a
general family of almost anti-self-dual to a family of anti-self-dual
connections. In fact, as first observed by Taubes in \cite{TauStable}, one
only needs $\|F_A^{+,g}\|_{L^{\sharp,2}(X)}$ to be uniformly small (or just
$\|F_A^{+,g}\|_{L^\sharp(X)}$ uniformly small if one adheres to the system
of norms used by Taubes in \cite[\S 5]{TauStable}). The first author is
especially grateful to Tom Mrowka for reminding him of this important
fact. See \S \ref{sec:Existence} here and
\cite[\S 5]{FeehanSlice}. Thus, contrary to the statement in the earlier
introduction, Theorems 3.4.10 and 3.4.17 in \cite{FrM} do not `implicitly
assume' that the family of spliced connections satisfy such a bound.
Hence, the final \S 9 of the previous version (which addressed the problem
that $\|F_A^{+,g}\|_{L^p}$ may not be uniformly bounded when $p>2$) was
unnecessary and has been omitted.

\subsection{Acknowledgments}
The first author is indebted to Tom Mrowka for his thoughtful criticism of
an earlier version of this article and, in particular, for reminding him of
a forgotten fact that some of the pre-conditions stated for deforming
almost anti-self-dual connections could be relaxed; this led to the
simplifications described above. He is grateful for the hospitality and
support of the many institutions where various parts of the present article
and its sequels were written, namely Harvard University, beginning in
Autumn 1993, Ohio State University (Columbus), the Max Planck Institut
f\"ur Mathematik (Bonn) --- where the principal new parts of the article
were written --- together with the Institute des Hautes Etudes
Scientifiques (Bures-sur-Yvette) and the Institute for Advanced Study
(Princeton), where the article was finally completed.  He would also like
to thank the National Science Foundation for their support during the
preparation of this article.

%

\section{The wall-crossing formula for Donaldson invariants}
\label{sec:WallCrossing}
Our first task is to give a precise statement of the conjectured shape
of the general wall-crossing formula. There are two slightly different
versions given in \cite{KotschickMorgan}: the one we aim to prove is a
stronger form of the second conjecture. It is only fair to admit
straightaway that since the advent of the Seiberg-Witten invariants,
the conjecture itself no longer appears to have a direct bearing on
smooth four-manifold topology, though it is still an important aspect
of Yang-Mills gauge theory.  Nevertheless, the conjecture does retain
considerable mathematical interest or, more precisely, its method of
proof because this is expected to provide a simpler illustration of
the proof of the corresponding conjecture for $\PU(2)$ monopoles
\cite{FLGeorgia}.

Our second task in this subsection is to describe the Kuranishi models
for essential singularities in parameterized moduli spaces
$M_\kappa^w(X,\gamma)$ when $b^+(X)=1$ and both the class $w\in
H^2(X;\ZZ)$ and path of metrics $\gamma:[-1,1]\to\Met(X)$ have been
chosen so as to eliminate all singularities except those due to
non-flat reducible connections at a critical metric
$g_0=\gamma(0)$. The properties of the model in turn dictate our
approach to gluing theory and the construction of models for
singularities due to reducibles occurring in all strata of the
compactification $\barM_\kappa^w(X,\gamma)$.

Reducible connections in $M_\kappa^w(X,g_0)$ occur in (possibly
obstructed) families of dimension $b^1(X)$. Links of
positive-dimensional families of singularities in moduli spaces of
$\PU(2)$ monopoles are constructed in \cite{FL2} but, as might be
expected, they are considerably more complicated than those of
isolated singularities.  Neighborhoods of isolated reducibles in lower
strata of $\barM_\kappa^w(X,\gamma)$ may be parameterized via Taubes'
approach to gluing, which uses spectral projections to construct an
equivariant, finite-rank normal bundle of the singularity in a
thickened moduli space which contains the singularity as a smooth
submanifold.  This is no longer possible in the case of
positive-dimensional families, due to spectral flow, and one must use
a combination of the link construction in \cite{FL2} and a gluing
method which leans closer to that of Donaldson
\cite{DonConn}, \cite{DK}. Gluing constructions of the latter kind lead to
further difficulties: see \cite{FL3}, \cite{FL4}.
However, the formulas for the wall-crossing terms have
only been predicted or proved modulo the Kotschick-Morgan conjecture,
when $b^1(X)=0$. Hence, for these reasons, we shall restrict our
attention to this case whenever necessary to simplify the discussion at
hand.

\subsection{Preliminaries}
\label{subsec:Notation}
To avoid ambiguity, we gather here the notational conventions we shall
employ throughout our work; we follow, whenever possible, the usual
practices of \cite{DK}, \cite{KMStructure} and, in a few instances,
those of \cite{FrM}.

Let $E$ be a Hermitian, rank-two vector bundle over $S^4$.
Let $\sA_\kappa=\sA_E=\sA_E(S^4)$ be the space of $L^2_k$
connections on the $\SO(3)$ bundle $\fg_E=\su(E)$ over $S^4$.  Let $\sG_E$
be the group of $L^2_{k+1}$ unitary gauge transformations of $E$ with
determinant one and let $\sG_E^s\subset\sG_E$ be the subgroup
consisting of gauge transformations $u\in\sG_E$ such that
$u(s)=\id_{E|_s}$.

Let $\sB_\kappa^s=\sB_E^s=\sB_E^s(S^4)$ be the
$\sG_E^s$-quotient space of $L^2_k$ connections on the bundle
$\fg_E$ over $S^4$. There is a map, defined by the choice $q\in \Fr(\fg_E)|_s$,
$$
\sB_E^s \simeq \sA_E\times_{\sG_E} \Fr(\fg_E)|_s, \qquad
[A]\mapsto [A,q],
$$
which is easily seen to be a homeomorphism. The existence of a
slice for the action of the based gauge group, $\sG_E^s$, has been
established by Groisser and Parker and it then follows that this
homeomorphism is in fact a diffeomorphism
\cite{GroisserParkerGeometryDefinite}, as there are no non-flat
reducible $\SO(3)$ connections over $S^4$. The group $\SO(3) =
\SU(2)/\{\pm\id\}$ acts freely on $\sB_E^s$ by the assignment
$[A]\mapsto [uAu^{-1}]$ or $[A,q]\mapsto [A,qu]$, if
$u\in\SU(2)$.

 Let $\Fr(\fg_E)$ denote the principal $\SO(3)$ bundle of
oriented, orthonormal frames for $\fg_E$. Let $\Fr(TX)$ denote the
principal $\SO(4)$ bundle of oriented, $g$-orthonormal frames for $TX$.

\subsection{The Kotschick-Morgan conjecture} 
Let $X$ be a closed, oriented, smooth four-manifold with
$b^+(X)=1$. We briefly recall the definition of Donaldson invariants
in this case \cite{DonHCobord}, \cite{FrQ}, \cite{Goettsche},
\cite{GoettscheZagier}, \cite{KotschickBPlus1}, \cite{KotschickMorgan}.
Following \cite[\S 1]{FrQ}, we call
$$
\Omega_X := \{\beta\in H^2(X;\RR):\beta^2 > 0\}/\RR^+
$$
the {\em positive cone\/} of $H^2(X;\RR)/\RR^+$, where we denote
$\RR^+=(0,\8)$.  Fix a class $w \in H^2(X;\ZZ)$ and a number $\kappa
\in \quarter\ZZ$ such that 
\begin{equation}
\label{eq:BasicChoiceOfPontrjaginStieffel}
-4\kappa = w^2 \pmod{\ZZ}
\quad\text{and}\quad
8\kappa-3(2-b^1)\geq 0.
\end{equation}
A {\em wall of type
$(w,\kappa)$\/} is a nonempty hyperplane in $\Omega_X$,
\begin{equation}
\label{eq:WallDefn}
W^\xi 
:= 
\{\beta \in \Omega_X: \beta\cdot\xi = 0\}
= \xi^\perp\cap\Omega_X,
\end{equation}
for some class $\xi \in H^2(X;\ZZ)$ with 
\begin{equation}
\label{eq:WallLabelDefn}
\xi = w \pmod{2}
\text{ and }
p \le \xi^2 < 0,
\end{equation}
where $p=-4\kappa$. Let $\ell = (\xi^2-p)/4$: observe 
that $0\leq\ell<\kappa$ and that $\ell$ is an integer.
\footnote{The case $\xi^2=0$ and $\ell=\kappa$ could occur if the background 
reducible connection were flat. However, this possibility is excluded
by choosing $w$ as suggested in \cite{MorganMrowkaPoly}; see Lemma
\ref{lem:MorganMrowka}.}  According to
\cite{DonHCobord} the set of $(w,\kappa)$-walls is locally finite.

The connected components of the complement in $\Omega_X$ of the walls
of type $(w,\kappa)$ are called the {\em chambers of type
$(w,\kappa)$}.  Require that $w\pmod{2}$ be such that no $\SO(3)$
bundle $V$ over $X$ with $w_2(V) = w\pmod{2}$ admits a flat
connection; see \cite{KotschickBPlus1}, \cite{MorganMrowkaPoly} for sufficient
criteria and how they may be achieved with no loss of generality
\cite{MorganMrowkaPoly} for the purposes of computing Donaldson invariants.
Then the Donaldson invariants of $X$ associated to $(w,\kappa)$ are
defined with respect to chambers of type $(w,\kappa)$. We shall write
$$
D_{X,\sC}^w(z), \qquad z\in\AAA(X),
$$
for the Donaldson invariants defined by the class $w\in H^2(X;\ZZ)$ and the
chamber $\sC$. Indeed, for a Riemannian metric $g$ on $X$, let
$\omega\in\Omega^{+,g}(X)$ be a $g$-self-dual, harmonic two-form (so
$d^{+,g}d^{*,g}\omega=0$) such that $\int_X\omega\wedge\omega = 1$. The form
is unique up to sign, so fix an $\omega$ and let it define an
orientation for $H^+(X;\RR)$. Let $(-1)\wedge\omega \in
\det(H^0(X;\RR)\oplus H^+(X;\RR))$ define the choice of homology
orientation. For generic $g$, the point $[\omega]\in H^2(X;\RR)$ does
not lie on a wall $W$. If $\sC\subset H^2(X;\RR)/\RR^+$ is the chamber
containing $[\omega]$ then, according to \cite[Theorem 3.2]{KotschickBPlus1},
\cite[Theorem 3.0.1]{KotschickMorgan}, one has 
$D_{X,g}^w(z) = D_{X,\sC}^w(z)$,
that is, the Donaldson invariants only depend on the chamber $\sC$
containing the period point $[\omega(g)]$ and not the metric $g$.

Recall that $M_\kappa^w(X,g)$ has the expected dimension
\begin{align*}
\dim M_\kappa^w(X,g)
&=
-2w^2 - 3(b^+-b_1+1)
\\
&= 
-2w^2 - 3(2-b_1).
\end{align*}
If $\deg z \neq -2w^2 - 3(2-b_1) \pmod{8}$
then $D^w_{X,\sC}(z) = 0$. Otherwise, find $\kappa\in\quarter\ZZ$ such that
$$
\deg z = 8\kappa - 3(2-b_1)
$$
and a Hermitian, rank-two vector bundle $E$ over $X$ such that
\begin{equation}
\label{eq:U2BundleCompatibleWithPontrjaginStieffel}
c_1(E) = w
\text{ and }
-\quarter p_1(\fg_E) = \langle c_2(E) - c_1(E)^2/4,[X]\rangle = \kappa.
\end{equation}
where, throughout this paper, we denote $\fg_E := \su(E)$.
We now define $D^w_{X,\sC}(z)$ following the recipe of \cite[\S
2]{KMStructure}. 

Suppose that $\sC_+$ and $\sC_-$ are separated by a single wall $W^\xi$;
there may be more than one class $\xi$ of type $(w,\kappa)$ defining
$W^\xi$. The difference
$$
\delta_{X,\sC_+,\sC_-}^w(z) := D_{X,\sC_+}^w(z) - D_{X,\sC_-}^w(z), \qquad
z\in\AAA(X), 
$$
defines a {\em transition formula\/}.
Based on known examples and their own work, Kotschick and Morgan proposed
the following two conjectures:

\begin{conj}
\label{conj:WeakKM}
The transition formula $\delta_{X,\sC_+,\sC_-}^w$ is a homotopy invariant
of the pair $(X,\xi)$; more precisely, if $\phi$ is an oriented homotopy
equivalence from $X'$ to $X$, then
$$
\delta_{X,\phi^*(\sC_+),\phi^*(\sC_-)}^{\phi^*w}
=
\phi^*\delta_{X,\sC_+,\sC_-}^w.
$$
\end{conj}

This is expected to be equivalent to the following apparently stronger
claim:

\begin{conj}
\label{conj:StrongKM}
The transition formula $\delta_{X,\sC_+,\sC_-}^w$ is a polynomial in
$\xi$ and the quadratic form $Q_X$, with coefficients depending only on
$\xi^2$, homotopy invariants of $X$, and universal constants. 
\end{conj}

The property of the transition formula which we aim to establish is,
in fact, stronger than that predicted by Conjecture
\ref{conj:StrongKM}. It is a remarkable result of G\"ottsche that
if one assumes Conjecture \ref{conj:WeakKM} then the difference
terms $\delta_{X,\sC_+,\sC_-}^w$ can be computed explicitly for any
simply-connected four-manifold. It seems very likely that the
condition $\pi_1(X)=1$ can be weakened to the case $b^1(X)=0$, as
G\"ottsche points out in \cite{Goettsche} that he uses this constraint
merely to remain consistent with the stated version of the
Kotschick-Morgan conjecture in \cite{KotschickMorgan} and the
Fintushel-Stern blow-up formula \cite{FSBlowUp}: see \cite[Corollary
5.1]{Goettsche}. G\"ottsche's argument only requires Conjecture
\ref{conj:WeakKM} in order to replace $X$ by
$\CC\PP^1\times\CC\PP^1$ in \cite[Lemma 4.10]{Goettsche}.

The first result in this direction is due to Donaldson, who gave a formula
in the case $w=0$ and $\kappa=1$; these were extended by Kotschick in
\cite{KotschickBPlus1}.  Using a gauge-theoretic approach, Yang \cite{Yang}
completely determined the transition formulas for $w=0$ and
$\ell=1$. This was further extended by Leness to the case $\ell=2$ and
arbitrary $w$ \cite{Leness}.

For a rational ruled complex surface $X$, the transition formulas for
$w=0$ and $2\leq\kappa\leq 4$ were determined by \cite{LiQin},
\cite{Mong}, \cite{Tyurin} using methods of algebraic geometry. Using
Thaddeus' method, Friedman and Qin completely determined the
transition formulas when $X$ is a rational surface with $-K_X$
effective \cite[p. 13]{FrQ}. Similar results were independently
obtained by Ellingsrud and G\"ottsche
\cite{EGVariation}, \cite{EGWall}. Wall-crossing formulae for the case
$b_1(X)>0$ have been computed by M{\~u}noz using the techniques of Friedman
and Qin \cite{MunozWallCrossing}.

In a quite different direction, G\"ottsche \cite{Goettsche} assumed the
validity of the Kotschick-Morgan conjecture and used this assumption to
explicitly compute the transition formulas, for arbitrary closed, oriented,
smooth four-manifolds $X$ with $b^1(X)=0$ and $b^+(X)=1$ in terms of
modular forms.

\subsection{Uhlenbeck compactness for parameterized moduli spaces}
\label{subsec:Compactness}
We recall the Uhlenbeck compactness results 
\cite{DK}, \cite{FU} for the moduli space of anti-self-dual $\SO(3)$
connections and their natural extension to the case of families of metrics.
Recall that the process of gluing up connections is essentially the inverse
of Uhlenbeck degeneration of connections.

\begin{defn}
\label{defn:UhlenbeckTop}
We say that a sequence of points $[A_\alpha,g_\alpha]$ in
$\sB_\kappa^w(X)\times\Met(X)$ {\em converges\/} to a point $[A,g,\bx]$
in $\sB_{\kappa-\ell}^w(X)\times\Met(X)\times\Sym^\ell(X)$, where
$0\leq\ell\leq\kappa$ is an integer,
if the following hold:
\begin{itemize}
\item The sequence of metrics $g_\alpha$ converges in $C^r$ to a metric $g$ on 
$T^*X$.
\item Let $E$, $E_0$ be rank-two Hermitian bundles over $X$, 
with determinant line bundle $w$ and induced $\SO(3)$ bundles, $\fg_E$ and
$\fg_{E_0}$, with first Pontrjagin numbers $-4\kappa$ and
$-4(\kappa-\ell)$, respectively.  There is a sequence of
$L^2_{k+1,\loc}(X\less\bx)$, determinant-one, unitary bundle isomorphisms
$u_\alpha:E|_{X\less\bx}\to E_0|_{X\less\bx}$ such that the sequence of
connections $u_\alpha(A_\alpha)$ on $\fg_{E_0}|_{X\less\bx}$ converges to
$A$ on $\fg_{E_0}|_{X\less\bx}$ in $L^2_{k,\loc}$ over $X\less\bx$, and
\item The sequence of measures  
$|F_{A_\alpha}|^2$ converges
in the weak-* topology on measures to $|F_A|^2 +
8\pi^2\sum_{x\in\bx}\delta(x)$, where $\delta(x)$ is the unit Dirac measure
centered at the point $x\in X$.
\end{itemize}
\end{defn}

Let $C\subset \Met(X)$ denote a compact subset of the Banach space of
$C^r$ Riemannian metrics on $T^*X$ ($r\geq 3$) and define the
corresponding parameterized moduli space by
$$
M_\kappa^w(X,C)
:=
\bigcup_{g\in C}M_\kappa^w(X,g)\times\{g\}
\subset
\sB_\kappa^w(X)\times \Met(X).
$$
Let $\barM_\kappa^w(X,C)$ be the closure of $M_\kappa^w(X,C)$ in the
space $\cup_{\ell=0}^{\lfloor\kappa\rfloor}M_{\kappa-\ell}\times\Sym^\ell(X)$
with respect to the induced Uhlenbeck topology, where
$\lfloor\kappa\rfloor$ is the greatest integer contained in $\kappa(E) :=
-\quarter p_1(\fg_E)$, where $-\quarter p_1(\fg_E) = 
\langle c_2(E)-\quarter c_1(E)^2,[X]\rangle$.
We call the intersection of $\barM_\kappa^w(X,C)$ with
$M_{\kappa-\ell}^w\times\Sym^\ell(X)$ a {\em lower-level\/} of the
compactification $\barM_\kappa^w(X,C)$ if $\ell>0$ and call
$M_\kappa^w(X,C)$ the {\em top\/} or {\em highest level\/}.

According to \cite[Theorem 4.4.3 \& Proposition 9.1.2]{DK},
$\barM_\kappa^w(X,C)$ is a compact, Haussdorf, second-countable space when
endowed with the induced Uhlenbeck topology. We call the intersection of
$\barM_\kappa^w(X,C)$ with $M_{\kappa-\ell}^w(X,C)\times\Sym^\ell(X)$ a
{\em lower-level\/} of the compactification $\barM_\kappa^w(X,C)$ if
$\ell>0$ and call $M_\kappa^w(X,C)$ the {\em top\/} or {\em highest
level\/}: in general, the spaces $M_{\kappa-\ell}^w(X,C)$ and
$\Sym^\ell(X)$ are at best smoothly stratified.

The following lemma gives a simple condition on $w$ which guarantees that there
will be no flat $\SO(3)$ connections in $M_\kappa^w(X,g)$, for any $g$.
Because it relies only on $w\pmod{2}$, the lemma applies
simultaneously to all levels of the Uhlenbeck compactification
$\barM_\kappa^w(X,g)$.
\footnote{According to \cite[p. 536]{FS1}, \cite[Proposition 4.1]{KotschickBPlus1}, 
if $V$ is an $\SO(3)$ bundle over a closed four-manifold $X$ such that
$w_2(V)$ is not pulled back from $H^2(K(\pi_1(X),1);\ZZ/2\ZZ)$, then
there are no flat connections on $V$. However, one sees more
directly that the Morgan-Mrowka criterion is satisfied via blow-ups.}

\begin{lem}
\cite[p. 226]{MorganMrowkaPoly}
\label{lem:MorganMrowka}
Let $X$ be a closed, oriented four-manifold. Then
\begin{itemize}
\item
If $E$ is a $\U(2)$ bundle over $X$ and $\beta\in H_2(X;\ZZ)$ is a
spherical class such that $\langle c_1(E),\beta\rangle\neq 0\pmod{2}$,
then no $\SO(3)$ bundle $V$ over $X$ with $w_2(V)=c_1(E)\pmod{2}$
admits a flat connection.
\item
In particular, if $E$ is any $\U(2)$ bundle over $X$ and $\hat E$ is
the $\U(2)$ bundle over $\hat X = X\#\overline{\CC\PP}^2$ with
$c_2(\hat E)=c_2(E)$ and $c_1(\hat E) = c_1(E)+e^*$, where
$e^*=\PD[e]$ and $e=[\overline{\CC\PP}^1]\in H_2(X;\ZZ)$, then no
$\SO(3)$ bundle $V$ over $\hat X$ with $w_2(V)=w_2(\su(\hat E))$ admits a
flat connection.
\footnote{Because the Donaldson 
invariants of $X\#\overline{\CC\PP}^2$ determine and are determined by
those of $X$, (see \cite{FSBlowUp}, \cite{MorganMrowkaPoly}), no
information about these invariants is lost by blowing up.}
\end{itemize}
\end{lem}

\begin{proof}
Let $S\subset X$ be an embedded two-sphere representing $\beta$.  If
$A$ were a flat $\SO(3)$ connection on $V$, then its restriction $A|_S$
would be flat and thus a trivial connection on a trivial $\SO(3)$
bundle $V|_S$, as $S$ is simply connected. But $\langle
w_2(V|_S),[S]\rangle = \langle c_1(E),\beta\rangle \neq 0\pmod{2}$,
yielding a contradiction. This proves the first assertion.

Replacing $X$ by $\hat X = X\#\overline{\cp}^2$ and $E$ by $\hat E$
with $c_2(\hat E) = c_2(E)$ and $c_1(\hat E) = c_1(E)+e^*$ ensures
that $\langle c_1(\hat E),e\rangle
\neq 0\pmod{2}$, as $c_1(E)\cdot e^* = 0$ and $e^*\cdot e^* = -1$. The second 
assertion follows.
\end{proof}

\subsection{Transversality for one-parameter families of moduli spaces}
Rather than consider a moduli space of anti-self-dual connections
defined by a single Riemannian metric, $g_0$, we shall need to
consider $C^\8$ paths of metrics $\gamma:[-1,1]\to \Met(X)$, where
$\Met(X)$ is the Banach manifold of $C^r$ metrics on $T^*X$ (with
$r\in\ZZ$ large). For this purpose, we shall need the following basic
transversality result \cite{DK}, \cite{FU}, \cite{KMStructure}.

\begin{thm}
\label{thm:Transversality}
Let $X$ be a closed, oriented, smooth four-manifold
and let $r\geq 3$ be an integer. Let $v\in H^2(X;\ZZ/2\ZZ)$ be a class
such that no $\SO(3)$ bundle with second Stiefel-Whitney class $v$
admits a flat connection. Then the following hold:
\begin{itemize}
\item
There is a first-category subset of the space of $C^r$ metrics such
that for all metrics $g$ in the complement of this subset, the
following hold: For any $\kappa\in\quarter\ZZ$ and $w\in H^2(X;\ZZ)$
with $-4\kappa=w^2\pmod{\ZZ}$, the moduli space
$M_\kappa^{w,*}(X,g)$ of irreducible
$g$-anti-self-dual $\SO(3)$ connections is
regular (as the zero-locus of the section $F^{+,g}(\,\cdot\,)$ over
$\sB_\kappa^{w,*}(X))$ and thus a $C^\8$ manifold of the expected
dimension $8\kappa - 3(1-b^1+b^+)$.
\item
When $b^+(X)>0$,
there is a first-category subset of the space of $C^r$ metrics such that
for all metrics $g$ in the complement of this subset, the following hold:
For any $\kappa\in\quarter\ZZ$ and $w\in H^2(X;\ZZ)$
with $-4\kappa=w^2\pmod{\ZZ}$, the moduli space
$M_\kappa^w(X,g)$ of $g$-anti-self-dual $\SO(3)$ connections
contains no reducible connections.
\item
Let $g_{-1}$, $g_1$ be metrics satisfying the conditions of the
preceding assertion and suppose $b^+(X)>0$. 
Then there is a first-category subset of the
space of $C^r$ paths $\gamma:[-1,1]\to\Met(X)$ with $\gamma(\pm
1)=g_{\pm 1}$ such that for all paths $\gamma$ in the complement of
this subset, the following hold.  For any $\kappa\in\quarter\ZZ$ and
$w\in H^2(X;\ZZ)$ with $-4\kappa=w^2\pmod{\ZZ}$, the parameterized
moduli space of $\gamma$-anti-self-dual $\SO(3)$ connections,
$$
M_\kappa^{w,*}(X,\gamma) 
:=
\{([A],t) \in \sB_\kappa^{w,*}(X)\times[-1,1]: F^{\gamma(t)}(A)=0\},
$$
is regular (as the zero-locus of the section
$F^{+,\gamma(\cdot)}(\,\cdot\,)$ over $\sB_\kappa^{w,*}(X))\times[-1,1]$
and thus a $C^\8$ manifold of the expected dimension. In particular,
$M_\kappa^{w,*}(X,\gamma)$ gives a $C^\8$ cobordism between the manifolds 
$M_\kappa^{w,*}(X,g_{-1})$ and $M_\kappa^{w,*}(X,g_1)$.
\end{itemize}
\end{thm}

\begin{proof}
All of the above assertions are standard, so we merely list their
provenance. The first assertion is proved as Corollary 4.3.18 in \cite{DK}
and Theorem 3.17 in \cite{FU}: the manifold $X$ is not required to be
simply connected, as assumed in
\cite[\S 4.3]{DK}, \cite[\S 3]{FU}, because flat connections are excluded by 
the condition on $w$ and twisted reducible connections are smooth
points by \cite[Corollary 2.5]{KMStructure}. Similarly, the third
assertion follows from \cite[Corollary 4.3.18]{DK}, \cite[Corollary
2.5]{KMStructure}, and the fact that flat connections are excluded.
The second assertion is proved as Corollary 4.3.15 in
\cite{DK}: the condition that $X$ be simply connected is not required.
\end{proof}

\subsection{Kuranishi models and reducible anti-self-dual connections}
\label{subsec:KuranishiReducibleASD}
It remains to consider the non-flat reducible connections in
$\barM_\kappa^w(X,\gamma)$, the only possible singular points of this
compact cobordism which cannot be excluded by perturbations or a
suitable choice of $w\pmod{2}$ or by dimension-counting when
considering the intersection of geometric representatives with
$\barM_\kappa^w(X,\gamma)$.  In this subsection we examine the
elliptic deformation complex for a reducible $[A_0]$ in
$M_\kappa^w(X,g_0)$, while in the following subsection we describe the
Kuranishi model for a neighborhood of a reducible $[A_0,g_0]$ in
$M_\kappa^w(X,\gamma)$.

In this section we classify the projectively anti-self-dual
connections on $E$ with non-trivial stabilizers in $\sG_E$.  If $A$ is
a connection on $E$, we write $A^{\ad}$ for the connection induced by
$A$ on the associated $\SO(3)$ bundle $\fg_E$ and denote its
stabilizer by $\Stab_A\subset\sG_E$.  Except when discussing bundle
splittings and reducible connections, our usual convention is to
consider connections on $\fg_E$ only, with their lifts to $E$ via the
fixed connection on $\det E$ being implicit when needed. We call $A$
{\em irreducible\/} if $\Stab_A = \{\pm 1\}$, the center of $\SU(2)$,
and {\em reducible\/} otherwise. The following classification is
standard; see, for example, Proposition 3.1 in \cite{FL2}.

\begin{lem} 
\label{lem:ClassificationOfStabilizers}
Let $X$ be a closed, oriented, smooth four-manifold with
Riemannian metric $g$. Let $E$ be a Hermitian
two-plane bundle over $X$ and suppose a unitary connection $A$ on $E$
represents a point $[A]\in M_E$. Then one of the following holds:
\begin{enumerate}
\item 
The connection $A$ is irreducible, with stabilizer $\Stab_A=\{\pm 1\}$.  
\item 
The connection $A$ is reducible, but not projectively flat:
the bundle $E$ splits as an ordered sum of line bundles, $E=L_1\oplus
L_2$, the connection $A$ has stabilizer $\Stab_A \cong S^1$ and $A =
A_1\oplus A_2$, where $A_1$ is a unitary connection on $L_1$ and $A_2
= A_e\otimes A_1^*$ is the corresponding connection on $L_2 = (\det
E)\otimes L_1^*$, with $A_e$ being the fixed connection on $\det E$.
\item 
The connection $A$ is projectively flat and reducible with respect to a
splitting $E = L_1\oplus L_2$, where $L_2 = (\det E)\otimes L_1^*$ and 
$c_1(E)-2c_1(L_1)\in H^2(X;\ZZ)$ is a torsion class.
\end{enumerate}
\end{lem}

For the second case in Lemma \ref{lem:ClassificationOfStabilizers}.
we remark that a projectively $g$-anti-self-dual connection $A$
on $E$ is reducible with respect to the splitting $E\cong L_1\oplus
L_2$ if and only if the induced $g$-anti-self-dual connection
$A^{\ad}$ on $\fg_E$ is reducible with respect to the splitting
$\fg_E\cong i\RR
\oplus \xi$, where $\xi\cong \det E\otimes L_1^{-2}$. Indeed, we have 
$$
A = \begin{pmatrix}A_1 & 0 \\ 0 & A_e\otimes A_1^*\end{pmatrix}
\quad\text{and}\quad
A^{\ad} = 0\oplus A_e\otimes 2A_1^*,
$$
with $F^{+,g}(A_e\otimes 2A_1^*) = 0$.

{}From Theorem \ref{thm:Transversality} we see that it remains to
consider the case of reducible connections in
$\barM_\kappa^w(X,\gamma)$. Suppose that a class $\xi\in H^2(X;\ZZ)$
defines a $(w,\kappa)$-wall as in definition \eqref{eq:WallDefn}
so $\xi=w\pmod{2}$ and $-4\kappa\leq
\xi^2 < 0$ as in \eqref{eq:WallLabelDefn}. 
Let $0\leq\ell<\kappa$ be the integer defined by $-4(\kappa-\ell) =
\xi^2$. Then there is a $\U(2)$ bundle $E$ and a splitting
$E=L\oplus w\otimes L^{-1}$ with $\xi=w\otimes
L^{-2}$, where $c_1(E)=w$ and $p_1(\fg_E) = -4(\kappa-\ell) =
\xi^2$. The corresponding $\SO(3)$ bundle $\fg_E\cong \ubarRR\oplus\xi$ 
admits a connection, reducible with respect to this splitting, given by the
trivial connection on $\ubarRR$ and a $g$-anti-self-dual connection on
the line bundle $\xi$ if and only if $c_1(\xi)\cdot[\omega_g]=0$,
where $\omega_g$ is a $g$-self-dual, harmonic two-form.

Recall that any (real) linear representation space $V$ of $S^1$ can be
decomposed as $V=V_\CC\oplus V_\RR$, where $S^1$ acts trivially on the
real space $V_\RR$ and by a direct sum of one-dimensional
representations on the complex space $V_\CC$. Furthermore, the only
linear irreducible representation which is free away from the origin
is the standard representation on $\CC$.  To determine the possible
dimensions for the cohomology of the elliptic deformation complex at a
reducible anti-self-dual $\SO(3)$ connection, we shall make use of the
following general observation of Freed and Uhlenbeck:

\begin{prop}
\label{prop:FUReducible}
\cite[Proposition 4.9]{FU}
Let $X$ be a closed, oriented four-manifold with Riemannian metric $g$ and
$b^1(X)=0$. Let $A$ be a reducible, unitary, projectively
anti-self-dual connection on a rank-two, Hermitian bundle $E$ over
$X$. If $H^\bullet_A$ is the cohomology of the deformation complex for
the induced reducible, anti-self-dual $\SO(3)$ connection on $\fg_E$,
then there are real isomorphisms
$$
H_A^0 \cong H^0(X;\RR) \cong \RR,
\quad
H_A^1 \cong \CC^q,
\quad\text{and}\quad
H_A^2 \cong \CC^p\oplus H^+(X;\RR),
$$
where $\Stab_A\cong S^1$ acts by the standard action on $\CC^q$,
$\CC^p$, and by the trivial action on $H^0(X;\RR)$, $H^+(X;\RR)$. Here,
$\Stab_A\cong \SO(2)$ denotes the stabilizer in $\sG_E$ of the
anti-self-dual $\SO(3)$ connection $A^{\ad}$ on $\fg_E$ and, as
usual, $\sG_E$ is the group of unitary automorphisms of $E$ with
determinant one.
\end{prop}

The decomposition given by Freed-Uhlenbeck can described more
explicitly by borrowing from \cite[pp. 69--70]{FU} and the
corresponding discussion for reducible $\PU(2)$ monopoles in
\cite{FL2}. Recall from \cite[\S 3.4]{FL2} that the elliptic
deformation complex for the $g_0$-anti-self-dual equation at a
reducible connection $A = A_1\oplus A_2$ on $E=L_1\oplus L_2$, where
$A_2 := A_e\otimes A_1^*$ and $L_2 := (\det E)\otimes L_1^{-1}$,
splits into a {\em tangential\/} and a {\em normal\/} deformation
complexes; the tangential complex corresponds to the trivial $\Stab_A
\cong S^1$ action. We can thus apply Lemma 3.17 in
\cite{FL2} to partially compute the cohomology groups $H_A^\bullet$ of the
deformation complex for $A$ on $\fg_E$,
$$
\begin{CD}
\Omega^0(\fg_E)
@>{d_A}>>
\Omega^1(\fg_E)
@>{d_A^+}>>
\Omega^+(\fg_E).
\end{CD}
$$
The rolled-up full elliptic deformation complex for $A$ has index
$-2p_1(\fg_E) - 3(1-b^1+b^+)$, that is $-2\xi^2-6$ for
$b^1(X)=0$ and $b^+(X)=1$, with
\begin{equation}
\label{eq:XiLabelsSO3Reduction}
\xi = \det E\otimes L_1^{-2} 
\quad\text{and}\quad 
\fg_E\cong i\RR\oplus \xi.
\end{equation}
When $A$ is reducible, as above, the
tangential and normal deformation complexes are
\begin{gather}
\label{eq:TangentialComplex}
\begin{CD}
\Omega^0(i\RR)
@>{d}>>
\Omega^1(i\RR)
@>{d^+}>>
\Omega^+(i\RR),
\end{CD}
\\
\label{eq:NormalComplex}
\begin{CD}
\Omega^0(\xi)
@>{d_{A_\xi}}>>
\Omega^1(\xi)
@>{d_{A_\xi}^+}>>
\Omega^+(\xi),
\end{CD}
\end{gather}
where $A_\xi = A_e\otimes 2A_1^*$.
Of course, the cohomology of the tangential deformation complex
\eqref{eq:TangentialComplex} is isomorphic to $H^\bullet(X;\RR)$, so
$H^0(X;\RR)\cong \RR$ while $H^1(X;\RR)=0$ and $H^+(X;\RR)\cong \RR$
by hypothesis; thus the index of its rolled-up complex is $-2$. The
condition that $H^1(X;\RR)=0$ ensures that reducibles do not occur in
obstructed 
positive-dimensional families and so are isolated points of their
moduli spaces. The fact that the reducibles are isolated if and only if
$H^1(X;\RR)=0$ can also be seen from \cite[Proposition 2.2.6]{DK}; see
\cite[Proposition 4.2.15]{DK} in addition.

Turning to the normal deformation complex \eqref{eq:NormalComplex}, we
recall from
\cite[Equation (3.45)]{FL2} that the (real) index of its rolled-up complex 
is given by 
\begin{equation}
\label{eq:DefnNxi}
2n := 2n_\xi = -2p_1(\fg_E)-(e(X)+\sigma(X)) = -2\xi^2-2(1-b^1+b^+). 
\end{equation}
For $b^1(X)=0$ and $b^+(X)=1$, this becomes
$2n = -2\xi^2-4$. Observe that the rolled-up normal deformation
complex has positive (real) index when $\xi^2\leq -3$, zero index
when $\xi^2 = -2$, and index $-2$ when $\xi^2=-1$. {}From condition
\eqref{eq:WallLabelDefn} we need only ever consider $\xi^2 \leq -1$ and 
thus $n_\xi \geq -1$ (when $b^1(X)=0$ and $b^+(X)=1$).

Proposition \ref{prop:FUReducible} implies that for
the normal deformation complex  \eqref{eq:NormalComplex}
we must have 
$$
H_{A_\xi}^0 = 0,
\quad
H_{A_\xi}^1 \cong \CC^{p + n},
\quad\text{and}\quad
H_{A_\xi}^2 \cong \CC^p,
$$
where $n = -\xi^2-2$ is the complex index of the rolled-up normal
complex. Hence, 
\begin{equation}
\label{eq:CohomEllipticCplxAtReducible}
H_A^0 \cong \RR,
\quad
H_A^1 \cong \CC^{p+n},
\quad\text{and}\quad
H_A^2 \cong \CC^p\oplus \RR.
\end{equation}
The preceding consequence of Proposition \ref{prop:FUReducible} can be
sharpened if the metric $g$ on $X$ is suitably generic, though this refinement
is not essential for our purposes.

\begin{thm}
\label{thm:GenericMetricReducible}
\cite[Theorem 4.19]{FU}
Continue the hypotheses and notation of Proposition
\ref{prop:FUReducible}, with $\xi$ a complex line bundle over $X$ as in
\eqref{eq:XiLabelsSO3Reduction} and
$n := -\xi^2-2\geq -1$, the index of the rolled-up normal deformation
complex. If, in addition, the metric $g$ on $X$ is {\em $\xi$-generic\/},
that is it satisfies $[\omega(g)]\in W_\xi$ but is otherwise generic, then
the following hold.
\begin{itemize}
\item
If $n\geq 0$, then $H_{A_\xi}^1 \cong \CC^n$ and
$H_{A_\xi}^2 \cong 0$. Thus,
$$
H_A^0 \cong \RR,
\quad
H_A^1 \cong \CC^n,
\quad\text{and}\quad
H_A^2 \cong \RR.
$$
\item
If $n = -1$, then $H_{A_\xi}^1 \cong 0$ and
$H_{A_\xi}^2 \cong \CC$. Thus,
$$
H_A^0 \cong \RR,
\quad
H_A^1 \cong 0,
\quad\text{and}\quad
H_A^2 \cong \CC\oplus\RR.
$$
\end{itemize}
\end{thm}

\begin{proof}
The argument of Freed-Uhlenbeck in \cite[pp. 69--73]{FU} leading to the
proof of \cite[Theorem 4.19]{FU} makes no use of the facts that $b^+(X)=0$,
that $X$ is simply-connected, or that $E$ is an $\SU(2)$ bundle with
$c_2(E)=1$. Observe that the variations of the metric preserve the
curvature, $F_{A_\xi}$.  The period points of these metrics thus remain in
the wall $W_\xi$.  Hence, the conclusions follow directly from \cite{FU}.
\end{proof}

\subsection{Kuranishi obstructions in the parameterized moduli space}
\label{subsec:ParameterizedKuranishi}
It remains to extend the remarks of the preceding section concerning
singularities in the moduli space $M_\kappa^w(X,g_0)$ to the case of the
one-parameter family of moduli spaces $M_\kappa^w(X,\gamma)$. 

The usual Kuranishi model for a neighborhood of $[A_0]\in
M_\kappa^w(X,g_0)$ results from considering the section $A\mapsto
P_+(g_0)F_A$ of the Hilbert bundle
$\sA_E(X)\times_{\sG_E}L^2_\ell(\Lambda^{+,g_0}\otimes\fg_E)$. If
$g_t$, $t\in [-1,1]$, is a smooth path of metrics on $T^*X$, we may
choose a smooth path $f_t \in \GL(T^*X)$ such that $g_t = f_t(g_0)$,
$t\in [-1,1]$; the transformations $f_t$ are uniquely determined by
$g_t$ and $g_0$ by requiring the $f_t$ to be symmetric.  Therefore, we
have
$$
P_+(g_t) = f_t^{-1}\circ P_+(g_0)\circ f_t:\Lambda^2\to \Lambda^{+,g_t}
$$
and so instead of considering the zero-locus of the section
$(A,t)\mapsto P_+(g_t)F_A$ of
$(\sA_E(X)\times[-1,1])\times_{\sG_E}L^2_\ell(\Lambda^{+,g_t}\otimes\fg_E)$,
a Hilbert bundle depending on $t$, we may equivalently consider the
zero-locus of the section $(A,t)\mapsto P_+(g_0)f_t(F_A)$ of the {\em
fixed\/} Hilbert bundle
$(\sA_E(X)\times[-1,1])\times_{\sG_E}L^2_\ell(\Lambda^{+,g_0}\otimes\fg_E)$.
The construction of the local Kuranishi model for a neighborhood of
$[A_0,g_0]\in M_\kappa^w(X,\gamma)$ then follows in the usual way
\cite{DK} to give:

\begin{lem}
\label{lem:ParameterizedKuranishi}
Let $E$ be a Hermitian two-plane bundle over a closed, oriented,
smooth four-manifold $X$ with $b^+(X)=1$ and $\gamma:[-1,1]\to\Met(X)$
a generic path of Riemannian metrics on $T^*X$.  Let $[A_0,g_0]$ be a
point in $M_E(\gamma)$ such that $\gamma(0)=g_0$ and $A_0$ is a
reducible unitary connection on $E$. Then $A_0$ has stabilizer
$\Stab_{A_0}\cong S^1$ in $\sG_E$ and there are
\begin{itemize}
\item
An open, $S^1$-invariant neighborhood $\sO=\sO_{A_0}$ of the origin in
$H^1_{A_0}$ and a positive constant $\eps=\eps_{A_0}$ together with a
smooth, $S^1$-equivariant embedding $\bga=\bga_{A_0,g_0}$,
$$
\bga:\sO\times(-\eps,\eps)\subset H^1_{A_0}\oplus\RR
\hookrightarrow
\sA_E,
$$
with $\bga(0,0) = (A_0,0)$ and
$(M_E\cap\bga(\sO\times(-\eps,\eps)))/S^1$ an open
neighborhood of $[A_0,g_0]$ in $M_E(\gamma)$; and
\item A smooth, $S^1$-equivariant map $\bchi = \bchi_{A_0,g_0}$,
$$
\bchi:\sO\times(-\eps,\eps)\subset H^1_{A_0}\oplus\RR
\to
H^2_{A_0}
$$
such that $\bga$ restricts to a
smoothly-stratified diffeomorphism from
$(\bchi^{-1}(0)\cap(\sO\times(-\eps,\eps)))/S^1$ onto
$(M_E\cap\bga(\sO\times(-\eps,\eps)))/S^1$.
\end{itemize}
The obstruction map can be defined by
$$
\bchi(a,t) := \Pi_{A_0,g_0}P_+(g_0)f_t(F(A_0+a)),
$$ 
while the existence of the embedding $\bga$ follows from the
Banach-space implicit function theorem.  Here, $\Pi_{A_0,g_0}$ is the
$L^2$-orthogonal projection from $L^2(\Lambda^{+,g_0}\otimes\fg_E)$
onto the span of the eigenvectors of $d_{A_0}^{+,g_0}d_{A_0}^{*,g_0}$
on $L^2(\Lambda^{+,g_0}\otimes\fg_E)$ with eigenvalue less than
$\half\mu$, where $\mu$ is the first positive eigenvalue of this
Laplacian.
\end{lem} 

The model described in Lemma \ref{lem:ParameterizedKuranishi} can be
simplified: we may remove the factor $H^+(X;\RR)$ from the
obstruction space at the cost of restricting the obstruction map
$\bchi$ to the preimage $\bchi_\RR^{-1}(0)\subset \RR\oplus\CC^n$,
where $\bchi_\RR$ is defined by \eqref{eq:RealKuranishi}.

\begin{lem}
\label{lem:SimpleParameterizedKuranishi}
Continue the notation and hypotheses of Lemma
\ref{lem:ParameterizedKuranishi}. Assume that $b^1(X)=0$ and that $A_0$ is
a reducible $\SO(3)$ connection with respect to the splitting
$\fg_E\cong i\RR\oplus\xi$. Then the following hold.
\begin{itemize}
\item
If $n_\xi\geq 1$, then there are
an open, $S^1$-invariant neighborhood $\sO$ of the origin in
$\CC^n$ and a positive constant $\eps$ together with a
smooth, $S^1$-equivariant embedding,
$$
\bga:\sO\times(-\eps,\eps)\subset \CC^{n_\xi}\oplus\RR
\hookrightarrow
\sA_E,
$$
with $\bga(0,0) = (A_0,0)$ and
$M_E(\gamma)\cap\bga(\sO/S^1\times(-\eps,\eps))$, an open
neighborhood of $[A_0,g_0]$ in $M_E(\gamma)$. 

The map $\bga$ restricts to a smoothly-stratified diffeomorphism
from $\bchi_\RR^{-1}(0)\cap(\sO/S^1\times(-\eps,\eps))$ onto
$M_E(\gamma)\cap\bga(\sO/S^1\times(-\eps,\eps))$. Moreover, there is
an $S^1$-equivariant diffeomorphism
$\bchi_\RR^{-1}(0)\cap(\sO\times(-\eps,\eps))\cong \CC^n$.
\item 
If $n_\xi\leq 0$, then $[A_0,g_0]$ is an isolated point of the parameterized
moduli space $M_E(\gamma)$, with vanishing obstruction when $n_\xi=0$ and
obstruction space $\CC$ when $n_\xi=-1$.
\end{itemize}
\end{lem} 

\begin{proof}
The proof can be extracted from a similar argument used to prove the
wall-crossing formula in \cite[Proposition 2.12]{DonHCobord}.  By
\eqref{eq:CohomEllipticCplxAtReducible},
Theorem \ref{thm:GenericMetricReducible}, and the fact that the metric
$g_0$ is $\xi$-generic, 
$H_{A_0}^1=\CC^{n_\xi}$ and $H^2_{A_0}=\RR$ when
$n_\xi\geq 0$, while $H_{A_0}^1=0$ and
$H^2_{A_0}=\RR\oplus\CC$ when $n_\xi=-1$. With respect to the
splitting $\Lambda^+\otimes\fg_E\cong
i\Lambda^+\oplus(\Lambda^+\otimes\xi)$, we may write $\bchi =
(\bchi_\RR,\bchi_\CC)$, where
\begin{equation}
\label{eq:RealKuranishi}
\begin{aligned}
\bchi_\RR(z,t) 
&= 
\langle \bchi(z,t),f_t(\omega_0)\rangle f_t(\omega_0), 
\\
&= 
\langle f_t(F(A_0+\bga(z,t))),f_t(\omega_0)\rangle f_t(\omega_0), 
\quad 
(z,t)\in H_{A_0}^1\oplus(-\eps,\eps).
\end{aligned}
\end{equation}
According to Lemma \ref{lem:LocalWallCrossing} 
we may assume without loss that the path
$\omega_t\equiv f_t(\omega_0)$ crosses the wall $W_\xi$ transversely, so
$(\partial\bchi_\RR/\partial t)(0,0) \neq 0$.  Hence, the map
$\pi_\CC:\bchi_\RR^{-1}(0)\cong \CC^{n_\xi}$ is an
$S^1$-equivariant diffeomorphism.
\end{proof}

\begin{rmk}
\label{rmk:SettingForObstructions}
Suppose $\xi\in H^2(X;\ZZ)$ defines a reduction $\fg_{E_0}\cong i\RR\oplus\xi$.
Because of Lemma \ref{lem:SimpleParameterizedKuranishi} we see that,
provided $b^1(X)=0$, we need only consider the following cases where the
family of reducibles in $M_{\kappa_0}^w(X,C)$ consists of either
\begin{itemize}
\item
A single point $[A_0,g_0]$ with vanishing obstruction, when $n_\xi\geq 0$,
or
\item
A single, isolated point $[A_0,g_0]$ with obstruction $H^2_{A_0,g_0}\cong \CC$,
when $n_\xi=-1$.
\end{itemize}
In the first case we would glue up an open neighborhood $\sU_0=\tilde\sU_0/S^1$
$M_{\kappa_0}^w(X,C)$ of the point $[A_0,g_0]$; in the second, $\sU_0$
could simply be the isolated point $[A_0,g_0]$ itself when gluing up. Here,
in either case, we can assume that
$\tilde\sU_0\subset\sA_{\kappa_0}^w(X)\times\Met(X)$ is contained in
the slice $(A_0+\Ker d_{A_0}^{*,g_0})\times C$, so $\tilde\sU_0/\sG_E =
\sU_0/S^1$. We shall make use of this simplification --- namely, that we
need only glue up anti-self-dual connections rather than {\em extended\/}
anti-self-dual connections even when gluing up non-smooth points --- when
defining the space of gluing data,
$\Gl_\kappa^w(X,\sU_0,\Sigma,\lambda_0)$, in \S
\ref{sec:Splicing}. 
\end{rmk}

\subsection{Paths of metrics and the period map}
There is a technical point with regard to the definition of the family
of background connections over $X$, in the presence of a path of
metrics, which requires some comment. Suppose $g_{\pm 1}$ are two
$C^\8$ metrics on $TX$ and $\omega_{\pm 1}$ are the corresponding
self-dual harmonic two-forms on $X$ (we assume $b^+(X)=1$). The
metrics $g_{\pm 1}$ thus define points $[\omega_{\pm 1}]$ in
$H^+(X;\RR)-\{0\}$, lying in the positive cone $\Omega_X$ of
$H^+(X;\RR)$ since if $\omega = \omega_{-1}$ or $\omega_1$ then
$$
\langle[\omega]\cup[\omega],[X]\rangle
=
\int_X\omega\wedge\omega
=
\int_X\omega\wedge*\omega
=
\|\omega\|_{L^2(X)}^2 > 0.
$$
In particular, the points $[\omega_{\pm 1}]$ lie in chambers $\sC_{\pm}$
of $\Omega_X$. If $\sC_+=\sC_-=\sC$, that is, these points lie in the {\em
same\/} chamber, there is a question of whether a path
$[\omega_t]\in\sC$, $t\in[-1,1]$, joining $[\omega_{-1}]$ to
$[\omega_1]$ can be lifted to a path of metrics $g_t$ in $\Met(X)$
joining $g_{-1}$ to $g_1$. 
For this purpose, one needs to consider the
`period map'
$$
\pi:\Met(X) \to \PP(H^2(X;\RR)) := H^2(X;\RR)/\RR^*,
\qquad
g \mapsto \pi(g) := \sH^+(g),
$$
where $\sH^\pm(g)$ are the spaces of self-dual and anti-self-dual
harmonic two-forms.  As remarked in
\cite[p. 434]{KotschickMorgan} it is not known whether two metrics whose 
period points lie in the same chamber can be connected by a path of
metrics whose period points lie in the same chamber. However,
following result of Donaldson will suffice for our needs:
\footnote{The statement of Proposition
\ref{prop:PeriodMapIsTransverseToWalls} is taken from
\cite[p. 336]{DonConn}: the period map defined there is precisely dual to the 
one considered in \cite{DK}.}

\begin{prop}
\label{prop:PeriodMapIsTransverseToWalls}
\cite[\S 6]{DonConn}, \cite[Proposition 4.3.14]{DK}
Let $X$ be a closed, oriented four-manifold and let $\GG(b^2,b^+)$ be
the real Grassmann manifold of $b^+(X)$-dimensional subspaces of
$H^2(X;\RR)$. Let $\UU\subset\GG(b^2,b^+)$ be the open subset 
of planes on which
the intersection form $Q$ of $X$ is positive definite.  Let
$\pi:\Met(X)\to\UU$ be the period map $g\mapsto \sH^+(g)$.
Let $\xi\in H^2(X;\ZZ)$ be a class such that $\xi^2<0$ and let
$W_\xi \subset \UU = \{K\in\UU: \xi\perp K\}$. Then $W_\xi$
is a $b^+(X)$-codimensional submanifold and the period map $\pi$ is
transverse to $W_\xi$.
\end{prop}

When $b^+(X)=1$ then $\GG(b^2,b^+)=\PP(H^2(X;\RR))$, $\UU = \Omega_X$,
and $W_\xi$ is a wall.
For the purposes of proving the wall-crossing formula, it is enough to
consider paths of metrics satisfying the conclusions of the following
lemma, which comprise the preamble to Theorem 5.2.2 in
\cite{KotschickMorgan}.

\begin{lem}
\label{lem:LocalWallCrossing}
Let $X$ be a closed, oriented, smooth four-manifold with $b^+(X)=1$
and suppose that $\sC_{\pm}$ are chambers in $\Omega_X$ separated by a
single geometric wall $W$.  Then there is a generic path of smooth metrics
$\gamma:[-1,1]\to\Met(X)$, $t\mapsto g_t$, with the following properties:
\begin{enumerate}
\item
Fix a class $w\in H^2(X;\ZZ)$ such that no $\SO(3)$ bundle $V$
over $X$ with $w_2(V)=w\pmod{2}$ admits a flat connection. Then for
every $\U(2)$ bundle $E$ with $c_1(E)=w$ and $\kappa \equiv -\quarter
p_1(\fg_E) > 0$, the parameterized moduli space
$M_\kappa^w(X,\gamma)$ contains no flat connections
and is a $C^\8$ manifold away from the reducible connections.
\item
For all $t\neq 0$, the compactification $\barM_\kappa^w(X,\gamma)$
contains no points whose background connection is reducible or flat.
\item
The self-dual two-form $\omega(g_0)$ lies in a single geometric
$(w,\kappa)$-wall, $W$, while $\omega(g_t)\in \sC_-$ for $t<0$ and
$\omega(g_t)\in\sC_+$ for $t>0$.
\item
The path of self-dual two-forms $\omega(g_t)$ defined by $\gamma$ crosses 
the wall $W$ transversely.
\end{enumerate}
\end{lem}

\begin{proof}
We may restrict our attention to an open
neighborhood $U\subset \PP(H^2(X;\RR))$ such that $U\subset\sC^-\cup
W\cup\sC^+$ and so we can express $U$ as a disjoint union $U =
U\cap\sC^-\sqcup U\cap W \sqcup U\cap\sC^+$. Furthermore, we choose
$U$ so that it only meets a single geometric wall $W$.  The wall $W$
is defined by some, not necessarily unique class $\xi\in
H^2(X;\ZZ)$.  According to Proposition
\ref{prop:PeriodMapIsTransverseToWalls}, the period map $\pi$ is
transverse to $W$ and so $\pi^{-1}(W)$ is a smooth, codimension-one
submanifold of $\Met(X)$. Plainly, we may choose a generic path
$\gamma:[-1,1]\to\Met(X)$ whose image is contained in the open subset
$\pi^{-1}(U)$ and which is transverse to $\pi^{-1}(W)$. Consequently,
the induced path $\omega(\gamma) \equiv \pi(\gamma)$ in $H^2(X;\RR)/\RR^+$
is transverse to $W$. Finally, we may choose $\gamma$ short enough
that $\omega(g_t)\in\sC_-$ for all $t<0$ and $\omega(g_t)\in\sC_+$ for
all $t>0$.
\end{proof}

\subsection{Links of reducibles and the computational problem}
It remains to state the computational problem which concerns the
present article and its sequels. Given the situation of Lemma
\ref{lem:LocalWallCrossing} one can prove (see \cite[Theorem
5.2.2]{KotschickMorgan}, though we shall give a different argument
elsewhere) the following local wall-crossing result:
$$
D_{X,g_{+1}}^w(z) - D_{X,g_{-1}}^w(z) 
= 
\sum_\xi (-1)^{(w-\xi)^2/4}\delta^w_\xi(z),
\quad z\in\AAA(X),
$$
where the sum is over all $\xi\in H^2(X;\ZZ)$ which define the
geometric $(w,\kappa)$-wall $W$ containing $\omega(g_0)$ and which satisfy
$\xi\cdot\omega(g_{-1}) < 0 <\xi\cdot\omega(g_1)$, with $\deg z
= 8\kappa-3(2-b^1)$. The terms $\delta^w_\xi(z)$ have the following properties:
\begin{itemize}
\item
$\delta^w_\xi$ is independent of where the path $\gamma$ crosses
the wall $W$,
\item
$\delta^w_\xi$ changes sign when the direction of crossing is reversed,
\item
$\delta^w_{-\xi} = (-1)^{\xi^2}\delta^w_\xi$.
\end{itemize}
It follows easily from these properties that $D_{X,g}^w(z) =
D_{X,g'}^w(z)$ for any pair of metrics $g$, $g'$ whose period points
lie in the same chamber $\sC$ of $\Omega_X$, so $D_{X,\sC}^w(z) =
D_{X,g}^w(z)$ is well-defined and the preceding wall-crossing formula
can be replaced by its global analogue:
$$
D_{X,\sC_+}^w(z) - D_{X,\sC_{-1}}^w(z) 
= 
\sum_\xi (-1)^{(w-\xi)^2/4}\delta^w_\xi(z),
\quad z\in\AAA(X),
$$
with $\sC_\pm$ in the same component of $\Omega_X$ and the sum is over
all $\xi\in H^2(X;\ZZ)$ which define geometric $(w,\kappa)$-walls
satisfying $\xi\cdot\sC_- < 0 <\xi\cdot\sC_+$, with
$\deg z = 8\kappa-3(2-b^1)$. 

Technicalities aside, each term $\delta^w_\xi(z)$ is defined by pairings
$$
\delta^w_\xi(z)
=
\langle\mu(z),[\bar\bL^w_{\kappa,\xi}]\rangle,
$$
where $\mu:H_i(X;\QQ)\to H^{4-i}(\sB_\kappa^w(X)\times\gamma;\QQ)$
and $\bL^w_{\kappa,\xi}$ is the ``link'' in $M_\kappa^{w,*}(X,\gamma)$ of the 
level containing the reducible $\SO(3)$ connection defined by $w,\xi,\kappa$,
$$
[A_0]\times\Sym^\ell(X) \subset \barM_\kappa^w(X,\gamma),
$$
while $\bar\bL^w_{\kappa,\xi}$ is its closure in $\barM_\kappa^w(X,\gamma)$.
Here, $\gamma$ is a path satisfying the conclusions of Lemma
\ref{lem:LocalWallCrossing} such that $A_0$ is a reducible
$g_0$-anti-self-dual $\SO(3)$ connection on an $\SO(3)$ bundle
$\fg_{E_0}$ with $w_2(\fg_{E_0})=w\pmod{2}$ and $-\quarter
p_1(\fg_{E_0})=\kappa-\ell$. The reduction of the $\SO(3)$ bundle
$\fg_{E_0}$ is defined by the class $\xi$: we have $E_0 = L\oplus
(w\otimes L^{-1})$ and $\fg_{E_0}\cong
\xi\oplus\ubarRR$, where $\xi = w\otimes L^{-2}$, according to our
standard convention, and thus $p_1(\fg_{E_0}) = \xi^2$. 

The heart of the problem then is to show, using a more precise version
of this ``definition'', that the $\delta^w_\xi$ depend at most on the
cohomology classes $w, \xi \in H^2(X;\ZZ)$, the homotopy type of $X$,
and the degree of $z$. The above computation has been carried out
directly when $\ell=0,1,2$ (see \cite{DonHCobord}, \cite{Leness},
\cite{KotschickBPlus1}, \cite{KotschickMorgan}, \cite{Yang}), which verifies
the conjecture in these cases, so our task is to extend this work so
we may conclude that the $\delta^w_\xi(z)$ have the desired invariance
properties even if an explicit formula cannot be deduced by purely
topological methods.


\section{Splicing connections}
\label{sec:Splicing}
We construct a {\em splicing map\/} $\bga'$ which will give a homeomorphism
from the total space of a gluing data bundle into an open neighborhood in
$\barsB_\kappa^w(X)$ of a level $M_{\kappa-\ell}^w(X)\times \Sigma$, where
$\Sigma\subset\Sym^\ell(X)$ is a smooth stratum.  In this section we keep
the Riemannian metric $g$ on $TX$ fixed, though we shall later allow the
metrics $g$ to vary in a precompact open subset of the space $\Met(X)$ of
all $C^\8$ Riemannian metrics on $TX$. When defining our spaces of
connections and gauge transformations, we shall always choose to work
$L^2_k$ connections modulo $L^2_{k+1}$ gauge transformations with $k\geq
4$: though only $k\geq 2$ is strictly necessary, we shall occasionally need
to consider estimates whose constants depend on the $C^0$ norm of
$F_A^+$. Thus, {\em for the remainder of this article and its sequels\/}
\cite{FLKM2}, \cite{FLKM3}, \cite{FLKM4}, we shall assume $k\geq 4$ in
order to make use of the Sobolev embedding $L^2_4\subset C^0$.

\subsection{Uhlenbeck neighborhoods}
\label{subsec:UhlenbeckCompactification}
In order to associate sequences of mass centers and scales to a sequence of
points $[A_\al]$ in $M_E(X,g)$ converging to a point $[A_0,\bx]$ in
$\barM_E(X,g)$, we need a suitable finite open cover of $\barM_E(X,g)$; for
simplicity, we take $C=\{g\}$.  Let $\varrho$ be the injectivity radius of
the Riemannian manifold $(X,g)$ and let $\dist_g(\cdot,\cdot)$ be the distance
function on $X$ defined by the metric $g$.

\begin{defn}
\label{defn:UhlenbeckNeighborhood}
  Suppose $[A_0,\bx]$ is a point in $\barM_E(X,g)$, where $A_0$ is a
  connection on a Hermitian, rank-two bundle $E_0$ over $X$ with
  $c_2(E_0)<c_2(E)$, $c_1(E_0)=c_1(E)$ and $\eps,\de,M$ are positive
  constants with $M\ge 1$, and $\eps<2\pi$, and $M\de\le
  \quarter\{1,\varrho,d_\bx\}$, where $d_\bx := \min_{i\ne
  j}\dist(x_i,x_j)$ is the minimum distance separating the distinct points
  $x_i$ of the multi-set $\bx = (x_1,\dots,x_m)$.  Assume $M\de$ is small
  enough that
$$
\int_{B(x_i,M\de)}|F_{A_0}|^2\,dV < \quarter\eps^2.
$$
Then $[A]\in M_E(X,g)$ lies in an
$(\eps,\de,M)$ {\em Uhlenbeck open neighborhood\/} $\sU$ of $[A_0,\bx]$ if the
following hold:
\begin{itemize}
\item There is an $L^2_{k+1,\loc}$ bundle map $u:E|_{X\less\{\bx\}}\to 
E_0|_{X\less\{\bx\}}$ such that the connections $u(A)$ and $A_0$
are $L^2_k$ close on $X\less B_\bx(\half\de)$:
$$
\|u(A) - A_0\|_{L^2_k(A_0;X\less B(\bx,\de/2))} < \eps;
$$
\item The measure $|F_A|^2$ is close to $|F_{A_0}|^2 +
8\pi^2\ka_\bx\de_\bx$:
\begin{align*}
\left|\int_{B(x_i,\de)}\left(|F_A|^2-|F_{A_0}|^2
- 8\pi^2\ka_i\de_{x_i}\right)\,dV\right| &< \quarter\eps^2, \\
\left|\int_{\Om(x_i;\de,M\de)}\left(|F_A|^2-|F_{A_0}|^2
\right)\,dV\right| &< \quarter\eps^2,
\end{align*}
for $i=1,\dots,m$ and
where $(x_1,\dots,x_m)$ is a representative for $\bx$ with multiplicities
$\ka_1,\dots,\ka_m$. 
\end{itemize}
\end{defn}

The constraints on the curvature $F_A$ in Definition
\ref{defn:UhlenbeckNeighborhood} imply that
\begin{align}
\left|8\pi^2\ka_i - \int_{B(x_i,\de)}|F_A|^2\,dV\right| &< \half\eps^2, 
\notag\\
\int_{\Om(x_i;\de,M\de)}|F_A|^2\,dV &< \half\eps^2, 
\label{eq:UhlenbeckCurvCond}\\
\left|8\pi^2\ka_i - \int_{B(x_i,M\de)}|F_A|^2\,dV\right| &< \eps^2.
\notag
\end{align}
Observe that no cutoff functions are required to smooth out the implied
characteristic functions of the balls 
with center $x_i$ and radii $\de$, $M\de$ above as Definition
\ref{defn:UhlenbeckNeighborhood} rules
out the possibility of the curvature $F_A$ concentrating on their
boundaries.

\subsection{Gluing data}
\label{subsec:GluingData}
In this section we describe the gluing-data bundles we shall need to
parameterize an open neighborhood in $\barM_\kappa^w(X,C)$ of a level
$(M_{\kappa-\ell}^w(X,C)\times\Sigma)\cap \barM_\kappa^w(X,C)$.

\subsubsection{Connections over the four-sphere}
\label{subsubsec:ConnectionsOnFourSphere}
We begin by defining the required gluing data associated with the
four-spheres: this will, essentially, define the fibers of our
gluing-data bundles.  Let $\kappa\ge 1$ be an integer and let $S^4$
have its standard round metric of radius one. Let $E$ be a
Hermitian, rank-two vector bundle over $S^4$ with $c_2(E) =
\kappa$ and let $\fg_E = \su(E)$ be the corresponding $\SO(3)$ 
bundle with $-\quarter p_1(\fg_E) = \kappa$.

A choice of frame $f$ in the principal
$\SO(4)$ frame bundle $\Fr(TS^4)$ for
$TS^4$, over the north pole $n\in S^4$, defines a conformal
diffeomorphism,
\begin{equation}
\label{eq:ConformalDiffeo}
\varphi_n:\RR^4 \to S^4\less\{s\},
\end{equation}
that is inverse to a stereographic projection from the south pole
$s\in S^4\subset \RR^5$; let $y(\,\cdot\,): S^4\less\{s\}\to\RR^4$ be
the corresponding coordinate chart.

\begin{defn}
\label{defn:CenterScale}
\cite[pp. 343--344]{TauFrame}
Continue the notation of the preceding paragraph.  The {\em center\/}
$z=z([A],f)\in\RR^4$ and the {\em scale\/} $\lambda=\lambda([A],f)$ of a
point $[A]\in\sB_E(S^4)=\sB_\kappa(S^4)$ are defined by
\begin{align}
z([A],f) 
&:= \frac{1}{8\pi^2\kappa}\int_{{\Bbb R}^4}y|F_A|^2\,d^4y, 
\label{eq:MassCenter}\\
\lambda([A],f)^2 
&:= \frac{1}{8\pi^2\kappa}\int_{{\Bbb R}^4} |y-z([A],f)|^2|F_A|^2\,d^4y. 
\label{eq:Scale}
\end{align}
A point $[A]$ is {\em centered\/} if $z([A],f)=0$ and $\lambda([A],f) = 1$.
An ideal point $[A,\bx]\in \sB_{\kappa-\ell}(S^4)\times\Sym^\ell(S^4)$ is
{\em centered\/} if $z([A,\bx],f)=0$ and $\lambda([A,\bx],f) = 1$, where
the center and scale are defined by replacing $|F_A|^2$ with
$|F_A|^2+8\pi^2\sum_{x\in\bx}\delta_x$ in equations \eqref{eq:MassCenter}
and \eqref{eq:Scale}.
\end{defn}

For any $(\lambda,z)\in\RR^+\times\RR^4$, we define a conformal
diffeomorphism of $\RR^4$ by
\begin{equation}
\label{eq:CenterScaleDiffeo}
f_{\lambda,z}:\RR^4\to\RR^4,\qquad y\mapsto (y-z)/\lambda.
\end{equation}
It is easy to see that the point $[(f_{\lambda,z}^{-1})^*A]$ is
centered if $z = z([A],f)$ and $\lambda = \lambda([A],f)$.  Via the
standard action of $\SO(4)\subset\GL(\RR^4)$ on $\RR^4$ and the
conformal diffeomorphism $\varphi_n^{-1}:S^4\less\{s\}\to\RR^4$, we
obtain a homomorphism from $\SO(4)$ into the group of conformal
diffeomorphisms of $S^4$ fixing the south pole $s$. The group $\SO(4)$
acts on the quotient space of connections based at the south pole,
$\sB_\kappa^s(S^4) \cong \sA_\kappa(S^4)\times_{\sG_E}\Fr(\fg_E)|_s$, by
\begin{equation}
\label{eq:BasedSphereModSpaceSO(4)Action}
\SO(4)\times \sB_\kappa^s(S^4) \to \sB_\kappa^s(S^4),
\quad
(u,[A,q]) \mapsto u\cdot[A,q] = [(u^{-1})^*A,q].
\end{equation}
The assignment $f\mapsto z([A],f)$ is $\SO(4)$-equivariant, that is,
$z([A],uf) = uz([A],f)$ for $u\in \SO(4)$, and so $\lambda[A] =
\lambda([A],f)$ is independent of the choice of $f \in \Fr(TS^4)|_n$. The
center and scale functions define a smooth map
$(z,\lambda):\sB_\kappa(S^4)\to\RR^4\times\RR^+$ and the preimage of $(0,1)\in
\RR^4\times\RR^+$, namely the codimension-five submanifold of {\em
centered connections\/}
\begin{equation}
\label{eq:QuotientSpaceCenteredConn}
\sB_\kappa^\diamond(S^4)
:=
\{[A]: z[A] = 0 \text{ and } \lambda[A] = 1\}
\subset
\sB_\kappa(S^4), 
\end{equation}
serves as a canonical slice for the action of $\RR^4\times\RR^+$ on
$\sB_\kappa(S^4)$; we simply write $z[A]$ when the frame $f$ for $TS^4|_n$ is
understood and employ the decoration $\diamond$ to indicate a subspace of
centered connections. Since the action of $\SO(4)$ preserves the north and
south poles of $S^4$, it acts on $\sB_\kappa^{\diamond,s}(S^4)$, and so
$\sB_\kappa^{\diamond,s}(S^4)$ admits an action of $\SO(3)\times\SO(4)$,
with $\SO(3)$ varying the frame of $\fg_E|_s$.

A simple {\em a priori\/} decay estimate for the curvature of
connections is provided by the Chebychev inequality:

\begin{lem}\label{lem:Chebychev}
Let $[A] \in \sB_\kappa(S^4)$ and suppose $z([A],f) = 0$. If
  $\lambda = \lambda[A]$ and $R>0$, then
$$
\frac{1}{8\pi^2\kappa} \int_{|y|\ge R\lambda}|F_A|^2\,d^4y \le R^{-2}.
$$
\end{lem}

\begin{pf}
From equation \eqref{eq:Scale} we see that
$$
\lambda^2 \ge \frac{1}{8\pi^2\kappa}\int_{|y|\ge R\lambda}|y|^2|F_A|^2\,d^4y
\ge 
\frac{R^2\lambda^2}{8\pi^2\kappa}\int_{|y|\ge R\lambda}|F_A|^2\,d^4y
$$
and this gives the required bound.
\end{pf}

An immediate consequence of the Chebychev inequality is that if $A$ is
a centered anti-self-dual connection over $S^4$, with second Chern
number $\kappa$ so its energy is equal to $8\pi^2\kappa$, then $A$
cannot bubble outside the ball $B(0,1/2\pi)\subset\RR^4$, as it has
most energy $4\pi^2$ in this region by Lemma \ref{lem:Chebychev}.

We let $M^\diamond_\kappa(S^4)\subset \sB_\kappa^\diamond(S^4)$ denote the
subspace of anti-self-dual connections, centered at the north pole, while
$M^s_\kappa(S^4)\subset \sB_\kappa^s(S^4)$ is the subspace of
anti-self-dual connections, framed at the south pole. Finally, we let
$M^{\diamond,s}_\kappa(S^4)\subset \sB_\kappa^{\diamond,s}(S^4)$ denote the
subspace of anti-self-dual connections, which framed at the south pole and
centered at the north pole. The definition of Uhlenbeck convergence
(Definition \ref{defn:UhlenbeckTop}) for sequences of points $[A_\alpha]$
in $\sB_\kappa(S^4)$ generalizes easily to the case of sequences
$[A_\alpha,q_\alpha]$ in $\sB_\kappa^s(S^4)$, when $\Fr(\fg_E|_s)\cong
\SO(3)$ is given its usual topology. Since Lemma \ref{lem:Chebychev}
implies that points in $M^{\diamond,s}_\kappa(S^4)$ cannot bubble near the
south pole, the proof of \cite[Theorem 4.4.3]{DK} adapts with no change to
show that $M_\kappa^{\diamond,s}(S^4)$ has compact closure
$\barM^{\diamond,s}_\kappa(S^4)$ in the space of framed, ideal
anti-self-dual connections,
$$
\bigcup_{\ell=0}^\kappa
\left(M^s_{\kappa-\ell}(S^4) \times \Sym^\ell(S^4) \right).
$$
Note that it is the ideal points $[A_0,\bx]$ in
$\barM^s_{\diamond,s}\kappa(S^4)$ which are centered in the sense of
Definition \ref{defn:CenterScale} and not the background connections $A_0$.
Note also that the points in $S^4=\RR^4\cup\{\8\}$ representing the
multisets $\bx$ will be constrained to lie in the ball $\bar B(0,1/2\pi)$
by Lemma \ref{lem:Chebychev}.

Observe that since the bundle map $u$ in the definition of an Uhlenbeck
neighborhood respects the $\SO(3)$ action on frames, there is a global
action of $\SO(3)$ on the space $\barM^{\diamond,s}_\kappa(S^4)$.  This
action is free on all strata except $[\Theta]\times\Sym^\kappa(S^4)$, where
the background connection is trivial.  In addition, the $\SO(4)$ action on
$M^{\diamond,s}_\kappa(S^4)$ extends over $\barM^{\diamond,s}_\kappa(S^4)$.

\subsubsection{Gluing data bundles}
\label{subsubsec:GluingDataBundles}
We now define bundles of gluing data,
$\Gl_\kappa^w(X,\sU_0,\Sigma,\lambda_0)$, associated to a choice of bundle
data $(\kappa,w)$, smooth stratum $\Sigma\subset\Sym^\ell(X)$, subset
$C\subset \GL(T^*X)$, open neighborhood $\sU_0\subset M_{\kappa_0}^w(X,C)$,
metric $g$ on $TX$, and positive constant
$\lambda_0$. Here, $C\subset\GL(T^*X)$ is a compact, connected, smooth
submanifold-with-boundary containing the identity: the pair $(C,g)$
defines a compact, connected, smooth submanifold-with-boundary of
$\Met(X)$ which we also denote, when no confusion can arise, by
$C$. Our definition is motivated by those of Friedman-Morgan
\cite[\S\S 3.4.2--3.4.4]{FrM}, Kotschick-Morgan
\cite[\S 4]{KotschickMorgan}, and Taubes \cite{TauFrame}, \cite{TauStable}.

Choose $\kappa\in\quarter\ZZ$ and $w\in H^2(X;\ZZ)$ according to
condition \eqref{eq:BasicChoiceOfPontrjaginStieffel} and let $E$ be 
the corresponding Hermitian, rank-two vector bundle over $X$ defined by
condition
\eqref{eq:U2BundleCompatibleWithPontrjaginStieffel}.  Let $\ell$ be an
integer such that $1\le\ell\le \lfloor\kappa\rfloor$ and let
$\Sigma\subset\Sym^\ell(X)$ be a smooth stratum, defined by an integer
partition $\kappa_1\ge\cdots\ge\kappa_m\ge 1$ of $\ell$ such that
$\kappa_1+\cdots+\kappa_m=\ell$. Let $E_0$ be the Hermitian, rank-two
vector bundle over $X$ with $\det E_0 = \det E$ and $c_2(E_0) =
c_2(E)-\ell$. Let $\Fr(TX)$ denote the principal $\SO(4)$ bundle of
oriented, $g$-orthonormal frames for $TX$.  Suppose $\tilde\sU_0\subset
\sA_{E_0}(X)\times C$ projects to an open subset $\sU_0\subset
M_{E_0}(X,C)$. Observe that the fibered product
\begin{equation}
\tilde\sU_0\times_{\sG_{E_0}}\prod_{i=1}^m 
\left(\Fr(\fg_{E_0})\times_X\Fr(TX)\right)
\to \sU_0\times \prod_{i=1}^m X
\label{eq:PreFiberedProductBundle}
\end{equation}
is a principal $\prod_{i=1}^m(\SO(3)\times\SO(4))$ bundle. Let
$\Delta\subset\prod_{i=1}^m X$ be the full diagonal, that is, the subset of
$m$-tuples of points in $X$ where some point occurs with multiplicity at
least two. Let $\fS = \fS_\Sigma = \fS(\kappa_1,\dots,\kappa_m)$ be the
largest subgroup of the symmetric group on $m$ symbols which preserves the
multiplicities.  The group $\fS(\kappa_1,\dots,\kappa_m)$ acts freely on
$\prod_{i=1}^m X-\Delta$ and we have an identification
$$
\left(\prod_{i=1}^m X-\Delta\right)/\fS(\kappa_1,\dots,\kappa_m)
=
\Sigma\subset\Sym^\ell(X).
$$
The group $\fS(\kappa_1,\dots,\kappa_m)$ acts on the fibered product bundle
\eqref{eq:PreFiberedProductBundle} to give a principal bundle
\begin{equation}
\bFr_{E_0}(X,\sU_0,\Sigma) \to \sU_0\times \Sigma,
\end{equation}
with total space $\bFr_\kappa^w(X,\sU_0,\Sigma) = \bFr_{E_0}(X,\sU_0,\Sigma)$,
where
\begin{equation}
\bFr_{E_0}(X,\sU_0,\Sigma)
:= 
\left(\tilde\sU_0\times_{\sG_{E_0}}\prod_{i=1}^m 
\left(\Fr(\fg_{E_0})\times_X\Fr(TX)\right)\right)|_{\prod_{i=1}^m X-\Delta}
\end{equation}
and structure group
\begin{equation}
\bG(\Sigma) 
= 
\bG(\kappa_1,\dots,\kappa_m) 
:= 
\left(\prod_{i=1}^m (\SO(3)\times\SO(4))\right)
\rtimes \fS(\kappa_1,\dots,\kappa_m). 
\end{equation}
Let $S^4$ have its standard round metric of radius one, let
$E_1,\dots,E_m$ be Hermitian, rank-two vector bundles over $S^4$ with
$c_2(E_j) = \kappa_j$ for $j=1,\dots,m$, and let $\fg_{E_j} \equiv
\su(E_j)$ be the corresponding $\SO(3)$ bundles with $-\quarter
p_1(\fg_{E_j}) = \kappa_j$.  Denote $\RR^+ := (0,\8)$. The smooth
manifold
\begin{equation}
\label{eq:ConeBundleFiber}
\sZ(\Sigma) 
= 
\sZ(\kappa_1,\dots,\kappa_m) 
:=
\prod_{i=1}^m \left(\tM_{E_i}^{\diamond}\times_{\sG_{E_i}}
\Fr(\fg_{E_i})|_s\times\RR^+\right)
\end{equation}
carries a natural action of $\bG(\kappa_1,\dots,\kappa_m)$ coming from
the action of $\SO(4)$ rotating $S^4$ and the action of $\SO(3)$ on
the choice of frame for $\fg_{E_i}|_s$ in each factor, and the action
of $\fS(\kappa_1,\dots,\kappa_m)$ permuting the factors whose
underlying bundles are isomorphic.  The bundle $\Fr(TX)$ has its usual
right $\SO(4)$ action given by
\begin{equation}
\label{eq:TangentFrameBundleAction}
\SO(4)\times\Fr(TX)\to \Fr(TX), \quad (u,f) \mapsto fu^{-1}.
\end{equation}
We form the fibered product $\Gl_\kappa^w(X,\sU_0,\Sigma) =
\Gl_E(X,\sU_0,\Sigma)$, where
\begin{equation}
\label{eq:FiberedProduct}
\Gl_E(X,\sU_0,\Sigma) 
:= 
\bFr_{E_0}(X,\sU_0,\Sigma)\times_{\bG(\Sigma)} \sZ(\Sigma),
\end{equation}
together with the natural projection map
\begin{equation}
\label{eq:ConeBundleOverStratum}
\pi:\Gl_E(X,\sU_0,\Sigma) \to \sU_0\times \Sigma.
\end{equation}
This is a locally trivial fiber bundle with structure group
$\bG(\kappa_1,\dots,\kappa_m)$. We let
$\sZ(\Sigma,\lambda_0)\subset\sZ(\Sigma)$ denote the open subspace
obtained by replacing $\RR^+$ by $(0,\lambda_0)$ in the definition
\eqref{eq:ConeBundleFiber} and let $\Gl_\kappa^w(X,\sU_0,\Sigma,\lambda_0) =
\Gl_E(X,\sU_0,\Sigma,\lambda_0)$ denote the corresponding gluing data bundle.

\subsubsection{Compactification of the gluing data}
\label{subsubsec:CompactGluingData}
The gluing-data bundle $\Gl_{E_0}(X,\sU_0,\Si)$ has two
sources of non-compactness, the non-compact fiber
$Z(\Si)$ and the non-compact base, $\Si$.  We now
describe a compactification of $\Gl_{E_0}(X,\sU_0,\Si)$. It is a trivial
matter to replace the precompact neighborhood $\sU_0\subset M_{E_0}(X,C)$
by its compactification $\bar\sU_0$, so we assume that $\sU_0$ is itself
compact. 

We first describe the closure of the fibers.
The definition of $\sZ(\Si)$ uses the space
$M^{\diamond,s}_{E_i}(S^4)$ and an dilation parameter
lying in $\RR^+=(0,\infty)$.  The dilation action
\eqref{eq:CenterScaleDiffeo} maps points in
$M^{\diamond,s}_{E_i}\times (0,\infty)$ to framed connections on
$\fg_{E_i}$.  This action extends naturally to
$\barM^{\diamond,s}_{E_i}(S^4)\times (0,\infty)$.  The Uhlenbeck limit of
these connections as the dilation parameter goes to zero is given by
$[\Theta,\bx]$ where $\Theta$ is a trivial connection and
$\bx\in\Sym^\kappa(S^4)$ is represented by $(n,n,\dots,n)$, where $n$ is the
north pole.  Thus, we define a ccompletion of
$$
\barM^{\diamond,s}_{E_i}(S^4)\times (0,\infty),
$$
by the cone $c(\barM^{\diamond,s}_{E_i}(S^4))$, 
where the cone point corresponds to zero dilation parameter. Then, let
$$
\bar{\sZ}(\Sigma)=\bar{\sZ}(\ka_1,\dots,\ka_m)
:=
\prod_{i=1}^m c(\barM^{\diamond,s}_{E_i}(S^4)).
$$
The fibers of the gluing
data bundle are then compactified if we work with the space:
\begin{equation}
\label{eq:GluingDataClosure}
\bar{\Gl}_E(X,\sU_0,\Sigma)
:= 
\bFr_{E_0}(X,\sU_0,\Sigma)\times_{\bG(\Sigma)} \bar{\sZ}(\Sigma).
\end{equation}
As the only difference between this space 
and $\Gl_E(X,\sU_0,\Sigma)$ is given by the change in
the fiber from $\sZ(\Si)$ to $\bar{\sZ}(\Si)$,
the projection \eqref{eq:ConeBundleOverStratum} extends to
$$
\pi: \bar{\Gl}_E(X,\sU_0,\Si)
\to \sU_0\times\Si,
$$
and this extension is still a locally trivial fiber bundle
with the same structure group.
We again let $\bar{\Gl}_E(X,\sU_0,\Si,\la_0)$
be the open subspace given by requiring all dilation
parameters to be less than $\la_0$.

We deal with the non-compact base, $\Si$, by replacing
$\bFr_{E_0}(X,\sU_0,\Si)$ with
\begin{equation}
\label{eq:ExtendedBundle}
\bFr_{E_0}(X,\sU_0,\bar\Si):=
\tilde\sU_0\times_{\sG_{E_0}}\prod_{i=1}^m 
\left(\Fr(\fg_{E_0})\times_X\Fr(TX)\right).
\end{equation}
The quotient $X^m/\fS_\Si$ can be identified with
$\bar\Si$, the closure of $\Si$ in $\Sym^\ell(X)$,
and so there is a projection
$$
\bFr_{E_0}(X,\sU_0,\bar\Si)
\to 
\sU_0\times\bar\Si.
$$
However, because the action of $\fS_\Si$ is not free
over $\bar\Si\less\Si$, this extension is not a fiber
bundle.  

The splicing and gluing maps are defined on a subspace of the following
extension of the gluing data bundle,
\begin{equation}
\label{eq:ExtendedGluingDataBundle}
\bar{\Gl}_E(X,\sU_0,\bar\Sigma)
:= 
\bFr_{E_0}(X,\sU_0,\bar\Sigma)\times_{\bG(\Sigma)} \bar{\sZ}(\Sigma),
\end{equation}
as we shall see in \S \ref{subsubsec:ExtensionSplicingGluing}.

\subsection{Splicing connections over $X$} 
\label{subsec:SplicingConstruction}
In this section we describe the {\em splicing construction\/} for
$\SO(3)$ connections on $\fg_E$ via a partition of unity on $X$ and
define the {\em splicing maps} $\bga'$ for connections.  Our
description is largely motivated by those of Taubes in 
\cite{TauSelfDual}, \cite{TauPath}, \cite{TauIndef}, \cite[\S
4]{TauFrame}, \cite{TauStable}. We keep the notation of \S
\ref{subsec:GluingData}.

\subsubsection{Cutting and splicing connections}
\label{subsubsec:CutAndSpliceConnections}
Fix a smooth stratum $\Sigma\subset\Sym^\ell(X)$, let $\bx\in\Sigma$, and
choose $r_0=r_0(\bx)>0$ to be one half the smaller of
\begin{itemize}
\item The minimum geodesic distance between distinct points $x_i$,
  $x_j$ of the representative $(x_1,\dots,x_m)$ of $\bx$, and
\item 
The injectivity radius of the Riemannian manifold $(X,g)$.
\end{itemize}
A choice of $L^2_k$ connection $A_0$ on $\fg_{E_0}$ and $\SO(3)$ frames $p_i\in
\Fr(\fg_{E_0})|_{x_i}$ define $L^2_{k+1}$ local sections of $\Fr(\fg_{E_0})$,
$$
\sigma_i = \sigma_i(A_0,p_i)
:B(x_i,r_0)\to \Fr(\fg_{E_0}), \quad i=1,\dots,m, 
$$
using parallel transport along radial geodesics emanating from $x_i$
together with smoothing for sections of $\fg_{E_0}$ using the heat kernel
$\exp(-td_{A_0}^*d_{A_0})$ for small $t$, as described in \cite[\S
A.1]{FL1}; see also the remarks in \cite[p. 177]{TauStable}. Let
\begin{equation}
\fg_{E_0}|_{B(x_i,r_0)}\simeq B(x_i,r_0)\times \SO(3)
\label{eq:BackroundBundleTriv}
\end{equation}
be the $L^2_{k+1}$ local trivializations defined by the sections $\sigma_i$. 

Define a smooth cutoff function $\chi_{x_0,\eps}:X\to [0,1]$ by setting
\begin{equation}
\label{eq:ChiCutoffFunctionDefn}
\chi_{x_0,\eps}(x) 
:= 
\chi(\dist(x,x_0)/\eps)
\end{equation}
where $\chi:\RR\to [0,1]$ is a smooth function such that $\chi(t)=1$
for $t\ge 1$ and $\chi(t)=0$ for $t\le 1/2$.  Thus, we have
$$
\chi_{x_0,\eps}(x) 
=
\begin{cases}
1 &\text{for } x\in X - B(x_0,\eps),
\\
0 &\text{for } x\in B(x_0,\eps/2).
\end{cases}
$$
We define a cut-off connection on the bundle $\fg_{E_0}$ over $X$ by
setting
\begin{equation}
\label{eq:CutOffBackgroundConnection}
\chi_{\bx,4\sqrt{\blambda}}A_0
:=
\begin{cases}
A_0 
&\text{over $X - \cup_{i=1}^m B(x_i,4\sqrt{\lambda_i})$}, 
\\
\Gamma + \chi_{x_i,4\sqrt{\lambda_i}}\sigma_i^*A_0
&\text{over $\Omega(x_i,2\sqrt{\lambda_i},4\sqrt{\lambda_i})$,} 
\\
\Gamma 
&\text{over $B(x_i,2\sqrt{\lambda_i})$,} 
\end{cases} 
\end{equation}
where $\Gamma$ denotes the product connection on $B(x_i,r_0)\times \SO(3)$,
while $\bx = (x_1,\dots,x_m)$ and $\blambda =
(\lambda_1,\dots,\lambda_m)$. 

Let $A_i$ be centered $L^2_k$ connections on the bundles $\fg_{E_i}$
over $S^4$. Let 
$$
\tau_i
=
\tau_i(A_i,q_i)
:S^4\less\{n\} \to \Fr(\fg_{E_i}), \quad i=1,\dots,m, 
$$
be the analogously-defined $L^2_{k+1}$ local sections induced by the
connections $A_i \in \sA_{E_i}$ and $\SO(3)$ frames $q_i \in
\Fr(\fg_{E_i})|_s$. Let
\begin{equation}
\fg_{E_i}|_{S^4\less\{n\}}\simeq S^4\less\{n\}\times \SO(3)
\label{eq:SphereBundleTriv}
\end{equation}
be the $L^2_{k+1}$ local trivializations defined by the $\tau_i$. Let
$\delta_\lambda:S^4\to S^4$ be the conformal diffeomorphism of $S^4$
defined by
\begin{equation}
\delta_\lambda := \varphi_n\circ\lambda\circ \varphi_n^{-1},
\end{equation}
where $\lambda:\RR^4\to\RR^4$ is the dilation given by $y\mapsto
y/\lambda$, so $\delta_\lambda$ is a conformal diffeomorphism of $S^4$
which fixes the north and south poles.  We define cut-off, rescaled
connections on the bundles $\fg_{E_i}$ over $S^4$, with mass centers at the
north pole, by setting
\begin{equation}
(1-\chi_{x_i,\sqrt{\lambda_i}/2})\delta_{\lambda_i}^*A_i
:=
\begin{cases}
\delta_{\lambda_i}^*A_i
&\text{over $\varphi_n(B(0,\quarter\sqrt{\lambda_i}))$}, 
\\
\Gamma + (1-\chi_{x_i,\sqrt{\lambda_i}/2})\tau_i^*\delta_{\lambda_i}^*A_i
&\text{over $\varphi_n(\Omega(0,\quarter\sqrt{\lambda_i},
\half\sqrt{\lambda_i}))$,} 
\\
\Gamma 
&\text{over $S^4 - \varphi_n(B(0,\half\sqrt{\lambda_i}))$,} 
\end{cases} 
\end{equation}
where $\Gamma$ denotes the product connection on $(S^4-\{n\})\times \SO(3)$.
Note that $\tau_i^*\delta_{\lambda_i}^*A_i
= \delta_{\lambda_i}^*\tau_i^*A_i$, since the section $\tau_i$ is
defined by parallel translation from the south pole via the connection $A_i$.

A choice of $\SO(4)$ frames $f_i\in\Fr(TX)|_{x_i}$ defines coordinate charts
$$
\varphi_i^{-1}:B(x_i,r_0)\subset X\to\RR^4, \quad i=1,\dots,m, 
$$ 
via the exponential maps $\exp_{f_i}:B(0,r_0)\subset TX|_{x_i}\to X$.  
The orientation-preserving diffeomorphism $\varphi_i\circ \varphi_n^{-1}$
identifies the annulus
$\varphi_n(\Omega(0,\half\sqrt{\lambda_i},2\sqrt{\lambda_i}))$ in $S^4$, 
$$
\Omega\left(0,\half\sqrt{\lambda_i},2\sqrt{\lambda_i}\right)
:=
\varphi_n
\left(\left\{x\in\RR^4:\half\sqrt{\lambda_i}<|x|<2\sqrt{\lambda_i}
\right\}\right)
\subset \RR^4
$$
with the annulus in $X$,
$$
\Omega\left(x_i,\half\sqrt{\lambda_i},2\sqrt{\lambda_i}\right)
:=
\left\{x\in X: \half\sqrt{\lambda_i} < \dist_g(x,x_i) < 2\sqrt{\lambda_i} 
\right\}
\subset X.
$$
We define a glued-up $\SO(3)$ bundle $\fg_E$ over $X$ by setting
\begin{equation}
\fg_E 
= 
\begin{cases}
\fg_{E_0} &\text{over $X - \cup_{i=1}^m B(x_i,\half\sqrt{\lambda_i})$,} \\
\fg_{E_i} &\text{over $B(x_i,2\sqrt{\lambda_i})$.} 
\end{cases}
\label{eq:DefnGluedUpSO(3)Bundle}
\end{equation}
The bundles $\fg_{E_0}$ and $\fg_{E_i}$ are identified over the annuli
$\Omega(x_i,\half\sqrt{\lambda_i},2\sqrt{\lambda_i}\})$ in $X$ via the
isomorphisms of $\SO(3)$ bundles defined by the 
orientation-preserving diffeomorphisms $\varphi_i\circ \varphi_n^{-1}$,
identifying the annuli
$\Omega(x_i,\half\sqrt{\lambda_i},2\sqrt{\lambda_i}\})$ with the
corresponding annuli in $S^4$ and the $\SO(3)$ bundle maps defined by the
trivializations \eqref{eq:BackroundBundleTriv} and
\eqref{eq:SphereBundleTriv}.   
Define a spliced $L^2_k$ connection $A$ on $\fg_E$ by setting
\begin{equation}
\label{eq:SplicedConnection}
A
:= 
\begin{cases}
  A_0 &\text{over $X - \cup_{i=1}^m B(x_i,4\sqrt{\lambda_i})$},
  \\
  \Gamma +
  \chi_{x_i,4\sqrt{\lambda_i}}A_0 +
  (1-\chi_{x_i,\sqrt{\lambda_i}/2})\delta_{\lambda_i}^*A_i &\text{over
    $\Omega(x_i,\quarter\sqrt{\lambda_i},4\sqrt{\lambda_i})$},
  \\
  \delta_{\lambda_i}^*A_i &\text{over $B(x_i,\quarter\sqrt{\lambda_i})$,
    $i=1,\dots,m$,}
\end{cases} 
\end{equation}
where the cut-off connections are defined as above; the bundle and
annulus identifications are understood. Also implicit in definition
\eqref{eq:SplicedConnection} are the following constraints on the
scales, given an $m$-tuple of distinct points, $(x_1,\dots,x_m)$,
\begin{equation}
\label{eq:ConstraintScales}
8\sqrt{\lambda_i} + 8\sqrt{\lambda_j} < \dist_g(x_i,x_j),
\quad i\neq j,
\end{equation}
representing a multiset $\bx\in\Sigma$.  The {\em splicing map\/},
$\bga'$, is defined on the open subset of the gluing-data bundle,
$\Gl^w_\kappa(X,\sU_0,\Sigma)$ in definition
\eqref{eq:FiberedProduct}, defined by
\begin{equation}
\label{eq:ConstrainedGluingDataBundle}
\Gl^{w,+}_\kappa(X,\sU_0,\Sigma)
:=
\{\pi^{-1}([A_0],\bx)\in \Gl^w_\kappa(X,\sU_0,\Sigma):
\text{$(\blambda,\bx)$ obeys \eqref{eq:ConstraintScales}} \}. 
\end{equation} 
If we replace $\Sigma$ in
$\Gl^w_\kappa(X,\sU_0,\Sigma)$ by a compact subset $K\Subset \Sigma$, so
that there is a uniform positive lower bound for $\dist_g(x_i,x_j)$,
$i\neq j$, for all $\bx \in K$, then we can plainly choose a positive
constant $\lambda_0 = \lambda_0(K,g)$ such that
$$
\Gl^w_\kappa(X,\sU_0,K,\lambda_0)
\subset
\Gl^{w,+}_\kappa(X,\sU_0,\Sigma),
$$
as the $m$-tuples of scales $\blambda = (\lambda_1,\dots,\lambda_m)$
present in $\Gl^w_\kappa(X,\sU_0,K,\lambda_0)$ satisfy the constraints
\eqref{eq:ConstraintScales}.

\subsubsection{Splicing maps}
\label{subsubsec:RefinedSplicingMap}
It remains to summarize what the splicing construction has achieved.
The construction, thus far, yields a splicing map giving a smooth
embedding
\begin{equation}
\label{eq:SplicingMap}
\bga':\Gl^{w,+}_\kappa(X,\sU_0,\Sigma)\to \sB_\kappa^{w,*}(X)\times\Met(X)
\end{equation}
The fact that the map is a smooth embedding is a fairly straightforward
consequence of the definitions; it will follow from the stronger result in
\cite{FLKM2} that the gluing map $\bga$ gives a smooth embedding of the
gluing data. Similarly, the fact that the image of $\bga'$ contains only
irreducible connections is also a fairly easy consequence of the
definitions; it is proved in \cite{FLKM2}.

\subsubsection{Extension of the splicing and gluing maps}
\label{subsubsec:ExtensionSplicingGluing}
It remains to define an extension of the 
splicing map $\bga'$ and gluing map $\bga$ from
$\Gl_\kappa^{w,+}(X,\sU_0,\Sigma,\lambda_0)$ to
$\bar{\Gl}_\kappa^{w,+}(X,\sU_0,\bar\Sigma,\lambda_0)$. Thus, given a point
$\bp$ in $\bar{\Gl}_\kappa^{w,+}(X,\sU_0,\bar\Sigma,\lambda_0)$, we obtain
a spliced, ideal connection $\bga'(\bp) =
[A,\bx]\in\sB_{\kappa-\ell}^w(X)\times\Sym^\ell(X)$ by the following
procedure:
\begin{itemize}
\item
The connection $A$ is obtained by splicing the background connection $A_0$
and the connections $A_i$ on $S^4$ for which $\lambda_i>0$, exactly as in
\S \ref{subsubsec:RefinedSplicingMap}.
\item
The multiset $\bx\in\Sym^\ell(X)$ is the union of 
\begin{itemize}
\item
All points $x_i\in
\pi_2(\bp)$, where $\pi_2$ is the projection from  
$\bar{\Gl}_\kappa^{w,+}(X,\sU_0,\bar\Sigma,\lambda_0)$ onto the second
factor, $\bar\Sigma$, for which the corresponding scale $\lambda_i$ is zero. 
\item
All points $\exp_{f_i}(\lambda_i\varphi_n^{-1}(y))$ such that
$\lambda_i\neq 0$, $y\in\by_i$, and $[A_i,\by_i,q_i]\in
\barM_{\kappa_i}^{\diamond,s}(S^4)$, where $\by_i$ is non-empty.
\end{itemize}
\end{itemize}
The requisite domain of the splicing and gluing maps
is given by $\bar{\Gl}_E^+(X,\sU_0,\bar\Sigma,\lambda_0)$,
where the $+$ indicates that the pairs $(\bx,\blambda)$ must satisfy the
constraints \eqref{eq:ConstraintScales} (with equality now allowed),
essentially mimicking the definition \eqref{eq:ConstrainedGluingDataBundle} of
$\Gl_E^+(X,\sU_0,\Sigma,\lambda_0)$. 

The splicing map \eqref{eq:SplicingMap} has an obvious extension to
a map
\begin{equation}
\label{eq:SplicingMapExtension}
\bga':\barGl^{w,+}_\kappa(X,\sU_0,\Sigma)\to 
\barsB_\kappa^{w,*}(X)\times\Met(X),
\end{equation}
and similarly for the gluing map $\bga$, which we define in \S
\ref{subsec:ProofGluingTheoremWoUhlenbeckContinuity}.  One of the principal
assertions of our main result, Theorem
\ref{thm:GluingTheorem}, is that our choice of extensions of these maps is
the right one, as the resulting extensions are shown to be continuous.


\section{Regularity of solutions to quasi-linear systems}
\label{sec:Regularity}
Let $A$ be a $C^\8$, orthogonal connection on the $\SO(3)$ bundle $\fg_E =
\su(E)$, where $E$ is a Hermitian, rank-two vector bundle over $X$. 
Based on the approach we gave in \cite{FL1} for the $\PU(2)$ monopole
equations, we give the regularity theory --- which we shall need both in
the present article and its sequels \cite{FLKM2}, \cite{FLKM3},
\cite{FLKM4} --- for solutions $a\in L^2_1(\Lambda^1\otimes\fg_E)$ to the
quasi-linear elliptic system consisting of a quasi-linear equation, with a
quadratic non-linearity, and the inhomogeneous Coulomb-gauge equation,
\begin{equation}
\label{eq:QuasiFirstASDCoulomb}
\begin{aligned}
d_A^+a + (a\wedge a)^+ &= w,
\\
d_A^*a &= \zeta,
\end{aligned}
\end{equation}
where $(\zeta,w)\in L^p_k(\fg_E) \oplus
L^p_k(\Lambda^+\otimes\fg_E)$. We use the isomorphism
$\ad:\fg_E\to\so(\fg_E)$ to view sections of $\so(\fg_E)$ as sections of
$\fg_E$. A systematic treatment is given in \cite[\S 3]{FL1} and we merely
summarize the main results here; some special cases appear in \cite{DK} and
\cite{FU}.  As in
\cite{TauSelfDual}, \cite{TauStable}, we use $a=d_A^*v$, where $v\in
L^2_2(X,\Lambda^+\otimes\fg_E)$, to rewrite the first-order quasi-linear
equation for $a$ in terms of $w$ in
\eqref{eq:QuasiFirstASDCoulomb} as
\begin{equation}
\label{eq:QuasiSecondASDCoulomb}
d_A^+d_A^*v + (d_A^*v\wedge d_A^*v)^+ = w,
\end{equation}
to give a second-order, quasi-linear elliptic equation for $v$ in terms of
$w$; the Coulomb-gauge equation reduces to $d_A^*d_A^*v = (F_A^+)^*v$.

Though our applications of the regularity results of this section will
often be to solutions $a$ to the anti-self-equation, $F^+(A+a)=0$, namely
equation \eqref{eq:QuasiFirstASDCoulomb} with $w=-F_A^+$, we shall also
need to consider the regularity theory for more general equations, such as
the various forms of the {\em extended anti-self-dual equation\/} (and
their solutions, defining {\em extended anti-self-dual connections\/})
introduced in \cite{DonPoly}, \cite{FrM}, \cite{TauIndef}, \cite{TauStable}
--- see equations \eqref{eq:QuickExtASDEqnForv},
\eqref{eq:LongSecOrderExtASDEqnForv}, and
\eqref{eq:RegLongSecOrderExtASDEqnForv}, for example. These more
general equations can be viewed as special cases of
\eqref{eq:QuasiFirstASDCoulomb} or \eqref{eq:QuasiSecondASDCoulomb}
obtained by a suitable choice of $w$, such as that given in equation
\eqref{eq:DefnOfwForExtASD}. Therefore, by considering the slightly more
general equations \eqref{eq:QuasiFirstASDCoulomb} and
\eqref{eq:QuasiSecondASDCoulomb}, we can address the regularity theory
simultaneously for all such equations.

\begin{rmk}
The precise dependence of the constants $C$ on the connection $A$ in this
section is not something we track carefully, as we are primarily interested
in local $L^2_k$ estimates of solutions $a$ to equations
\eqref{eq:QuasiFirstASDCoulomb} or of solutions $v$ to equations
\eqref{eq:QuasiSecondASDCoulomb} when $k\geq 2$ or $k\geq 3$. For this
range of the Sobolev indices, the constants no longer depend on simply
$\|F_A\|_{L^2(X)}$ or $\|F_A^+\|_{L^\8(X)}$, for example. Rather, they
depend on higher-order covariant derivatives of the connection $A$. In 
general, even when $A$ is not anti-self-dual, we can always bound
$L^2_{k+1}$ norms of local connection one-forms in Coulomb gauge by
$L^2_{k,A}$ norms of the curvature, following the method of proof for
\cite[Lemma 2.3.11]{DK}. If $A$ is anti-self-dual, the dependence on $A$
reduces to the $L^2$ norm of $F_A$. Thus, throughout this section,
constants depending on $A$ can be gauge-invariantly given as constants
depending on $L^2_{k,A}$ norms of $F_A$. 
\end{rmk} 

\subsection{Regularity for $L^2_1$ solutions to the inhomogeneous
Coulomb-gauge and quasi-linear equations}
\label{subsec:L2_1InhomoGlobal}
We apply the more general results of \cite[\S 3]{FL1} --- concerning
regularity theory for the $\PU(2)$ monopole equations --- to systems
\eqref{eq:QuasiFirstASDCoulomb} and \eqref{eq:QuasiSecondASDCoulomb}.
Briefly, an $L^2_1$ solution $a$ to the anti-self-dual and Coulomb-gauge
equations \eqref{eq:QuasiFirstASDCoulomb}, with an $L^2_k$ inhomogeneous
term (with $k \ge 1$) is in $L^2_{k+1}$. Thus, if the inhomogeneous term is
in $C^\8$ then $a$ is in $C^\8$.

Fix an exponent $1\le p\le\8$ and an integer $k\ge 0$ such that
$(k+1)p>4$. Let $A$ be an orthogonal $L^p_k$ connection on $\fg_E$.
For any $a$ in $C^\8(\Lambda^1\otimes\fg_E)$, with $1\leq q\leq\8$
an exponent and $\ell\geq$ an integer such that
$L^p_k\subset L^q_\ell$ (see \cite[Theorem 5.4]{Adams}), set
$$
\|a\|_{L^q_{\ell,A}(X,g)} 
:= 
\left(\sum_{i=0}^\ell\|\cov_A^i a\|_{L^q(X,g)}^q\right)^{1/q},
$$
and let $L^p_k(\Lambda^1\otimes\fg_E)$ be the completion of
$\Om^1(\fg_E)$ with respect to this norm. Define the Sobolev norms and
Banach-space completions of the spaces $\Omega^0(\fg_E)$ and
$\Omega^+(\fg_E)$ in the analogous manner.

\begin{prop}
\label{prop:L2_1InhomoReg}
(see \cite[Proposition 3.2]{FL1})
Let $X$ be a closed, oriented, Riemannian four-manifold with metric $g$
and let $E$ be a Hermitian, rank-two bundle
over $X$.  Let $A$ be a $C^\8$ connection on the $C^\8$ bundle
$\fg_E$ over $X$ and let $2\le p < 4$.  Then there are
positive constants $\eps=\eps(A,p)$ and $C=C(A,p)$ with
the following significance. Suppose that $a\in
L^2_1(X,\La^1\otimes\fg_E)$ is a
solution to the elliptic system
\eqref{eq:QuasiFirstASDCoulomb} over $X$, where
$(\zeta,w)$ is in $L^p$.  If $\|a\|_{L^4(X)}<\eps$ then
$a$ is in $L^p_1$ and
$$
\|a\|_{L^p_{1,A}(X)} \le C\left(\|(\zeta,w)\|_{L^p(X)}
+ \|a\|_{L^2(X)}\right).
$$
\end{prop}

\begin{prop}
\label{prop:Lp_1InhomoReg}
(see \cite[Proposition 3.3]{FL1})
Continue the notation of Proposition \ref{prop:L2_1InhomoReg}.
Let $k\ge 1$ be an integer and let $2<p<\8$. Let
$A$ be a $C^\8$ connection on the $C^\8$ bundle $\fg_E$
over $X$. Suppose that $a\in L^p_1(X,\La^1\otimes\fg_E)$ is a solution
to the elliptic system 
\eqref{eq:QuasiFirstASDCoulomb} over $X$, where $(\zeta,w)$ is in $L^2_k$. Then
$a$ is in $L^2_{k+1}$ and there is a universal polynomial $Q_k(x,y)$, with
positive real coefficients, depending at most on $A,k$, such that
$Q_k(0,0)=0$ and
$$
\|a\|_{L^2_{k+1,A}(X)} 
\le Q_k\left(\|(\zeta,w)\|_{L^2_{k,A}(X)},
\|a\|_{L^p_{1,A}(X)}\right).
$$
In particular, if $(\zeta,w)$ is in $C^\8$ then $a$ is in
$C^\8$ and if $(\zeta,w)=0$, then
$$
\|a\|_{L^2_{k+1,A}(X)} \le C\|a\|_{L^p_{1,A}(X)}.
$$
\end{prop}

By combining Propositions \ref{prop:L2_1InhomoReg} and
\ref{prop:Lp_1InhomoReg} we obtain the desired regularity result for
$L^2_1$ solutions to the inhomogeneous Coulomb-gauge and anti-self-dual
equations:

\begin{cor}
\label{cor:L2_1InhomoReg}
(see \cite[Corollary 3.4]{FL1})
Continue the notation of Proposition \ref{prop:L2_1InhomoReg}.
Let $A$ be a $C^\8$ connection on the $C^\8$
bundle $\fg_E$ over $X$.  Then there is a positive constant
$\eps=\eps(A)$ such that the following hold. 
Suppose that $a\in 
L^2_1(X,\La^1\otimes\fg_E)$, is a solution to the elliptic system
\eqref{eq:QuasiFirstASDCoulomb} over $X$, where $(\zeta,w)$ is in $L^2_k$ and 
$\|a\|_{L^4(X)}<\eps$ and $k\ge 0$ is an integer.
Then $a$ is in $L^2_{k+1}$ and 
there is a universal polynomial $Q_k(x,y)$, with positive real coefficients,
depending at most on $A,k$, such that $Q_k(0,0)=0$ and
$$
\|a\|_{L^2_{k+1,A}(X)} 
\le Q_k\left(\|(\zeta,w)\|_{L^2_{k,A}(X)},
\|a\|_{L^2(X)}\right).
$$
In particular, if $(\zeta,w)$ is in $C^\8$ then $a$ is in
$C^\8$ and if $(\zeta,w)=0$, then
$$
\|a\|_{L^2_{k+1,A}(X)} \le C\|a\|_{L^2(X)}.
$$
\end{cor}

\subsection{Regularity for $L^2_2$ solutions to the inhomogeneous,
second-order quasi-linear equation}
\label{subsec:L2_2SecondOrderInhomoGlobal}
It remains to consider the regularity properties of the solution $v$
to the second-order equation \eqref{eq:QuasiSecondASDCoulomb}, where
$a=d_A^{+,*}v$ and thus
$$
d_A^+d_A^{+,*}v = d_A^+a.
$$
The Laplacian $d_A^+d_A^{+,*}$ is elliptic, with $C^\8$ coefficients, and
thus standard regularity theory implies that if $a\in
L^p_k(X,\La^1\otimes\fg_E)$ then $d_A^+a\in
L^p_{k-1}(X,\La^1\otimes\fg_E)$ and so $v\in
L^p_{k+1}(X,\La^1\otimes\fg_E)$. 

It will prove convenient to consider the following generalization of the
second-order system \eqref{eq:QuasiSecondASDCoulomb}:
\begin{equation}
\label{eq:GeneralSecondOrderASD}
d_A^+d_A^*v + (d_A^*v\wedge d_A^*v)^+ + \{\alpha,d_A^*v\} = w,
\end{equation}
where $v\in L^2_2(X,\Lambda^+\otimes\fg_E)$, $w\in
L^p_k(X,\Lambda^+\otimes\fg_E)$, and $\alpha \in
C^\8(X,\Lambda^1\otimes\fg_E)$. Here,
$$
\{\, ,\, \}:\Gamma(\Lambda^1\otimes\fg_E)\times
\Gamma(\Lambda^1\otimes\fg_E)
\to \Gamma(\Lambda^+\otimes\fg_E),
$$
denotes the linear map defined by
$$
\{\alpha,\beta\} := (\alpha\wedge\beta)^+ + (\beta\wedge\alpha)^+.
$$
Then the proofs of the regularity results of
\S \ref{subsec:L2_1InhomoGlobal} adapt with essentially no change to give
the following regularity results for 
equation \eqref{eq:GeneralSecondOrderASD}.

\begin{prop}
\label{prop:L2_2SecondOrderInhomoReg}
(see \cite[Proposition 3.2]{FL1})
Let $X$ be a closed, oriented, Riemannian four-manifold with metric $g$
and let $E$ be a Hermitian, rank-two bundle
over $X$.  Let $A$ be a $C^\8$ connection on the $C^\8$ bundle
$\fg_E$ over $X$ and let $2\le p < 4$.  Suppose $\alpha \in
C^\8(X,\Lambda^1\otimes\fg_E)$. Then there are
positive constants $\eps=\eps(A,\alpha,p)$ and $C=C(A,\alpha,p)$ with
the following significance. Suppose that $v\in
L^2_2(X,\La^+\otimes\fg_E)$ is a
solution to the elliptic system
\eqref{eq:GeneralSecondOrderASD} over $X$, where
$w$ is in $L^p(X,\La^+\otimes\fg_E)$.  
Further, if $\|d_A^*v\|_{L^4(X)}<\eps$ then $v$ is in $L^p_2$ and
$$
\|v\|_{L^p_{2,A}(X)} \le C\left(\|w\|_{L^p(X)}
+ \|v\|_{L^2(X)}\right).
$$
\end{prop}

\begin{prop}
\label{prop:Lp_2SecondOrderInhomoReg}
(see \cite[Proposition 3.3]{FL1})
Continue the notation of Proposition \ref{prop:L2_1InhomoReg}.
Let $k\ge 1$ be an integer and let $2<p<\8$. Let
$A$ be a $C^\8$ connection on the $C^\8$ bundle $\fg_E$
over $X$. Suppose $\alpha \in
C^\8(X,\Lambda^1\otimes\fg_E)$. Suppose that $v\in
L^p_2(X,\La^+\otimes\fg_E)$ is a
solution to the elliptic system
\eqref{eq:GeneralSecondOrderASD} over $X$, where $w$ is in $L^2_k$. Then
$v$ is in $L^2_{k+2}$ and there is a universal polynomial $Q_k(x,y)$, with
positive real coefficients, depending at most on $A,\alpha,k$, such that
$Q_k(0,0)=0$ and
$$
\|v\|_{L^2_{k+2,A}(X)} 
\le Q_k\left(\|w\|_{L^2_{k,A}(X)},
\|v\|_{L^p_{2,A}(X)}\right).
$$
In particular, if $w$ is in $C^\8$ then $v$ is in
$C^\8$ and if $w=0$, then
$$
\|v\|_{L^2_{k+2,A}(X)} \le C\|v\|_{L^p_{2,A}(X)}.
$$
\end{prop}

By combining Propositions \ref{prop:L2_2SecondOrderInhomoReg} and
\ref{prop:Lp_2SecondOrderInhomoReg} we obtain the desired regularity result for
$L^2_2$ solutions to the inhomogeneous general second-order anti-self-dual
equation:

\begin{cor}
\label{cor:L2_2SecondOrderInhomoReg}
(see \cite[Corollary 3.4]{FL1})
Continue the notation of Proposition \ref{prop:L2_2SecondOrderInhomoReg}.
Let $A$ be a $C^\8$ connection on the $C^\8$
bundle $\fg_E$ over $X$.  Suppose $\alpha \in
C^\8(X,\Lambda^1\otimes\fg_E)$. Then there is a positive constant
$\eps=\eps(A,\alpha)$ such that the following hold. 
Suppose that $a\in 
L^2_1(X,\La^1\otimes\fg_E)$, is a solution to the elliptic system
\eqref{eq:GeneralSecondOrderASD} over $X$, where $w$ is in $L^2_k$ and 
$\|d_A^*v\|_{L^4(X)}<\eps$ and $k\ge 0$ is an integer. 
Then $v$ is in $L^2_{k+2}$ and 
there is a universal polynomial $Q_k(x,y)$, with positive real coefficients,
depending at most on $A,\alpha,k$, such that $Q_k(0,0)=0$ and
$$
\|v\|_{L^2_{k+2,A}(X)} 
\le Q_k\left(\|w\|_{L^2_{k,A}(X)},
\|v\|_{L^2(X)}\right).
$$
In particular, if $w$ is in $C^\8$ then $v$ is in
$C^\8$ and if $w=0$, then
$$
\|v\|_{L^2_{k+2,A}(X)} \le C\|v\|_{L^2(X)}.
$$
\end{cor}

\subsection{Local regularity and interior estimates for $L^2_1$ solutions
to the inhomogeneous Coulomb-gauge and quasi-linear equations}
\label{subsec:L2_1InhomoLocal} 
In this section we specialize the results of \S 
\ref{subsec:L2_1InhomoGlobal} to the case where the reference connection is
a $C^\8$ flat connection, so $A=\Ga$ on the bundle
$\fg_E$ over an open subset $\Om\subset X$.

We continue to assume that $X$ is a closed, oriented four-manifold with
$C^\8$ Riemannian metric $g$ and Hermitian, rank-two bundle $E$ extending
that on $\Om\subset X$.  We then have the following local versions of
Propositions \ref{prop:L2_1InhomoReg} and \ref{prop:Lp_1InhomoReg} and
Corollary \ref{cor:L2_1InhomoReg}:

\begin{prop}
\label{prop:L2_1InhomoRegLocal}
(see \cite[Proposition 3.9]{FL1})
Continue the notation of the preceding paragraph. Let $\Om'\Subset\Om$ be a
precompact open subset and let $2\le p < 4$.  Then there are positive
constants $\eps=\eps(\Om,p)$ and $C=C(\Om',\Om,p)$ with the following
significance. Suppose that $a$ is an $L^2_1(\Om)$ solution to the
elliptic system \eqref{eq:QuasiFirstASDCoulomb} over $\Om$, with $A = \Ga$
and where $(\zeta,w)$ is in $L^p(\Om)$. If
$\|a\|_{L^4(\Om)}<\eps$ then $a$ is in $L^p_1(\Om')$ and
$$
\|a\|_{L^p_{1,\Ga}(\Om')} 
\le C\left(\|(\zeta,w)\|_{L^p(\Om)} + \|a\|_{L^2(\Om)}\right).
$$
\end{prop}

\begin{prop}
\label{prop:Lp_1InhomoRegLocal}
(see \cite[Proposition 3.10]{FL1})
Continue the notation of Proposition \ref{prop:L2_1InhomoRegLocal}.
Let $k\ge 1$ be an integer, and let $2<p<\8$. Suppose that
$a$ is an $L^p_1(\Om)$ solution to the elliptic system
\eqref{eq:QuasiFirstASDCoulomb} over $\Om$ with $A=\Ga$, 
where $(\zeta,w)$ is in $L^2_k(\Om)$. Then $a$ is in
$L^2_{k+1}(\Om')$ and there is a universal polynomial $Q_k(x,y)$, with
positive real coefficients, depending at most on
$k$, $\Om'$, $\Om$, such that
$Q_k(0,0)=0$ and
$$
\|a\|_{L^2_{k+1,\Ga}(\Om')} 
\le Q_k\left(\|(\zeta,w)\|_{L^2_{k,\Ga}(\Om)},
\|a\|_{L^p_{1,\Ga}(\Om)}\right).
$$
If $(\zeta,w)$ is in $C^\8(\Om)$ then $a$ is in $C^\8(\Om')$
and if $(\zeta,w)=0$, then
$$
\|a\|_{L^2_{k+1,\Ga}(\Om')} 
\le C\|a\|_{L^p_{1,\Ga}(\Om)}.
$$
\end{prop}

\begin{cor}
\label{cor:L2_1InhomoRegLocal}
(see \cite[Corollary 3.11]{FL1})
Continue the notation of Proposition \ref{prop:Lp_1InhomoRegLocal}.  Then
there is a positive constant $\eps=\eps(\Om)$ with the following
significance. Suppose that $a$ is an $L^2_1(\Om)$ solution to the elliptic
system \eqref{eq:QuasiFirstASDCoulomb} over $\Om$ with
$A=\Ga$, where $(\zeta,w)$ is in $L^2_k(\Om)$ and
$\|a\|_{L^4(\Om)}<\eps$.  Then $a$ is in $L^2_{k+1}(\Om')$ and there is a
universal polynomial $Q_k(x,y)$, with positive real coefficients, depending
at most on $k$, $\Om'$, $\Om$, such that $Q_k(0,0)=0$ and
$$
\|a\|_{L^2_{k+1,\Ga}(\Om')} 
\le Q_k\left(\|(\zeta,w)\|_{L^2_{k,\Ga}(\Om)},
\|a\|_{L^2(\Om)}\right).
$$
If $(\zeta,w)$ is in 
$C^\8(\Om)$ then $a$ is in $C^\8(\Om')$ and if
$(\zeta,w)=0$, then
$$
\|a\|_{L^2_{k+1,\Ga}(\Om')} \le C\|a\|_{L^2(\Om)}.
$$
\end{cor}

Corollary \ref{cor:L2_1InhomoRegLocal} thus yields a sharp local
elliptic regularity result for $L^2_1$ anti-self-dual connections $A$
which are given to us in Coulomb gauge relative to $\Ga$.

\begin{prop}
\label{prop:L2_1CoulMonoRegLocal}
(see \cite[Proposition 3.12]{FL1})
Continue the notation of Corollary \ref{cor:L2_1InhomoRegLocal}.  Then
there is a positive constant $\eps=\eps(\Om)$ and, if $k\ge 1$ is an
integer, there is a positive constant $C=C(\Om',\Om,k)$ with the following
significance. Suppose that $A$ is an $L^2_1$ solution to the
$g$-anti-self-dual equation over $\Om$, which is in Coulomb gauge over
$\Om$ relative to $\Ga$, so $d_\Ga^*(A-\Ga) = 0$, and obeys
$\|A-\Ga\|_{L^4(\Om)}<\eps$.  Then $A-\Ga$ is in $C^\8(\Om')$ and
for any $k\ge 1$,
$$
\|A-\Ga\|_{L^2_{k,\Ga}(\Om')} 
\le C\|A-\Ga\|_{L^2(\Om)}.
$$
\end{prop}

\subsection{Local regularity and interior estimates for $L^2_2$ solutions
to the general, inhomogeneous, second-order quasi-linear equation}
\label{subsec:L2_2SecondOrderInhomoLocal} 
In this section we specialize the results of \S 
\ref{subsec:L2_2SecondOrderInhomoGlobal} to the case of an open subset
$\Om\subset X$. 

We continue to assume that $X$ is a closed, oriented four-manifold with
$C^\8$ Riemannian metric $g$ and Hermitian, rank-two bundle $E$ extending
that on $\Om\subset X$.  We then have the following local versions of
Propositions \ref{prop:L2_2SecondOrderInhomoReg} and
\ref{prop:Lp_2SecondOrderInhomoReg} and Corollary
\ref{cor:L2_2SecondOrderInhomoReg}:

\begin{prop}
\label{prop:L2_2SecondOrderInhomoRegLocal}
(see \cite[Proposition 3.9]{FL1}) Continue the notation of the preceding
paragraph.  Suppose $\alpha \in C^\8(X,\Lambda^1\otimes\fg_E)$. Let
$\Om'\Subset\Om$ be a precompact open subset and let $2\le p < 4$.  Then
there are positive constants $\eps=\eps(A,\alpha,\Om,p)$ and
$C=C(A,\alpha,\Om',\Om,p)$ with the following significance. Suppose that
$v$ is an $L^2_2(\Om)$ solution to the elliptic system
\eqref{eq:GeneralSecondOrderASD} over $\Om$,  where
$w$ is in $L^p(\Om)$. If $\|d_A^*v\|_{L^4(\Om)}<\eps$ then $v$ is in
$L^p_2(\Om')$ and
$$
\|v\|_{L^p_{2,A}(\Om')} 
\le C\left(\|w\|_{L^p(\Om)} + \|v\|_{L^2(\Om)}\right).
$$
\end{prop}

\begin{prop}
\label{prop:Lp_2SecondOrderInhomoRegLocal}
(see \cite[Proposition 3.10]{FL1})
Continue the notation of Proposition \ref{prop:L2_2SecondOrderInhomoRegLocal}.
Let $k\ge 1$ be an integer, and let $2<p<\8$. Suppose that
$v$ is an $L^p_2(\Om)$ solution to the elliptic system
\eqref{eq:GeneralSecondOrderASD} over $\Om$, 
where $w$ is in $L^2_k(\Om)$. Then $v$ is in
$L^2_{k+2}(\Om')$ and there is a universal polynomial $Q_k(x,y)$, with
positive real coefficients, depending at most on
$A$, $\alpha$, $k$, $\Om'$, $\Om$, such that
$Q_k(0,0)=0$ and
$$
\|v\|_{L^2_{k+2,A}(\Om')} 
\le Q_k\left(\|w\|_{L^2_{k,A}(\Om)},
\|v\|_{L^p_{2,A}(\Om)}\right).
$$
If $w$ is in $C^\8(\Om)$ then $v$ is in $C^\8(\Om')$
and if $w=0$, then
$$
\|v\|_{L^2_{k+2,A}(\Om')} 
\le C\|v\|_{L^p_{2,A}(\Om)}.
$$
\end{prop}

\begin{cor}
\label{cor:L2_2SecondOrderInhomoRegLocal}
(see \cite[Corollary 3.11]{FL1}) Continue the notation of Proposition
\ref{prop:Lp_2SecondOrderInhomoRegLocal}.  Then there is a positive
constant $\eps=\eps(A,\alpha,\Om)$ with the following significance. Suppose
that $v$ is an $L^2_2(\Om)$ solution to the elliptic system
\eqref{eq:GeneralSecondOrderASD} over $\Om$, where $w$
is in $L^2_k(\Om)$ and $\|d_A^*v\|_{L^4(\Om)}<\eps$.  Then $v$ is in
$L^2_{k+2}(\Om')$ and there is a universal polynomial $Q_k(x,y)$, with
positive real coefficients, depending at most on $A$, $\alpha$, $k$,
$\Om'$, $\Om$, such that $Q_k(0,0)=0$ and
$$
\|v\|_{L^2_{k+2,A}(\Om')} 
\le Q_k\left(\|w\|_{L^2_{k,A}(\Om)},
\|v\|_{L^2(\Om)}\right).
$$
If $w$ is in 
$C^\8(\Om)$ then $v$ is in $C^\8(\Om')$ and if
$w=0$, then
$$
\|v\|_{L^2_{k+2,A}(\Om')} \le C\|v\|_{L^2(\Om)}.
$$
\end{cor}

\subsection{Estimates for extended anti-self-dual connections in a good local
gauge}\label{subsec:Uhlenbeck} 
It remains to combine the local regularity results and estimates of \S
\ref{subsec:L2_1InhomoLocal}, for extended anti-self-dual connections $A$ where
the connection $A$ is assumed to be in Coulomb gauge relative to the
product $\SO(3)$ connection $\Ga$, with Uhlenbeck's local, Coulomb
gauge-fixing theorem. We then obtain regularity results and estimates for
anti-self-dual connections $A$ with small curvature $F_A$ (see
Theorem 2.3.8 and Proposition 4.4.10 in \cite{DK}).

In order to apply Corollary \ref{cor:L2_1InhomoRegLocal} we need
Uhlenbeck's Coulomb gauge-fixing result \cite[Theorem 2.1 \& Corollary
2.2]{UhlLp}.  Let $B$ (respectively, $\barB$) be the open (respectively,
closed) unit ball centered at the origin in $\RR^4$ and let $G$ be a
compact Lie group. In order to provide universal constants we assume
$\RR^4$ has its standard metric, though the results of this subsection
naturally hold for any $C^\8$ Riemannian metric, with comparable constants
for metrics which are suitably close.

\begin{thm}\label{thm:CoulombBallGauge}
There are positive constants $c$ and $\eps$ with the following
significance. If $2\le p < 4$ is a constant and $A\in
L^p_1(B,\La^1\otimes\fg)\cap L^p_1(\rd B,\La^1\otimes\fg)$ is a
connection matrix whose curvature satisfies $\|F_A\|_{L^p(B)} < \eps$,
then there is a gauge transformation 
$u\in L^p_2(B, G)\cap  L^p_2(\rd B, G)$ such that $u(A)
:= uAu^{-1} - (du)u^{-1}$ satisfies
\begin{align}
d^*u(A) &= 0\text{ on }B, \tag{1}\\
{\textstyle{\frac{\rd}{\rd r}}}\lrcorner u(A) &= 0\text{ on }\rd B, \tag{2}\\
\|u(A)\|_{L^p_1(B)} &\le c\|F_A\|_{L^p(B)}.\tag{3}
\end{align}
If $A$ is in $L^p_k(B)$, for $k\ge 2$, then $u$ is in $L^p_{k+1}(B)$. The gauge
transformation $u$ is unique up to multiplication by a constant element of $G$.
\end{thm}

\begin{rmk}
\begin{enumerate}
\item The $L^p_2$ gauge transformation $u$ is obtained by solving the
second-order quasi-linear elliptic equation $d^*A^u=0$, that is
$$
d^*du = -*(duu^{-1}\wedge *du) -d^*A u + *(*A\wedge du) 
- *(duu^{-1}\wedge *Au),
$$
given $A$ in $L^p_1$. For any open subset $U\Subset B$, 
there is a positive constant $\eps_0(U)$
such that if $A$ in $L^p_2(B)$ satisfies $\|A\|_{L^p_2(B)}<\eps_0$ 
then $u$ is in $L^p_3(U)$ by the implicit function theorem.
In particular, if $A$ is in
$C^\8(B)$ then standard elliptic theory and the Sobolev multiplication
theorem imply that $u$ is in $C^\8(U)$ \cite[Proposition 2.3.4]{DK}.
\item Uhlenbeck's methods \cite{UhlLp} imply that the gauge transformation 
$u_A$ depends analytically on the connection $A$.
\item If $G$ is abelian then the requirement that $\|F_A\|_{L^p(B)} < \eps$ can be
omitted. 
\end{enumerate}
\end{rmk}

It is often useful to rephrase Theorem \ref{thm:CoulombBallGauge} in two
other slightly different ways. Suppose $A$ is an $L^2_k$ connection on a
$C^\8$ principal $G$ bundle $P$ over $B$ with $k\ge 2$ and $\|F_A\|_{L^2(B)} <
\eps$. Then the assertions of Theorem \ref{thm:CoulombBallGauge} are 
equivalent to each of the following:
\begin{itemize}
\item There is an $L^2_{k+1}$ trivialization $\tau:P\to B\times G$ such
that (i) $d_\Gamma^*(\tau(A)-\Gamma)=0$, where $\Gamma$ is the product
connection on $B\times G$, (ii) $\frac{\rd}{\rd r}
\lrcorner (\tau(A)-\Gamma) = 0$, and (iii)
$\|(\tau(A)-\Gamma)\|_{L^2_1(B)} \le c\|F_A\|_{L^2(B)}$.

\item There is an $L^2_{k+1}$ flat connection $\Gamma$ on $P$ such
that (i) $d_\Gamma^*(A-\Gamma)=0$, (ii) $\frac{\rd}{\rd r}
\lrcorner (A-\Gamma) = 0$, and (iii)
$\|(A-\Gamma)\|_{L^2_1(B)} \le c\|F_A\|_{L^2(B)}$, and
an $L^2_{k+1}$ trivialization $P|_B\simeq B\times G$
taking $\Gamma$ to the product connection.
\end{itemize}

We can now combine Theorem \ref{thm:CoulombBallGauge} with 
Corollary \ref{cor:L2_1InhomoRegLocal} or Proposition
\ref{prop:L2_1CoulMonoRegLocal} to give (see
\cite[Theorem 2.3.8]{DK}, \cite[Proposition 8.3]{FU}): 

\begin{thm}
\label{thm:CoulombASDBall}
Let $B\subset\RR^4$ be the open unit ball with center at the origin, let
$U\Subset B$ be an open subset, and let $\Ga$ be the product connection on
$B\times\SO(3)$.  Then there is a positive constant $\eps$ and if $k\ge
1$ is an integer, there is a positive constant $C(k,U)$ with the
following significance.  Suppose that $A$ is an $L^2_1$ solution to the
{\em extended anti-self-dual connection\/} on $\fg_E$ over $B$, namely
$$
F_A^+ = w,
$$
for some $w \in L^2_{k-1}(B,\Lambda^+\otimes\fg_E)$, and that the curvature
of the $\SO(3)$ connection matrix $A$ obeys
$$
\|F_A\|_{L^2(B)} < \eps.  
$$
Then there is an $L^2_{k+1}$ gauge transformation $u:B\to\SU(2)$ such that
$u(A)-\Ga$ is in $L^2_k(U)$ with $d^*(u(A)-\Ga)=0$ over $B$ and
$$
\|u(A)-\Ga\|_{L^2_{k,\Gamma}(U)} 
\le 
Q_k\left(\|u(w)\|_{L^2_{k-1,\Gamma}(B)}, \|F_A\|_{L^2(B)}\right).
$$
If $u(w)$ is in $C^\8(B)$, then $u(A)-\Ga$ is in $C^\8(U)$. If
$w = 0$, so that $A$ is an anti-self-dual connection, then
$$
\|u(A)-\Ga\|_{L^2_{k,\Gamma}(U)} \le C\|F_A\|_{L^2(B)}.
$$
\end{thm}

\begin{rmk}
  The $L^2_k$ estimates for $A$ obtained from Theorem
  \ref{thm:CoulombASDBall} are uniform over all metrics in a small
  $C^k$ neighborhood of a fixed metric \cite[p. 120]{FU}.
\end{rmk}

Again, it is often useful to rephrase Theorem \ref{thm:CoulombASDBall} in
the two other slightly different ways. Suppose $\ell\ge 1$ and that $A$ is
an $L^2_k$ anti-self-dual connection on $\fg_E$ over the unit ball
$B\subset\RR^4$ with $k\ge \max\{2,\ell\}$, $\|F_A\|_{L^2(B)}<\eps$ and
$U\Subset B$. Then the assertions of Theorem \ref{thm:CoulombASDBall} are
equivalent to each of the following:
\begin{itemize}
\item There is a $C^\8$ trivialization $\tau:\fg_E|_B\to
B\times \su(2)$ and a $L^2_{k+1}$
determinant-one, unitary bundle automorphism $u$ of
$E|_B$ such that, with respect to the product connection $\Gamma$ on
$B\times \su(2)$, we have
(i) $d_\Gamma^*(\tau u(A)-\Gamma)=0$, and (ii)
$\|(\tau u(A)-\Gamma,\tau u\Phi)\|_{L^2_{\ell,\Gamma}(U)} 
\le C\|F_A\|_{L^2(B)}$.

\item There is an $L^2_{k+1}$ flat connection $\Gamma$ on $\fg_E|_B$ such
that (i) $d_\Gamma^*(A-\Gamma)=0$, and (ii)
$\|A-\Gamma\|_{L^2_{\ell,\Gamma}(U)} \le c\|F_A\|_{L^2(B)}$,
and an $L^2_{k+1}$ trivialization $\fg_E|_B\simeq B\times \su(2)$
taking $\Gamma$ to the product connection.
\end{itemize}

We will need interior estimates for anti-self-dual connections in a
good local gauge over more general simply-connected
regions than the open balls
considered in Theorem \ref{thm:CoulombASDBall}. Specifically, recall
that a domain $\Om\subset X$ is {\em strongly simply-connected\/} if it has
an open covering by balls $D_1,\dots,D_m$ (not necessarily geodesic) such
that for $1\le r\le m$ the intersection $D_r\cap (D_1\cup\cdots\cup
D_{r-1})$ is connected. We recall (see \cite[Proposition 2.2.3]{DK} or
\cite[Proposition I.2.6]{Kobayashi}):

\begin{prop}
\label{prop:FlatProduct}
If $\Gamma$ is a $C^\8$ flat connection on a principal $G$ bundle
$P$ over a simply-connected manifold $\Om$, then there is a $C^\8$
isomorphism $P\simeq \Om\times G$ taking $\Gamma$ to
the product connection on $\Om\times G$. 
\end{prop}

More generally, if $A$ is $C^\8$ connection on a $G$ bundle $P$ over 
a simply-connected manifold-with-boundary $\barOm=\Om\cup\rd\Om$
with $L^p$-small curvature (with $p>2$), then Uhlenbeck's theorem implies
that $A$ is $L^p_2$-gauge equivalent to a connection which is $L^p_1$-close
to an $L^p_1$ flat connection on $P$ (see \cite[Corollary 4.3]{UhlChern} or
\cite[p. 163]{DK}). We recall the following {\em a priori\/}
interior estimate (see also \cite[Proposition 3.18]{FL1} and
\cite[Proposition 4.4.10]{DK}).

\begin{prop}
\label{prop:MonopoleGoodGaugeIntEst}
Let $X$ be a closed, oriented, Riemannian four-manifold and let $\Om\subset
X$ be a strongly simply-connected open subset. Then there is a positive
constant $\eps(\Om)$ with the following significance. For $\Om'\Subset \Om$
a precompact open subset and an integer $\ell\ge 1$, there is a constant
$C(\ell,\Om',\Om)$ such that the following holds. Suppose $A$ is an $L^2_k$
{\em extended anti-self-dual connection\/} on $\fg_E$ over $\Om$, namely
$$
F_A^+ = w,
$$
for some $w \in L^2_{k-1}(\Omega,\Lambda^+\otimes\fg_E)$,  with $k\ge
\max\{2,\ell\}$ such that
$$
\|F_A\|_{L^2(\Om)} < \eps.
$$
Then there is an $L^2_{k+1}$ flat connection $\Gamma$ on $\fg_E|_{\Om'}$ such
that 
$$
\|A-\Ga\|_{L^2_{\ell,\Ga}(\Om')}
\le 
Q_k\left(\|w\|_{L^2_{k-1,\Gamma}(\Om)}, \|F_A\|_{L^2(\Om)}\right),
$$
and an $L^2_{k+1}$ trivialization $\fg_E|_{\Om'}\simeq \Om'\times \su(2)$
taking $\Gamma$ to the product connection.
If $w$ is in $C^\8(\Om)$, then $A$, $\Ga$ are in $C^\8(\Om')$. If
$w = 0$, so that $A$ is an anti-self-dual connection, then
$$
\|A-\Ga\|_{L^2_{\ell,\Ga}(\Om')} \le C\|F_A\|_{L^2(\Om)}.
$$
\end{prop}

For some local patching arguments over simply-connected open regions
$Y\subset X$ it is very useful to be able to lift gauge
transformations in $\sG_{\su(E)}(Y)$ to gauge transformations in
$\sG_E(Y)$. The following result tells us that this is always possible
when $Y$ is simply connected; it is an extension of Theorem IV.3.1 in
\cite{MorganGTNotes} from $\SU(2)$ to $\U(2)$ bundles.

\begin{prop}
\label{prop:GaugeGroupSequence}
\cite[Proposition 2.16]{FL1}
Let $E$ be a rank-two, Hermitian bundle over
a connected manifold $X$. Then there is an exact sequence
\footnote{It can also be shown that the map 
$\sG_{\su(E)}\rightarrow H^1(X;\ZZ/2\ZZ)$ is surjective.}
$$
1\rightarrow
\{\pm \id_E\}
\rightarrow
\sG_E
\rightarrow
\sG_{\su(E)}
\rightarrow
H^1(X;\ZZ/2\ZZ).
$$
\end{prop}

\subsection{Regularity of Taubes' solutions to 
extended anti-self-dual equations}
\label{subsec:TaubesRegularity}
Finally, we come to one of the main results of this section.

\begin{prop}
\label{prop:RegularityTaubesSolution}
Let $X$ be a closed, oriented, $C^\8$ four-manifold with 
$C^\8$ Riemannian metric $g$
and let $E$ be a Hermitian, rank-two bundle
over $X$.  Let $A$ be a $C^\8$ connection on the $C^\8$ bundle
$\fg_E$ over $X$. Suppose that $a\in
L^2_1(X,\La^1\otimes\fg_E)$ is a
solution to an extended anti-self-dual equation over $X$,
$$
F^{+,g}(A+a) = w,
$$
where $a=d_A^{+,*}v$, for some $v\in (C^0\cap L^2_2)(X,\La^+\otimes\fg_E)$,
and $w\in C^\8(X,\La^+\otimes\fg_E)$. Then $a$ is in
$C^\8(X,\La^1\otimes\fg_E)$ and $v\in C^\8(X,\La^+\otimes\fg_E)$.
\end{prop}

\begin{rmk}
\label{rmk:AnalysisForTopologists}
The reader may wonder why we cannot apply Corollary \ref{cor:L2_1InhomoReg}
to the resulting elliptic system \eqref{eq:GluingExtASDforRegularityProof}
and deduce straightaway that the solution $a$ is in $C^\8$.  The
difficulty, of course, is that the hypotheses of Corollary
\ref{cor:L2_1InhomoReg} are not necessarily satisfied due to the possibly
unfavorable dependence of the constant $\eps(A)$ on the connection $A$. For
example, by combining results of \S \ref{sec:Decay} and \S
\ref{sec:Existence}, we can show that $\|a\|_{L^4(X)} \leq C\lambda$ when
$A_\lambda$ is a spliced connection (produced by the alogorithm of \S
\ref{sec:Splicing}), $\lambda$ is a parameter which can be made arbitrarily
small, and $C$ is an essentially universal constant, depending at most on
$\|F_{A_\lambda}^+\|_{L^\8(X)}$ and $\|F_{A_{\lambda}}\|_{L^2}$: in
particular, $C$ is independent of $\lambda$. Now, as we make $\lambda$
smaller to try to satisfy the constraint of Corollary
\ref{cor:L2_1InhomoReg}, we will in general find that $\eps(A_\lambda)$ 
also becomes smaller (as it will, for example, if $\eps(A_\lambda)$ depends
on an $L^p_\ell$ norm of $F_{A_\lambda}$ with $p>2$, $\ell=0$ or $\ell\geq
1$, $p=2$). Thus, as we do not necessarily have a lower bound for $\eps(A)$
which is uniform with respect to $A$, we are left `chasing our tail' and
Corollary \ref{cor:L2_1InhomoReg} will not apply.
\end{rmk}

\begin{proof}[Proof of Proposition \ref{prop:RegularityTaubesSolution}]
Since $a = d_A^{+,*}v$, the extended anti-self-dual and Coulomb-gauge equations
become
\begin{equation}
\begin{aligned}
\label{eq:GluingExtASDforRegularityProof}
d_A^+a + (a\wedge a)^+ &= w-F_A^+,
\\
d_A^*a &= (F_A^+)^*v,
\end{aligned}
\end{equation}
where $F_A^+ \in C^\8(\La^+\otimes\fg_E)$ and $(F_A^+)^*v \in (C^0\cap
L^2_2)(\La^+\otimes\fg_E)$. Observe that
$$
(\zeta,w') := ((F_A^+)^*v,w-F_A^+) 
\in 
L^2_2(\fg_E)\oplus L^2_2(\La^+\otimes\fg_E).
$$
{}From Remark \ref{rmk:AnalysisForTopologists} we note that we {\em
cannot\/} take the apparently obvious route and apply Corollary
\ref{cor:L2_1InhomoReg}. 
Instead, we fix a point $x\in X$ and a geodesic ball $B(x,\delta)$ centered
at $x$, fix a $C^\8$ trivialization $\fg_E|_{B(x,\delta)}\simeq
B(x,\delta)\times\su(2)$, and let $\Gamma$ be the resulting product
connection. Then
\begin{equation}
\label{eq:LocalASDRegularityEqn}
\begin{aligned}
d^+b + (b\wedge b)^+ &= w-F_\Gamma^+ = w,
\\
d^*b &= *((A-\Gamma)\wedge *a) + d^*(A-\Gamma) + (F_A^+)^*v,
\end{aligned}
\end{equation}
where $b := A-\Gamma+a \in L^2_1(B(x,\delta),\La^1\otimes\fg_E)$ and
the first equation is the expansion of $F^+(\Gamma+b)=w$. In
particular, the right-hand side of equation \eqref{eq:LocalASDRegularityEqn} is
in $L^2_1(B(x,\delta))$.

In Corollary \ref{cor:L2_1InhomoRegLocal}, let $\Omega = B\subset \RR^4$ be
the unit ball, let $\eps = \eps(\Omega)$ be the corresponding positive
constant, and let $\Omega' = B'$ be a slightly smaller concentric ball.
Since $a\in L^4(X,\La^1\otimes\fg_E)$, we know that $\|a\|_{L^4(U)}\to 0$
as $\Vol(U)\to 0$, where $U$ is any measurable subset of $X$. In
particular, for a small enough ball $B(x,\delta)$, we have
$\|b\|_{L^4(B(x,\delta))} < \eps$. We now rescale the metric on the ball
$B(x,\delta)$ to give a ball $\tilde B$ with unit radius, center $x$, and
metric which is approximately flat. Note that the equations
\eqref{eq:LocalASDRegularityEqn} are scale-equivariant and that
$$
\|b\|_{L^4(\tilde B)} 
=
\|b\|_{L^4(B(x,\delta))} < \eps,
$$
by the scale-invariance of the $L^4$ norm on one-forms. Corollary
\ref{cor:L2_1InhomoRegLocal} now implies that $b$ --- and thus $a$ --- is
in $L^2_2$ on the slightly smaller rescaled ball $\tilde B'$. Since
$d_A^+d_A^{+,*}v = d_A^+a$ is now in $L^2_1(\tilde B')$, we have $v\in
L^2_3(\tilde B')$ by standard elliptic regularity theory for a Laplacian,
$d_A^+d_A^{+,*}$, with $C^\8$ coefficients \cite{Hormander}.
Hence, the right-hand side of equation
\eqref{eq:LocalASDRegularityEqn} is in $L^2_2(\tilde B')$.  Repeating the
above argument shows that $a$ is in $C^\8$ on the ball $\tilde
B_{1/2}\subset \tilde B$ with center $x$
and radius $\half$. Hence, $a$ is in $C^\8$ on $B(x,\delta/2)$ and so, as
the point $x\in X$ was arbitrary, we see that $a \in
C^\8(X,\La^1\otimes\fg_E)$ and $v\in C^\8(X,\La^+\otimes\fg_E)$.
\end{proof}


\section{Decay estimates for curvature and connections in
  Coulomb and radial gauges}
\label{sec:Decay}
Our splicing and unsplicing constructions for connections requires
estimates for local connection one-forms in suitable gauges, usually
Coulomb or radial, over small balls, over their complements, and over
annuli. The most useful estimates bound Sobolev or H\"older norms of local
connection one-forms in terms of the $L^2$ norms of their curvatures and
the radii defining the given regions, with estimate constants which are
essentially universal. These are the kinds of estimates we derive in this
section. Some of these bounds can be derived from the curvature decay
estimates due, in various forms, to Donaldson, Groisser-Parker, R\aa de,
Taubes, and Uhlenbeck. Perhaps the most useful and powerful statement is
due to Donaldson, Groisser-Parker and R\aa de; see the discussion in
\cite{FLKM2}.  However, we prefer to use methods which are as simple as
possible. Furthermore, our chosen methods have the the following
advantages:
\begin{itemize}
\item
They extend easily to the case of solutions to {\em extended\/}  
anti-self-dual equations, such as equations \eqref{eq:ASD-A} or 
\eqref{eq:QuickASDEqnForv}. 
\item
They extend easily to the case of $\PU(2)$ {\em monopoles\/}
\cite{FL3}.
\end{itemize}
In contrast, the proofs of the elegant, pointwise decay estimates of Donaldson,
Groisser-Parker, and R\aa de --- for the curvature of anti-self-dual or
Yang-Mills connections --- requires a delicate application of the
Chern-Simons functional and considerably more work is needed to extend them
to the case of $\PU(2)$ monopoles, though decay estimates of the desired
type for Seiberg-Witten monopoles are described in \cite{MST}. It is likely
that their estimates can be extended to cover the case of extended
anti-self-dual connections, but again some work is required to show this.
The disadvantage of our approach is that we must separately consider the
cases of (i) balls, (ii) complements of balls in $\RR^4$, and (iii) annuli,
while the decay estimate theorem of Donaldson, Groisser-Parker, and R\aa de
would give a more uniform treatment of these different regions.

In all cases, the estimates of this section continue to hold in
situations where the Riemannian metric is uniformly $C^r$ close enough
to the standard metric. Thus, for the sake of notational simplicity, we
shall assume that our connections are defined on regions in $\RR^4$
equipped with its standard Riemannian metric $g$. We let $B(r)$ denote
the open ball $\{x\in\RR^4:|x|<r\}$ of radius $r$ and center at the
origin and let $\barB(r)$ denote the closed ball $\{x\in\RR^4:|x|\le r\}$.

\subsection{Estimates for anti-self-dual connections over balls}
\label{subsec:EstimateConnctionBall}
The basic case is that of ball in $\RR^4$: estimates for connections
over the complements of balls and over annuli follow from this.
Theorem \ref{thm:CoulombASDBall} yields an interior $L^\8$
bound for the curvature $F_A$ of an anti-self-dual connection over a ball
in terms of the $L^2$ norm of $F_A$; compare Proposition 8.3 in
\cite{FU} and Lemma 9.1 in \cite{TauFrame}. 

\begin{lem}
\label{lem:C0ASDFAFixedBall}
There are positive constants $\eps$ and, given 
a constant $1\leq p\leq \8$ and an integer $k\geq 0$, a
positive constant $c(k,p)$ with the following significance.  Let $G$ be
a compact Lie group. Suppose $r_1>0$ is a constant and that $A$ is an
anti-self-dual connection on the product bundle $P= B(r_1)\times G$,
where $B(r_1)$ is a ball centered at the origin in $\RR^4$, with its
Euclidean metric, such that
$$
\|F_A\|_{L^2(B(r_1))} < \eps,
$$
then
$$
\|\cov_\Gamma^kF_A\|_{L^p(B(r_1/2))} 
\leq c(k,p)r_1^{(4/p)-k-2}\|F_A\|_{L^2(B(r_1))}.
$$
Moreover, if $\Gamma$ is the product connection 
on $\RR^4\times G$ and $d_\Gamma^*(A-\Gamma)=0$, 
then 
$$
\|\cov_\Gamma^k(A-\Gamma)\|_{L^p(B(r_1/2))} 
\leq c(k,p)r_1^{(4/p)-k-1}\|F_A\|_{L^2(B(r_1))}.
$$
Finally, the estimates for $F_A$ also hold when $\cov_\Gamma$ is
replaced by $\cov_A$.
\end{lem}

\begin{pf}
  We define a rescaled metric $\tg$ on $\RR^4$ by setting
  $\tg_{\mu\nu} = r_1^{-2}g_{\mu\nu}$, so $\tg^{\mu\nu} =
  r_1^2g^{\mu\nu}$ and $B_g(r_1) = B_\tg(1)$. Then,
\begin{equation}
\label{eq:RescalingCovDerivConnCurvPtwiseNorms}
\begin{aligned}
|\cov_\Gamma^k(A-\Gamma)|_g
&= r_1^{-k-1}|\cov_\Gamma^k(A-\Gamma)|_\tg,
\\
|\cov_A^k F_A|_g
&= r_1^{-k-2}|\cov_A^k F_A|_\tg.
\end{aligned}
\end{equation}
Since $\|F_A\|_{L^2(B(1),\tg)}<\eps$ and $\eps$ is sufficiently small,
Theorem \ref{thm:CoulombBallGauge} provides a gauge transformation $u$
over $B_\tg(1)$ such that $d_\Gamma^*(u(A)-\Gamma)=0$, so without loss
of generality we may assume that $A$ is already in Coulomb gauge, as
we suppose in the second assertion of the Lemma. Theorem
\ref{thm:CoulombASDBall}, noting that its proof extends from $p=2$ to
the general case in a straightforward manner, then yields 
$$
\|A-\Gamma\|_{L^p_{k,\Gamma}(B(1/2),\tg)} 
\le c(k,p)\|F_A\|_{L^2(B(1),\tg)},
$$
and the desired Sobolev $L^p_k$ estimates follow, since 
\begin{align*}
\|\cov_\Gamma^k(A-\Gamma)\|_{L^p(B(r_1/2),g)}
&= r_1^{(4/p)-k-1}\|\cov_\Gamma^k(A-\Gamma)|_{L^2(B(1/2),\tg)},
\\
\|\cov_A^k F_A\|_{L^p(B(r_1/2),g)}
&= r_1^{(4/p)-k-2}\|\cov_A^k F_A\|_{L^2(B(1/2),\tg)}.
\end{align*}
This completes the proof.
\end{pf}

In our treatment of `ungluing' in \cite{FLKM2} we shall need $L^2_1$
estimates for covariant derivatives of radial-gauge connection
one-forms, while $L^p$ estimates are needed throughout our article.
Therefore, suppose that the connection one-form $A-\Gamma$ is in radial
gauge with respect to the origin, $x^\mu A_\mu(x) = 0$, so \cite{UhlRem}
\begin{equation}
\label{eq:RadialGaugeConnInTermsCurv}
A_\nu(x) = \int_0^1 tx^\mu(F_A)_{\mu\nu}(ty)\,dt,
\end{equation}
where $A_\mu = (A-\Gamma)_\mu$. We then have the following estimates
for the local connection one-forms and their derivatives (compare
Taubes' Lemmas 9.1 and 9.2 in \cite{TauFrame}):

\begin{lem}
\label{lem:SobolevEstimateRadialGaugeConnOneForm}
Continue the hypotheses and notation of Lemma
\ref{lem:C0ASDFAFixedBall}. Suppose that the connection one-form
$A-\Gamma$ is in radial gauge with respect to the origin. Then, for
any $1\leq p\leq \8$,
$$
\|\cov_\Gamma^k(A-\Gamma)\|_{L^p(B(r_1/2))}
\leq
c(k,p)r_1^{(4/p)-k-1}\|F_A\|_{L^2(B(r_1))}.
$$
\end{lem}

\begin{proof}
Denote $A-\Gamma$ simply by $A$.
We confine our attention to the case $k=0, 1$ (which are the only ones we
shall need) and leave the general case to the reader. {}From
equation \eqref{eq:RadialGaugeConnInTermsCurv} we have
\begin{align}
\label{eq:PointwiseCovRadialConnEst}
|A|(x) &\leq |x|\cdot |F_A|(x), \quad x\in \RR^4,
\\
\notag
|\cov_\Gamma A|(x) 
&\leq |x|\cdot|\cov_\Gamma F_A|(x) + |F_A|(x).
\end{align}
Lemma \ref{lem:C0ASDFAFixedBall} then yields the following Sobolev
estimates. First,
$$
\|A\|_{L^p(B(r_1/2)}
\leq cr_1^{(4/p)-1}\|F_A\|_{L^2(B(r_1))},
$$
which gives the desired bound when $k=0$. For $k=1$, we have
\begin{align*}
\|\cov_\Gamma A\|_{L^p(B(r_1/2)}
&\leq r_1\|\cov_\Gamma F_A\|_{L^p(B(r_1/2))} + \|F_A\|_{L^p(B(r_1/2))}
\\
&\leq cr_1^{4/p}(r_1\cdot r_1^{-3} + r_1^{-2})\|F_A\|_{L^2(B(r_1))}
\\
&= cr_1^{(4/p)-2}\|F_A\|_{L^2(B(r_1))}.
\end{align*}
The cases $k\geq 2$ are similar.
\end{proof}

\begin{rmk}
  The estimates in Lemmas \ref{lem:C0ASDFAFixedBall} and
  \ref{lem:SobolevEstimateRadialGaugeConnOneForm} agree in the case
  $k=0$ for $F_A$ with that obtained by integrating the decay estimate
  of Donaldson et al. when $r_0=0$. The estimates for
  $F_A$ and $A-\Gamma$ agree with those obtained by integrating the
  pointwise decay estimates of \cite[Lemma 9.1]{TauFrame}.
\end{rmk}

\subsection{Estimates for anti-self-dual connections over the
  complement of a ball} 
\label{subsec:EstimateConnctionBallComplement}
We now consider the case where $A$ is an anti-self-dual connection
with small curvature on the complement of a ball in $\RR^4$. 

Our first task is to note that given a $C^\8$ anti-self-dual
connection $A$ on $\RR^4-B(r_0)$, with finite energy, there is no loss
in generality --- thanks to the removable singularities theorem of K.
Uhlenbeck \cite[Theorem 4.4.12]{DK}, \cite[Theorem 4.1]{UhlRem} --- if
we assume that $A$ is a $C^\8$ connection on $S^4-\varphi_n(B(r_0))$,
where the origin of $\RR^4$ is identified with the north pole of
$S^4$.

\begin{thm}
\label{thm:RemSing}
Let $A$ be a $C^\8$ anti-self-dual connection on a $G$ bundle $P$ over the
punctured unit ball $B\less\{0\}$ in $\RR^4$, with its Euclidean metric. If
$$
\|F_A\|_{L^2(B)} < \8,
$$
there is a $C^\8$ anti-self-dual connection $\barA$ on a $C^\8$ bundle
$\barP$ over $B$ and a $C^\8$ bundle isomorphism $u:\barP|_{B\less\{0\}}\to
P|_{B\less\{0\}}$ such that $u(\barA) = A$ over $B\less\{0\}$.
\end{thm}

We now consider the problem of estimating the $L^p_k$ norms of
$F_A$ and $A-\Gamma$ (in radial gauge) on $\RR^4-B(r_0)$. 

\begin{lem}
\label{lem:C0ASDFAFixedBallComplement}
There are positive constants $\eps$ and, given an integer $k\geq 0$
and a constant $1\leq p\leq\8$, a
positive constant $c(k,p)$ with the following significance.  Let $G$ be
a compact Lie group. Suppose $r_0>0$ is a constant and that $A$ is an
anti-self-dual connection on the product bundle $P= (\RR^4-B(r_0))\times G$,
where $B(r_0)$ is a ball centered at the origin in $\RR^4$, with its
Euclidean metric, such that
$$
\|F_A\|_{L^2(\RR^4-B(r_0))} < \eps,
$$
then, provided $k+2-4/p\geq 0$,
$$
\|\cov_A^kF_A\|_{L^p(\RR^4-B(2r_0))} 
\leq c(k,p)r_0^{(4/p)-k-2}\|F_A\|_{L^2(\RR^4-B(r_0))}.
$$
Moreover, if $\Gamma$ is the product connection on
$(\RR^4-B(r_0))\times G$ and $A-\Gamma$ is in radial gauge with
respect to $\8$ then, for $k+1-4/p\geq 0$,
$$
\|\cov_\Gamma^k(A-\Gamma)\|_{L^p(\RR^4-B(2r_0))} 
\leq c(k,p)r_0^{(4/p)-k-1}\|F_A\|_{L^2(\RR^4-B(r_0))}.
$$
\end{lem}

\begin{proof}
  Identifying $\RR^4$ with the quaternions $\HH$, we recall that $y =
  \iota^*(y) = \varphi_s^{-1}\circ\varphi_n(x) = x^{-1}$ for
  $x\in\RR^4-\{0\}$, where $\varphi_s:S^4-\{n\}\to\RR^4$ and
  $\varphi_n:S^4-\{s\}\to\RR^4$ are stereographic projections.  The
  inversion map is conformal, so $\iota^*A$ is an anti-self-dual
  connection on the punctured ball $\barB(1/r_0)-\{0\}$, with curvature
  bound
$$
\|F_{\iota^*A}\|_{L^2(B(1/r_0))} 
= \|F_A\|_{L^2(\RR^4-B(r_0))} < \eps.
$$
As $\|F_{\iota^*A}\|_{L^2(B(1/r_0))} < \8$,
Theorem \ref{thm:RemSing} implies that there is a $C^\8$ gauge transformation
$u$ of $\barB(1/r_0)\times G$ such that the connection one-form
$u(\iota^*A)$ extends to a $C^\8$ connection over the ball
$B(1/r_0)$. Thus, we may consider $A$ to be a $C^\8$ connection on the
product bundle $P = (S^4-\varphi_n(B(r_0)))\times G$, with
product connection $\Gamma$.
 
{}From the identity $|dx| = |y|^{-2}|dy|$, we obtain
\begin{align*}
|\varphi_n^*\cov_A^k F_A|(x)
&=
|y|^{-2(k+2)}\cdot|\varphi_s^*\cov_A^k F_A|(y),
\\
|\varphi_n^*\cov_\Gamma^k(A-\Gamma)|(x)
&=
|y|^{-2(k+1)}\cdot|\varphi_s^*\cov_\Gamma^k(A-\Gamma)|(y).
\end{align*}
Combining these observations and noting that $d^4x =
|y|^{-8}d^4y$, we see that (for $2(k+2)-8/p\geq 0$ and
$2(k+1)-8/p\geq 0$, respectively)
\begin{align*}
\|\varphi_n^*\cov_A^kF_A\|_{L^p(\RR^4-B(2r_0))} 
&= \||y|^{2(k+2)-8/p}\varphi_s^*\cov_A^kF_A\|_{L^p(B(1/2r_0))}
\\
&\leq (2r_0)^{-2(k+2)+8/p}\|\varphi_s^*\cov_A^kF_A\|_{L^p(B(1/2r_0))},
\\
\|\varphi_n^*\cov_\Gamma^k(A-\Gamma)\|_{L^p(\RR^4-B(2r_0))}
&=
\||y|^{2(k+1)-8/p}\varphi_s^*\cov_\Gamma^k(A-\Gamma)\|_{L^p(B(1/2r_0))}
\\
&\leq
(2r_0)^{-2(k+1)+8/p}\|\varphi_s^*\cov_\Gamma^k(A-\Gamma)\|_{L^p(B(1/2r_0))}.
\end{align*}
Applying Lemmas \ref{lem:C0ASDFAFixedBall} and
\ref{lem:SobolevEstimateRadialGaugeConnOneForm} to the connection
$\varphi_s^*A$ on $B(1/r_0)$, in radial gauge with respect to the
origin, we see that
\begin{align*}
\|\varphi_s^*\cov_A^kF_A\|_{L^p(B(1/2r_0))} 
&\leq c(k,p)r_0^{-(4/p)+k+2}\|F_A\|_{L^2(B(1/r_0))},
\\
\|\varphi_s^*\cov_\Gamma^k(A-\Gamma)\|_{L^p(B(1/2r_0))} 
&\leq c(k,p)r_0^{-(4/p)+k+1}\|F_A\|_{L^2(B(1/r_0))}.
\end{align*}
Combining the preceding two pairs of estimates then yields the desired
bounds.
\end{proof}

\begin{rmk}
The estimates in Lemma \ref{lem:C0ASDFAFixedBallComplement} agree in
the case $k=0$ for $F_A$ with that obtained by integrating the decay
estimate of Donaldson et al. when $r_1=\8$. The estimates
for $F_A$ and $A-\Gamma$ agree with those obtained by integrating the
pointwise decay estimates of \cite[Lemma 9.2]{TauFrame}.
\end{rmk}

\subsection{Estimates for anti-self-dual connections over annuli}
\label{subsec:EstimateConnctionAnnulus}
The case of an annulus requires a slightly more sophisticated approach
and a complete treatment --- where both inner and outer radius can be
specified independently --- requires the decay estimate of
Donaldson et al.  For positive constants $r_0<r_1$,
let
$$
\Om(r_0,r_1) := \{x\in\RR^4:r_0<|x|<r_1\}.
$$
We then have the following consequence of Proposition
\ref{prop:MonopoleGoodGaugeIntEst}.

\begin{lem}
\label{lem:ASDBasicAnnGauge}
There are positive constants $\eps$ and, given an integer $k\geq 0$
and a constant $1\leq p\leq\8$, a
positive constant $c(k,p)$ with the following significance.  Let $G$ be
a compact Lie group. Suppose $r_1>0$ is a constant and that $A$ is an $C^\8$
anti-self-dual connection on a $G$ bundle $P$ over the annulus
$\Omega = \Omega(\quarter r_1,4r_1)$ in $\RR^4$, with its
Euclidean metric, such that
$$
\|F_A\|_{L^2(\Omega(\quarter r_1,4r_1))} < \eps.
$$
Then there is an $C^\8$ flat connection on $P|_{\Omega'}$ such that
\begin{align*}
\|\cov_\Gamma^k(A-\Gamma)\|_{L^p(\Omega(\half r_1,2r_1))} 
\leq c(k,p)r_1^{(4/p)-k-1}\|F_A\|_{L^2(\Omega(\quarter r_1,4r_1))},
\\
\|\cov_\Gamma^kF_A\|_{L^p(\Omega(\half r_1,2r_1))} 
\leq c(k,p)r_1^{(4/p)-k-2}\|F_A\|_{L^2(\Omega(\quarter r_1,4r_1))},
\end{align*}
and a $C^\8$ bundle isomorphism $u:P|_{\Omega'}\to\Omega'\times G$
taking $\Gamma$ to the product connection. Finally, if $\Gamma$ is the
product connection on $\Omega\times G$, $P = \Omega\times G$, and
$A-\Gamma$ is in radial gauge on $\Omega$ with respect to the origin,
then the preceding Sobolev  bounds for $A-\Gamma$ continue to hold.
\end{lem}

\begin{pf}
  We proceed much as in the proof of Lemma \ref{lem:C0ASDFAFixedBall}.
  Define a rescaled metric $\tg$ on $\RR^4$ by setting $\tg_{\mu\nu} =
  r_1^{-2}g_{\mu\nu}$, so $\tg^{\mu\nu} = r_1^2g^{\mu\nu}$ and
  $\Omega_g(\quarter r_1,4r_1) = \Omega_\tg(\quarter,4)$. 
Since $\|F_A\|_{L^2(\Omega(\quarter,4),\tg)}<\eps$ and $\eps$ is
sufficiently small, Proposition \ref{prop:MonopoleGoodGaugeIntEst}
provides a $C^\8$ flat connection $\Gamma$ and a $C^\8$ bundle
isomorphism $P|_{\Omega'}\to \Omega'\times G$, taking $\Gamma$ to the
product connection, such that (noting that its proof extends from
$p=2$ to the general case in a straightforward manner)
\begin{align*}
\|A-\Gamma\|_{L^p_{k,\Gamma}(\Omega_\tg(\half,2))} 
\leq c(k,p)\|F_A\|_{L^2(\Omega_\tg(\quarter,4))},
\\
\|F_A\|_{L^p_{k,\Gamma}(\Omega_\tg(\half,2))} 
\leq c(k,p)\|F_A\|_{L^2(\Omega_\tg(\quarter,4))}.
\end{align*}
The rescaling behavior of pointwise norms
\eqref{eq:RescalingCovDerivConnCurvPtwiseNorms} implies that
\begin{align*}
\|\cov_\Gamma^k(A-\Gamma)\|_{L^p(\Omega(\half r_1,2r_1),g)}
&= r_1^{(4/p)-k-1}\|\cov_\Gamma^k(A-\Gamma)|_{L^2(\Omega(\half,2),\tg)},
\\
\|\cov_A^k F_A\|_{L^p(\Omega(\half r_1,2r_1),g)}
&= r_1^{(4/p)-k-2}\|\cov_A^k F_A\|_{L^2(\Omega(\half,2),\tg)}.
\end{align*}
Therefore, the desired Sobolev $L^p_k$ estimates follow by combining
the preceding two pairs of estimates.

Lastly, if $A-\Gamma$ is in radial gauge with respect to the origin,
that is $x^\mu A_\mu(x)=0$ for all $x\in\Omega$, then 
\eqref{eq:RadialGaugeConnInTermsCurv} yields
$$
A_\nu(x) = \int_{r_1/4}^{|x|} t\theta^\mu(F_A)_{\mu\nu}(t\theta)\,dt,
\quad x \in \Omega,
$$
where $\theta = x/|x|\in S^3$. The rest of the proof now mimics
that of Lemma \ref{lem:SobolevEstimateRadialGaugeConnOneForm}.
\end{pf}

\subsection{Estimates for cut-off functions}
\label{subsec:EstimateCutoff}
For any $N>1$ and $\la>0$, Lemma 7.2.10 in \cite{DK} provides a
clever construction of a cut-off function $\be$ on $\RR^4$ such that $\be(x)
= 0$ for $|x|\le N^{-1}\sqrt{\la}$, $\be(x) = 1$ for $|x|\ge
N\sqrt{\la}$, and satisfying the crucial estimate
$$
\|\cov\be\|_{L^4} \le c(\log N)^{-3/4},
$$
where $c$ is a constant independent of $N$ and $\la$. We will need the
following refinement of this result; for the definition and properties of
the $L^\sharp$ and related families of Sobolev norms, see \cite[\S
4]{FeehanSlice}.

\begin{lem}
\label{lem:dBetaEst}
There is a positive constant $c$ such that the following holds. For any
$N>4$ and 
$\la>0$, there is a $C^\8$ cut-off function $\be=\be_{N,\la}$ on $\RR^4$
such that 
$$
\be(x) = 
\begin{cases}
1 &\text{if }|x|\ge \half\sqrt{\la}, \\
0 &\text{if }|x|\le N^{-1}\sqrt{\la}, 
\end{cases}
$$
and satisfying the following estimates:
\begin{align}
\|\cov\be\|_{L^4} &\le c(\log N)^{-3/4}, \tag{1}\\
\|\cov^2\be\|_{L^2} &\le c(\log N)^{-1/2}, \tag{2}\\
\|\cov\be\|_{L^{2\sharp}} &\le c(\log N)^{-1/2}, \tag{3}
\\
\|\cov\be\|_{L^2} &\le c(\log N)^{-1}\sqrt{\lambda}, \tag{4}
\\
\|\cov^2\be\|_{L^{4/3}} &\le c(\log N)^{-1}\sqrt{\lambda}. \tag{5}
\end{align}
\end{lem}

\begin{pf}
We fix a $C^\8$ cut-off function $\ka:\RR\to [0,1]$ such that
$\ka(t) = 1$ if $t\ge 1$ and $\ka(t)=0$ if $t\le 0$.
Now define a $C^\8$ cut-off function $\al=\al_N:\RR\to [0,1]$, depending on
the parameter $N$, by setting
$$
\al(t) := \ka\left(\frac{\log N + t}{\log N - \log 2}\right).
$$
Therefore, 
$$
\alpha(t)
=
\begin{cases}
1 &\text{if }t\ge -\log 2, \\
0 &\text{if }t\le -\log N,
\end{cases}
$$
and there is a constant $c$ independent of $N$ such that
$$
\left|\frac{d\al}{dt}\right| \le \frac{c}{\log N} \quad\text{and}\quad
\left|\frac{d^2\al}{dt^2}\right| \le \frac{c}{(\log N)^2}.
$$
Now define the cut-off function $\be=\be_{N,\la}:\RR\to [0,1]$, depending on the
parameters $N$ and $\la$, by setting
$$
\be(x) := \al(\log|x| - \log\sqrt{\la}),\qquad x\in \RR^4,
$$
and observe that $\be$ satisfies the following pointwise bounds on $\RR^4$:
\begin{equation}
\label{eq:PointwiseBetaEstimates}
|\cov\be(x)| \le \frac{c}{r\log N} \quad\text{and}\quad
|\cov^2\be(x)| \le \frac{c}{r^2\log N},
\end{equation}
where $r=|x|$. Consequently, using $dV = r^3\,drd\theta$,
\begin{align*}
\int_{\RR^4}|\cov\be|^4\,dV 
&\le \frac{c}{(\log N)^4}\int_{N^{-1}\sqrt{\la}}^{\half\sqrt{\la}}\frac{dr}{r} 
\le c\frac{\log N - \log 2}{(\log N)^4} \le \frac{c}{(\log N)^3},
\end{align*}
and so this gives (1). Next, by the
pointwise estimate \eqref{eq:PointwiseBetaEstimates} we have
\begin{align*}
\int_{\RR^4}|\cov^2\be|^2\,dV 
&\le \frac{c}{(\log N)^2}\int_{N^{-1}\sqrt{\la}}^{\half\sqrt{\la}}\frac{dr}{r} 
\le c\frac{\log N - \log 2}{(\log N)^2} \le \frac{c}{\log N},
\end{align*}
which gives (2). Lastly, \cite[Lemma 5.2]{FeehanSlice} implies that
$$
\|\cov\be\|_{L^{2\sharp}} \le c(\|\cov^2\be\|_{L^2} + \|\cov\be\|_{L^4}),
$$
and so (3) follows from (1) and (2). 
Next, the pointwise estimate \eqref{eq:PointwiseBetaEstimates} yields
\begin{align*}
\int_{\RR^4}|\cov\be|^2\,dV 
&\le \frac{c}{(\log N)^2}\int_{N^{-1}\sqrt{\la}}^{\half\sqrt{\la}}r dr
\le c\frac{\lambda(1-N^{-2})}{(\log N)^2} 
\le \frac{c\lambda}{(\log N)^2},
\end{align*}
and so this gives (4). Finally, the pointwise estimate
\eqref{eq:PointwiseBetaEstimates} yields
\begin{align*}
\int_{\RR^4}|\cov\be|^{4/3}\,dV 
&\le \frac{c}{(\log N)^{4/3}}\int_{N^{-1}\sqrt{\la}}^{\half\sqrt{\la}}
r^{1/3} dr
\le c\frac{\lambda^{2/3}(1-N^{-4/3})}{(\log N)^{4/3}} 
\le \frac{c\lambda^{2/3}}{(\log N)^{4/3}},
\end{align*}
and so this gives (5). 
\end{pf}

\subsection{Universal estimate for the \boldmath{$L^{\sharp}$} norm of the
the curvature of an anti-self-dual connection}
\label{subsec:UniversalEstimateSharpNormSDCurv}
We need to establish the existence of a universal
estimate for the $L^{\sharp,2}$ norm of the curvature of a
$g$-anti-self-dual connection.

\begin{lem}
\label{lem:UniversalLSharp2BoundFA}
Let $(X,g)$ be a $C^\8$, closed, oriented, Riemannian four-manifold
and let $\kappa\geq 0$ be an integer. Then there is a constant
$c=c(X,g,\kappa)$ such that $\|F_A\|_{L^{\sharp,2}(X,g)} \leq c$ for
all $[A] \in M_\kappa^w(X,g)$.
\end{lem}

\begin{proof}
  By hypothesis and the usual Chern-Weil identity, we have
  $\|F_A\|_{L^2(X,g)} = \sqrt{8\pi^2\kappa}$, for any $[A] \in
  M_\kappa^w(X,g)$. The moduli space $M_\kappa^w(X,g)$ has a finite
  covering by Uhlenbeck neighborhoods of the type given in Definition
  \ref{defn:UhlenbeckNeighborhood}, so it suffices to consider one such
  neighborhood $\sU$. Following the notation of Definition
  \ref{defn:UhlenbeckNeighborhood}, the compact subset $Y := X-\cup_{i=1}^m
  B(x_i,\delta)$ of $X-\bx$ has a finite covering by $N$ balls $B(y_j,r)$
  such that $B(y_j,4r) \subset Y' := X-\cup_{i=1}^m B(x_i,\half\delta)$,
  where $N = N(X,g,\kappa,\delta)$ and $r=r(\sU)$ is chosen so
  $\|F_A\|_{L^2(B(y_j,4r),g)}$ is small enough that the usual elliptic
  estimates hold for Yang-Mills connections (see, for example, \cite[Theorem
  2.1]{UhlRem}, \cite[Theorem 2.3.8]{DK}, \cite[Proposition 8.3]{FU}, or
  Theorem \ref{thm:CoulombASDBall} in the present article):
\begin{equation}
\label{eq:LocalConnEstBj}
\|A-\Gamma_j\|_{L^3_{1,\Gamma}(B(y_j,2r),g)} 
\leq 
c\|F_A\|_{L^2(B(y_j,4r),g)},  
\end{equation}
where $c=c(g,r)$ and $\Gamma_j$ is a flat connection over $B(y_j,3r)$ such that
$d_{\Gamma_j}^*(A-\Gamma_j) = 0$. Thus,
\begin{align*}
\|F_A\|_{L^\sharp(B(y_j,r),g)} 
&\leq 
c\|F_A\|_{L^3(B(y_j,2r),g)}
\quad\text{(Sobolev embedding \cite[Lemma 4.1]{FeehanSlice})}
\\
&\leq 
c\|A-\Gamma_j\|_{L^3_{\Gamma_j,1}(B(y_j,2r),g)}
\\
&\leq 
c\|F_A\|_{L^2(B(y_j,4r),g)}
\quad\text{(by estimate \eqref{eq:LocalConnEstBj})}
\end{align*}
Hence, combining the preceding estimates over each ball $B(y_j,r)$
covering $Y$, and noting that $B(y_j,4r)$ is contained in $Y'$, 
for $j=1,\dots,N$, we see that 
$$
\|F_A\|_{L^\sharp(Y,g)} \leq c\|F_A\|_{L^2(Y',g)},
$$
where $c=c(g,\sU)$.

We now consider each of the balls $B(x_i,2\delta)$ in turn and
estimate $\|F_A\|_{L^\sharp(B(x_i,\delta),g)}$. The $L^\sharp$ norm on
two-forms is invariant under rescalings of the metric, so we may
assume without loss of generality that $\delta=1$.  We can now apply
the usual bubble-tree compactness argument \cite[\S
4]{FeehanGeometry}, \cite[\S 5]{TauFrame} to conclude that there is an
Uhlenbeck neighborhood $\sU_i$ of the type given in Definition
\ref{defn:UhlenbeckNeighborhood} with $X$ replaced by
$B(x_i,2\delta)$, noting that no bubbling occurs in the annuli
$\Omega(x_i,\threequarter\delta,2\delta)$.  
We repeat the preceding argument to bound the
$L^\sharp$ norm of $F_A$ on a precompact subset of $B(x_i,2\delta)$
away from the (at most
$\kappa_i$ and at least two if $\kappa_i\geq 2$) points of curvature
concentration. This process terminates after at most $\kappa$
rescalings and so yields the desired $L^\sharp$ bound on $F_A$.
\end{proof}

We use the preceding universal bound on $\|F_A\|_{L^{\sharp,2}(X,g)}$
in the next subsection.

\subsection{Estimate for the self-dual curvature of cut-off
  connections} 
\label{subsec:EstimateSDCurvCutoffConn}
An important application of our local estimates for anti-self-dual
connections is to bound the self-dual curvature of approximately
anti-self-dual connections produced by splicing.

We briefly recall some elementary facts concerning metrics,
connections, and rescaling. Suppose, to begin, that $g$ is a Euclidean
metric on $\RR^4$ and so there are coordinates $\{x^\mu\}$ on $\RR^4$
with respect to which $g$ has the standard form, $g(x) =
\delta_{\mu\nu}dx^\mu\otimes dx^\nu = dx^\mu\otimes dx^\mu$. Given a
positive constant $\lambda$, let $\tg = \lambda^{-2}g$ denote the
rescaled metric, so $\tg(x) = \lambda^{-2} dx^\mu\otimes dx^\mu$. With
respect to the coordinates $y=x/\lambda$, the metric $\tg$ becomes
$\tg(y) = dy^\mu\otimes dy^\mu$. Analogous remarks apply to metrics
which are not necessarily Euclidean.

Let $g$ be a metric on $\RR^4$ and let $\lambda[A,g]$ be the scale of a
unitary connection $A$ over $\RR^4$ with energy $8\pi^2\kappa$ and center of
mass at the origin,
$$
\lambda[A,g]^2
:=
\frac{1}{8\pi^2\kappa}
\int_{\RR^4}\dist_g(0,\cdot)^2|F_A|_g^2\,dV_g.
$$
With respect to the rescaled metric, $\tg=\lambda^{-2}g$, we see that
$$
\lambda[A,\tg]^2
:=
\frac{1}{8\pi^2\kappa}
\int_{\RR^4}\dist_\tg(0,\cdot)^2|F_A|_\tg^2\,dV_\tg = 1,
$$
since $\dist_\tg(0,\cdot) = \lambda\dist_g(0,\cdot)$ and
$|F_A|_\tg^2\,dV_\tg = |F_A|_g^2\,dV_g$.

\begin{prop}
\label{prop:CutoffConnSelfDualCurvEst}
If $1\leq p\leq \8$, there is a positive constant $c(g,p)$ with the following
significance.  Let $A$ be a spliced connection on $\fg_E$, 
produced by the splicing construction of \S \ref{sec:Splicing}
(see Definition \eqref{eq:SplicedConnection}). Then, there is a constant
$\lambda_0 = \lambda_0(g,\kappa)$ such that for all $\lambda <\lambda_0$,
where $\lambda := \max_{1\leq i\leq m}\lambda_i$, such that
\begin{align*}
\|F^{+,g}_A\|_{L^p(X,g)} 
&\leq c\lambda^{2/p}
+ c\|F_{A_0}^{+,g}\|_{L^p(X,g)} 
+ c\sum_{i=1}^m\|F^{+,g_0}_{A_i}\|_{L^p(S^4,g_0)}.
\end{align*}
If $A_0$ is $g$-anti-self-dual over $X$ and the $A_i$
are $g_0$-anti-self-dual over $S^4$, then
\begin{equation}
\label{eq:MainLpEstimateSelfDualCurv}
\|F^{+,g}_A\|_{L^{\sharp,2}(X,g)} \leq c\lambda.
\end{equation}
\end{prop}

\begin{pf}
Let $\eps$ be a positive constant which is less than or equal to the
constants of Lemmas \ref{lem:C0ASDFAFixedBall} and
\ref{lem:C0ASDFAFixedBallComplement}, and let $0<\varrho\leq 1$ be a
constant. We choose $\varrho$ small enough that the background connection,
$A_0$, satisfies $\sup_{x\in X}\|F_A\|_{L^2(B(x,\varrho))} < \eps$.  Also
require that $\varrho$ be small enough that $\varrho\sqrt{8\pi^2\kappa_i} <
\eps$, for all $i\geq 1$. Choose $\lambda_0$ small enough that
$\sqrt{\lambda_0}\leq\varrho/16$.

The Chebychev inequality (see Lemma \ref{lem:Chebychev}) implies that
$$
\|F_{A_i}\|_{L^2(\RR^4-B(0,1/\varrho),\lambda_i^{-2}\delta)}
=
\|F_{A_i}\|_{L^2(\RR^4-B(0,\lambda_i/\varrho),\delta)}
\leq
\varrho\sqrt{8\pi^2\kappa_i} < \eps,
$$
since $\varrho\sqrt{8\pi^2\kappa_i} < \eps$, by definition of
$\varrho$, and as
$$
\lambda[A_i,\delta]^2
=
\frac{1}{8\pi^2\kappa}
\int_{\RR^4}|x|^2|F_{A_i}|_{\delta}^2\,d^4x = \lambda_i^2.
$$
Recall from definition \eqref{eq:SplicedConnection} that
\begin{equation}
\label{eq:SplicedConnectionToo}
A
:= 
\begin{cases}
  A_0 &\text{on $X - \cup_{i=1}^m B(x_i,4\sqrt{\lambda_i})$},
  \\
  \Gamma + \chi_{x_i,4\sqrt{\lambda_i}}\sigma_i^*A_0 &\text{on
    $\Omega(x_i,2\sqrt{\lambda_i},4\sqrt{\lambda_i})$},
  \\
  \Gamma &\text{on
    $\Omega(x_i,\sqrt{\lambda_i}/2,2\sqrt{\lambda_i})$},
  \\
    \Gamma + (1-\chi_{x_i,\sqrt{\lambda_i}/2})
    (\varphi_n\circ\varphi_i^{-1})^*\tau_i^*A_i &\text{on
    $\Omega(x_i,\sqrt{\lambda_i}/4,\sqrt{\lambda_i}/2)$},
  \\
    (\varphi_n\circ\varphi_i^{-1})^*A_i 
    &\text{on $B(x_i,\sqrt{\lambda_i}/4)$},
\end{cases} 
\end{equation}
where the cut-off functions are defined in equation
\eqref{eq:ChiCutoffFunctionDefn}.  
Over the region $X - \cup_{i=1}^m B(x_i,4\sqrt{\lambda_i})$ we have
$F_A = F_{A_0}$, while over the annuli
$\Omega(x_i,\sqrt{\lambda_i}/2,2\sqrt{\lambda_i})$ we have $F_A =
0$, and over the balls $B(x_i,\sqrt{\lambda_i}/4)$ we have $F_A =
F_{A_i}$. Thus, it remains to consider the two annuli
where cutting off occurs and the balls $B(x_i,\sqrt{\lambda_i}/4)$.

The hypotheses
imply that $F_{A_0}^{+,g}=0$ on each ball $B(x_i,\varrho)$. 
Over the outer annulus $\Omega(x_i,2\sqrt{\lambda_i},4\sqrt{\lambda_i})$ 
we see that \eqref{eq:SplicedConnectionToo} yields
\begin{equation}
\label{eq:ExpandCurvOuterAnnulus}
F_A = \chi_iF_{A_0} + d\chi_i\wedge\si_i^*A_0 
+ (\chi_i^2-\chi_i)\si_i^*A_0\wedge \si_i^*A_0,
\end{equation}
where we set $\chi_i = \chi_{x_i,4\sqrt{\lambda_i}}$ for convenience.
According to Lemma \ref{lem:C0ASDFAFixedBall}, we have the $L^\8$
curvature estimate and by Lemma
\ref{lem:SobolevEstimateRadialGaugeConnOneForm}, since the connection
one-form $\sigma_i^*A_0$ is in radial gauge with respect to the point
$x_i$, the $L^p$ connection one-form estimate,
\begin{align*}
\|F_{A_0}\|_{L^\8(B(x_i,\varrho/2))}
&\leq c\varrho^{-2}\|F_{A_0}\|_{L^2(B(x_i,\varrho))},
\\
\|\sigma_i^*A_0\|_{L^p(B(x_i,4\sqrt{\lambda_i}))}
&\leq
c\lambda_i^{(2/p)-1/2}\|F_{A_0}\|_{L^2(B(x_i,8\sqrt{\lambda_i}))}.
\end{align*}
The lemmas are applicable because of the choice of $\eps$. Because
$$
\|F_{A_0}\|_{L^2(B(x_i,8\sqrt{\lambda_i}))} 
\leq
c\lambda_i\|F_{A_0}\|_{L^\8(B(x_i,8\sqrt{\lambda_i}))}, 
$$
we can combine the preceding two estimates to get the $L^p$ estimate, for
any $1\leq p\leq \8$,
\begin{equation}
\label{eq:LpEstimateBackroundConnOneFormBall}
\|\sigma_i^*A_0\|_{L^p(B(x_i,4\sqrt{\lambda_i}))}
\leq
c\lambda_i^{(2/p)+1/2}\varrho^{-2}\|F_{A_0}\|_{L^2(B(x_i,\varrho))},
\quad
\sqrt{\lambda_i} \leq \varrho/16.
\end{equation}
Given $1\leq p \leq \8$, choose $4/3\leq q\leq \8$ by setting $1/p=1/4+1/q$.
{}From the expression \eqref{eq:ExpandCurvOuterAnnulus} for $F_A$ over
the annulus $\Omega(x_i,2\sqrt{\lambda_i},4\sqrt{\lambda_i})$ and the integral
estimates \eqref{eq:LpEstimateBackroundConnOneFormBall} for
$\sigma_i^*A_0$ and that of Lemma \ref{lem:dBetaEst} for $d\chi_i$, we
obtain an estimate for $F_A^{+,g}$,
\begin{align}
\label{eq:LpEstimateSDCurvOuterAnnulus}
\|F_A^{+,g}\|_{L^p(\Omega(x_i,2\sqrt{\lambda_i},4\sqrt{\lambda_i}))} 
&\leq 
\|F_{A_0}^{+,g}\|_{L^p} 
+ \|d\chi_i\|_{L^4}
\|\si_i^*A_0\|_{L^q(\Omega(x_i,2\sqrt{\lambda_i},4\sqrt{\lambda_i}))}
\\
\notag
&\quad
+ \|\si_i^*A_0\|_{L^{2p}(\Omega(x_i,2\sqrt{\lambda_i},4\sqrt{\lambda_i}))}^2
\\
\notag
&\leq \|F_{A_0}^{+,g}\|_{L^p} 
+ c\varrho^{-2}(\lambda_i^{2/p} 
+ \lambda_i^{(2/p)+1})\|F_{A_0}\|_{L^2(B(x_i,\varrho))}.
\end{align}
The preceding derivation of an $L^p$ estimate for $F_A^{+,g}$ over
$\Omega(x_i,2\sqrt{\lambda_i},4\sqrt{\lambda_i})$ can easily be adapted to
the case of an $L^\sharp$ estimate, using Lemmas \ref{lem:C0ASDFAFixedBall}
and \ref{lem:dBetaEst} and the embedding $L^2_1\subset L^{2\sharp}$ of
\cite[Lemma 4.1]{FeehanSlice},
\begin{equation}
\label{eq:LSharpEstimateSDCurvOuterAnnulus}
\|F_A^{+,g}\|_{L^\sharp(\Omega(x_i,2\sqrt{\lambda_i},4\sqrt{\lambda_i}))} 
\leq 
\|F_{A_0}^{+,g}\|_{L^\sharp} 
+ c\varrho^{-2}(\lambda_i + \lambda_i^2)\|F_{A_0}\|_{L^2(B(x_i,\varrho))}.
\end{equation}
This completes the required estimates for $F_A^{+,g}$ over the outer
annuli.

Next we turn to the estimate of $F_A^{+,g}$ over the inner annulus
$\Omega(x_i,\sqrt{\lambda_i}/4,\sqrt{\lambda_i}/2)$. The hypotheses imply
that $F_{A_i}^{+,\delta}=0$ on $\RR^4$ or $F_{A_i}^{+,g_0}=0$ on $S^4$,
where $\delta$ is the Euclidean metric on $\RR^4$ and $g_0$ is the
round metric of radius one on $S^4$. This time, our
expression \eqref{eq:SplicedConnectionToo} for the connection $A$ yields
\begin{equation}
\label{eq:ExpandCurvInnerAnnulus}
F_A = \chi_iF(\tau_i^*A_i) +
 d\chi_i\wedge \tau_i^*A_i 
+ (\chi_i^2-\chi_i)\tau_i^*A_i 
\wedge \tau_i^*A_i,
\end{equation}
where we now abbreviate $\chi_i = 1-\chi_{x_i,\sqrt{\lambda_i}/2}$ and
$\tau_i^*A_i = (\varphi_n\circ\varphi_i^{-1})^*\tau_i^*A_i$. Let $\tg_i$
denote the round metric on $S^4$ induced by pulling back the rescaled metric
$\lambda_i^{-2}\delta$ from $\RR^4$ via the diffeomorphism
$\varphi_n^{-1}:S^4-\{s\}\to\RR^4$.

According to Lemma \ref{lem:C0ASDFAFixedBall}, we have the $L^\8$
curvature estimate and by Lemma
\ref{lem:C0ASDFAFixedBallComplement}, since the connection
one-form $\tau_i^*A_i$ is in radial gauge with respect to the south pole
$s$, the $L^p$ connection one-form estimate,
\begin{align*}
\|F_{A_i}\|_{L^\8(\varphi_s(B(\varrho/2)),\tg_i)}
&\leq c\varrho^{-2}\|F_{A_i}\|_{L^2(\varphi_s(B(\varrho)),\tg_i)},
\\
\|\tau_i^*A_i\|_{L^p(\varphi_n(\RR^4-B(\sqrt{\lambda_i}/4)),g_0)}
&\leq
c\lambda_i^{(2/p)-1/2}
\|F_{A_i}\|_{L^2(\varphi_n(\RR^4-B(x_i,\sqrt{\lambda_i}/8)),g_0)}.
\end{align*}
Again, the lemmas are applicable because of the choice of $\eps$. Since
\begin{align*}
\|F_{A_i}\|_{L^2(\varphi_s(B(\varrho)),\tg_i)}
&=
\|F_{A_i}\|_{L^2(\varphi_n(\RR^4-B(1/\varrho)),\tg_i)}
=
\|F_{A_i}\|_{L^2(\varphi_n(\RR^4-B(\lambda_i/\varrho)),g_0)},
\\
\|F_{A_i}\|_{L^2(\varphi_s(B(\sqrt{\lambda_i}/8)),\tg_i)} 
&\leq
c\lambda_i\|F_{A_i}\|_{L^\8(\varphi_s(B(\sqrt{\lambda_i}/8)),\tg_i)},
\\
\|F_{A_i}\|_{L^2(\varphi_n(\RR^4-B(\sqrt{\lambda_i}/8)),g_0)}
&=
\|F_{A_i}\|_{L^2(\varphi_n(\RR^4-B(1/8\sqrt{\lambda_i})),\tg_i)}
=
\|F_{A_i}\|_{L^2(\varphi_s(B(8\sqrt{\lambda_i})),\tg_i)},
\end{align*}
we can combine the two preceding sets of estimates to get the $L^p$
bound, for any $1\leq p\leq \8$,
\begin{equation}
\label{eq:LpEstimateSphereConnOneFormBallComp}
\|\tau_i^*A_i\|_{L^p(\varphi_n(\RR^4-B(\sqrt{\lambda_i}/4)),g_0)}
\leq
c\lambda_i^{(2/p)+1/2}\varrho^{-2}
\|F_{A_i}\|_{L^2(\varphi_s(B(\varrho)),\tg_i)},
\quad
\sqrt{\lambda_i} \leq \varrho/16.
\end{equation}
{}From the expression \eqref{eq:ExpandCurvInnerAnnulus} for $F_A$ over
the inner
annulus $\Omega(x_i,\sqrt{\lambda_i}/4,\sqrt{\lambda_i}/2)$ and the integral
estimates \eqref{eq:LpEstimateSphereConnOneFormBallComp} for
$\varphi_n^*\tau_i^*A_i$ 
and that of Lemma \ref{lem:dBetaEst} for $d\chi_i$, we
obtain an estimate for $F_A^{+,g}$,
\begin{align}
\label{eq:LpEstimateSDCurvInnerAnnulus}
&\|F_A^{+,g}\|_{L^p(\Omega(x_i,\sqrt{\lambda_i}/4,\sqrt{\lambda_i}/2))} 
\\
\notag
&\leq 
\|F^{+,g}_{A_i}\|_{L^p(\Omega(x_i,\sqrt{\lambda_i}/4,\sqrt{\lambda_i}/2),g)} 
+ \|d\chi_i\|_{L^4}
\|\tau_i^*A_i\|_{L^q(\Omega(n,\sqrt{\lambda_i}/4,\sqrt{\lambda_i}/2),g_0)}
\\
\notag
&\quad
+ \|\tau_i^*A_i
\|_{L^{2p}(\Omega(n,\sqrt{\lambda_i}/4,\sqrt{\lambda_i}/2),g_0)}^2
\\
\notag
&\leq 
\|F^{+,g}_{A_i}\|_{L^p(\Omega(x_i,\sqrt{\lambda_i}/4,\sqrt{\lambda_i}/2),g)} 
+ c\varrho^{-2}(\lambda_i^{2/p} 
+ \lambda_i^{(2/p)+1})\|F_{A_i}\|_{L^2(\Omega(n,\lambda_i/\varrho,\8),g_0)}.
\end{align}
Again, the preceding derivation of an $L^p$ estimate for $F_A^{+,g}$ over
$\Omega(x_i,\sqrt{\lambda_i}/4,\sqrt{\lambda_i}/2)$ can easily be adapted to
the case of an $L^\sharp$ estimate, using Lemmas \ref{lem:C0ASDFAFixedBall}
and \ref{lem:dBetaEst} and the embedding $L^2_1\subset L^{2\sharp}$ of
\cite[Lemma 4.1]{FeehanSlice},
\begin{equation}
\label{eq:LSharpEstimateSDCurvInnerAnnulus}
\begin{aligned}
\|F_A^{+,g}\|_{L^\sharp(\Omega(x_i,\sqrt{\lambda_i}/4,\sqrt{\lambda_i}/2))} 
&\leq 
\|F^{+,g}_{A_i}
\|_{L^\sharp(\Omega(x_i,\sqrt{\lambda_i}/4,\sqrt{\lambda_i}/2),g)} 
\\
&\quad + c\varrho^{-2}(\lambda_i 
+ \lambda_i^2)\|F_{A_i}\|_{L^2(\Omega(n,\lambda_i/\varrho,\8),g_0)}.
\end{aligned}
\end{equation}
This completes the required estimate for $F_A^{+,g}$ over the inner
annulus, modulo the $L^p$ and $L^\sharp$ estimates for $F^{+,g}_{A_i}$ over
$\Omega(x_i,\sqrt{\lambda_i}/4,\sqrt{\lambda_i}/2)$, which we consider below.

Finally, we turn to the balls $B(x_i,\sqrt{\lambda_i}/4)$. We use geodesic
normal coordinates centered at $x_i$ to attach the connections $A_i$,
defined on $S^4 = (TX)_{x_i}\cup\{\8\}$ with its standard metric,
$g_0$. Just as in \cite[Equation (8.20)]{TauSelfDual}, using
$F_{A_i}^{+,g_0}=0$ and
$$
F_{A_i}^{+,g} = \half(1+*_g)F_{A_i} = \half(*_g-*_{g_0})F_{A_i}
$$
gives, via $g_{\mu\nu} = \delta_{\mu\nu} + O(r^2)$,
\begin{align*}
\|F_{A_i}^{+,g}\|_{L^\sharp(B(x_i,\sqrt{\lambda_i}/4))}
&\leq 
c\|*_g-*_{g_0}\|\cdot\|F_{A_i}\||_{L^\sharp(B(x_i,\sqrt{\lambda_i}/4))}
&\leq 
c\lambda_i\|F_{A_i}\||_{L^\sharp(B(x_i,\sqrt{\lambda_i}/4))}.
\end{align*}
According to Lemma \ref{lem:UniversalLSharp2BoundFA}, we have
$$
\|F_{A_i}\|_{L^\sharp(S^4,g_0)}
\leq 
c\|F_{A_i}\|_{L^2(S^4,g_0)}.
$$
In the same vein, for the $L^p$ estimate for
$F^{+,g}_{A_i}$ over $\Omega(x_i,\sqrt{\lambda_i}/4,\sqrt{\lambda_i}/2)$,
we have
$$
\|F_{A_i}^{+,g}\|_{L^p(\Omega(x_i,\sqrt{\lambda_i}/4,\sqrt{\lambda_i}/2),g)} 
\leq
c\lambda_i
\|F_{A_i}\|_{L^p(\Omega(x_i,\sqrt{\lambda_i}/4,\sqrt{\lambda_i}/2),g)}.
$$ 
The required bounds for $F^{+,g}_A$ now follow by combining the estimates
\eqref{eq:LpEstimateSDCurvOuterAnnulus} and
\eqref{eq:LpEstimateSDCurvInnerAnnulus} for $F_A^{+,g}$ over the annuli
where cutting off occurs, together with the estimates for
$F_A^{+,g}$ on the complement of these annuli.
\end{pf}


\section{Eigenvalue estimates for Laplacians on self-dual two-forms}
\label{sec:Eigenvalue2}
We need to analyze the behavior of the small eigenvalues of the Laplacian
$d_A^+d_A^{+,*}$ on $L^2(X,\Lambda^+\otimes\fg_E)$ as the point $[A]$ in
$\sB_E(X)$ converges, in the Uhlenbeck topology, to a point $([A_0],\bx)$
in $\sB_{E_0}(X)\times\Sym^\ell(X)$.  Our strategy is modelled on the
method used by Taubes to describe the small eigenvalues of $d_A^+d_A^{+,*}$
\cite[Lemma 4.8]{TauIndef}, where $A$ is an approximately anti-self-dual
connection constructed by splicing $|\bx|\leq \ell$ one-instantons onto the
background product connection $\Ga$. However, the upper and lower bounds
for the small eigenvalues we obtain here are more precise and more general
than those of
\cite{TauIndef} or the analogous arguments are used by Mrowka in the proof of
Theorem 6.1.0.0.4 \cite{MrowkaThesis} in order to establish his `main
eigenvalue estimate'. Our derivation in \cite{FLKM2} of the
eigenvalue bounds for the Laplacian $d^*_Ad_A$ on $L^2(X,\fg_E)$ is
closely modelled on the arguments given here for the Laplacian
$d_A^+d_A^{+,*}$ on $L^2(X,\Lambda^+\otimes\fg_E)$.

In the present article, our main application of the eigenvalue bounds of
this section will be to deriving an estimate for the partial right inverse
$P_{A,\mu} = d_A^{+,*}(d_A^+d_A^{+,*})^{-1}\Pi_{A,\mu}^\perp$ for
$d_A^+$, as given in Proposition \ref{prop:EstPA0mu}.  We have the
following analogue of Proposition 4.6 in \cite{TauIndef}.

\begin{thm}
\label{thm:H2SmallEval}
Let $M$ be a positive constant and let $(X,g)$ be a closed, oriented,
Riemannian four-manifold.  Then there are positive constants $c=c(g)$,
$\eps=\eps(g)$, $C = C(g,M)$, and small enough $\lambda_0 = \lambda_0(C)$
such that for all $\lambda \in (0,\lambda_0]$, the following holds.  Let
$E$, $E_0$ be Hermitian, rank-two vector bundles over $X$ with $\det E =
\det E_0$, $c_2(E)=c_2(E_0)+\ell$, for some integer $\ell\geq 0$, and
suppose $-\quarter p_1(\fg_E)\leq M$. Let $[A_0]\subset\sB_{E_0}(X)$ be a
point such that $\|F_{A_0}^{+,g}\|_{L^{\sharp,2}(X,g)} < \eps$.  Suppose
$\bx\in\Sym^\ell(X)$ and define a precompact open subset of $X\less\bx$ by
setting $U:= X-\cup_{x\in\bx}\bar B(x,4\sqrt{\lambda})$.  Let $A$ be a
connection representing a point in the open subset of $\sB_E(X)$ defined by
the following constraints:
\begin{itemize}
\item $\|F_A^{+,g}\|_{L^{\sharp,2}(X,g)} < \eps$,
\item $\|A-A_0\|_{L^2(U)} < c\sqrt{\lambda}$ and
$\|A-A_0\|_{L^2_{1,A_0}(U)} < M$, (noting that
$E|_{X\less\bx} \cong E_0|_{X\less\bx}$).
\end{itemize}
Let $\nu_2[A_0]$ be the {\em least positive eigenvalue\/} of the Laplacian
$d_{A_0}^+d_{A_0}^{+,*}$ on $L^2(\Lambda^+\otimes\fg_{E_0})$ and suppose
$\max\{\nu_2[A_0],\nu_2[A_0]^{-1}\}
\leq M$. Assume that the kernel of
$d_{A_0}^+d_{A_0}^{+,*}$ has dimension $n\leq M$. Let $\{\mu_i[A]\}_{i\geq
1}$ denote the eigenvalues, repeated according to their multiplicity, of
the Laplacian $d_A^+d_A^{+,*}$ on $L^2(\Lambda^+\otimes\fg_E)$. Then
\begin{align*}
\mu_i[A] &\leq
\begin{cases}
C\lambda, &\text{if }i=1,\dots,n, \\
\nu_2[A_0] + C\sqrt{\lambda}, &\text{if }i=n+1,
\end{cases} 
\\
\mu_i[A] &\geq \nu_2[A_0]- C\sqrt{\la},  \quad \text{if }i\geq n+1.
\end{align*}
\end{thm}

\begin{rmk}
\begin{enumerate}
\item
Recall that by \cite[\S 2.1.4]{DK}
$$
\kappa 
= -\quarter p_1(\fg_E) 
= \frac{1}{8\pi^2}(\|F_A^-\|_{L^2}^2 - \|F_A^+\|_{L^2}^2)
= \frac{1}{8\pi^2}(\|F_A\|_{L^2}^2 - 2\|F_A^+\|_{L^2}^2),
$$
and therefore,
$$
\|F_A\|_{L^2}^2 \leq 2\eps^2 + 8\pi^2\kappa.
$$
So, we have a universal energy bound for all points in
$\sB_E(X)$ with $\|F_A^{+,g}\|_{L^2(X,g)} < \eps$.
\item
If a point $[A]$ in $\sB_E(X)$ lies in an Uhlenbeck neighborhood (see
Definition \ref{defn:UhlenbeckNeighborhood}) of a point $([A_0],\bx)$ in
$\sB_{E_0}(X)\times\Sym^\ell(X)$, then the hypotheses of Theorem
\ref{thm:H2SmallEval} are fulfilled.
\end{enumerate} 
\end{rmk}

\begin{cor}
\label{cor:H2SmallEval}
Let $(X,g)$ be a closed, oriented, Riemannian four-manifold. Suppose
$[A_\alpha]$ is a sequence of points in $\sB_E(X)$ which converges, in the
Uhlenbeck topology, to a point
$([A_0],\bx)\in\sB_{E_0}(X)\times\Sym^\ell(X)$ and suppose the kernel of
$d_{A_0}^+d_{A_0}^{+,*}$ has dimension $n$.  Then the first $n$ eigenvalues
(counted with multiplicity) of $d_{A_\alpha}^+d_{A_\alpha}^{+,*}$ converge
to zero, while the $(n+1)$-st eigenvalue of
$d_{A_\alpha}^+d_{A_\alpha}^{+,*}$ converges to $\nu_2[A_0]$, the least
positive eigenvalue of $d_{A_0}^+d_{A_0}^{+,*}$.
\end{cor}

\begin{rmk}
\label{rmk:H2SmallEvalFamilies}
The statements of Theorem \ref{thm:H2SmallEval} and Corollary
\ref{cor:H2SmallEval} extend without difficulty from the case of a single,
fixed metric $g$ and connection $[A_0]$ to the case of a family
$\sU_0\subset \sB_E(X)\times\Met(X)$ satisfying the constraints
\begin{itemize}
\item
$\Spec d_{A_0}^{+,g}d_{A_0}^{+,g,*} \subset [0,t)\cup (2\mu,\8)$,
for some positive constants $t$, $\mu$ such that $t \leq \mu/2$ and
all $[A_0,g]\in\sU_0$,
\item
The number of eigenvalues, counted with multiplicity, of
$d_{A_0}^{+,g}d_{A_0}^{+,g,*}$ which are less than $\mu$ is constant,
say $n$, as $[A_0,g]$ varies over $\sU_0$, and
\item
$\|F_{A_0}^{+,g}\|_{L^{\sharp,2}(X,g)} < \eps$.
\end{itemize}
Thus, for example, in Corollary
\ref{cor:H2SmallEval} we can replace $g$ by a sequence $\{g_\alpha\}$ of smooth
Riemannian metrics on $X$ converging in $C^r$, $r\geq 1$, to a metric
$g_0$. In the statements of Theorem \ref{thm:H2SmallEval} and Corollary
\ref{cor:H2SmallEval}, we replace $\nu_2[A_0]$ by $\mu$. Then, in the
proofs below, we replace $\nu_2[A_0]$ by $\mu$,
$\Ker d_{A_0}^{+,g}d_{A_0}^{+,g,*}$ by
$\Ran\Pi_{A_0,\mu}$ and $\Pi_{A_0}$ by $\Pi_{A_0,\mu}$, the projection from
$L^2(\Lambda^+\otimes\fg_{E_0})$ onto the span of the eigenvectors of
$d_{A_0}^{+,g}d_{A_0}^{+,g,*}$ with eigenvalue less than $\mu$. Thus, in
the statement of Theorem \ref{thm:H2SmallEval}, the upper bounds are
replaced by $C(t+\lambda)$ and $\mu + C(\sqrt{t}+\sqrt{\lambda})$,
respectively, while the lower bound is replaced by
$\mu-C(\sqrt{t}+\sqrt{\lambda})$, for all $t < t_0(C)$ and $\lambda <
\lambda_0(C)$. Similar modifications apply to Corollary
\ref{cor:H2SmallEval}, noting that $\lambda$ is replaced by a sequence
$\{\lambda_\alpha\}$ converging to zero.
\end{rmk}

The proof of Theorem \ref{thm:H2SmallEval} takes up the remainder
of this section.  While the argument is lengthy, the basic strategy is
reasonably straightforward:
\begin{itemize}
\item
Find upper bounds for the first $n+1$ eigenvalues of $d_A^+d_A^{+,*}$;
\item
Construct an approximately $L^2$-orthonormal basis for $\Ker
d_{A_0}^+d_{A_0}^{+,*}$ by cutting off the eigenvectors corresponding to
the first $n$ eigenvalues of $d_A^+d_A^{+,*}$;
\item
Compute a lower bound for the $(n+1)$-st eigenvalue of
$d_A^+d_A^{+,*}$. 
\end{itemize}
We begin in the next subsection by disposing of a few technical lemmas.

\subsection{Technical preliminaries}
\label{subsec:Eigen2TechPrelim}
We begin with some elementary lemmas before proceeding to 
the different stages of the proof of Theorem
\ref{thm:H2SmallEval}.  The following Kato-Sobolev inequality is
well-known; see Lemma 4.6 in \cite{TauSelfDual} for the case $p=2$.

\begin{lem}
\label{lem:Kato}
Let $X$ be a $C^\8$, closed, oriented, Riemannian, four-manifold 
and let $1\le p<4$. Then there is a positive constant $c$ with the
following significance.  Let $A$ be an orthogonal $L^p_1$ connection
on a Riemannian vector bundle $V$ over $X$. Then, for any $a\in
L^p_1(V)$, we have
$$
\|a\|_{L^q(X)} \le c\|a\|_{L^p_{1,A}(X)},
$$
when $1/p=1/q+1/4$.
\end{lem}

\begin{pf}
The estimate follows from the pointwise Kato inequality $|\cov|a||\le
|\cov_Aa|$ (see \cite[Equation (6.20)]{FU} and the Sobolev inequality
$$
\|f\|_{L^q} \le c(\|\cov f\|_{L^q} + \|f\|_p),
$$
for any $f\in L^p_1(X)$ \cite[Theorem 5.4]{Adams}.
\end{pf}

We shall need the Bochner-Weitzenb\"ock formulas,
\begin{align}
d_Ad_A^* + 2d_A^{+,*}d_A^+ 
&= \cov_A^*\cov_A + \{\Ric,\cdot\}
- 2\{F_A^-,\cdot\}, \label{eq:BW1} \\
2d_A^+d_A^{+,*}& = \cov_A^*\cov_A - 2\{W^+,\cdot\} + \frac{R}{3} +
\{F_A^+,\cdot\}, \label{eq:BW+}
\end{align}
for the Laplacians on $\Gamma(\Lambda^1\otimes\fg_E)$ and
$\Gamma(\Lambda^+\otimes\fg_E)$, respectively; see \cite[p. 94]{FU}.
We give $L^2_{1,A}$ estimates for $v$ in terms of the $L^2$ norm of
$d_A^{+,*}$; similar estimates are given in \cite[Lemma
5.2]{TauSelfDual} and in \cite[Appendix A]{TauIndef}.

\begin{lem}
\label{lem:L21AEstv}
Let $X$ be a $C^\8$, closed, oriented, Riemannian four-manifold.  Then
there are positive constants $c$ and $\eps=\eps(c)$
with the following significance. Let $A$
be an $L^2_4$ orthogonal connection on an $\SO(3)$ bundle $\fg_E$ over
$X$ such that $\|F_A^+\|_{L^2}<\eps$.
Then, for all $v\in L^{4/3}_2(\Lambda^+\otimes\fg_E)$,
\begin{align}
\tag{1}
\|v\|_{L^2_{1,A}(X)} 
&\le 
c(1+\|F_A^+\|_{L^2(X)})^{1/2}(\|d_A^{+,*}v\|_{L^2(X)}+\|v\|_{L^2(X)}),
\\
\tag{2}
\|d_A^{+,*}v\|_{L^2(X)} 
&\leq 
(1+\|F_A^+\|_{L^2(X)})^{1/2}\|d_A^+d_A^{+,*}v\|_{L^{4/3}(X)}.
\end{align}
\end{lem}

\begin{pf}
  The Bochner-Weitzenbock formula \eqref{eq:BW+} and integration
  by parts yields
\begin{align*}
\|\cov_Av\|_{L^2}^2 &= (\cov_A^*\cov_Av,v)_{L^2} \\
&= 2(d_A^+d_A^{+,*}v,v)_{L^2} + (\{W^+,v\},v)_{L^2} 
-\third(Rv,v)_{L^2} - (\{F_A^+,v\},v)_{L^2} \\
&\le 2\|d_A^{+,*}v\|_{L^2}^2 + c\|v\|_{L^2}^2
+ c\|F_A^+\|_{L^2}\|v\|_{L^4}^2.
\end{align*}
Apply Kato's inequality (Lemma \ref{lem:Kato}) to give $\|v\|_{L^4} \leq
c\|v\|_{L^2_{1,A}}$ and rearrangement with $\eps$ chosen small enough that
$c\|F_A^+\|_{L^2} \leq 1/2$ to give the first assertion.
The second assertion follows by
integrating by parts, H\"older's inequality, and Kato's inequality
(Lemma \ref{lem:Kato}),
$$
\|d_A^{+,*}v\|_{L^2}^2 
=
(d_A^+d_A^{+,*}v,v)_{L^2}
\leq
\|d_A^+d_A^{+,*}v\|_{L^{4/3}}\|v\|_{L^4}
\leq
c\|d_A^+d_A^{+,*}v\|_{L^{4/3}}\|v\|_{L^2_{1,A}},
$$
together with the first inequality.
\end{pf}

The following consequence of \cite[Lemma 5.9]{FeehanSlice} is the
principal tool we use to obtain the required eigenvector estimates.

\begin{lem}
\label{lem:LInftyEstDe2AEvec}
Let $X$ be a $C^\8$, closed, oriented, Riemannian four-manifold.  Then
there are positive constants $c$ and $\eps=\eps(c)$ with the following
significance. Let $A$ be an $L^2_4$ orthogonal connection on an $\SO(3)$
bundle $\fg_E$ over $X$ such that $\|F_A^+\|_{L^2}<\eps$.  Let
$\eta=\eta(A)\in L^2_3(\Lambda^+\otimes\fg_E)$ be an eigenvector of
$d_A^+d_A^{+,*}$ with eigenvalue $\mu[A]$. Then
$$
\|\eta\|_{L^\8(X)} + \|\eta\|_{L^2_{2,A}(X)}  
\le C[A]\|\eta\|_{L^2(X)},
$$
where $C[A] := C([A],\mu[A])$ with $\mu=\mu[A]$ in
$$
C([A],\mu)
:=
c(1+\|F_A\|_{L^2(X)})(1+\|F_A^+\|_{L^2(X)})^{1/2}
(1+\mu)(1+\sqrt{\mu}).
$$
More generally, if $\eta\in L^2(\Lambda^+\otimes\fg_E)$ lies in the
span of the eigenvectors of $d_A^+d_A^{+,*}$ with eigenvalues less
than or equal to a positive constant $\mu$, then
$$
\|\eta\|_{L^\8(X)} + \|\eta\|_{L^2_{2,A}(X)}
\le 
C([A],\mu)\|\eta\|_{L^2(X)}.
$$
\end{lem}

\begin{pf}
By \cite[Lemma 5.9]{FeehanSlice} (or Lemma
\ref{lem:LinftyL22CovLapEstv} here),
\begin{align*}
&\|\eta\|_{L^\8} + \|\eta\|_{L^2_{2,A}(X)}
\\
&\quad\le 
c(1+\|F_A\|_{L^2})(\|\eta\|_{L^2} + \|d_A^+d_A^{+,*}\eta\|_{L^\sharp}) 
\\
&\quad = c(1+\|F_A\|_{L^2})(1+\mu[A])\|\eta\|_{L^{\sharp,2}},
\end{align*}
using $d_A^+d_A^{+,*}\eta = \mu[A]\eta$ in the second line.
On the other hand, 
\begin{align*}
\|\eta\|_{L^{\sharp,2}} 
&\le c\|\eta\|_{L^4} \quad\text{(by \cite[Lemma 4.1]{FeehanSlice})}\\
&\le c\|\eta\|_{L^2_{1,A}} \quad\text{(by Lemma \ref{lem:Kato})}\\
&\le c(1+\|F_A^+\|_{L^2})^{1/2}
(\|\eta\|_{L^2} + \|d_A^{+,*}\eta\|_{L^2}) 
\quad\text{(by Lemma \ref{lem:L21AEstv})} \\
&= c(1+\|F_A^+\|_{L^2})^{1/2}
(1+\sqrt{\mu[A]})\|\eta\|_{L^2},
\end{align*}
using $\|d_A^{+,*}\eta\|_{L^2} = 
\sqrt{\mu[A]}\|\eta\|_{L^2}$ to obtain the final equality.
Combining the preceding estimates completes the proof of the first inequality.

Next suppose that $\eta\in L^2(\Lambda^+\otimes\fg_E)$ lies in the
span of the $n$ $L^2$-orthonormal eigenvectors $\eta_i(A)$ of
$d_A^+d_A^{+,*}$ with eigenvalues $\mu_i[A] \leq \mu$. Since $\eta =
\sum_{i=1}^n(\eta,\eta_i)_{L^2}\eta_i$
and denoting $C'=(1+\|F_A\|_{L^2})(1+\|F_A^+\|_{L^2})^{1/2}$, we have
\begin{align*}
\|\eta\|_{L^\8}
&\leq
\sum_{i=1}^n|(\eta,\eta_i)_{L^2}|\cdot\|\eta_i\|_{L^\8}
\\
&\leq 
C'\sum_{i=1}^n|(\eta,\eta_i)_{L^2}|
\cdot(1+\mu_i)(1+\sqrt{\mu_i})\|\eta_i\|_{L^2}
\quad\text{(by the special case)},
\\
&\leq C'(1+\mu)(1+\sqrt{\mu})
\left(\sum_{i=1}^n(\eta,\eta_i)_{L^2}^2\right)^{1/2}
\\
&= C'(1+\mu)(1+\sqrt{\mu})\|\eta\|_{L^2},
\end{align*}
where we use the inequality $(\sum_{i=1}^n a_i)^2 \leq
n\sum_{i=1}^na_i^2$, for $a_i\in\RR$, in the second last line. 
This completes the proof.
\end{pf}

The other basic tool we shall need is Courant's min-max principle for
eigenvalues (see, for example, \cite[pp. 16-17]{Chavel}).

\begin{lem}
\label{lem:MaxMin}
  Let $T$ be a non-negative, self-adjoint operator on a Hilbert space,
  $\fH$, such that $(\id+T)^{-1}$ is compact. Let
  $0\leq\mu_1\leq\mu_2\leq\cdots$ be the eigenvalues of $T$, repeated
  according to their multiplicity. Then, for
  any subspace $V_n\subset\fH$ with $\dim V_n = n$,
$$
\mu_{n+1} \geq \inf_{\eta\in V_n^\perp}\frac{(\eta,T\eta)}{\|\eta\|^2},
$$
where $V_n^\perp$ is the orthogonal complement of $V_n$ in $\fH$.
If $V_n$ is the span of the eigenvectors of $T$ corresponding to the
eigenvalues $\lambda_1,\dots,\lambda_n$, then equality holds above.
\end{lem}

\subsection{Good sequences of cut-off functions}
\label{subsec:GoodSeqCutoffFunctions}
In many of our arguments involving Uhlenbeck continuity of eigenvalues,
eigenvectors, gluing maps, gluing-map differentials and so on, we shall need
suitable sequences of cut-off functions to compare connections, and objects
depending on them, which live on different bundles but which become
isomorphic when restricted to an open subset $X\less\bx$, for some $\bx \in
\Sym^\bullet(X)$. Since the construction of such sequences will arise many
times, it is useful to extract the essential properties we shall need once
and for all in the following definition-lemma:

\begin{lem}
\label{lem:GoodCutoffSequence}
Let $X$ be a smooth Riemannian four-manifold and let
$\bx\in\Sym^\bullet(X)$. Then there is a sequence of cut-off functions
$\{\chi_n\}_{n=1}^\8$ such that the following hold:
\begin{itemize}
\item
$\cup_{n=1}^\8 U_n = X\less\bx$, where $U_n\Subset X\less\bx$ is
the open subset given by the interior of $\supp\chi_n$,
\item
$\chi_m\leq\chi_n$ for all $m\leq n$, so $U_m\subset U_n$ and $\supp\chi_m
\subset \supp\chi_n$,
\item 
$\supp\chi_n\Subset X\less\bx$,
\item
$\lim_{n\to\8}(\|\cov\chi_n\|_{L^2_1} + \|\cov^2\chi_n\|_{L^\sharp}) = 0$,
\end{itemize}
We call any sequence of cut-off functions
satisfying the above properties an {\em $\bx$-good sequence\/}.
\end{lem}

\begin{proof}
We may assume without loss that the multiset $\bx\in\Sym^\bullet(X)$ is
represented by a single point $x_0\in X$. Recall from Lemma
\ref{lem:dBetaEst} that we can define a cut-off function
$\beta_{N,\lambda}:X\to\RR$ with the property that
$$
\beta_{N,\lambda}(x)
=
\begin{cases}
1 &\text{if }\dist(x,x_0)\geq \half\sqrt{\lambda}, \\
0 &\text{if }\dist(x,x_0)\leq N^{-1}\sqrt{\lambda}, 
\end{cases}
$$
and satisfying the estimates
\begin{align*}
\|\cov\beta_{N,\lambda}\|_{L^4} &\le c(\log N)^{-3/4}, 
\\
\|\cov^2\beta_{N,\lambda}\|_{L^2} &\le c(\log N)^{-1/2}, 
\\
\|\cov\beta_{N,\lambda}\|_{L^{2\sharp}} &\le c(\log N)^{-1/2}, 
\\
\|\cov\beta_{N,\lambda}\|_{L^2} &\le c(\log N)^{-1}\sqrt{\lambda}.
\end{align*}
If we choose $N=n$ and $\lambda = n^{-4}$, for example, so $\chi_n :=
\beta_{n,n^{-4}}$ and $U_n = X\less \bar B(x_0,\half n^{-2})$, then it is 
immediate that the sequences $\{\chi_n\}_{n=1}^\8$ and $\{U_n\}_{n=1}^\8$
have the desired properties.
\end{proof}

\subsection{An upper bound for the small eigenvalues}
\label{subsec:Eigen2UpperBoundSmall}
We turn to the proof of Theorem \ref{thm:H2SmallEval}.  We define an
``$\bx$-good'' cut-off function $\chi_0$ on $X$ 
(see Lemma \ref{lem:GoodCutoffSequence}) by setting $\chi_0 =
\chi_{x_i,8\sqrt{\lambda_i}}$ (see equation \eqref{eq:ChiCutoffFunctionDefn}
for its definition) on each ball $B(x_i,8\sqrt{\lambda_i})$, so $\chi_0$ is
equal to zero on the balls $B(x_i,4\sqrt{\la_i})$, and setting $\chi_0=1$
on the complement of these balls, $X-\cup_{i=1}^m
B(x_i,8\sqrt{\lambda_i})$. We take
$$
U := X -\bigcup_{i=1}^mB(x_i,4\sqrt{\lambda_i}),
$$
so $U\Subset X\less\bx$ and $U$ is the interior of $\supp\chi_0$.

Though it suffices to consider sequences of $\bx$-good cut-off functions
only satisfying the broad properties of Lemma \ref{lem:GoodCutoffSequence},
we assume for simplicity that $\chi_0$ is defined as in the preceding
paragraph and observe that $\supp\chi_0\Subset X\less\bx$. 
We may assume without loss that $E|_{X\less\bx} = E_0|_{X\less\bx}$ and
thus, over $X\less\bx$, we may write
\begin{equation}
\label{eq:DiffConn}
a := A - A_0 \in L^2_{k,\loc}(X\less\bx,\Lambda^1\otimes\fg_{E_0}),
\end{equation}
where the section $a$ may also be viewed as an element of
$L^2_{k,\loc}(X\less\bx,\Lambda^1\otimes\fg_E)$. Our initial task is to
find an upper bound for the first $n+1$ eigenvalues of $d_A^+d_A^{+,*}$.

\begin{lem}
\label{lem:H2SmallEvalUpperBound}
Continue the hypotheses and notation of Proposition
\ref{thm:H2SmallEval}.  Let $\nu_2[A_0]$ be the least positive
eigenvalue of $d_{A_0}^+d_{A_0}^{+,*}$ and let
$\mu_1[A]\le\cdots\le\mu_{n+1}[A]$ be the first $n+1$ eigenvalues of
$d_A^+d_A^{+,*}$. Then there are positive constants 
$$
C = C(\|F_{A_0}\|_{L^2},\nu_2[A_0]) 
\quad\text{and}\quad
\lambda_0 = \lambda_0(C)
$$ 
such that for all $\lambda \in (0,\lambda_0]$,
\begin{equation}
\mu_i[A] \le
\begin{cases}
C\la, &\text{if }i=1,\dots,n, \\
\nu_2[A_0] + C\sqrt{\la}, &\text{if }i=n+1.
\end{cases} 
\label{eq:H2DeAEvalEst}
\end{equation}
\end{lem}

\begin{proof}
Let $\{\xi_i(A_0)\}_{i=1}^{n}\subset L^2(\Lambda^+\otimes\fg_{E_0})$ be an 
$L^2$-orthonormal basis for $\Ker d_{A_0}^+d_{A_0}^{+,*}$ and let
$\xi_{n+1}[A_0]\in(\Ker d_{A_0}^+d_{A_0}^{+,*})^\perp$ be an $L^2$ unit
eigenvector corresponding to the least positive eigenvalue $\nu_2[A_0]$ of
$d_{A_0}^+d_{A_0}^{+,*}$. 
We begin by constructing a set of approximate eigenvectors,
$\{\xi_i'\}_{i=1}^{n+1}$, for the first $n+1$ eigenvalues of $d_A^+d_A^{+,*}$
by setting $\xi_i':=\chi_0\xi_i$ for each $i\in\{1,\dots,n+1\}$.  

\begin{claim}
\label{claim:H2DeAEvectorEst}
Continue the above notation. Then there are positive constants 
$$
C = C(\|F_{A_0}\|_{L^2},\nu_2[A_0])
\quad\text{and}\quad
\lambda_0 = \lambda_0(C)
$$ 
such that for all $\lambda\in (0,\lambda_0]$,
\begin{align}
\tag{1}
\|d_A^{+,*}\xi_i'\|_{L^2(X)} &\le 
\begin{cases}
C\sqrt{\la}, &\text{if }1\le i\le n, \\
\sqrt{\nu_2[A_0]}+C\sqrt{\la}, &\text{if }i= n+1,
\end{cases}\\
\tag{2}
\left|(\xi_i',\xi_j')_{L^2(X)} - \de_{ij}\right|
&\le C\la^2, \quad 1\le i,j\le n+1.
\end{align}
\end{claim}

\begin{proof}
To verify the first assertion in
Claim \ref{claim:H2DeAEvectorEst}, note that $A=A_0+a$ on $X\less\bx$ and 
$d_A^{+,*} = *d_A$ on $\Gamma(\Lambda^+\otimes\fg_E)$, so
\begin{equation}
\begin{aligned}
\label{eq:ExpanddA+*OnCutoffEvectorXi}
d_A^{+,*}\xi_i' 
&= 
*d_A\xi_i'
= 
*d_{A_0}(\chi_0\xi_i) + *(a\wedge\chi_0\xi_i) 
\\
&= 
*(d\chi_0\wedge\xi_i) + \chi_0 d_{A_0}^{+,*}\xi_i + *(\chi_0a\wedge\xi_i),
\end{aligned}
\end{equation}
while
$$
\|d_{A_0}^{+,*}\xi_i\|_{L^2} 
=
\begin{cases}
0, &\text{if $1\leq i\leq n$},
\\
\sqrt{\nu_2[A_0]}\|\xi_i\|_{L^2}, &\text{if $i = n+1$}.
\end{cases}
$$
By Lemma \ref{lem:dBetaEst} (yielding $\|d\chi_0\|_{L^2} \leq
c\sqrt{\lambda}$), Lemma \ref{lem:LInftyEstDe2AEvec}, the preceding
eigenvalue estimate, and equation
\eqref{eq:ExpanddA+*OnCutoffEvectorXi}, 
\begin{align*}
\|d_A^{+,*}\xi_i'\|_{L^2} 
&\le (\|d\chi_0\|_{L^2}+ \|\chi_0a\|_{L^2})\|\xi_i\|_{L^\8} 
+ \|d_{A_0}^{+,*}\xi_i\|_{L^2} 
\\
&\le c\sqrt{\la}\|\xi_i\|_{L^\8} + \|d_{A_0}^{+,*}\xi_i\|_{L^2} 
\quad\text{(by Lemma \ref{lem:dBetaEst} and hypothesis)}
\\
&\le C\sqrt{\la}\|\xi_i\|_{L^2} 
+ \|d_{A_0}^{+,*}\xi_i\|_{L^2} 
\quad\text{(by Lemma \ref{lem:LInftyEstDe2AEvec})}
\\
&= C\sqrt{\la}\|\xi_i\|_{L^2} 
+ \sqrt{\nu_2[A_0]}\|\xi_i\|_{L^2}, 
\quad 1\le i\le n+1,
\end{align*}
and so the first assertion in Claim \ref{claim:H2DeAEvectorEst} follows.
Using $\xi_i' = \chi_0\xi = \xi_i + (\chi_0-1)\xi_i$, 
$$
(\xi_i',\xi_j')_{L^2}
= (\xi_i,\xi_j)_{L^2}
+ 2((\chi_0-1)\xi_i,\xi_j)_{L^2}
+ ((\chi_0-1)\xi_i,(\chi_0-1)\xi_j)_{L^2}.
$$
But $(\xi_i,\xi_j)_{L^2}=\de_{ij}$, while
the second two terms on the right above are bounded by
\begin{align*}
&\left|2((\chi_0-1)\xi_i,\xi_j)_{L^2}
+ ((\chi_0-1)\xi_i,(\chi_0-1)\xi_j)_{L^2}\right| \\
&\qquad\le c\la^2\|\xi_i\|_{L^\8}\|\xi_j\|_{L^\8} \le C\la^2,
\end{align*}
courtesy of Lemma \ref{lem:LInftyEstDe2AEvec}.
This gives the second assertion in Claim \ref{claim:H2DeAEvectorEst}.
\end{proof}

Now $\|\xi_i'\|_{L^2}\ge 1-C\la$ by the second assertion of
Claim \ref{claim:H2DeAEvectorEst} and so its first assertion yields
\begin{align*}
\|d_A^{+,*}\xi_i'\|_{L^2} 
&\le C\sqrt{\la}(1-C\la)^{-1}\|\xi_i'\|_{L^2}
\\
&\le C\sqrt{\la}\|\xi_i'\|_{L^2}, \quad 1\le i\le n
\quad\text{(choosing $\lambda<C/2$),}
\\
\|d_A^{+,*}\xi_{n+1}'\|_{L^2} 
&\le \left(\sqrt{\nu_2[A_0]}+C\sqrt{\la}\right)\|\xi_{n+1}'\|_{L^2}.
\end{align*}
Let $V_i:=[\xi_1',\dots,\xi_i']$ denote the span of the vectors
$\xi_1',\dots,\xi_i'$ in $L^2(\Lambda^+\otimes\fg_E)$ for
$i=1,\dots,n+1$.  By applying Gram-Schmidt orthonormalization to the basis
$\{\xi_1',\dots,\xi_i'\}$ for $V_i$, the preceding estimates imply that
\begin{equation}
\|d_A^{+,*}\eta\|_{L^2} \le 
\begin{cases}
C\sqrt{\la}\|\eta\|_{L^2}, &\text{if }\eta\in V_n, \\
\left(\sqrt{\nu_2[A_0]}+C\sqrt{\la}\right)\|\eta\|_{L^2}, &\text{if }\eta\in
V_{n+1}. 
\end{cases}
\label{eq:dA+*SmallEvalEst}
\end{equation}
Recall from the statement of Theorem \ref{thm:H2SmallEval} 
that $\mu_1[A]\le\cdots\le\mu_{n+1}[A]$
are the first $n+1$ eigenvalues of $d_A^+d_A^{+,*}$; let
$\eta_1(A),\dots,\eta_{n+1}(A)$ be the corresponding $L^2$-orthonormal
eigenvectors. Let $[\eta_1,\dots,\eta_i]$ denote the span of the
eigenvectors $\eta_1,\dots,\eta_i$ in $L^2(\Lambda^+\otimes\fg_E)$ and
let $[\eta_1,\dots,\eta_i]^\perp$ be its $L^2$-orthogonal complement
in $L^2(\Lambda^+\otimes\fg_E)$. Note that $\dim
([\eta_1,\dots,\eta_i]^\perp\cap V_{n}) \ge n-i$. By the estimates
\eqref{eq:dA+*SmallEvalEst} we have
$$
\mu_i[A] = \inf_{\eta\in[\eta_1,\dots,\eta_{i-1}]^\perp}
\frac{\|d_A^{+,*}\eta\|_{L^2}^2}{\|\eta\|_{L^2}^2} 
\le \inf_{\eta\in[\eta_1,\dots,\eta_{i-1}]^\perp\cap V_{n}}
\frac{\|d_A^{+,*}\eta\|_{L^2}^2}{\|\eta\|_{L^2}^2} 
\le C\la, 
$$
for $i=1,\dots,n$, which gives the first inequality in the stated
eigenvalue bounds
\eqref{eq:H2DeAEvalEst}. Therefore, $d_A^+d_A^{+,*}$ has at least $n$ 
eigenvalues which are less than or equal to $C\la$.  Similarly, the estimates
\eqref{eq:dA+*SmallEvalEst} yield
$$
\mu_{n+1}[A] = \inf_{\eta\in[\eta_1,\dots,\eta_{n}]^\perp}
\frac{\|d_A^{+,*}\eta\|_{L^2}^2}{\|\eta\|_{L^2}^2} 
\le \inf_{\eta\in[\eta_1,\dots,\eta_{n}]^\perp\cap V_{n+1}}
\frac{\|d_A^{+,*}\eta\|_{L^2}^2}{\|\eta\|_{L^2}^2} 
\le \left(\sqrt{\nu_2[A_0]} + C\sqrt{\la}\right)^2,
$$
completing the proof of the second inequality in the eigenvalue bounds
\eqref{eq:H2DeAEvalEst}.
\end{proof}

\subsection{An approximately orthonormal basis of eigenvectors}
\label{subsec:ApproxOrthoBasisEigenvectors}
We next construct an approximately $L^2$-orthonormal basis for $\Ker
d_{A_0}^+d_{A_0}^{+,*}$ using the eigenvectors corresponding to the
first $n$ eigenvalues of $d_A^+d_A^{+,*}$. Let $\Pi_{A_0}$ be the
$L^2$-orthogonal projection from $L^2(\Lambda^+\otimes\fg_{E_0})$ onto
$\Ker d_{A_0}^+d_{A_0}^{+,*}$ and let $\Pi_{A_0}^\perp := 1-\Pi_{A_0}$
be the $L^2$-orthogonal projection from
$L^2(\Lambda^+\otimes\fg_{E_0})$ onto $(\Ker
d_{A_0}^+d_{A_0}^{+,*})^\perp = \Ran d_{A_0}^+d_{A_0}^{+,*}$.

Recall from the proof of Lemma \ref{lem:H2SmallEvalUpperBound} that
$\{\eta_1[A],\dots,\eta_{n+1}[A]\}$ 
were defined to be the $L^2$-orthonormal
eigenvectors associated with the eigenvalues
$\mu_1[A]\le\cdots\le\mu_{n+1}[A]$ of $d_A^+d_A^{+,*}$ and construct
an approximately $L^2$-orthonormal subset
$\{\eta_i'\}_{i=1}^{n+1}\subset L^2(\Lambda^+\otimes\fg_{E_0})$ by
setting $\eta_i':=\chi_0\eta_i$ for $i=1,\dots,n+1$.

\begin{lem}
\label{lem:H2SmallEvalApproxOrthonBasis}
Continue the hypotheses and notation of Theorem
\ref{thm:H2SmallEval}. Then
  there are positive constants 
$$
C = C(\|F_A\|_{L^2},\|F_{A_0}\|_{L^2},\nu_2[A_0])
\quad\text{and}\quad
\lambda_0 = \lambda_0(C) 
$$
such that for all $\lambda \in
(0,\lambda_0]$,
\begin{gather}
\tag{1}
\|\Pi_{A_0}^\perp\eta_i'\|_{L^2}
\leq C\sqrt{\la},
\\
\tag{2}
\left|\left(\Pi_{A_0}\eta_i',\Pi_{A_0}\eta_j'\right)_{L^2(X)}
-\de_{ij}\right| \le C\la, \quad 1\le i,j\le n,
\\
\tag{3}
\|\Pi_{A_0}\eta_{n+1}'\|_{L^2(X)} \le C\sqrt{\la}
\quad\text{and}\quad
\|\Pi_{A_0}^\perp\eta_{n+1}'\|_{L^2(X)} \ge 1-C\sqrt{\la}. 
\end{gather}
\end{lem}

\begin{proof}
The following estimates are the key to the proof of Lemma
\ref{lem:H2SmallEvalApproxOrthonBasis}.

\begin{claim}
\label{claim:H2DeA0EvectorEst} 
Continue the above notation. Then there are positive constants 
$$
C = C(\|F_A\|_{L^2},\|F_{A_0}\|_{L^2},\nu_2[A_0])
\quad\text{and}\quad
\lambda_0 = \lambda_0(C)
$$
such that for all $\lambda \in (0,\lambda_0]$,
\begin{align}
\tag{1}
\|d_{A_0}^{+,*}\eta_i'\|_{L^2(X)} &\le C\sqrt{\la}, \quad 1\le i\le n,\\
\tag{2}
\|d_{A_0}^+d_{A_0}^{+,*}\eta_i'\|_{L^{4/3}(X)} 
&\le C\sqrt{\la}, \quad 1\le i\le n, 
\\
\tag{3}
\left|(\eta_i',\eta_j')_{L^2(X)}-\de_{ij}\right|
&\le C\la^2, \quad 1\le i,j\le n+1.
\end{align}
\end{claim}

\begin{proof}
To verify the first assertion in Claim \ref{claim:H2DeA0EvectorEst}, note that
$d_{A_0}^{+,*} = *d_{A_0}$ 
on $\Gamma(\Lambda^+\otimes\fg_{E_0})$, and so
\begin{equation}
\begin{aligned}
\label{eq:ExpanddA+*OnCutoffEvectorEta}
d_{A_0}^{+,*}\eta_i' 
&=
d_A^{+,*}(\chi_0\eta_i) - *(a\wedge \chi_0\eta_i)
\\
&= 
*(d\chi_0\wedge\eta_i) + \chi_0 d_A^{+,*}\eta_i  - *(\chi_0a\wedge \eta_i),
\quad i=1,\dots,n+1,
\end{aligned}
\end{equation}
where we use the fact that $A=A_0+a$ on $X\less\bx$. Thus, by the estimates in
\eqref{eq:H2DeAEvalEst} (see Lemma \ref{lem:H2SmallEvalUpperBound}) and Lemmas
\ref{lem:dBetaEst} and \ref{lem:LInftyEstDe2AEvec}, our hypotheses, and the 
identity \eqref{eq:ExpanddA+*OnCutoffEvectorEta}, 
\begin{align*}
\|d_{A_0}^{+,*}\eta_i'\|_{L^2} 
&\le (\|d\chi_0\|_{L^2}+\|\chi_0a\|_{L^2})\|\eta_i\|_{L^\8} 
+ \|d_A^{+,*}\eta_i\|_{L^2} \\
&\le c(\sqrt{\la}\|d\chi_0\|_{L^4}
+\|\chi_0a\|_{L^2})\|\eta_i\|_{L^\8} + \sqrt{\mu_i[A]} 
\\
&\le C\sqrt{\la}, \quad 1\le i\le n,
\end{align*}
using $\|d_A^{+,*}\eta_i\|_{L^2} = \sqrt{\mu_i[A]}\|\eta_i\|_{L^2}$;
this proves the first assertion in Claim
\ref{claim:H2DeA0EvectorEst}. Next, observe that schematically
(and again using the fact that $A=A_0+a$ on $X\less\bx$)
\begin{align*}
d_{A_0}^+d_{A_0}^{+,*}\eta_i' 
&= 
d_{A-a}^+(*(d\chi_0\wedge\eta_i) + \chi_0 d_A^{+,*}\eta_i
- *(\chi_0a\wedge\eta_i))
\\
&=
(\cov^2\chi_0+\chi_0\cov_{A_0}a+\chi_0a\otimes a)\otimes\eta_i
\\
&\quad+ \cov\chi_0\otimes (a\otimes\eta_i + \cov_A\eta_i)
+ \chi_0 d_A^+d_A^{+,*}\eta_i,
\end{align*}
and so, using Lemmas \ref{lem:dBetaEst}, \ref{lem:Kato}, 
\ref{lem:LInftyEstDe2AEvec}, and \ref{lem:H2SmallEvalUpperBound}, we get
\begin{align*}
\|d_{A_0}^+d_{A_0}^{+,*}\eta_i'\|_{L^{4/3}}
&\le \left(\|\cov^2\chi_0\|_{L^{4/3}} 
+ \|\chi_0\cov_{A_0}a\|_{L^{4/3}} \right.
\\
&\quad\left. + \|\cov\chi_0\|_{L^2}\|a\|_{L^4(\supp\chi_0)}
+ \|\sqrt{\chi_0}a\|_{L^2}\|\sqrt{\chi_0}a\|_{L^4}\right)
\|\eta_i\|_{L^\8} 
\\
&\quad + c\|\cov\chi_0\|_{L^2}\|\cov_A\eta_i\|_{L^4}
+ \|d_A^+d_A^{+,*}\eta_i\|_{L^2} \\
&\le C\sqrt{\la}\|\eta_i\|_{L^\8\cap L^4_{1,A}} + \mu_i[A]
\\
&\le C\sqrt{\la}, \quad 1\le i\le n. 
\end{align*}
Hence, the second assertion in Claim \ref{claim:H2DeA0EvectorEst}
follows from Lemma \ref{lem:dBetaEst} and our hypotheses. 
Using $\eta_i' = \chi_0\eta_i
= \eta_i + (\chi_0-1)\eta_i$, 
$$
(\eta_i',\eta_j')_{L^2}
= (\eta_i,\eta_j)_{L^2}
+ 2((\chi_0-1)\eta_i,\eta_j)_{L^2}
+ ((\chi_0-1)\eta_i,(\chi_0-1)\eta_j)_{L^2}.
$$
Now $(\eta_i,\eta_j)_{L^2}=\de_{ij}$, while
the second two terms on the right above are bounded by
\begin{align*}
&\left|2((\chi_0-1)\eta_i,\eta_j)_{L^2}
+ ((\chi_0-1)\eta_i,(\chi_0-1)\eta_j)_{L^2}\right| \\
&\qquad\le c\la^2\|\eta_i\|_{L^\8}\|\eta_j\|_{L^\8} \le C\la^2,
\end{align*}
and this gives the third assertion in Claim \ref{claim:H2DeA0EvectorEst}. 
\end{proof}

By Lemmas \ref{lem:Kato} and \ref{lem:L21AEstv}, integration by parts, and
Claim \ref{claim:H2DeA0EvectorEst} we have, for any 
$\xi\in (\Ker d_{A_0}^+d_{A_0}^{+,*})^\perp$,
\begin{align*}
(\eta_i',d_{A_0}^+d_{A_0}^{+,*}\xi)_{L^2} 
&= (d_{A_0}^+d_{A_0}^{+,*}\eta_i',\xi)_{L^2} 
\le \|d_{A_0}^+d_{A_0}^{+,*}\eta_i'\|_{L^{4/3}}\|\xi\|_{L^4} \\
&\le C\sqrt{\la}\left(\|\xi\|_{L^2} + \|d_{A_0}^{+,*}\xi\|_{L^2}\right) \\
&\le C\sqrt{\la}\left(\|\xi\|_{L^2} 
+ \|d_{A_0}^+d_{A_0}^{+,*}\xi\|_{L^2}\right) \\
&\le C(\nu_0[A_0]^{-1}+1)\sqrt{\la}\|d_{A_0}^+d_{A_0}^{+,*}\xi\|_{L^2},
\end{align*}
and therefore we obtain
\begin{equation}
(\eta_i',d_{A_0}^+d_{A_0}^{+,*}\xi)_{L^2} 
\le C\sqrt{\la}\|d_{A_0}^+d_{A_0}^{+,*}\xi\|_{L^2},
\quad \xi\in (\Ker d_{A_0}^+d_{A_0}^{+,*})^\perp, \quad i=1,\dots,n.
\label{eq:H2SmallEvecAngleEst}
\end{equation}
Since $\Ran \Pi_{A_0}^\perp = \Ran d_{A_0}^+d_{A_0}^{+,*}$, we may
choose $\rho_i\in L^2_2(\Lambda^+\otimes\fg_{E_0})$ such that
$\Pi_{A_0}^\perp\eta_i'=d_{A_0}^+d_{A_0}^{+,*}\rho_i$. If $\theta_i$
denotes the angle between $\eta_i'$ and $\Pi_{A_0}^\perp\eta_i'$ in
$L^2(\Lambda^+\otimes\fg_{E_0})$, then
$$
\|\Pi_{A_0}^\perp\eta_i'\|_{L^2}
= \|\eta_i'\|_{L^2}\cos\theta_i
= \frac{\left(\eta_i',d_{A_0}^+d_{A_0}^{+,*}\rho_i\right)_{L^2}}{
\|d_{A_0}^+d_{A_0}^{+,*}\rho_i\|_{L^2}}.
$$
Together with equation \eqref{eq:H2SmallEvecAngleEst}, this yields
for $i=1,\dots, n$, 
\begin{equation}
\|\Pi_{A_0}^\perp\eta_i'\|_{L^2}
\leq
C\sqrt{\la}.
\label{eq:H2PerpProjEigenvecEst} 
\end{equation}
By Claim \ref{claim:H2DeA0EvectorEst} we have, for $1\le i,j\le n$,
\begin{align*}
\left|\left(\Pi_{A_0}\eta_i',\Pi_{A_0}\eta_j'\right)_{L^2}
-\de_{ij}\right|
&\le \left|\left(\Pi_{A_0}\eta_i',\Pi_{A_0}\eta_j'\right)_{L^2}
-\left(\eta_i',\eta_j'\right)_{L^2}\right|
+ \left|\left(\eta_i',\eta_j'\right)_{L^2}-\de_{ij}\right| \\
&\le \left|\left(\Pi_{A_0}^\perp\eta_i',
\Pi_{A_0}^\perp\eta_j'\right)_{L^2}\right| + C\la^2 \\
&\le \|\Pi_{A_0}^\perp\eta_i'\|_{L^2}
\|\Pi_{A_0}^\perp\eta_j'\|_{L^2} + C\la^2 \\
&\le C(\la + \la^2) 
\quad\text{(from equation \eqref{eq:H2PerpProjEigenvecEst})},
\end{align*}
and so
\begin{equation}
\left|\left(\Pi_{A_0}\eta_i',\Pi_{A_0}\eta_j'\right)_{L^2}
-\de_{ij}\right| \le C\la, \quad 1\le i,j\le n.
\label{eq:H2KerDeA0ApproxBasis}
\end{equation}
Thus, for small enough $\la$, the vectors
$\{\Pi_{A_0}\eta_i'\}_{i=1}^{n}$ form an approximately $L^2$-orthonormal
basis for $\Ker d_{A_0}^+d_{A_0}^{+,*}$.

Recall from the proof of Lemma \ref{lem:H2SmallEvalUpperBound} that
$\{\xi_i(A_0)\}_{i=1}^{n}\subset L^2(\Lambda^+\otimes\fg_{E_0})$
denotes an ordered, $L^2$-orthonormal basis for $\Ker
d_{A_0}^+d_{A_0}^{+,*}$.  We can therefore assume without loss of
generality that $\{\xi_i\}_{i=1}^{n}$ is the basis for $\Ker
d_{A_0}^+d_{A_0}^{+,*}$ obtained by applying Gram-Schmidt
orthonormalization to the basis $\{\Pi_{A_0}\eta_i'\}_{i=1}^{n}$:
\begin{equation}
\label{eq:H2GramSchmidt}
\tilde\xi_\ell 
= \Pi_{A_0}\eta_\ell' 
- \sum_{i=1}^{\ell-1}(\Pi_{A_0}\eta_\ell',\xi_i)_{L^2}\xi_i
\quad\text{and}\quad
\xi_\ell = \frac{\tilde\xi_\ell}{\|\tilde\xi_\ell\|_{L^2}},
\quad 1\leq \ell \leq n.
\end{equation}
By Claim \ref{claim:H2DeA0EvectorEst} and the orthogonal-complement
estimates \eqref{eq:H2PerpProjEigenvecEst} we have
\begin{align}
\label{eq:H2PerpPropjLastEigenvectorEst}
|(\eta_{n+1}',\Pi_{A_0}\eta_i')_{L^2}|
&=
|(\eta_{n+1}',\eta_i')_{L^2}-(\eta_{n+1}',\Pi_{A_0}^\perp\eta_i')_{L^2}|
\\
\notag
&\leq
C\lambda^2 + C(1+\lambda)\sqrt{\lambda}
\leq 
C\sqrt{\lambda}.
\end{align}
Consequently, by the Gram-Schmidt equations \eqref{eq:H2GramSchmidt}
together with the inner-product estimates
\eqref{eq:H2PerpPropjLastEigenvectorEst}, and the fact that
$$
\Pi_{A_0}\eta_{n+1}' 
= \sum_{i=1}^{n}\left(\eta_{n+1}',\xi_i\right)_{L^2}\xi_i,
$$
we have
\begin{equation}
\|\Pi_{A_0}\eta_{n+1}'\|_{L^2} \le C\sqrt{\lambda} \quad\text{and}\quad
\|\Pi_{A_0}^\perp\eta_{n+1}'\|_{L^2} \ge 1-C\sqrt{\lambda}.
\label{eq:H2KerDeA0ApproxPerpVec}
\end{equation}
Thus, inequalities \eqref{eq:H2PerpProjEigenvecEst},
\eqref{eq:H2KerDeA0ApproxBasis}, and
\eqref{eq:H2KerDeA0ApproxPerpVec} give the required estimates for
$\{\Pi_{A_0}\eta_i'\}_{i=1}^{n+1}$ and completes the proof of Lemma
\ref{lem:H2SmallEvalApproxOrthonBasis}.
\end{proof}

\subsection{A lower bound for the first non-small eigenvalue}
\label{subsec:Eigen2LowerBoundFirstNonSmallEigenvalue}
Lastly, we compute a lower bound for the $(n+1)$-st eigenvalue of
$d_A^+d_A^{+,*}$. Recall that $\eta_{n+1}=\eta_{n+1}(A)$ is the
$L^2$-unit eigenvector associated to the $(n+1)$-st eigenvalue
$\mu_{n+1}[A]$ of $d_A^+d_A^{+,*}$.

\begin{lem}
\label{lem:H2SmallEvalLowerBound}
Continue the above notation. Then there are positive constants 
$$
C = C(\|F_A\|_{L^2},\|F_{A_0}\|_{L^2},\nu_2[A_0])
\quad\text{and}\quad
\lambda_0 = \lambda_0(C)
$$ 
such that for all $\lambda \in (0,\lambda_0]$,
$$
\mu_{n+1}[A] \ge \nu_2[A_0] - C\sqrt{\la}.
$$
\end{lem}

\begin{pf}
Recall that
$\eta_{n+1}'=\chi_0\eta_{n+1}\in\Gamma(\Lambda^+\otimes\fg_{E_0})$.
Therefore, via the identity \eqref{eq:ExpanddA+*OnCutoffEvectorEta},
Lemma \ref{lem:dBetaEst}, our hypotheses, and the fact that
$\|d_A^{+,*}\eta_{n+1}\|_{L^2} =
\sqrt{\mu_{n+1}[A]}\|\eta_{n+1}\|_{L^2}$, we get
\begin{align*}
\|d_{A_0}^{+,*}\eta_{n+1}'\|_{L^2} 
&\le \left(\|d\chi_0\|_{L^2}+\|\chi_0a\|_{L^2}\right)\|\eta_{n+1}\|_{L^\8}
+ \|d_A^{+,*}\eta_{n+1}\|_{L^2} 
\\
&\le c\sqrt{\la}\|\eta_{n+1}\|_{L^\8}
+ \sqrt{\mu_{n+1}[A]}\|\eta_{n+1}\|_{L^2}.
\end{align*}
According to Lemma \ref{lem:LInftyEstDe2AEvec} we have
\begin{align*}
\|\eta_{n+1}\|_{L^\8}
&\leq c(1+\|F_A\|_{L^2})
(1+\mu_{n+1}[A])(1+\sqrt{\mu_{n+1}[A]})\|\eta_{n+1}\|_{L^2}
\\
&\leq C\|\eta_{n+1}\|_{L^2},
\end{align*}
using the upper bound for $\mu_{n+1}[A]$ in Lemma
\ref{lem:H2SmallEvalUpperBound}. Hence, we obtain
$$
\|d_{A_0}^{+,*}\eta_{n+1}'\|_{L^2} 
\leq C\sqrt{\la} + \sqrt{\mu_{n+1}[A]}.
$$
In the other direction, we have the eigenvalue estimate
$\|d_{A_0}^{+,*} v\|_{L^2} \geq \sqrt{\nu_2[A_0]}\|v\|_{L^2}$ for all
$v\in(\Ker d_{A_0}^+d_{A_0}^{+,*})^\perp
\subset L^2(\Lambda^+\otimes\fg_{E_0})$, and so by the orthogonal projection 
estimates \eqref{eq:H2KerDeA0ApproxPerpVec},  
$$
\|d_{A_0}^{+,*}\eta_{n+1}'\|_{L^2}
\ge 
\sqrt{\nu_2[A_0]}\|\Pi_{A_0}^\perp\eta_{n+1}'\|_{L^2} 
\ge 
\sqrt{\nu_2[A_0]}(1 - C\la).
$$
Combining these upper and lower $L^2$ bounds for
$d_{A_0}^{+,*}\eta_{n+1}'$ yields
$$
\sqrt{\nu_2[A_0]}\left(1 - C\sqrt{\la}\right)
\le C\sqrt{\la} + \sqrt{\mu_{n+1}[A]},
$$
and this gives the lower bound on $\mu_{n+1}[A]$ for small enough $\la$.
\end{pf}

\begin{proof}[Proof of Theorem \ref{thm:H2SmallEval}]
  Combine Lemmas \ref{lem:H2SmallEvalUpperBound} and
  \ref{lem:H2SmallEvalLowerBound}.
\end{proof}

\subsection{Variation of eigenvalues with the connection}
\label{subec:Eigenvalue2Variation}
Our final task is to estimate the variation in the eigenvalues
$\mu_i[A]$ of $d_A^+d_A^{+,*}$ with the connection $A$.

\begin{lem}
\label{lem:PerturbationSecondLaplacianEigenvalues}
  Let $X$ be a $C^\8$, closed, oriented, Riemannian four-manifold.
  Then there is a positive constant $c$ with the following
  significance.  Let $E$ be a rank-two, Hermitian vector bundle over
  $X$, let $A$ be an orthogonal $L^2_4$ connection on $\fg_E$, and let
  $0\leq \mu_1[A] \leq \mu_2[A] \leq \cdots$ be the eigenvalues of
  $d_A^+d_A^{+,*}$ on $L^2(\Lambda^+\otimes\fg_E)$, repeated according to their
  multiplicity.  Then, for any integer $i\geq 1$ and $a\in
  L^2_4(\Lambda^1\otimes\fg_E)$,
\begin{align}
\tag{1}
\mu_i[A+a] 
&\le \mu_i[A](1 + C\|a\|_{L^4(X)} + C\|a\|_{L^4(X)}^2) 
+ C\|a\|_{L^4(X)}(1+\|a\|_{L^4(X)}), \\
\tag{2}
\mu_i[A] 
&\le \mu_i[A+a](1 + C\|a\|_{L^4(X)} + C\|a\|_{L^4(X)}^2) 
+ C\|a\|_{L^4(X)}(1+\|a\|_{L^4(X)}),
\end{align}
where $C[A] := c(1+\|F_A^+\|_{L^2})$.
\end{lem}

\begin{pf}
Recall that by Lemma \ref{lem:MaxMin}
$$
\mu_i[A+a] 
= 
\inf_{\eta\in[\eta_1[A+a],\dots,\eta_{i-1}[A+a]]^\perp}
\frac{\|d_{A+a}^{+,*}\eta\|_{L^2}^2}{\|\eta\|_{L^2}^2},
\quad i\in\NN,
$$
where the infimum is taken over all $\eta\in
L^2_1(\Lambda^+\otimes\fg_E)$ which are $L^2$ orthogonal to the span
of $\eta_1[A+a],\dots,\eta_{i-1}[A+a]$ in
$L^2(\Lambda^+\otimes\fg_E)$, the first $i-1$ eigenvectors of   
$d_{A+a}^+d_{A+a}^{+,*}$ on $L^2(\Lambda^+\otimes\fg_E)$. 
Since $d_{A+a}^{+,*}\eta =
d_A^{+,*}\eta + *(a\wedge\eta)$ we have, for any $\eta \in
L^2_1(\Lambda^+\otimes\fg_E)$,
\begin{align*}
\|d_{A+a}^{+,*}\eta\|_{L^2}^2 
&= \|d_A^{+,*}\eta\|_{L^2}^2 + \|*(a\wedge\eta)\|_{L^2}^2 
+ 2(d_A^{+,*}\eta,*(a\wedge\eta))_{L^2} \\
&\le \|d_A^{+,*}\eta\|_{L^2}^2(1 + \|a\|_{L^4})
+ (1 + \|a\|_{L^4})\|a\|_{L^4}\|\eta\|_{L^4}^2.
\end{align*}
But Lemmas \ref{lem:Kato} and \ref{lem:L21AEstv} imply that
$$
\|\eta\|_{L^4} \le c\|\eta\|_{L^2_{1,A}} \le c(1+\|F_A^+\|_{L^2})^{1/2}
(\|d_A^{+,*}\eta\|_{L^2}+\|\eta\|_{L^2}), 
$$
for any $\eta\in L^2_1(\Lambda^+\otimes\fg_E)$ and so
\begin{align*}
\|d_{A+a}^{+,*}\eta\|_{L^2}^2 
&\le \|d_A^{+,*}\eta\|_{L^2}^2(1 + \|a\|_{L^4})
+ C[A](1 + \|a\|_{L^4})\|a\|_{L^4}
(\|d_A^{+,*}\eta\|_{L^2} + \|\eta\|_{L^2})^2 \\
&\le \|d_A^{+,*}\eta\|_{L^2}^2(1 + C[A]\|a\|_{L^4}+ C[A]\|a\|_{L^4}^2)
+ C[A](1 + \|a\|_{L^4})\|a\|_{L^4}\|\eta\|_{L^2}^2,
\end{align*}
where $C[A] = c(1+\|F_A^+\|_{L^2})$. 
By Lemma \ref{lem:MaxMin} and the preceding estimate we see that
\begin{align*}
\mu_i[A+a] &= \inf_{\eta\in[\eta_1[A+a],\dots,\eta_{i-1}[A+a]]^\perp}
\frac{\|d_{A+a}^{+,*}\eta\|_{L^2}^2}{\|\eta\|_{L^2}^2} \\
&\leq \inf_{\eta\in[\eta_1[A+a],\dots,\eta_{i-1}[A+a]]^\perp}
\frac{\|d_A^{+,*}\eta\|_{L^2}^2}{\|\eta\|_{L^2}^2}
(1 + C[A]\|a\|_{L^4}+ C[A]\|a\|_{L^4}^2)
\\
&\quad + C[A](1 + \|a\|_{L^4})\|a\|_{L^4} \\
&\le \mu_i[A](1 + C[A]\|a\|_{L^4}+ C[A]\|a\|_{L^4}^2) 
+ C[A](1 + \|a\|_{L^4})\|a\|_{L^4},
\end{align*}
and this gives the upper bound for $\mu_i[A+a]$ in terms of
$\mu_i[A]$. Interchanging the roles of $A$ and $A+a$ gives the upper bound
for $\mu_i[A]$ in terms of $\mu_i[A+a]$.
\end{pf}


\section{Existence of gluing maps and obstruction sections}
\label{sec:Existence}
In this section we prove the assertions of Theorem \ref{thm:GluingTheorem}
concerning the existence and smoothness of the gluing map $\bga$ and
obstruction section $\bchi$, but not the assertions concerning the
existence of continuous extensions of these maps to the Uhlenbeck
compactification of the space of gluing data,
$\bar{\Gl}_\kappa^{w,+}(X,\sU_0,\Sigma,\lambda_0)$; the latter Uhlenbeck
continuity assertions are taken up in \S \ref{sec:Continuity}.

\subsection{The extended anti-self-dual equations}
\label{subsec:ExtendedASD}
As we saw in \S \ref{subsec:ParameterizedKuranishi}, we need to allow for
the presence of cokernel obstructions to deformation in the parameterized
moduli space of background connections contained in the gluing data and,
thus, small-eigenvalue obstructions to gluing. In this subsection we set up
a version of the `extended anti-self-dual equations' which we shall use to
take account of such obstructions. Though similar setups are described in
\cite{DonConn}, \cite{DonPoly}, \cite{DK}, \cite{FrM},
\cite{MrowkaThesis}, \cite{TauIndef}, and \cite{TauStable}, they are not
identical so we shall describe our version in some detail.

Let $A$ be an orthogonal connection on an $\SO(3)$ bundle $\fg_E=\su(E)$
over $X$, where $E$ is a Hermitian, rank-two vector bundle over $X$.  Our
first task is to obtain existence and global estimates for solutions $a\in
L^2_1(\Lambda^1\otimes\fg_E)$ to the quasi-linear equation,
\begin{equation}
\label{eq:ASD-A}
d_A^+a + (a\wedge a)^+ = w \quad\text{and}\quad
a=d_A^{+,*}v, 
\end{equation}
where $v \in (C^0\cap L^{\sharp,2}_2)(\Lambda^+\otimes\fg_E)$ and $w\in
L^{\sharp,2}(\Lambda^+\otimes\fg_E)$. (For the definition and properties of
the $L^\sharp$ and related families of Sobolev norms, see \cite[\S
4]{FeehanSlice}.)  As in \cite{TauSelfDual},
\cite{TauIndef}, \cite{TauStable}, we use $a=d_A^{+,*}v$ to rewrite the
first-order quasi-linear equation
\eqref{eq:ASD-A} for $a$ in terms of $w$ as
\begin{equation}
\label{eq:SecondOrderASD}
d_A^+d_A^{+,*}v + (d_A^{+,*}v\wedge d_A^{+,*}v)^+ = w, 
\end{equation}
to give a second-order, quasi-linear elliptic equation for $v$ in terms of $w$.
For the reader familiar with the construction of gluing maps in \cite{DK},
the current section gives a generalization of the treatment in \cite[\S
7.2.2 \& 7.2.3]{DK}, where existence of solutions (in
certain cases) to the anti-self-dual equation is established.

As in \cite{TauIndef}, \cite{TauStable}, we allow for the
possibility that the Laplacian $d_A^+d_A^{+,*}$ has small
eigenvalues. {}From \S \ref{sec:Eigenvalue2} we know that this occurs
when $[A]$ is close, in the Uhlenbeck topology, to a point
$[A_0,\bx]\in\sB_\kappa^w(X)\times\Sym^\ell(X)$, where
$d_{A_0}^+d_{A_0}^{+,*}$ has a positive-dimensional kernel. Thus, in
place of the anti-self-dual equation,
\begin{equation}
\label{eq:QuickASDEqnForv}
F^+(A+d_A^{+,*}v) = 0,
\end{equation}
for $v \in L^2_{k+1}(\Lambda^+\otimes\fg_E)$, we 
follow \cite{TauIndef}, \cite{TauStable} and consider the pair of equations
\begin{align}
\label{eq:QuickExtASDEqnForv}
\Pi_{A,\mu}^\perp\cdot F^+(A+d_A^{+,*}v) = 0,
\\
\label{eq:QuickFinitePartASDEqnForv}
\Pi_{A,\mu}\cdot F^+(A+d_A^{+,*}v) = 0,
\end{align}
for solutions $v \in L^2_{k+1}(\Lambda^+\otimes\fg_E)$ satisfying the
constraint $\Pi_{A,\mu}v=0$; equation \eqref{eq:QuickExtASDEqnForv} is
called the {\em extended anti-self-dual equation\/}. Here,
$\Pi_{A,\mu}$ is the $L^2$-orthogonal projection from
$L^2(\Lambda^+\otimes\fg_E)$ onto the finite-dimensional subspace
spanned by the eigenvectors of $d_A^+d_A^{+,*}$ with eigenvalues in
$[0,\mu]$, where $\mu$ is a positive constant (an upper bound for the
small eigenvalues), and $\Pi_{A,\mu}^\perp := \id - \Pi_{A,\mu}$.
Because of the small-eigenvalue problem, we shall only explicitly solve
equation \eqref{eq:QuickExtASDEqnForv}; the true solutions to 
equation \eqref{eq:QuickASDEqnForv} are then
simply the solutions to equation \eqref{eq:QuickExtASDEqnForv} which also
satisfy equation \eqref{eq:QuickFinitePartASDEqnForv}.

Equation \eqref{eq:QuickExtASDEqnForv} may be rewritten as 
\begin{equation}
\label{eq:LongSecOrderExtASDEqnForv}
\Pi_{A,\mu}^\perp(d_A^+d_A^{+,*}v + (d_A^{+,*}v\wedge d_A^{+,*}v)^+)
= - \Pi_{A,\mu}^\perp F_A^+,
\end{equation}
where again $v \in L^2_{k+1}(\Lambda^+\otimes\fg_E)$ satisfies the
constraint $\Pi_{A,\mu}v=0$.  For the purposes of the regularity
theory of \S \ref{sec:Regularity}, it is useful to rewrite equation
\eqref{eq:LongSecOrderExtASDEqnForv} as
\begin{equation}
\label{eq:RegLongSecOrderExtASDEqnForv}
d_A^+d_A^{+,*}v + (d_A^{+,*}v\wedge d_A^{+,*}v)^+ 
=
\Pi_{A,\mu}(d_A^+d_A^{+,*}v + (d_A^{+,*}v\wedge d_A^{+,*}v)^+)
- \Pi_{A,\mu}^\perp F_A^+.
\end{equation}
On the other hand, for the purposes of applying Banach-space fixed
theory to solve the non-linear equation
\eqref{eq:LongSecOrderExtASDEqnForv}, it is convenient to write
$$
v = G_{A,\mu}\xi, 
\quad\text{for some } \xi \in L^2_{k-1}(\Lambda^+\otimes\fg_E),
$$
where $G_{A,\mu}:L^2_{k-1}(\Lambda^+\otimes\fg_E) \to
L^2_{k+1}(\Lambda^+\otimes\fg_E)$ is the Green's operator for the Laplacian
$$
d_A^+d_A^{+,*}: \Pi_{A,\mu}^\perp\cdot L^2_{k+1}(\Lambda^+\otimes\fg_E)
\to \Pi_{A,\mu}^\perp\cdot L^2_{k-1}(\Lambda^+\otimes\fg_E).
$$
Thus, $d_A^+d_A^{+,*}G_{A,\mu} = \Pi_{A,\mu}^\perp$ and $\Pi_{A,\mu}v
= \Pi_{A,\mu}G_{A,\mu}\xi = 0$ for any $\xi \in
L^2_{k-1}(\Lambda^+\otimes\fg_E)$. Therefore, it suffices to solve for
$\xi$ such that
\begin{equation}
\label{eq:IntegralASDEqnForxi}
\xi + (d_A^{+,*}G_{A,\mu}\xi\wedge d_A^{+,*}G_{A,\mu}\xi)^+ = - F_A^+,
\end{equation}
because, noting that $\Pi_{A,\mu}^\perp\xi =
\Pi_{A,\mu}^\perp d_A^+d_A^{+,*}G_{A,\mu}\xi$, the resulting 
$v = G_{A,\mu}\xi$ solves equation \eqref{eq:LongSecOrderExtASDEqnForv},
namely
$$ 
\Pi_{A,\mu}^\perp
(\xi + (d_A^{+,*}G_{A,\mu}\xi\wedge d_A^{+,*}G_{A,\mu}\xi)^+ + F_A^+) = 0,
$$
{\em as well as} the constraint $\Pi_{A,\mu}v = 0$.
 
Noting that $\Pi_{A,\mu}d_A^+d_A^{+,*}G_{A,\mu}\xi = 0$, the
solution $\xi$ to \eqref{eq:IntegralASDEqnForxi} solves equation
\eqref{eq:QuickFinitePartASDEqnForv}, namely
$$
\Pi_{A,\mu}
((d_A^{+,*}G_{A,\mu}\xi\wedge d_A^{+,*}G_{A,\mu}\xi)^+ + F_A^+) = 0,
$$
{\em if and
only if\/}
\begin{equation}
\label{eq:IntegralFiniteASDEqnForxi}
\Pi_{A,\mu}\xi = 0.
\end{equation}
For applications of Banach-space fixed-point theory, the integral form
\eqref{eq:IntegralASDEqnForxi} of the extended anti-self-dual equation 
is a little more convenient, while for the purposes of regularity theory
and localization, we use the ``differential'' form
\eqref{eq:RegLongSecOrderExtASDEqnForv}. 

\subsection{Global estimates: linear theory} 
\label{subsec:GlobalEstLinear}
The results we need are treated in detail in \cite{FeehanSlice}, so we just
summarize the main conclusions we need here; the $L^\sharp$ and related
families of Sobolev norms are defined in \cite[\S 4]{FeehanSlice}.

\begin{lem}
\label{lem:LinftyL22CovLapEstv}
\cite[Lemma 5.9]{FeehanSlice}
Let $X$ be a closed, oriented four-manifold with metric $g$. Then there are
positive constants $c$ and $\eps=\eps(c)$ with the following significance.
Let $E$ be a Hermitian, rank-two vector bundle over $X$ and let $A$ be an
orthogonal $L^2_4$ connection on $\fg_E$ with curvature $F_A$, such that
$\|F_A\|_{L^{\sharp,2}(X)} <\eps$.  Then the following estimate holds for
any $v\in L^{\sharp,2}_2(\Lambda^+\otimes\fg_E)$:
\begin{align*}
\|v\|_{L^2_{2,A}(X)} + \|v\|_{C^0(X)}
\le c(1+\|F_A\|_{L^2(X)})
(\|d_A^+d_A^{+,*}v\|_{L^{\sharp,2}(X)} + \|v\|_{L^2(X)}).
\end{align*}
\end{lem}

\begin{cor}
\label{cor:L21AEstdA*v}
\cite[Corollary 5.10]{FeehanSlice}
Continue the hypotheses of Lemma \ref{lem:LinftyL22CovLapEstv}. Then: 
\begin{align*}
\|d_A^{+,*}v\|_{L^2_{1,A}(X)}
\le c(1+\|F_A\|_{L^2(X)})
(\|d_A^+d_A^{+,*}v\|_{L^{\sharp,2}(X)} + \|v\|_{L^2(X)}).
\end{align*}
\end{cor}

Note that if $a\in L^2_k(\Lambda^1\otimes\fg_E)$ is $L^2$-orthogonal to
$\Ker d_A^+$, so that $a=d_A^{+,*}v$ for some $v\in
L^2_{k+1}(\Lambda^+\otimes\fg_E)$, then the estimate of Corollary
\ref{cor:L21AEstdA*v} can be written in the more familiar form (see
\cite[Eq. (5.10)]{FeehanSlice})
\begin{equation}
\label{eq:UniformFirstOrderEllipticEst}
\|a\|_{L^2_{1,A}(X)}
\le
c(1+\|F_A\|_{L^2(X)})
(\|d_A^+a\|_{L^{\sharp,2}(X)} + \mu^{-1/2}\|a\|_{L^2(X)}),
\end{equation}
where $\mu = \nu_2[A]>0$ and
we make use of the eigenvalue estimate $\|v\|_{L^2} \le
\mu^{-1/2}\|d_A^{+,*}v\|_{L^2}$; the term $d_A^+a$ above can be
replaced by $(d_A^{+,*}d_A^+)a$ without changing the estimate constants.
More generally, if $\mu$ is any positive constant and $\Pi_{A,\mu}v=0$,
then the estimate \eqref{eq:UniformFirstOrderEllipticEst} continues to hold.

Fix a positive constant $\mu$.  Given an $L^2_4$ orthogonal connection $A$
on $\fg_E$, we define a partial right inverse
$$
P_{A,\mu} := d_A^{+,*}\left(d_A^+d_A^{+,*}\right)^{-1}\cdot \Pi_{A,\mu}^\perp:
L^{\sharp,2}(\Lambda^+\otimes\fg_E) \to L^2_1(\Lambda^1\otimes\fg_E),
$$
to the operator
$$
d_A^+:L^2_1(\Lambda^1\otimes\fg_E)\to L^{\sharp,2}(\Lambda^+\otimes\fg_E).
$$
We have the following estimates for $P_{A,\mu}$:

\begin{prop}
\label{prop:EstPA0mu}
Let $X$ be a closed, oriented, smooth four-manifold with Riemannian metric
$g$.  Then is a positive constant $\eps$ with the following
significance. Let $A$ be an $L^2_4$ orthogonal connection on an $\SO(3)$
bundle $\fg_E$ over $X$ such that $\|F_A^+\|_{L^{\sharp,2}}<\eps$.  Let
$\mu$ be a positive constant.  Then there is a positive constant
$C(\|F_A\|_{L^2},\mu)$ such that for all $\xi\in
L^{\sharp,2}(\Lambda^+\otimes\fg_E)$,
\begin{align}
\|P_{A,\mu}\xi\|_{L^2} &\le C\|\xi\|_{L^{4/3}}, \tag{1} \\
\|P_{A,\mu}\xi\|_{L^2_{1,A}} &\le C\|\xi\|_{L^{\sharp,2}}. \tag{2}
\end{align}
\end{prop}

\begin{pf}
We set $v=G_{A,\mu}\xi$, so that $\Pi_{A,\mu}v = 0$ with
$$
d_A^+d_A^{+,*}v=\Pi_{A,\mu}^\perp\xi 
\quad\text{and}\quad 
a=d_A^{+,*}v=P_{A,\mu}\xi.
$$ 
We first consider the $L^2$ estimate for $a$.  We integrate by parts to get
\begin{align*}
\|a\|_{L^2}^2 
&= 
\|d_A^{+,*}v\|_{L^2}^2 = (d_A^+d_A^{+,*}v,v)_{L^2}
= (\Pi_{A,\mu}^\perp\xi,v)_{L^2} = (\xi,v)_{L^2} 
\\
&\le \|\xi\|_{L^{4/3}}\|v\|_{L^4} 
\le c\|\xi\|_{L^{4/3}}\|v\|_{L^2_{1,A}}.
\end{align*}
{}From Lemma \ref{lem:L21AEstv} and the eigenvalue estimate $\|v\|_{L^2} \leq
\mu^{-1/2}\|d_A^{+,*}v\|_{L^2}$, we have
$$
\|v\|_{L^2_{1,A}} \le C(\|d_A^{+,*}v\|_{L^2} + \|v\|_{L^2})
\le C(\|d_A^{+,*}v\|_{L^2} + \mu^{-1/2}\|d_A^{+,*}v\|_{L^2}) 
\le C\|d_A^{+,*}v\|_{L^2}.
$$
Therefore,
$$
\|d_A^{+,*}v\|_{L^2}^2 \le C\|\xi\|_{L^{4/3}}\|d_A^{+,*}v\|_{L^2},
$$
and so
$$
\|d_A^{+,*}v\|_{L^2} \le C\|\xi\|_{L^{4/3}},
$$
which gives the desired $L^2$ estimate for $a=d_A^{+,*}v$. We can also
combine the preceding estimates, noting that we can choose to either retain or
drop the projection $\Pi_{A,\mu}^\perp$ as above, to give
\begin{equation}
\label{eq:L21EstForvInTermsL43Data}
\begin{aligned}
\|v\|_{L^2_{1,A}} &\le C\|\Pi_{A,\mu}^\perp\xi\|_{L^{4/3}},
\\
\|v\|_{L^2_{1,A}} &\le C\|\xi\|_{L^{4/3}},
\end{aligned}
\end{equation}
estimates we shall later need.

We next derive the $L^2_{1,A}$ estimate for $a$. 
Since $a = d_A^{+,*}v$, then
$$
\|a\|_{L^2_{1,A}} \le \|d_A^{+,*}v\|_{L^2_{1,A}} .
$$
Using integration by parts, we see that
\begin{align*}
\|d_A^{+,*}v\|_{L^2}^2 
&= (d_A^{+,*}v,d_A^{+,*}v)_{L^2} = (d_A^+d_A^{+,*}v,v)_{L^2}
\\
&\le \|d_A^+d_A^{+,*}v\|_{L^2}\|v\|_{L^2} 
\le \half(\|d_A^+d_A^{+,*}v\|_{L^2}^2 + \|v\|_{L^2}^2),
\end{align*}
and so
\begin{equation}
\|d_A^{+,*}v\|_{L^2} 
\leq 
{\frac{1}{\sqrt{2}}}(\|d_A^+d_A^{+,*}v\|_{L^2} + \|v\|_{L^2}). 
\label{eq:EasyEstdA*v}
\end{equation}
By combining Corollary \ref{cor:L21AEstdA*v} with equation
\eqref{eq:EasyEstdA*v} and noting that
$d_A^+d_A^{+,*}v=\Pi_{A,\mu}^\perp\xi = \xi - \Pi_{A,\mu}\xi$, we obtain
\begin{align*}
\|a\|_{L^2_{1,A}} &\le C(\|d_A^+d_A^{+,*}v\|_{L^{\sharp,2}} + \|v\|_{L^2}) \\
&= C(\|d_A^+d_A^{+,*}v\|_{L^{\sharp,2}} +
\mu^{-1}\|d_A^+d_A^{+,*}v\|_{L^2}) \\
&\le C\|\Pi_{A,\mu}^\perp\xi\|_{L^{\sharp,2}}
\\
&\le C\|\xi\|_{L^{\sharp,2}},
\quad\text{(by Lemma \ref{lem:LInftyEstDe2AEvec}
and \cite[Lemma 4.1]{FeehanSlice})}
\end{align*}
which gives the required $L^2_{1,A}$ estimate for $a$. Of course, the
preceding argument also yields
\begin{equation}
\begin{aligned}
\label{eq:C0L22EstForvInTermsLSharp2Data}
\|v\|_{C^0\cap L^2_{2,A}} &\le C\|\Pi_{A,\mu}^\perp\xi\|_{L^{\sharp,2}},
\\
\|v\|_{C^0\cap L^2_{2,A}} &\le C\|\xi\|_{L^{\sharp,2}},
\end{aligned}
\end{equation}
estimates we shall later need.
\end{pf}

\subsection{Global estimates: non-linear theory}
\label{subsec:GlobalEstimatesNonLinear} 
We now extend our global estimates in \S \ref{subsec:GlobalEstLinear} for
the linear operator $d_A^+d_A^{+,*}$ to estimates for the non-linear
operator $d_A^+d_A^{+,*}+(d_A^{+,*}(\cdot)\wedge d_A^{+,*}(\cdot))^+$. 
By adapting the arguments of
\cite{TauSelfDual}, \cite{TauFrame}, \cite{TauStable} 
we first obtain global $C^0\cap L^2_{2,A}$ estimates for solutions $v$ to
equation
\eqref{eq:SecondOrderASD} and thus global $L^2_{1,A}$ estimates for
solutions $a=d_A^{+,*}v$ to equation \eqref{eq:ASD-A}. These global
estimates have the property that the constant $C$ is {\em uniform\/} with
respect to the point $[A]\in\sB_E(X)$ in the sense that $C$ depends on $A$
only through $\|F_A^+\|_{L^{\sharp,2}}$, $\|F_A\|_{L^2}$, and $\mu$, a
choice of positive, small-eigenvalue cutoff constant for the Laplacian
$d_A^+d_A^{+,*}$ on $L^2(\Lambda^+\otimes\fg_E)$.  We will make frequent
use of the following {\em a priori\/} estimate for $L^2_2$ solutions to
equation \eqref{eq:SecondOrderASD}.

\begin{lem}
\label{lem:GlobalASDEst}
Let $X$ be a closed, oriented, smooth four-manifold with Riemannian metric
$g$, and let $M, \La$ be positive constants.  Then there are positive
constants $\eps, C$ such that the following holds. Let $A$ be an $L^2_4$
connection on $\fg_E$ with $\|F_A^+\|_{L^{\sharp,2}} < \eps$,
$\|F_A\|_{L^2}\le M$, and let $\mu$ be a positive constant such that
$\La^{-1} \le \mu \le \La$.  Let $v\in (C^0\cap L^2_2)(\Lambda^+\otimes\fg_E)$
and $w\in L^{\sharp,2}(\Lambda^+\otimes\fg_E)$ solve
\eqref{eq:SecondOrderASD}, with $\Pi_{A,\mu}v=0$, and suppose that
$\|v\|_{L^2_{2,A}(X)}<\eps$.  Then
\begin{align}
\|v\|_{L^2_{1,A}(X)} &\le C\|w\|_{L^{4/3}(X)}, \tag{1}\\
\|v\|_{L^2_{2,A}(X)} + \|v\|_{L^\8(X)} &\le C\|w\|_{L^{\sharp,2}(X)}. \tag{2}
\end{align}
\end{lem}

\begin{pf}
Note that estimate \eqref{eq:L21EstForvInTermsL43Data} yields
$$
\|v\|_{L^2_{1,A}} \le C\|d_A^+d_A^{+,*} v\|_{L^{4/3}},
$$
while from equation \eqref{eq:SecondOrderASD} and the Sobolev multiplications
$L^2_{1,A}\times L^2\to L^{4/3}$ (with constant depending only on the
metric) we have
\begin{align*}
\|d_A^+d_A^{+,*} v\|_{L^{4/3}}
&\le \|(d_A^{+,*}v\wedge d_A^{+,*}v)^+\|_{L^{4/3}} + \|w\|_{L^{4/3}} \\
&\le c\|d_A^{+,*}v\|_{L^2_{1,A}}\|d_A^{+,*}v\|_{L^2} + \|w\|_{L^{4/3}},
\end{align*}
and so,
$$
\|v\|_{L^2_{1,A}} 
\le C\left(\|v\|_{L^2_{2,A}}\|v\|_{L^2_{1,A}} + \|w\|_{L^{4/3}}\right).
$$
Thus, taking $\eps\le 1/(2C)$, the first estimate follows by
rearrangement. 

H\"older's inequality and 
our elliptic estimate in Lemma \ref{lem:LinftyL22CovLapEstv} 
for the Laplacian $d_A^+d_A^{+,*}$ yields
$$
\|v\|_{L^2_{2,A}} + \|v\|_{L^\8} 
\le C\|d_A^+d_A^{+,*} v\|_{L^{\sharp,2}} + \|v\|_{L^2},
$$
while from equation \eqref{eq:SecondOrderASD} and the Sobolev
multiplications $L^2_{1,A}\times L^2_{1,A}\to L^{\sharp,2}$ (again with
constant depending only on the metric) we have
\begin{align*}
\|d_A^+d_A^{+,*} v\|_{L^{\sharp,2}}
&\le \|(d_A^{+,*}v\wedge d_A^{+,*}v)^+\|_{L^{\sharp,2}} + \|w\|_{L^{\sharp,2}} \\
&\le c\|d_A^{+,*}v\|_{L^2_{1,A}}^2 + \|w\|_{L^{\sharp,2}},
\end{align*}
and so,
$$
\|v\|_{L^2_{2,A}} + \|v\|_{L^\8}
\le C\left(\|v\|_{L^2_{2,A}}^2 + \|w\|_{L^{\sharp,2}}\right).
$$
Thus, taking $\eps\le 1/2C$, the second estimate follows by
rearrangement. 
\end{pf} 

\subsection{Existence of solutions to the extended anti-self-dual equation}
\label{subsec:ExistSolnExtASD}
We now turn to the question of existence and uniqueness of solutions to 
equation \eqref{eq:QuickExtASDEqnForv} or one of its equivalent forms,
namely equations \eqref{eq:LongSecOrderExtASDEqnForv} or
\eqref{eq:RegLongSecOrderExtASDEqnForv}.
We first recall the following elementary fixed-point result \cite[Lemma
7.2.23]{DK}.

\begin{lem}\label{lem:FixedPoint}
Let $q:\fB\to\fB$ be a continuous map on a Banach space $\fB$ with $q(0)=0$
and  
$$
\|q(x_1)-q(x_2)\|\le K\left(\|x_1\|+\|x_2\|\right)\|x_1-x_2\|
$$
for some positive constant $K$ and all $x_1,x_2$ in $\fB$. Then for each
$y$ in $\fB$ with $\|y\|<1/(10K)$ there is a unique $x$ in $\fB$ such that
$\|x\|\le 1/(5K)$ and $x+q(x)=y$.
\end{lem}

With the aid of Lemma \ref{lem:FixedPoint}, we obtain existence and
uniqueness for solutions to equation \eqref{eq:QuickExtASDEqnForv}. 

\begin{prop}
\label{prop:ASDExist}
Let $(X,g)$ be a closed, oriented, $C^\8$ Riemannian four-manifold, and let
$M,\La$ be positive constants.  Then there are positive constants $\eps, C$
such that the following holds. Let $E$ be a Hermitian, rank-two vector
bundle and let $A$ be an $L^2_4$ orthogonal connection on $\fg_E$ with
$\|F_A^+\|_{L^{\sharp,2}(X)}\le M$, $\|F_A\|_{L^2(X)}\le M$, and let $\mu$
be a positive constant such that $\La^{-1} \le \mu
\le \La$. Then there is a unique solution
$v\in (C^0\cap L^2_2)(\Lambda^+\otimes\fg_E)$ to equation
\eqref{eq:QuickExtASDEqnForv} such that $\Pi_{A,\mu}v=0$ and
\begin{align}
\|v\|_{L^2_{1,A}(X)} &\le C\|F_A^+\|_{L^{4/3}(X)}, \tag{1}\\
\|v\|_{L^2_{2,A}(X)} + \|v\|_{L^\8(X)} 
&\le C\|F_A^+\|_{L^{\sharp,2}(X)}. \tag{2}
\end{align}
Moreover, if $A$ is an $L^2_l$ connection, for any $l\geq 4$, then the
solution $v$ is contained in $L^2_{l+1}(\Lambda^+\otimes\fg_E)$.
\end{prop}

\begin{pf}
We try to solve equation \eqref{eq:LongSecOrderExtASDEqnForv} for solutions
of the form $v=G_{A,\mu}\xi$, where $\xi\in
L^{\sharp,2}(\Lambda^+\otimes\fg_E)$ and $G_{A,\mu}$ is the Green's
operator for the Laplacian $d_A^+d_A^{+,*}$ on
$\Ran\Pi_{A,\mu}^\perp$. Therefore, as explained in \S
\ref{subsec:ExtendedASD}, we seek solutions $\xi$ to equation
\eqref{eq:IntegralASDEqnForxi}.  We now apply Lemma
\ref{lem:FixedPoint} to equation \eqref{eq:IntegralASDEqnForxi}, with
$\fB=L^{\sharp,2}(\Lambda^+\otimes\fg_E)$, choosing $y=-F_A^+$ and $q(\xi)
= (d_A^{+,*}G_{A,\mu}\xi\wedge d_A^{+,*}G_{A,\mu}\xi)^+$, so our goal is to
solve
$$
\xi + q(\xi) = -F_A^+\quad\text{on }L^{\sharp,2}(\Lambda^+\otimes\fg_E). 
$$
For any $\xi_1,\xi_2\in L^{\sharp,2}(\Lambda^+\otimes\fg_E)$, 
the estimate in Proposition \ref{prop:EstPA0mu} yields
$$
\|d_A^{+,*}G_{A,\mu}\xi_i\|_{L^2_{1,A}} \le C\|\xi_i\|_{L^{\sharp,2}}.
$$
Therefore, applying the preceding estimates, H\"older's inequality, Lemma
\ref{lem:Kato}, and the $L^\sharp$-family of embedding and multiplication
results of \cite[Lemmas 4.1 \& 4.3]{FeehanSlice}, we obtain
$$
\|q(\xi_1)-q(\xi_2)\|_{L^{\sharp,2}}
\le K_2\left(\|\xi_1\|_{L^{\sharp,2}} + \|\xi_2\|_{L^{\sharp,2}}\right) 
\|\xi_1-\xi_2\|_{L^{\sharp,2}},
$$
where $K_2 = cZ_2$ with $c$ a universal constant and $Z_2$ the constant of 
Proposition \ref{prop:EstPA0mu}. 

Thus, provided $\eps \le (10K_2)^{-1}$, we have $\|F_A^+\|_{L^{\sharp,2}} <
(10K_2)^{-1}$, and Lemma \ref{lem:FixedPoint} implies that there is a
unique solution $\xi\in L^{\sharp,2}(\Lambda^+\otimes\fg_E)$ to equation 
\eqref{eq:IntegralASDEqnForxi} 
such that $\|\xi\|_{L^{\sharp,2}}\le (5K_2)^{-1}$.
The $L^{4/3}$ and $L^{\sharp,2}$ estimates for $\xi$,
\begin{equation}
\label{eq:LpEstForXi}
\|\xi\|_{L^{4/3}} \le 2\|F_A^+\|_{L^{4/3}} \quad\text{and}\quad
\|\xi\|_{L^{\sharp,2}} \le 2\|F_A^+\|_{L^{\sharp,2}},
\end{equation}
follow from Proposition \ref{prop:EstPA0mu} and a
standard rearrangement argument.
To obtain the stated $L^2_{1,A}$ and $L^2_{2,A}$ estimates for $v$, observe
that 
\begin{equation}
\label{eq:LaplacianOfvWithEigenvalCutoff}
d_A^+d_A^{+,*}v = d_A^+d_A^{+,*}G_{A,\mu}\xi = \Pi_{A,\mu}^\perp\xi,
\end{equation}
so that the $L^2_{1,A}$ bound is obtained from
\begin{align*}
\|v\|_{L^2_{1,A}} &\leq C\|\xi\|_{L^{4/3}}
\quad\text{(by estimate \eqref{eq:L21EstForvInTermsL43Data})}
\\
&\le C\|F_A^+\|_{L^{4/3}}, 
\quad\text{(by estimate \eqref{eq:LpEstForXi}).}
\end{align*}
The $C^0\cap L^2_{2,A}$ estimate similarly follows, as
\begin{align*}
\|v\|_{C^0\cap L^2_{2,A}} &\leq C\|\xi\|_{L^{\sharp,2}}
\quad\text{(by estimate \eqref{eq:C0L22EstForvInTermsLSharp2Data})}
\\
&\le C\|F_A^+\|_{L^{\sharp,2}} 
\quad\text{(by estimate \eqref{eq:LpEstForXi}).}
\end{align*}
The final regularity assertion follows from
Proposition \ref{prop:RegularityTaubesSolution}, our regularity result for
gluing solutions to the (extended) anti-self-dual equation.
\end{pf}

\subsection{Proof of main theorem, without Uhlenbeck continuity assertions}
\label{subsec:ProofGluingTheoremWoUhlenbeckContinuity}
At this stage, we are ready to complete the proof of at least part of
Theorem \ref{thm:GluingTheorem}, namely the existence of gluing and
obstruction maps, but not the Uhlenbeck continuity assertions for these
maps, which requires considerably more work and is taken up in \S
\ref{sec:Continuity}. On the other hand, the proof of this first part of
Theorem \ref{thm:GluingTheorem} is an essentially formal, almost immediate
consequence of Proposition \ref{prop:ASDExist}.

\begin{proof}[Proof of Theorem \ref{thm:GluingTheorem}, without Uhlenbeck
continuity assertions] 
According to Proposition \ref{prop:ASDExist} there are well-defined,
$\sG_E$-equivariant maps $\bdelta$ and $\bvarphi$ from the open subset of
$\sA_\kappa^w(X)\times C$, 
\begin{equation}
\begin{aligned}
\label{eq:IntrinsicDomainTaubesGluing}
\sA_\kappa^w(X,C;\mu,\eps)
&:=
\{(A,g)\in\sA_\kappa^w(X)\times C: 
\\
&\qquad\Spec(d_A^{+,g}d_A^{+,g,*}) \subset [0,\half\mu)\cup (\mu,\8),
\\
&\qquad\text{ and } \|F_A^{+,g}\|_{L^{\sharp,2}} < \eps\}
\\
\sB_\kappa^w(X,C;\mu,\eps) &:= \sA_\kappa^w(X,C;\mu,\eps)/\sG_E,
\end{aligned}
\end{equation}
given by 
\begin{align*}
&\bdelta:(A,g)\mapsto (A+d_A^{+,g,*}v,g),
\\
&\bvarphi:(A,g)\mapsto \Pi_{A,\mu}F^{+,g}(A+d_A^{+,g,*}v), 
\end{align*}
into $\sA_\kappa^w(X)\times C$ and $L^2(X,\Lambda^{+,g}\otimes\fg_E)$,
respectively. These descend to a map and a section,
\begin{equation}
\begin{aligned}
\label{eq:TaubesDeformationMapnSection}
&\bdelta:\sB_\kappa^w(X,C;\mu,\eps) \to \sB_\kappa^w(X)\times C,
\\
&\bvarphi:\sB_\kappa^w(X,C;\mu,\eps) \to \fV_\mu,
\end{aligned}
\end{equation}
where $\fV_\mu\to \sB_\kappa^w(X,C;\mu)$ is the 
continuous, finite-rank, pseudovector bundle 
$$
\fV_\mu
:= 
\{(A,g,w): (A,g)\in\sA_\kappa^w(X,C;\mu) 
\text{ and }w\in\Ran\Pi_{A,g,\mu}\}/\sG_E,
$$
where we recall that $\Ran\Pi_{A,g,\mu}\subset
L^2_{k-1}(X,\Lambda^{+,g}\otimes\fg_E)$. Here, the spaces
$\sA_\kappa^w(X,C;\mu)$ and $\sB_\kappa^w(X,C;\mu)$ are defined by simply
omitting the conditions in definition
\eqref{eq:IntrinsicDomainTaubesGluing} involving $\eps$.
The map $\bdelta$ and section $\bvarphi$ have the desired property that
$$
\bdelta(\bvarphi^{-1}(0)) \subset M_\kappa^w(X,C),
$$
and are smooth respect to the quotient $L^2_k$ topology upon restriction to
$\sB_\kappa^{w,*}(X)$: smoothness follows by the same argument as used in
the proof of smoothness in the case of the Banach-space inverse function
theorem
\cite{AMR}, where the existence of an inverse is established via a
contraction-mapping argument, just as we employ here.  

The gluing map, gluing-data obstruction bundle, and gluing-data obstruction
section are then defined by
\begin{align*}
\bga &= \bdelta\circ\bga':\Gl_\kappa^{w,+}(X,\sU_0,\Sigma,\lambda_0) 
\to \sB_\kappa^w(X)\times C,
\\
\bchi &= \bvarphi\circ\bga:\Gl_\kappa^{w,+}(X,\sU_0,\Sigma,\lambda_0)\to\Xi,
\\
\Xi &= \bga^*\fV_\mu,
\end{align*}
where the open neighborhood $\sU_0\subset M_{\kappa_0}^w(X,C)$ is chosen
according to Remark \ref{rmk:SettingForObstructions}.  The estimates in \S
\ref{sec:Decay} and \S \ref{sec:Eigenvalue2} ensure that the image of
$\bga'$, when $\bga'$ is restricted to
$\Gl_\kappa^{w,+}(X,\sU_0,\Sigma,\lambda_0)$, is contained in
$\sB_\kappa^{w,*}(X,C;\mu,\eps)$: in particular, that the {\em a priori}
estimates on $F_A^{+,g}$ are obeyed (by Proposition
\ref{prop:CutoffConnSelfDualCurvEst}) and that the spectrum of
$d_{A_0}^{+,g}d_{A_0}^{+,*,g}$ obeys the required open constraints for a
suitable $\mu$ (by Theorem \ref{thm:H2SmallEval}).  We can choose $\mu$ to
be the infimum of $\half\nu_2[A_0,g]$ (the least positive eigenvalue of the
Laplacian $d_{A_0}^{+,g}d_{A_0}^{+,*,g}$), as
$[A_0,g]\in\sB_{\kappa_0}^w(X)\times C$ ranges over the space of
gauge-equivalence classes of background connections (and metrics in $C$)
contained in $\Gl_\kappa^{w,+}(X,\sU_0,\Sigma,\lambda_0)$.
\end{proof}

\begin{rmk}
The rank of $\fV_\mu$ varies over $\sB_\kappa^w(X,C;\mu,\eps)$, so that
is why we only refer to it as being a `pseudovector bundle'. However, the
construction of $\Gl_\kappa^{w,+}(X,\sU_0,\Sigma,\lambda_0)$ (via the
choice of family $\sU_0$ of background connections $A_0$) in \S
\ref{sec:Splicing} and Theorem
\ref{thm:H2SmallEval} ensure that its image under $\bga'$ is
contained in an open subset of $\sB_\kappa^w(X,C;\mu,\eps)$ over which
$\fV_\mu$ has constant rank and is thus a vector bundle when restricted to
this open subset.
\end{rmk}


\section{Continuity of gluing maps and obstruction sections}
\label{sec:Continuity}
Our goal in this section is to prove the assertions in Theorem
\ref{thm:GluingTheorem} concerning the Uhlenbeck continuity of the gluing
map $\bga$ and continuity of the Kuranishi obstruction map $\bchi$.

As we briefly discussed in \S \ref{subsec:Introduction}, the statements of
Theorems 3.4.10, 3.4.17, and 3.4.18 in
\cite{FrM} describe a gluing result similar to but less general than
our main gluing theorem for anti-self-dual connections (see \S
\ref{subsec:Introduction} here and \cite[\S 1]{FLKM2}). While proofs of 
\cite[Theorems 3.4.10, 3.4.17, \& 3.4.18]{FrM}
are attributed by Friedman and Morgan \cite{FrM} to Taubes (specifically,
\cite{TauSelfDual} and \cite{TauIndef}), it is nonetheless fair to
point out that several important assertions of \cite[Theorems 3.4.10,
3.4.17, \& 3.4.18]{FrM} are not proved in \cite{TauSelfDual},
\cite{TauIndef} (or their sequels \cite{TauFrame}, \cite{TauStable}).
In the preceding four papers of Taubes, the gluing map $\bga$ is shown ---
either directly or by routine arguments or modifications of existing
arguments --- to have the following properties:
\begin{itemize}
\item
The gluing construction gives a map $\bga$ from
$\Gl_\kappa^w(X,\sU_0,K,\lambda_0)$ into $\barM_\kappa^w(X,C)$, which is
smooth upon restriction to smooth strata of the gluing-data bundle
$\Gl_\kappa^w(X,\sU_0,K,\lambda_0)$. (Here, $K\Subset\Sigma$ is a
precompact subset of a smooth stratum $\Sigma$ of the symmetric product
$\Sym^\ell(X)$ and $C\subset\Met(X)$ is a compact subset.)
\item
Points in $M_\kappa^w(X,C)$ are shown to lie in the image of certain
iterated gluing maps, whose definition is significantly different from that
of $\bga$.
\end{itemize}
In the versions stated above we have assumed, for simplicity, that the
``background connections'' $A_0$ in $\Gl_\kappa^w(X,\sU_0,K,\lambda_0)$ have
non-vanishing $H^2_{A_0}$ as the issue of obstructions to gluing --- though
a complication --- is not germane to the main concern of this section,
namely Uhlenbeck continuity. 

Contrary to the assertions of \cite{FrM}, it is not shown that
$\bga$ is actually a homeomorphism from
$\bar{\Gl}_\kappa^w(X,\bar\sU_0,K,\lambda_0)$ onto an open subset of
$\barM_\kappa^w(X,C)$, nor is it shown to be a diffeomorphism upon
restriction to smooth strata. Though desirable, these properties of
$\bga$ are not obvious and do not follow by routine modifications of
arguments contained in his papers or those, for example, of Donaldson
\cite{DonConn}, \cite{DK} or Mrowka
\cite{MrowkaThesis}.  

Our goal in this section is to complete the proof of one assertion of
\cite{FrM}, namely that the gluing maps extend continuously to the 
Uhlenbeck closure of the gluing data,
$\bar{\Gl}_\kappa^{w,+}(X,\bar\sU_0,\bar\Sigma,\lambda_0)$, as discussed in \S
\ref{subsec:Introduction}. Note that this continuity
assertion is considerably stronger than that of
\cite[Theorem 3.4.10]{FrM}, where $\Sigma$ has been replaced by a
precompact open subset $K\subset\Sigma$ and so the gluing sites are not
allowed to approach each other.  Our proof will also show that the gluing
maps are continuous with respect to bubble-tree limits. A special case of
the continuity assertion is treated by Donaldson and Kronheimer in \cite[\S
7.2.8]{DK}.

To show that our gluing maps are open, it suffices to consider their
restrictions to smooth strata: openness then follows from dimension
counting and the fact that the gluing maps are local diffeomorphisms
onto their images after restriction to smooth strata. These and other
properties are proved in \cite{FLKM2}.

\subsection{Statement of continuity results}
\label{subsec:ContinuityPrelim}
Let $A$, $A_0$ be connections on $\SO(3)$ bundles $\fg_E$, $\fg_{E_0}$,
respectively, over $X$ with $c_1(E)=c_1(E_0)$ and
$c_2(E)-c_2(E_0)=:\ell>0$. Though we make no specific use of this fact for
now, the case we have in mind is where
the point $[A]\in\sB_E(X,g)$ lies in an
Uhlenbeck neighborhood of the point $[A_0,\bx]\in
\sB_{E_0}(X,g)\times\Sym^{\ell}(X)$ (see Definition
\ref{defn:UhlenbeckNeighborhood}). 
Let $a\in L^2_k(X,\Lambda^1\otimes\fg_E)$ and $a_0\in
L^2_k(X,\Lambda^1\otimes\fg_{E_0})$ be solutions to the quasi-linear
equations
\eqref{eq:ASD-A} on the bundles $\fg_E$ and $\fg_{E_0}$, respectively,
\begin{align}
d_A^+a + (a\wedge a)^+ &= w \quad\text{and}\quad
a=d_A^{+,*}v, \label{eq:FirstOrderQuasiLinearA}\\
d_{A_0}^+a_0 + (a_0\wedge a_0)^+ &= w_0\quad\text{and}\quad
a_0=d_{A_0}^{+,*}v_0, \label{eq:FirstOrderQuasiLinearA0}
\end{align}
where $v,w\in\Gamma(X,\Lambda^+\otimes\fg_E)$ and
$v_0,w_0\in\Gamma(X,\Lambda^+\otimes\fg_{E_0})$. The use of $w$, $w_0$ on
the right-hand-sides of equations \eqref{eq:FirstOrderQuasiLinearA} and
\eqref{eq:FirstOrderQuasiLinearA0} allows us to give a common solution to a
wider range of Uhlenbeck continuity problems than would be possible if we
simply chose $w=-F_A^+$ or $w_0=F_{A_0}^+$.

We wish to compare the solutions $a$ and $a_0$ to equations
\eqref{eq:FirstOrderQuasiLinearA} and \eqref{eq:FirstOrderQuasiLinearA0}
when the connection $A$ is close to $(A_0,\bx)$ in the Uhlenbeck topology and
$w$ is close to $w_0$ in $L^2_{k-1,\loc}$ on $X\less\bx$. The purpose of this
section is to show that $a$ then approaches $a_0$ in $L^2_{k,\loc}$ on
$X\less\bx$; more precise bounds are given below.

As earlier we use $a(A):=d_A^{+,*}v$ and $a_0=a(A_0):=d_{A_0}^{+,*}v_0$ to
rewrite the first-order quasi-linear equations
\eqref{eq:FirstOrderQuasiLinearA} and
\eqref{eq:FirstOrderQuasiLinearA0} for $a$ and $a_0$ as
\begin{align}
\label{eq:SecondOrderQuasiLinearA}
d_A^+d_A^{+,*}v + (d_A^{+,*}v\wedge d_A^{+,*}v)^+ &= w, 
\\
\label{eq:SecondOrderQuasiLinearA0}
d_{A_0}^+d_{A_0}^{+,*}v_0 + (d_{A_0}^{+,*}v_0\wedge d_{A_0}^{+,*}v_0)^+ &=
w_0,  
\end{align}
to give a pair of second-order, quasi-linear elliptic equations for
$v$ and $v_0$. In our applications, the terms $w$ and $w_0$ will
always depend on $A$ and $A_0$ in such a way as to preserve gauge
invariance. For example, when the Laplacian $d_A^+d_A^{+,*}$ has no
small eigenvalues, we set $w(A):=-F_A^+$ in equation
\eqref{eq:FirstOrderQuasiLinearA} and so we recover the usual
anti-self-dual equation. If the Laplacian $d_A^+d_A^{+,*}$ does have
small eigenvalues, we require that $\Pi_{A,\mu}v=0$ and set
\begin{equation}
\label{eq:DefnOfwForExtASD}
w(A):=\Pi_{A,\mu}((d_A^{+,*}v\wedge d_A^{+,*}v)^+ + F_A^+) 
- F_A^+,
\end{equation}
and so recover the extended anti-self-dual equation 
 \eqref{eq:QuickExtASDEqnForv} or, in detail,
\eqref{eq:LongSecOrderExtASDEqnForv}, 
where $\Pi_{A,\mu}$ is the $L^2$-orthogonal projection from
$L^2(\Lambda^+\otimes\fg_E)$ onto the finite-dimensional subspace spanned
by the eigenvectors of $d_A^+d_A^{+,*}$ with eigenvalues in $[0,\mu]$,
where $\mu$ is a positive constant (an upper bound for the small
eigenvalues), usually $\half\nu_2[A_0]$: see the small-eigenvalue estimates
supplied by Theorem \ref{thm:H2SmallEval}.  Similarly for the
connection $A_0$.

We can now state our main results.

\begin{thm}
\label{thm:UhlContSecOrderQuasiLin}
Let $X$ be a closed, oriented, smooth four-manifold $X$ with Riemannian
metric $g$. Then there is a positive constant $\eps$ such that the
following holds.  Let $E$, $E_0$ be Hermitian vector bundles over $X$, such
that $c_1(E_0)=c_2(E)$ and $c_2(E_0)<c_2(E)$. Let $M<\8$ be a constant, let
$A_0$ be an orthogonal $L^2_k$ connection on $\fg_{E_0}$ and suppose
$\{A_\alpha\}$ is a sequence of $L^2_k$ connections on $\fg_E$ such that
\begin{equation}
\label{eq:UniformF+BoundOnConnSequence}
\sup_\alpha\|F_{A_\alpha}\|_{L^2(X)} \leq M
\quad\text{and}\quad
\|F_{A_0}^+\|_{L^{\sharp,2}(X)}, \quad
\sup_\alpha\|F_{A_\alpha}^+\|_{L^{\sharp,2}(X)} < \eps,
\end{equation}
and which converges to $(A_0,\bx)\in\sA_{E_0}\times\Sym^\ell(X)$ in
the sense that
\begin{itemize}
\item
The sequence of connections $A_\alpha$ converges to $A_0$ in
$L^2_{k,A_0,\loc}$ on $X\less\bx$, and
\item
The sequence of measures $|F_{A_\alpha}|^2$ converges to
$|F_{A_0}|^2+8\pi^2\sum_{x\in\bx}\delta_x$ on $X$.
\end{itemize}
We assume that $E|_{X\less\bx} = E_0|_{X\less\bx}$.
Suppose $v_0 \in L^2_{k+1}(\Lambda^+\otimes\fg_{E_0})$ is the solution
to the second-order quasi-linear equation
\eqref{eq:SecondOrderQuasiLinearA} on $\fg_{E_0}$, with right-hand side $w_0$,
while $v_\alpha \in L^2_{k+1}(\Lambda^+\otimes\fg_E)$ are
the solutions to the second-order quasi-linear equation
\eqref{eq:SecondOrderQuasiLinearA} on $\fg_E$, with right-hand side $w_\alpha$.
Suppose that the sequence $\{w_\alpha\}$ converges to $w_0$ in
$L^2_{k-1,A_0,\loc}(X\less\bx)$.  Then the following hold:
\begin{itemize}
\item
The sequence of solutions $v_\alpha$ to equation
\eqref{eq:SecondOrderQuasiLinearA} converges to the solution $v_0$ to
equation \eqref{eq:SecondOrderQuasiLinearA0} in $L^2_{k+1,A_0,\loc}$ on
$X\less\bx$,
\item
The sequence of solutions $a_\alpha = d_{A_\alpha}^{+,*}v_\alpha$ to
equation \eqref{eq:FirstOrderQuasiLinearA} converges to the solution
$a_0 = d_{A_0}^{+,*}v_0$ to \eqref{eq:FirstOrderQuasiLinearA0} in
$L^2_{k,A_0,\loc}$ on $X\less\bx$,
\item
The sequence of connections $A_\alpha+a_\alpha$ converges to $A_0+a_0$ in
$L^2_{k,A_0,\loc}$ on $X\less\bx$, and
\item
The sequence of measures $|F(A_\alpha+a_\alpha)|^2$
converges to $|F(A_0+a_0)|^2+8\pi^2\sum_{x\in\bx}\delta_x$ on $X$.
\end{itemize}
\end{thm}

The preceding result has the following consequence when applied to the
extended anti-self-dual equation:

\begin{thm}
\label{thm:UhlContExtSecOrderASD}
Continue the hypotheses of Theorem \ref{thm:UhlContSecOrderQuasiLin}, but
now suppose $v_0 \in
\Pi_{{A_0},\mu}^\perp\cdot L^2_{k+1}(\Lambda^+\otimes\fg_{E_0})$ is
the solution to the extended second-order anti-self-dual equation
\eqref{eq:SecondOrderQuasiLinearA0} on $\fg_{E_0}$,
while $v_\alpha \in \Pi_{{A_\alpha},\mu}^\perp\cdot
L^2_{k+1}(\Lambda^+\otimes\fg_E)$ are the solutions to the extended
second-order anti-self-dual equation
\eqref{eq:SecondOrderQuasiLinearA} on $\fg_E$, where we now have
\begin{align*}
w_\alpha := w(A_\alpha) = \Pi_{{A_\alpha},\mu}(
(d_{A_\alpha}^{+,*}v_\alpha\wedge d_{A_\alpha}^{+,*}v_\alpha)^+ 
+ F_{A_\alpha}^+) - F_{A_\alpha}^+,
\\
w_0 := w(A_0) = \Pi_{A_0,\mu}(
(d_{A_0}^{+,*}v_0\wedge d_{A_0}^{+,*}v_0)^+ + F_{A_0}^+) - F_{A_0}^+.
\end{align*}
Then the conclusions of Theorem \ref{thm:UhlContSecOrderQuasiLin} hold and,
moreover,
\begin{itemize}
\item
The sequence of terms $w_\alpha$ converges to the term $w_0$ in
$L^2_{k-1,A_0,\loc}$ on $X\less\bx$.
\end{itemize}
\end{thm}

Theorem \ref{thm:UhlContExtSecOrderASD} has the following important
consequence.

\begin{proof}[Proof of Uhlenbeck continuity assertions of Theorem
\ref{thm:GluingTheorem}]
The assertions in Theorem
\ref{thm:GluingTheorem} concerning the Uhlenbeck continuity of the gluing
map $\bga$ and continuity of the Kuranishi obstruction map $\bchi$ follow
immediately from the definition of the splicing map $\bga'$ and Theorem
\ref{thm:UhlContExtSecOrderASD}. 
\end{proof}

The remainder of this section is taken up with the proofs of Theorems
\ref{thm:UhlContSecOrderQuasiLin} and \ref{thm:UhlContExtSecOrderASD}.

\subsection{Degeneration of solutions to the second-order quasi-linear
equation: linear \boldmath{$L^2_2$} theory} 
\label{subsec:DegenerationLinearTheory}
Prior to analyzing the asymptotic behavior of
solutions to the full anti-self-dual equation, it is useful to first
consider the asymptotic behavior of solutions 
$v\in L^2_{k+1}(X,\Lambda^+\otimes\fg_E)$ to the linear part of the
second order, elliptic anti-self-dual equation \eqref{eq:SecondOrderASD}, so
\begin{equation}
d_A^+d_A^{+,*}v = w, \label{eq:LinearSecondOrderASD}
\end{equation}
where $w\in L^2_{k-1}(X,\Lambda^+\otimes\fg_E)$ and the point
$[A]\in\sB_E(X,g)$ converges to the point $[A_0,\bx]\in
\sB_{E_0}(X,g)\times\Sym^{\ell}(X)$.

We need to compare the solutions $v$ to equation
\eqref{eq:LinearSecondOrderASD} as $A$ converges to $A_0$ and $w$ converges
to $w_0$ in $L^2_{k-1,A_0,\loc}$ on $X\less\bx$. For this purpose we use an
elementary Leibnitz formula for $d_A^+d_A^{+,*}$. We recall from
\eqref{eq:BW+} that $d_A^+d_A^{+,*} = \half\cov_A^*\cov_A +
\text{zeroth-order terms}$ on $L^2_{k-1}(X,\Lambda^+\otimes\fg_E)$ and so
we have, schematically, for any function $\psi_0\in C^\8(X)$,
\begin{equation}
\label{eq:LaplacianLeibnitz}
d_A^+d_A^{+,*}(\psi_0 v) 
= 
\cov^2\psi_0 \otimes v + \cov\psi_0\otimes\cov_A v + \psi_0 d_A^+d_A^{+,*} v.
\end{equation}
The `excision' argument which allows us to compare solutions to the
anti-self-dual equation on different bundles $E$ and $E_0$ relies on
a suitable choice of cutoff functions. 

Over $X\less\bx$ we may write $A = A_0 + (A-A_0)$ on 
$E|_{X\less\bx} = E_0|_{X\less\bx}$. Though we shall ultimately require a
pair of sequences of $\bx$-good cutoff functions satisfying the properties
of Lemma \ref{lem:GoodCutoffSequence}, it suffices for now to choose
smooth cutoff functions $\chi_0$ and $\psi_0$ on $X$ such that
\begin{itemize}
\item
$\supp\chi_0\Subset X\less\bx$, $\supp\psi_0\Subset X\less\bx$, 
\item
$\chi_0=1$ on $\supp\psi_0$,
\end{itemize}
and define open subsets of $X$ by
$$
\Om := \{x\in X: (d\psi_0)_x\ne 0\}
\quad\text{and}\quad
U := \{x\in X: \psi_0(x) > 0\}.
$$
We use the $\chi$'s to cut off and compare connections and the $\psi$'s to
cut off and compare solutions to the anti-self-dual equation on different
bundles. Noting that $E|_{X\less\bx} = E_0|_{X\less\bx}$ and setting
\begin{equation}
\label{eq:DefnOfAminusA0Cutoff}
b := \chi_0(A - A_0) \in L^2_k(X,\Lambda^1\otimes\fg_{E_0}),
\end{equation}
we have $A = A_0 + b$ over $\supp\psi_0\subset X$ since $\chi_0 = 1$ on
$\supp\psi_0$; the section $b$ may also be viewed as an element of
$L^2_k(X,\Lambda^1\otimes\fg_E)$.

The following comparison of solutions to the linear part of the
anti-self-dual equation is the basic prototype we later use to compare
solutions to the full anti-self-dual equation on different bundles $\fg_E$
and $\fg_{E_0}$.

\begin{lem}
\label{lem:LaplacianDegEst}
Let $(X,g)$ be a closed, oriented, smooth, Riemannian four-manifold. Then
there is a positive constant $\eps$ with the following significance.
Let $A$, $A_0$ be connections on $\SO(3)$ bundles $\fg_E$, $\fg_{E_0}$,
respectively, with $\|F_{A_0}\|_{L^{\sharp,2}(X)} < \eps$,
and assume that $\Ker d_{A_0}^+d_{A_0}^{+,*}=0$.  Suppose
$U\subset X$ is an open subset such that $E|_U = E_0|_U$. Let $\psi_0$,
$\chi_0$ be cutoff functions on $X$ such that $\supp\chi_0\subset U$,
$\supp\psi_0\subset U$, and $\chi_0=1$ on $\supp\psi_0$.  Given $w\in
L^2_{k-1}(X,\Lambda^+\otimes\fg_E)$ and $w_0\in
L^2_{k-1}(X,\Lambda^+\otimes\fg_{E_0})$, suppose $v\in
L^2_{k+1}(X,\Lambda^+\otimes\fg_E)$ and $v_0\in
L^2_{k+1}(X,\Lambda^+\otimes\fg_{E_0})$ satisfy
$$
d_A^+d_A^{+,*}v = w \quad\text{and}\quad d_{A_0}^+d_{A_0}^{+,*}v_0 = w_0.
$$
Then, 
\begin{align*}
&\|\psi_0(v - v_0)\|_{L^\8(X)} + \|\psi_0(v - v_0)\|_{L^2_{2,A_0}(X)} 
\\
&\leq 
C\|\psi_0(w - w_0)\|_{L^{\sharp,2}} 
\\
&\quad + C(\|\cov^2\psi_0\|_{L^{\sharp,2}}+ \|b\|_{L^{2\sharp,4}}
+\|\cov_{A_0} b\|_{L^{\sharp,2}} + \|b\|_{L^{2\sharp,4}}^2
\\
&\quad + \|\cov\psi_0\|_{L^{2\sharp,4}}\|b\|_{L^{2\sharp,4}}
+\|\cov\psi_0\|_{L^{2\sharp,4}})
\\
&\qquad \times (\|v\|_{L^\8} + \|\cov_A v\|_{L^{2\sharp,4}}
+ \|v_0\|_{L^\8}+ \|\cov_{A_0} v_0\|_{L^{2\sharp,4}}),
\end{align*}
where $b := \chi_0(A - A_0) \in L^2_k(X,\Lambda^1\otimes\fg_{E_0})$, and
$C(\nu_2[A_0],\|F_{A_0}\|_{L^2})$ is a constant, and
$\nu_2[A_0]$ is the least positive eigenvalue of the Laplacian
$d_{A_0}^+d_{A_0}^{+,*}$.
\end{lem}

\begin{pf}
We first compute a suitable $L^{\sharp,2}$ estimate for
$d_{A_0}^+d_{A_0}^{+,*}\psi_0(v-v_0)$ and then apply our elliptic estimates
for $d_{A_0}^+d_{A_0}^{+,*}$ given in Lemma \ref{lem:LinftyL22CovLapEstv}.
We have $A = A_0 + b$ over $\supp\psi_0\subset X$, since $\chi_0 = 1$ on
$\supp\psi_0$, and thus on $\supp\psi_0$, we see that
\begin{align*}
d_{A_0}^+d_{A_0}^{+,*}v 
&= -d_{A_0}^+*d_{A_0}v \\
&= -d_A^+*d_{A_0}v + (b\wedge *d_{A_0}v)^+ \\
&= -d_A^+*d_Av + d_A^+*(b\wedge v)
+ (b\wedge *d_A v)^+ 
- (b\wedge *(b\wedge v))^+,
\end{align*}
and so
\begin{equation}
\label{eq:ExpandLaplacianConnection}
d_{A_0}^+d_{A_0}^{+,*}v 
= 
d_A^+d_A^{+,*}v + d_A^+*(b\wedge v)
+ (b\wedge *d_A v)^+ - (b\wedge *(b\wedge v))^+.
\end{equation}
According to our Leibnitz formula \eqref{eq:LaplacianLeibnitz} we have
\begin{equation}
\label{eq:ExpandLaplacianCutoffDiff}
\begin{aligned}
d_{A_0}^+d_{A_0}^{+,*}\psi_0(v-v_0)
&= \cov^2\psi_0\otimes (v-v_0) 
+ \cov\psi_0\otimes\cov_{A_0}(v-v_0) 
\\
&\quad + \psi_0d_{A_0}^+d_{A_0}^{+,*}(v-v_0).
\end{aligned}
\end{equation}
Therefore, combining \eqref{eq:ExpandLaplacianConnection} and
\eqref{eq:ExpandLaplacianCutoffDiff}, we obtain
\begin{equation}
\label{eq:LaplacianOnCutoffDiff}
\begin{aligned}
d_{A_0}^+d_{A_0}^{+,*}\psi_0(v-v_0) 
&= 
\cov^2\psi_0\otimes(v-v_0) + \cov\psi_0\otimes\cov_{A_0}(v-v_0)
\\
&\quad + \psi_0(d_A^+d_A^{+,*}v - d_{A_0}^+d_{A_0}^{+,*}v_0) 
+ \psi_0d_A^+*(b\wedge (v - v_0))
\\
&\quad 
+ \psi_0(b\wedge *d_A (v - v_0))^+ - \psi_0(b\wedge *(b\wedge (v - v_0)))^+.
\end{aligned}
\end{equation}
Consequently, as $d_A^+d_A^{+,*}v = w$ and $d_{A_0}^+d_{A_0}^{+,*}v_0 =
w_0$, we have the schematic expression, employing tensor products to
denote all bilinear products and using
$\cov_{A_0}(v-v_0) = \cov_Av + b\otimes v - \cov_{A_0}v_0$ in
the identity \eqref{eq:LaplacianOnCutoffDiff},  
\begin{equation}
\label{eq:LapaceDiffSoln}
\begin{aligned}
d_{A_0}^+d_{A_0}^{+,*}\psi_0(v - v_0) 
&= \psi_0(w - w_0) + \cov^2\psi_0\otimes (v-v_0)
\\
&\quad + \cov\psi_0\otimes (\cov_Av - \cov_{A_0}v_0 - b\otimes v_0)) 
\\
&\quad + \psi_0\cov_{A_0} b\otimes(v-v_0)          
\\
&\quad + \psi_0b\otimes (\cov_Av - \cov_{A_0}v_0)
+ \psi_0 b\otimes b\otimes(v+v_0).  
\end{aligned}
\end{equation}
We need an $L^{\sharp,2}$ bound on $d_{A_0}^+d_{A_0}^{+,*}\psi_0(v - v_0)$
in terms of the $L^{\sharp,2}$ norm of $\psi_0(w - w_0)$ and the lower order
terms in the expression \eqref{eq:LapaceDiffSoln}
for $d_{A_0}^+d_{A_0}^{+,*}\psi_0(v - v_0)$:
\begin{align}
\notag
&\|d_{A_0}^+d_{A_0}^{+,*}\psi_0(v - v_0)\|_{L^{\sharp,2}}
\\
\label{eq:Laplacianof CutoffDifference}
&\leq 
\|\psi_0(w - w_0)\|_{L^{\sharp,2}} 
+ (\|\cov^2\psi_0\|_{L^{\sharp,2}}+ \|b\|_{L^{2\sharp,4}}
+\|\cov_{A_0} b\|_{L^{\sharp,2}} + \|b\|_{L^{2\sharp,4}}^2
\\
\notag
&\qquad + \|\cov\psi_0\|_{L^{2\sharp,4}}\|b\|_{L^{2\sharp,4}})
(\|v\|_{L^\8}+\|v_0\|_{L^\8}) 
\\
\notag
&\quad + (\|\cov\psi_0\|_{L^{2\sharp,4}}+ \|b\|_{L^{2\sharp,4}})
(\|\cov_A v\|_{L^{2\sharp,4}}
+ \|\cov_{A_0} v_0\|_{L^{2\sharp,4}}). 
\end{align}
{}From Lemma \ref{lem:LinftyL22CovLapEstv} and the eigenvalue bound, 
$$
\|\psi_0(v-v_0)\|_{L^2} 
\leq 
\nu_2[A_0]^{-1}\|d_{A_0}^+d_{A_0}^{+,*}\psi_0(v-v_0)\|_{L^2},
$$
we have the elliptic estimate
\begin{equation}
\label{eq:EllipticDifferenceEst}
\|\psi_0(v-v_0)\|_{L^\8} + \|\psi_0(v-v_0)\|_{L^2_{2,A_0}}
\le
C\|d_{A_0}^+d_{A_0}^{+,*}\psi_0(v-v_0)\|_{L^{\sharp,2}},
\end{equation}
for a constant $C(\nu_2[A_0],\|F_{A_0}\|_{L^2})$.
Combining the bound \eqref{eq:EllipticDifferenceEst} with that of
\eqref{eq:Laplacianof CutoffDifference} then yields the 
desired estimate.
\end{pf}

\begin{rmk}
Lemma \ref{lem:LaplacianDegEst} allows us to compare solutions $v\in
L^2_{k+1}(X,\Lambda^+\otimes\fg_E)$ to the second order, linear equation
$d_A^+d_A^{+,*}v = w$ with solutions $v_0\in
L^2_{k+1}(X,\Lambda^+\otimes\fg_{E_0})$ to $d_{A_0}^+d_{A_0}^{+,*}v_0 = w_0$
when the point $[A]$ lies in an Uhlenbeck neighborhood of a point
$[A_0,\bx]$. Otherwise said, it allows us to compare the second-order
Green's operator for $d_A^+d_A^{+,*}$ with the Green's operator for
$d_{A_0}^+d_{A_0}^{+,*}$ as the point $[A]$ approaches $[A_0,\bx]$.
\end{rmk} 

It is a simpler matter to bound the $L^2_1$ norm of the difference between
$v$ and $v_0$ over compact subsets of $X\less\bx$:

\begin{lem}
\label{lem:L21LaplacianDegEst}
Continue the notation and hypotheses of
of Lemma \ref{lem:LaplacianDegEst}. Then, 
\begin{align*}
&\|\psi_0(v - v_0)\|_{L^2_{1,A_0}(X)} 
\\
&\leq C(\|\psi_0(w - w_0)\|_{L^{4/3}} +
\|\cov^2\psi_0\|_{L^{4/3}}(\|v\|_{L^\8}+\|v_0\|_{L^\8})) 
\\
&\quad 
+ C\|\cov\psi_0\|_{L^2}(\|\cov_Av\|_{L^4}
+ \|\cov_{A_0}v_0\|_{L^4} + \|b\|_{L^4}\|v_0\|_{L^\8})
\\
&\quad + C\|\cov_{A_0} b\|_{L^{4/3}}(\|\psi_0(v-v_0)\|_{L^\8}          
\\
&\quad + C(\|b\|_{L^2}(\|\cov_Av\|_{L^4} + \|\cov_{A_0}v_0)\|_{L^4}
+ \|b\|_{L^2}\|b\|_{L^4}(\|v\|_{L^\8}+\|v_0\|_{L^\8})).  
\end{align*}
where $b := \chi_0(A - A_0) \in L^2_k(X,\Lambda^1\otimes\fg_{E_0})$, and
$C(\nu_2[A_0],\|F_{A_0}\|_{L^2})$ is a constant, and
$\nu_2[A_0]$ is the least positive eigenvalue of the Laplacian
$d_{A_0}^+d_{A_0}^{+,*}$.
\end{lem}

\begin{proof}
According to the elliptic estimates of Lemma \ref{lem:L21AEstv}, we have
$$
\|\psi_0(v-v_0)\|_{L^2_{1,A_0}}
\leq 
C\|d_{A_0}^+d_{A_0}^{+,*}\psi_0(v-v_0)\|_{L^{4/3}}.
$$ 
On the other hand, our schematic expression \eqref{eq:LapaceDiffSoln} for 
$d_{A_0}^+d_{A_0}^{+,*}\psi_0(v-v_0)$ yields
\begin{align*}
\|d_{A_0}^+d_{A_0}^{+,*}\psi_0(v - v_0)\|_{L^{4/3}} 
&\leq \|\psi_0(w - w_0)\|_{L^{4/3}} +
\|\cov^2\psi_0\|_{L^{4/3}}(\|v\|_{L^\8}+\|v_0\|_{L^\8}) 
\\
&\quad 
+ \|\cov\psi_0\|_{L^2}(\|\cov_Av\|_{L^4}
+ \|\cov_{A_0}v_0\|_{L^4} + \|b\|_{L^4}\|v_0\|_{L^\8}) 
\\
&\quad + \|\cov_{A_0} b\|_{L^{4/3}}\|\psi_0(v-v_0)\|_{L^\8}          
\\
&\quad + \|b\|_{L^2}(\|\cov_Av\|_{L^4} + \|\cov_{A_0}v_0)\|_{L^4}
\\
&\quad + \|b\|_{L^2}\|b\|_{L^4}(\|v\|_{L^\8}+\|v_0\|_{L^\8}).  
\end{align*}
Combining the preceding estimates gives the desired bound.
\end{proof}

It remains to give a bound for the difference $\psi_0(a-a_0)$.

\begin{lem}
\label{lem:FirstOrderASDDegEst}
Continue the notation and hypotheses of
of Lemma \ref{lem:LaplacianDegEst}. Given $w\in
L^2_{k-1}(X,\Lambda^+\otimes\fg_E)$ and $w_0\in
L^2_{k-1}(X,\Lambda^+\otimes\fg_{E_0})$, suppose $v\in
L^2_{k+1}(X,\Lambda^+\otimes\fg_E)$ and $v_0\in
L^2_{k+1}(X,\Lambda^+\otimes\fg_{E_0})$. Let
$a(A):=d_A^{+,*}v\in L^2_k(X,\Lambda^1\otimes\fg_E)$ and
$a_0:=a(A_0)=d_{A_0}^{+,*}v_0\in L^2_k(X,\Lambda^1\otimes\fg_{E_0})$. Then
\begin{align*}
\|\psi_0(a - a_0)\|_{L^2_{1,A_0}} 
&\le \|\psi_0(v - v_0)\|_{L^2_{2,A_0}} 
\\
&\quad + (\|b\|_{L^2_{1,A_0}}
+ \|b\|_{L^4}\|\cov\psi_0\|_{L^4} + \|\cov\psi_0\|_{L^2_1})
(\|v_0\|_{L^\8} + \|v\|_{L^\8})
\\ 
&\quad + (\|\cov\psi_0\|_{L^4} + \|b\|_{L^4})
(\|\cov_{A_0}v\|_{L^4} + \|\cov_Av\|_{L^4} + \|b\|_{L^4}\|v\|_{L^\8}).
\end{align*}
where $b = \chi_0(A - A_0) \in L^2_k(X,\Lambda^1\otimes\fg_{E_0})$. 
\end{lem}

\begin{pf}
Schematically, we have 
\begin{align}
\notag
\psi_0(a - a_0) 
&= \psi_0(\cov_Av - \cov_{A_0}v_0) 
\\
\label{eq:DiffSolnAsdDiffBundles}
&= \cov_A\psi_0(v-v_0) - \cov\psi_0\otimes(v-v_0)
\\
\notag
&= \cov_{A_0}\psi_0(v-v_0) 
+ b\otimes\psi_0(v-v_0) - \cov\psi_0\otimes(v-v_0).
\end{align}
Thus, we can bound the $L^2_{1,A_0}$ norm of $\psi_0(a - a_0)$ in terms of 
the $L^2_{1,A_0}$ norm of $\cov_{A_0}(\psi_0(v - v_0))$ and the lower-order
terms in the expression \eqref{eq:DiffSolnAsdDiffBundles}
for $\psi_0(a - a_0)$. For the $L^2$ norm we have
\begin{equation}
\label{eq:L2EstDiffSolnAsdDiffBundles}
\begin{aligned}
\|\psi_0(a - a_0)\|_{L^2} 
&\leq  \|\cov_{A_0}\psi_0(v - v_0)\|_{L^2}
\\
&\quad + (\|b\|_{L^2} + \|\cov\psi_0\|_{L^2})(\|v\|_{L^\8} + \|v_0\|_{L^\8}). 
\end{aligned}
\end{equation}
To compute the $L^2$ norm of $\cov_{A_0}\psi_0(a - a_0)$, note that,
schematically, the identity \eqref{eq:DiffSolnAsdDiffBundles} yields
\begin{align*}
\cov_{A_0}\psi_0(a - a_0) 
&=  \cov_{A_0}^2\psi_0(v - v_0)
+ \cov_{A_0}b\otimes \psi_0(v - v_0) 
\\
&\quad + b\otimes \cov \psi_0\otimes (v - v_0) 
+ b\otimes\psi_0\cov_{A_0}(v - v_0) 
\\
&\quad - \cov^2 \psi_0\otimes (v - v_0) 
- \cov \psi_0\otimes\cov_{A_0}(v - v_0),
\end{align*}
with $\cov_{A_0}v = \cov_Av - b\otimes v$. Therefore,
\begin{align*}
\|\cov_{A_0}\psi_0(a - a_0)\|_{L^2} 
&\le \|\cov_{A_0}^2\psi_0(v - v_0)\|_{L^2}
\\
&\quad + (\|\cov_{A_0} b\|_{L^2}
+ \|b\|_{L^4}\|\cov\psi_0\|_{L^4} + \|\cov^2\psi_0\|_{L^2})
(\|v_0\|_{L^\8} + \|v\|_{L^\8})
\\ 
&\quad + (\|\cov\psi_0\|_{L^4} + \|b\|_{L^4})
(\|\cov_{A_0}v\|_{L^4} + \|\cov_Av\|_{L^4} + \|b\|_{L^4}\|v\|_{L^\8}).
\end{align*}
Combining this estimate with the $L^2$ bound
\eqref{eq:L2EstDiffSolnAsdDiffBundles} then
gives the desired result.
\end{pf}

Again, the $L^2$ bound is simpler:

\begin{lem}
\label{lem:L2FirstOrderASDDegEst}
Continue the notation of Lemma \ref{lem:FirstOrderASDDegEst}. Then,
\begin{align*}
\|\psi_0(a - a_0)\|_{L^2(X)} 
&\le \|\psi_0(v - v_0)\|_{L^2_{1,A_0}} 
\\
&\quad + (\|b\|_{L^2} + \|\cov\psi_0\|_{L^2})\|\chi_0(v-v_0)\|_{L^\8}. 
\end{align*}
\end{lem}

\begin{proof}
  {}From the schematic expression \eqref{eq:DiffSolnAsdDiffBundles}
  for $\psi_0(a - a_0)$ we see that
\begin{align*}
\|\psi_0(a - a_0)\|_{L^2} 
&\leq  \|\cov_{A_0}\psi_0(v - v_0)\|_{L^2}
\\
&\quad + (\|b\|_{L^2} + \|\cov\psi_0\|_{L^2})\|\chi_0(v-v_0)\|_{L^\8}. 
\end{align*}
This completes the proof.
\end{proof}

It remains to bound the Sobolev norms of $b$ with estimates for  Sobolev
norms of $A-A_0$ and $\chi_0$.

\begin{lem}
\label{lem:L21BoundForbInTermsOfA-A0}
Continue the above notation. Then the following estimates hold,
\begin{align}
\tag{1}
\|b\|_{L^{2\sharp,4}(X)}
&\leq c\|b\|_{L^2_{1,A_0}(X)},
\\
\tag{2}
\|b\|_{L^2_{1,A_0}(X)}
&\leq (1+\|\cov\chi_0\|_{L^2(X)})\|A-A_0\|_{L^2_{1,A_0}(U)},
\\
\tag{3}
\|\cov_{A_0}b\|_{L^{\sharp,2}(X)}
&\leq
C(1+\|\cov\chi_0\|_{L^4(X)})\|A-A_0\|_{L^2_{2,A_0}(U)},
\end{align}
where $c$ depends only on the metric $g$, while $C$ depends on $U$ and $g$.
\end{lem}

\begin{proof}
Assertion 1 follows immediately from 
\cite[Lemma 4.1]{FeehanSlice}. Next,
using $\cov_{A_0}\chi_0(A-A_0)=\cov\chi_0\otimes(A-A_0) +
\chi_0\cov_{A_0}(A-A_0)$, we see that
$$
\|b\|_{L^2_{1,A_0}(X)}
\leq
(1+\|\cov\chi_0\|_{L^2(X)})\|A-A_0\|_{L^2(U)} +
\|\cov_{A_0}(A-A_0)\|_{L^2(U)}, 
$$
as required for Assertion 2. Lastly, for any $2<p<4$ and
$4<q<\8$ for which $1/p=1/4+1/q$, the Sobolev embeddings and H\"older
inequality yield
\begin{align*}
\|\cov_{A_0}b\|_{L^\sharp(X)}
&\leq c\|\cov_{A_0}b\|_{L^p(X)}
\\
&\leq c\|\cov\chi_0\|_{L^4(X)}\|A-A_0\|_{L^q(U)} +
\|\cov_{A_0}(A-A_0)\|_{L^p(U)}, 
\\
&\leq C\|\cov\chi_0\|_{L^4(X)}\|A-A_0\|_{L^2_{2,A_0}(U)} +
C\|\cov_{A_0}(A-A_0)\|_{L^2_{1,A_0}(U)}, 
\end{align*}
which completes the proof of Assertion 3.
\end{proof}

\begin{rmk}
\label{rmk:NoKernelAssumption}
Although we assumed $\Ker d_{A_0}^+d_{A_0}^{+,*}=0$ in Lemmas
\ref{lem:LaplacianDegEst} and \ref{lem:L21LaplacianDegEst}, we did
not need to explicitly need this assumption in Lemmas
\ref{lem:FirstOrderASDDegEst} and \ref{lem:L21BoundForbInTermsOfA-A0}. 
\end{rmk}

\subsection{Degeneration of solutions to the second-order quasi-linear
equation: non-linear \boldmath{$L^2_2$} theory} 
\label{subsec:ContinuityNonLinearL22}
We now consider the difference of the
non-linear terms in the anti-self-dual equation and, in conjunction with
Lemma \ref{lem:LaplacianDegEst}, we obtain the following comparison of
solutions to the anti-self-dual equation.

\begin{prop}
\label{prop:SecondOrderASDDegEst}
Let $(X,g)$ be a closed, oriented, smooth, Riemannian four-manifold. Then
there is a positive constant $\eps$ with the following significance.
Let $A$, $A_0$ be connections on $\SO(3)$ bundles $\fg_E$, $\fg_{E_0}$,
respectively, over $(X,g)$, with $\|F_{A_0}^+\|_{L^{2,\sharp}(X)} < \eps$.
Then there are positive constants $C=C(\nu_2[A_0],\|F_{A_0}\|_{L^2})$ and
$\delta=\delta(C)$, where $\nu_2[A_0]$ is the least positive eigenvalue of
the Laplacian $d_{A_0}^+d_{A_0}^{+,*}$, such that the following holds.
Suppose $U\subset X$ is an open subset such that $E|_U = E_0|_U$. Let
$\psi_0$, $\chi_0$ be cutoff functions on $X$ such that $\supp\chi_0\subset
U$, $\supp\psi_0\subset U$, and $\chi_0=1$ on $\supp\psi_0$. Given $w\in
L^2_{k-1}(X,\Lambda^+\otimes\fg_E)$ and $w_0\in
L^2_{k-1}(X,\Lambda^+\otimes\fg_{E_0})$, suppose $v\in
L^2_{k+1}(X,\Lambda^+\otimes\fg_E)$ and $v_0\in
L^2_{k+1}(X,\Lambda^+\otimes\fg_{E_0})$ satisfy equations
\eqref{eq:SecondOrderQuasiLinearA} and
\eqref{eq:SecondOrderQuasiLinearA0}, respectively. If
\begin{equation}
\label{eq:PropSecondOrderASDDegEstCondition}
\max\left\{\|b\|_{L^{2\sharp,4}}, \|\cov\psi_0\|_{L^{2\sharp,4}},
\|\cov_{A_0}v_0\|_{L^{2\sharp,4}}, \|\cov_Av\|_{L^{2\sharp,4}}\right\}
< 
\delta,
\end{equation}
where $b = \chi_0(A - A_0) \in L^2_k(X,\Lambda^1\otimes\fg_{E_0})$, then 
\begin{equation}
\begin{aligned}
\label{eq:SecondOrderASDDegEst}
&\|\psi_0(v - v_0)\|_{L^\8\cap L^2_{2,A_0}(X)} 
\\
&\leq
C(\|\psi_0(w - w_0)\|_{L^{\sharp,2}} 
+ \|b\|_{L^{2\sharp,4}}\|v\|_{L^\8}\|\cov_Av\|_{L^{2\sharp,4}} 
+ \|b\|_{L^{2\sharp,4}}^2\|v\|_{L^\8}^2)
\\
&\quad + C\|\cov_{A_0}\psi_0(v-v_0)\|_{L^{2\sharp,4}}
(\|\cov_{A_0}v_0\|_{L^{2\sharp,4}}
+ \|\cov_{A_0}v\|_{L^{2\sharp,4}}
+ \|b\|_{L^{2\sharp,4}}\|v\|_{L^\8})
\\
&\quad + (\|\cov^2\psi_0\|_{L^{\sharp,2}} 
+ \|\cov\psi_0\|_{L^4(X)} + \|b\|_{L^4} 
+ \|\cov_{A_0}b\|_{L^{\sharp,2}})
\\
&\quad\times(\|v\|_{L^\8\cap L^2_{2,A}} + \|v_0\|_{L^\8\cap L^2_{2,A_0}}).
\end{aligned}
\end{equation}
\end{prop}

\begin{pf}
{}From equations \eqref{eq:SecondOrderQuasiLinearA} and 
\eqref{eq:SecondOrderQuasiLinearA0} we have
\begin{equation}
\label{eq:DiffSecondOrderASDEquations}
\begin{aligned}
d_A^+d_A^{+,*}v &= -(d_A^{+,*}v\wedge d_A^{+,*}v)^+ + w =: \tilde w, \\
d_{A_0}^+d_{A_0}^{+,*}v_0 &= -(d_{A_0}^{+,*}v_0\wedge d_{A_0}^{+,*}v_0)^+ 
+ w_0 =: \tilde w_0.
\end{aligned}
\end{equation}
Therefore, applying Lemma \ref{lem:LaplacianDegEst} and the embedding
$L^2_1\subset L^{2\sharp,4}$ \cite[Lemma 4.1]{FeehanSlice} to the pair of
equations \eqref{eq:DiffSecondOrderASDEquations} yields the estimate (when
$\|b\|_{L^2_{1,A_0}}\leq 1$)
\begin{equation}
\begin{aligned}
\label{eq:DiffSecondOrderSolnEstimate1}
\|\psi_0(v - v_0)\|_{L^\8\cap L^2_{2,A_0}} 
&\le 
C\|\psi_0(\tilde w - \tilde w_0)\|_{L^{\sharp,2}}
\\
&\quad + C\left(\|\cov^2\psi_0\|_{L^{\sharp,2}} + \|\cov\psi_0\|_{L^4}
+ \|b\|_{L^4} + \|\cov_{A_0}b\|_{L^{\sharp,2}}\right)
\\
&\quad\times(\|v\|_{L^\8\cap L^2_{2,A}} + \|v_0\|_{L^\8\cap L^2_{2,A_0}}),
\end{aligned}
\end{equation}
where $b$ and the constant $C$ are as in the statement of the proposition;
this gives the first and last terms on the right of inequality
\eqref{eq:SecondOrderASDDegEst}. 

It remains to bound the term $\psi_0(\tilde w - \tilde w_0)$ (giving the
remaining terms in the bound \eqref{eq:SecondOrderASDDegEst}): this
difference may be written schematically as
\begin{equation}
\label{eq:DiffSecondOrderSoln1}
\psi_0(\tilde w - \tilde w_0)
= 
\psi_0(w - w_0) + \psi_0\left(\cov_Av\otimes\cov_Av 
- \cov_{A_0}v_0\otimes\cov_{A_0}v_0\right).
\end{equation}
On $\supp\psi_0$ we have
\begin{align}
\label{eq:DiffSecondOrderSoln2}
\cov_Av\otimes\cov_Av
&= 
\cov_{A_0}v\otimes\cov_{A_0}v 
+ b\otimes v\otimes\cov_{A_0}v 
\\
\notag
&\quad + \cov_{A_0}v\otimes b\otimes v
+  b\otimes v\otimes  b\otimes v.
\end{align}
Now, using $v=(v-v_0)+v_0$ on $\supp\psi_0$ we observe that
\begin{align}
\label{eq:DiffSecondOrderSoln3}
\cov_{A_0}v\otimes\cov_{A_0}v
&=
\cov_{A_0}(v-v_0)\otimes\cov_{A_0}(v-v_0)
+ \cov_{A_0}v_0\otimes\cov_{A_0}(v-v_0)
\\
\notag
&\quad + \cov_{A_0}(v-v_0)\otimes\cov_{A_0}v_0
+ \cov_{A_0}v_0\otimes\cov_{A_0}v_0.
\end{align}
Lastly, note that 
\begin{align}
\notag
\psi_0\cov_{A_0}(v-v_0)\otimes\cov_{A_0}v_0
&=
\cov_{A_0}\psi_0(v-v_0)\otimes\cov_{A_0}v_0
\\
\notag
&\quad - \cov\psi_0\otimes(v-v_0)\otimes\cov_{A_0}v_0
\\
\label{eq:DiffSecondOrderSoln4}
\psi_0\cov_{A_0}v_0\otimes\cov_{A_0}(v-v_0)
&=
\cov_{A_0}v_0\otimes\cov_{A_0}\psi_0(v-v_0)
\\
\notag
&\quad - \cov_{A_0}v_0\otimes\cov\psi_0\otimes(v-v_0)
\\
\notag
\psi_0\cov_{A_0}(v-v_0)\otimes\cov_{A_0}(v-v_0)
&=
\cov_{A_0}\psi_0(v-v_0)\otimes\cov_{A_0}(v-v_0)
\\
\notag
&\quad - \cov\psi_0\otimes(v-v_0)\otimes\cov_{A_0}(v-v_0).
\end{align}
Combining equations \eqref{eq:DiffSecondOrderSoln1},
\eqref{eq:DiffSecondOrderSoln2}, \eqref{eq:DiffSecondOrderSoln3}, and
\eqref{eq:DiffSecondOrderSoln4},  
we finally arrive at a schematic equation for the difference:
\begin{align}
\notag
&\psi_0(\tilde w - \tilde w_0)
\\
\notag
&= \psi_0(w - w_0) 
\\
\label{eq:DiffSecondOrderSolnCombo}
&\quad + \psi_0(b\otimes v\otimes\cov_{A_0}v 
+ \cov_{A_0}v\otimes b\otimes v
+  b\otimes v\otimes  b\otimes v)
\\
\notag
&\quad + \cov_{A_0}\psi_0(v-v_0)\otimes\cov_{A_0}v_0
- \cov\psi_0\otimes(v-v_0)\otimes\cov_{A_0}v_0
\\
\notag
&\quad + \cov_{A_0}v_0\otimes\cov_{A_0}\psi_0(v-v_0)
- \cov_{A_0}v_0\otimes\cov\psi_0\otimes(v-v_0)
\\
\notag
&\quad + \cov_{A_0}\psi_0(v-v_0)\otimes\cov_{A_0}(v-v_0)
- \cov\psi_0\otimes(v-v_0)\otimes\cov_{A_0}(v-v_0).
\end{align}
Hence, from equation \eqref{eq:DiffSecondOrderSolnCombo}, H\"older's
inequality, and the Sobolev multiplication $L^{2\sharp,4}\otimes
L^{2\sharp,4} \to L^{\sharp,2}$ \cite[Lemma 4.3]{FeehanSlice},
we have a bound
\begin{align}
\notag
&\|\psi_0(\tilde w - \tilde w_0)\|_{L^{\sharp,2}}
\\
\label{eq:DiffSecondOrderSolnEstimate2}
&\leq C\left(\|\psi_0(w - w_0)\|_{L^{\sharp,2}} 
+ \|b\|_{L^{2\sharp,4}}\|v\|_{L^\8}\|\cov_Av\|_{L^{2\sharp,4}} 
+ \|b\|_{L^{2\sharp,4}}^2\|v\|_{L^\8}^2\right.
\\
\notag
&\quad + \|\cov_{A_0}\psi_0(v-v_0)\|_{L^{2\sharp,4}}
\left(\|\cov_{A_0}v_0\|_{L^{2\sharp,4}}
+ \|\cov_{A_0}v\|_{L^{2\sharp,4}}
+ \|b\|_{L^{2\sharp,4}}\|v\|_{L^\8}\right)
\\
\notag
&\quad\left. + \|\cov\psi_0\|_{L^{2\sharp,4}}\|v-v_0\|_{L^\8}
\left(\|\cov_{A_0}v_0\|_{L^{2\sharp,4}}
+ \|\cov_Av\|_{L^{2\sharp,4}}
+ \|b\|_{L^{2\sharp,4}}\|v\|_{L^\8}\right)\right).
\end{align}
For small enough $\delta$, our hypotheses 
\eqref{eq:PropSecondOrderASDDegEstCondition} implies that the terms
$\|b\|_{L^{2\sharp,4}}$, $\|\cov\psi_0\|_{L^{2\sharp,4}}$,
$\|\cov_{A_0}v_0\|_{L^{2\sharp,4}}$, and $\|\cov_Av\|_{L^{2\sharp,4}}$ are
small enough to allow us to rearrange the bound
\eqref{eq:DiffSecondOrderSolnEstimate1} and,
together with the estimate \eqref{eq:DiffSecondOrderSolnEstimate2}, yield
the desired bound \eqref{eq:SecondOrderASDDegEst}.
\end{pf}

\begin{rmk}
Just as in \S \ref{subsec:DegenerationLinearTheory}, Lemma
\ref{lem:FirstOrderASDDegEst} provides an $L^2_{1,A_0}$ bound for the
difference $\psi_0(a-a_0)$, while Lemma
\ref{lem:L21BoundForbInTermsOfA-A0} gives a bound for the Sobolev
norms of $b$ in terms of estimates on the Sobolev norms of $A-A_0$ and
$\chi_0$.
\end{rmk}

\begin{rmk}
To motivate the last paragraph of the proof of Proposition
\ref{prop:SecondOrderASDDegEst}, we note that in our application
$A$ will be close to $A_0$ in $L^2_{k,A_0,\loc}$ on $X\less\bx$, so the
term $\|b\|_{L^{2\sharp,4}}$ will be close to zero. We can also choose
$\psi_0$ such that $\|\cov\psi_0\|_{L^{2\sharp,4}}$ is small.  Since $v$
and $v_0$ solve the second-order anti-self-dual equations, we may assume
that the terms $\|\cov_{A_0}v_0\|_{L^{2\sharp,4}}$ and
$\|\cov_Av\|_{L^{2\sharp,4}}$ are small, though not necessarily converging
to zero if $A$ converges to $A_0$ on on $X\less\bx$.
\end{rmk}

\subsection{Degeneration of solutions to the second-order extended 
anti-self-dual equation: \boldmath{$L^2_2$} theory}
\label{subsec:ExtASDContinuityL22}
The new difficulty here is that we need to estimate
$\|\psi_0(w-w_0)\|_{L^{\sharp,2}(X)}$ where $w$, $w_0$ are now
explicitly defined by equation \eqref{eq:DefnOfwForExtASD}, so
\begin{equation}
\label{eq:DefnOfww_0ForExtASDCompared}
\begin{aligned}
w(A) &:= \Pi_{A,\mu}((d_A^{+,*}v\wedge d_A^{+,*}v)^+ + F_A^+) 
- F_A^+,
\\
w_0 &:= w(A_0) = \Pi_{A_0,\mu}((d_{A_0}^{+,*}v_0\wedge d_{A_0}^{+,*}v_0)^+ 
+ F_{A_0}^+) - F_{A_0}^+, 
\end{aligned}
\end{equation}
and where $v$ and $v_0$ must satisfy $\Pi_{A_0,\mu}v_0=0$ and
$\Pi_{A,\mu}v=0$.  Typically, we choose $\mu=\half\nu_2[A_0]$.

The obvious problem is that our resulting quasi-linear equation
\eqref{eq:SecondOrderQuasiLinearA} is no longer local, that is, given
purely by a partial differential equation, as we have included spectral
projections. However, by virtue of our comparison in \S
\ref{sec:Eigenvalue2} of the eigenvectors of $d_A^+d_A^{+,*}$ with the
approximate eigenvectors defined by cutting off those of
$d_{A_0}^+d_{A_0}^{+,*}$, we can replace our new non-local quasi-linear
equations with `essentially local' ones plus error terms which converge to
zero if $A$ converges to $(A_0,\bx)$ in the Uhlenbeck sense. Given this
broad outline, we proceed to the details.

Our task is to compare the solution $v$ to equation
\eqref{eq:SecondOrderQuasiLinearA} with the solution $v_0$ to the
corresponding equation \eqref{eq:SecondOrderQuasiLinearA0}, where $w$ and
$w_0$ are defined by \eqref{eq:DefnOfww_0ForExtASDCompared}.
Hence, the proof proceeds in much the same way as that of Proposition
\ref{prop:SecondOrderASDDegEst}, the main additional step being that in 
order to estimate the $L^{\sharp,2}$ norms of the cut-off differences,
$\psi_0(w-w_0)$, we need to also compare the $L^2$-orthogonal
projections $\Pi_{A_0,\mu}$ and $\Pi_{A,\mu}$.

\subsubsection{Preliminary estimates}
With the aid of two elementary lemmas and our eigenvector estimates
from \S \ref{sec:Eigenvalue2}, we can effectively reduce the problem
of estimating the difference between solutions to the extended
anti-self-dual equation to the simpler case of the ordinary anti-self-dual
equation.

We begin with the more specific problem of estimating $v-v_0$ on the subset
$\{x\in X:\psi_0(x)=1\}$.

\begin{lem}
\label{lem:2ndOrderExtASDBasicComp}
Let $X$ be a closed, oriented, smooth, Riemannian four-manifold.  Then
there is a positive constant $\eps$ with the following significance.  Let
$A$, $A_0$ be $L^2_4$ connections on $\SO(3)$ bundles $\fg_E$, $\fg_{E_0}$
over $X$, respectively, with $\|F_{A_0}^+\|_{L^{\sharp,2}(X)} < \eps$.
Suppose $U\subset X$ is an open subset such that $E|_U = E_0|_U$. Let
$\psi_0$, $\chi_0$ be cutoff functions on $X$ such that $\supp\chi_0\subset
U$, $\supp\psi_0\subset U$, and $\chi_0=1$ on $\supp\psi_0$.  Let $v\in
L^2_3(\Lambda^+\otimes\fg_E)$, $v_0\in L^2_3(\Lambda^+\otimes\fg_{E_0})$,
and suppose that $\Pi_{A,\mu}v=0$ and $\Pi_{A_0,\mu}v_0=0$. Then,
\begin{align*}
\|\psi_0(v-v_0)\|_{L^\8\cap L^2_{2,A_0}(X)}
&\leq
C\|d_{A_0}^+d_{A_0}^{+,*}\psi_0(v-v_0)\|_{L^{\sharp,2}(X)}
\\
&\quad + C(\|1-\psi_0\|_{L^2(X)}(\|v_0\|_{L^\8(X)} + \|v\|_{L^\8(X)})
\\
&\quad + \|(\Pi_{A,\mu}-\Pi_{A_0,\mu})\psi_0v\|_{L^2(X)}).
\end{align*}
where $C=C(\mu,\|F_{A_0}\|_{L^2})$.
\end{lem}

\begin{proof}
Our estimate for $\|\psi_0(v-v_0)\|_{L^\8\cap L^2_{2,A_0}}$ is
initiated by
\begin{equation}
\begin{aligned}
\label{eq:2ndOrderExtASDBasicComparison2}
\|\psi_0(v-v_0)\|_{L^\8\cap L^2_{2,A_0}}
&\leq
\|\Pi_{A_0,\mu}^\perp\psi_0(v-v_0)\|_{L^\8\cap L^2_{2,A_0}}
\\
&\quad +
\|\Pi_{A_0,\mu}\psi_0(v-v_0)\|_{L^\8\cap L^2_{2,A_0}}.
\end{aligned}
\end{equation}
The last term on the right 
in \eqref{eq:2ndOrderExtASDBasicComparison2} is `small' because
$\Pi_{A,\mu}v=0$ and $\Pi_{A_0,\mu}v_0=0$. Indeed,
\begin{align*}
\|\Pi_{A_0,\mu}\psi_0(v-v_0)\|_{L^\8\cap L^2_{2,A_0}}
&\leq
C\|\Pi_{A_0,\mu}\psi_0(v-v_0)\|_{L^2}
\quad\text{(by Lemma \ref{lem:LInftyEstDe2AEvec})}
\\
&\leq
C(\|\Pi_{A_0,\mu}v_0\|_{L^2}
+ \|\Pi_{A_0,\mu}(1-\psi_0)v_0\|_{L^2}
\\
&\quad + \|\Pi_{A,\mu}\psi_0v\|_{L^2}
+ \|(\Pi_{A,\mu}-\Pi_{A_0,\mu})\psi_0v\|_{L^2})
\\
&\leq
C(\|\Pi_{A_0,\mu}(1-\psi_0)v_0\|_{L^2}+ \|\Pi_{A,\mu}v\|_{L^2}
\\
&\quad + \|\Pi_{A,\mu}(1-\psi_0)v\|_{L^2}
+ \|(\Pi_{A,\mu}-\Pi_{A_0,\mu})\psi_0v\|_{L^2}),
\end{align*}
and thus
\begin{equation}
\begin{aligned}
\label{eq:ASDtoExtASDIntermedEst1}
\|\Pi_{A_0,\mu}\psi_0(v-v_0)\|_{L^\8\cap L^2_{2,A_0}}
&\leq
C(\|1-\psi_0\|_{L^2}(\|v_0\|_{L^\8} + \|v\|_{L^\8})
\\
&\quad + \|(\Pi_{A,\mu}-\Pi_{A_0,\mu})\psi_0v\|_{L^2}).
\end{aligned}
\end{equation}
We estimate the second-last term on the right in
\eqref{eq:2ndOrderExtASDBasicComparison2} by
\begin{equation}
\begin{aligned}
\label{eq:ASDtoExtASDIntermedEst2}
&\|\Pi_{A_0,\mu}^\perp\psi_0(v-v_0)\|_{L^\8\cap L^2_{2,A_0}}
\\
&=
\|G_{A_0,\mu}d_{A_0}^+d_{A_0}^{+,*}\psi_0(v-v_0)\|_{L^\8\cap L^2_{2,A_0}}
\\
&\leq
C\|\Pi_{A_0,\mu}^\perp d_{A_0}^+d_{A_0}^{+,*}\psi_0(v-v_0)\|_{L^{\sharp,2}}
\quad
\text{ by Lemma \ref{lem:LinftyL22CovLapEstv}}
\\
&\leq
C\|d_{A_0}^+d_{A_0}^{+,*}\psi_0(v-v_0)\|_{L^{\sharp,2}}
+
C\|\Pi_{A_0,\mu}d_{A_0}^+d_{A_0}^{+,*}\psi_0(v-v_0)\|_{L^{\sharp,2}}
\\
&\leq
C\|d_{A_0}^+d_{A_0}^{+,*}\psi_0(v-v_0)\|_{L^{\sharp,2}}
+
C\|\Pi_{A_0,\mu}d_{A_0}^+d_{A_0}^{+,*}\psi_0(v-v_0)\|_{L^2}
\\
&\leq
C\|d_{A_0}^+d_{A_0}^{+,*}\psi_0(v-v_0)\|_{L^{\sharp,2}}
+
C\mu\|\Pi_{A_0,\mu}\psi_0(v-v_0)\|_{L^2}.
\end{aligned}
\end{equation}
The last term on the right in \eqref{eq:ASDtoExtASDIntermedEst2} can
be estimated via inequality \eqref{eq:ASDtoExtASDIntermedEst1}.  The
conclusion follows by combining inequalities
\eqref{eq:2ndOrderExtASDBasicComparison2},
\eqref{eq:ASDtoExtASDIntermedEst1}, and
\eqref{eq:ASDtoExtASDIntermedEst2}.
\end{proof}

Suppose $u\in\Gamma(\Lambda^+\otimes\fg_E)$ and
$u_0\in\Gamma(\Lambda^+\otimes\fg_{E_0})$.  {}From equations
\eqref{eq:SecondOrderQuasiLinearA}, \eqref{eq:SecondOrderQuasiLinearA0}, 
\eqref{eq:DefnOfww_0ForExtASDCompared} and Proposition 
\ref{prop:SecondOrderASDDegEst}, we see that we shall need at 
least $L^{\sharp,2}(X)$ estimates for terms of the general form
$\psi_0(\Pi_{A,\mu}u-\Pi_{A_0,\mu}u_0)$.

\begin{lem}
\label{lem:2ndOrderExtASDQuadDataComp}
Continue the hypotheses and notation of Lemma
\ref{lem:2ndOrderExtASDBasicComp} and the preceding paragraph. Then, for
any $1\leq\ell\leq k$,
\begin{align*}
\|\psi_0(\Pi_{A,\mu}u-\Pi_{A_0,\mu}u_0)\|_{L^{\sharp,2}(X)}
&\leq C(\|\psi_0(u-u_0)\|_{L^2(X)}
+ \|\Pi_{A,\mu}(1-\psi_0)u\|_{L^2(X)}
\\
&\quad + \|\Pi_{A_0,\mu}(1-\psi_0)u_0\|_{L^2(X)}
\\
&\quad + \|(\Pi_{A,\mu}-\Pi_{A_0,\mu})\psi_0u\|_{L^2(X)}),
\end{align*}
where $C=C(\mu,\|F_{A_0}\|_{L^2})$.
\end{lem}

\begin{proof}
Observe that, via the embedding $L^4(X)\subset L^\sharp(X)$ of \cite[Lemma
4.1]{FeehanSlice} and $L^2_{1,A_0}(X)\subset L^4(X)$ of Lemma \ref{lem:Kato},
\begin{align*}
\|\psi_0(\Pi_{A,\mu}u-\Pi_{A_0,\mu}u_0)\|_{L^{\sharp,2}}
&\leq 
\|\Pi_{A,\mu}u-\Pi_{A_0,\mu}u_0\|_{L^{\sharp,2}}
\\
&\leq C\|\Pi_{A,\mu}u-\Pi_{A_0,\mu}u_0\|_{L^2}
\quad\text{(by Lemma \ref{lem:LInftyEstDe2AEvec}, above embeddings)}
\\
&\leq C(\|\Pi_{A,\mu}\psi_0u-\Pi_{A_0,\mu}\psi_0u_0\|_{L^2}
+
\|\Pi_{A,\mu}(1-\psi_0)u\|_{L^2}
\\
&\quad + \|\Pi_{A_0,\mu}(1-\psi_0)u_0\|_{L^2})
\\
&\leq
C(\|\Pi_{A_0,\mu}\psi_0(u-u_0)\|_{L^2}
+
\|(\Pi_{A,\mu}-\Pi_{A_0,\mu})\psi_0u\|_{L^2}
\\
&\quad + \|\Pi_{A,\mu}(1-\psi_0)u\|_{L^2}
+ \|\Pi_{A_0,\mu}(1-\psi_0)u_0\|_{L^2}),
\end{align*}
and the conclusion follows.
\end{proof}

\begin{rmk}
Recall from \S \ref{sec:Eigenvalue2} that we can approximate $\Pi_{A,\mu}$
by $L^2$-orthogonal projection onto a space of vectors in
$\Gamma(\Lambda^+\otimes\fg_E)$ with support contained in $\{x\in
X:\psi_0(x)=1\}$; similarly for $\Pi_{A_0,\mu}$.
\end{rmk}

With the aid of Proposition \ref{prop:SecondOrderASDDegEst} 
and its proof, we obtain the
following comparison result for solutions to the extended
anti-self-dual equations.

\begin{prop}
\label{prop:SecondOrderExtASDDegEst}
Continue the hypotheses and notation of Lemma
\ref{lem:2ndOrderExtASDBasicComp}. Suppose that $v\in
L^2_3(X,\Lambda^+\otimes\fg_E)$ and $v_0\in
L^2_3(X,\Lambda^+\otimes\fg_{E_0})$ satisfy equations
\eqref{eq:SecondOrderQuasiLinearA} and 
\eqref{eq:SecondOrderQuasiLinearA0}, respectively, with $w$, $w_0$ defined by 
equations \eqref{eq:DefnOfww_0ForExtASDCompared}. Then
\begin{equation}
\begin{aligned}
\label{eq:SecondOrderExtASDDegEst}
&\|\psi_0(v - v_0)\|_{L^\8(X)} + \|\psi_0(v - v_0)\|_{L^2_{2,A_0}(X)}
\\
&\leq C\|\psi_0(F_A^+-F_{A_0}^+)\|_{L^{\sharp,2}} 
\\
&\quad + C(\|1-\psi_0\|_{L^2(X)}(\|v_0\|_{L^\8(X)} + \|v\|_{L^\8(X)})
+ \|(\Pi_{A,\mu}-\Pi_{A_0,\mu})\psi_0v\|_{L^2(X)})
\\
&\quad + C(\|\Pi_{A,\mu}(1-\psi_0)u\|_{L^2(X)}
+ \|\Pi_{A_0,\mu}(1-\psi_0)u_0\|_{L^2(X)}
\\
&\quad + \|(\Pi_{A,\mu}-\Pi_{A_0,\mu})\psi_0u\|_{L^2(X)})
\\
&\quad + \text{\em All terms on right
of inequality \eqref{eq:SecondOrderASDDegEst}
in Proposition \ref{prop:SecondOrderASDDegEst}},
\end{aligned}
\end{equation}
where $u$, $u_0$ are defined by \eqref{eq:ChoiceOfuu0},
$b = \chi_0(A - A_0) \in L^2_4(X,\Lambda^1\otimes\fg_{E_0})$ and 
$C=C(\mu,\|F_{A_0}\|_{L^2})$ is a constant.
\end{prop}

\begin{proof}
According to Lemma \ref{lem:2ndOrderExtASDBasicComp}, we have
$$
\|\psi_0(v-v_0)\|_{L^\8\cap L^2_{2,A_0}(X)}
\leq
C\|d_{A_0}^+d_{A_0}^{+,*}\psi_0(v-v_0)\|_{L^{\sharp,2}(X)}
+
\text{Small terms},
$$ 
where the `small terms' in question (those which tend to zero if
$A$ converges to $(A_0,\bx)$ in the Uhlenbeck topology) have been
added to the right-hand side of \eqref{eq:SecondOrderExtASDDegEst}
above. To estimate the term
$\|d_{A_0}^+d_{A_0}^{+,*}\psi_0(v-v_0)\|_{L^{\sharp,2}(X)}$, we apply
Proposition \ref{prop:SecondOrderASDDegEst} to reduce the problem to
one of estimating the term $\|\psi_0(w-w_0)\|_{L^{\sharp,2}}$, with
$w$, $w_0$ defined by equations
\eqref{eq:DefnOfww_0ForExtASDCompared}, again modulo small terms which
have been added to the right-hand side of inequality
\eqref{eq:SecondOrderExtASDDegEst}.

To bound the term $\|\psi_0(w-w_0)\|_{L^{\sharp,2}}$, we apply Lemma 
\ref{lem:2ndOrderExtASDQuadDataComp}, where we choose
\begin{equation}
\begin{aligned}
\label{eq:ChoiceOfuu0}
u(A) &:= (d_A^{+,*}v\wedge d_A^{+,*}v)^+ + F_A^+,
\\
u_0 &:= u(A_0) = (d_{A_0}^{+,*}v_0\wedge d_{A_0}^{+,*}v_0)^+ + F_{A_0}^+.
\end{aligned}
\end{equation}
The other term we need to estimate is
$\|\psi_0(F_A^+-F_{A_0}^+)\|_{L^{\sharp,2}(X)}$, which requires no
further comment. {}From Lemma \ref{lem:2ndOrderExtASDQuadDataComp} we
see that it suffices to estimate the term $\|\psi_0(u-u_0)\|_{L^2}$,
as again the remaining terms from Lemma
\ref{lem:2ndOrderExtASDQuadDataComp} are small and have been added to
the right-hand side of inequality
\eqref{eq:SecondOrderExtASDDegEst}. Using the definition of $u$, $u_0$ in 
\eqref{eq:ChoiceOfuu0} we see from the proof of Proposition 
\ref{prop:SecondOrderASDDegEst} that, starting from the definition of 
$\tilde w$, $\tilde w_0$ in \eqref{eq:DiffSecondOrderASDEquations} and
using rearrangement, the term $\|\psi_0(u-u_0)\|_{L^2}$ is bounded by
terms already appearing on the right-hand side of inequality
\eqref{eq:SecondOrderExtASDDegEst}.
\end{proof}

Proposition \ref{prop:SecondOrderExtASDDegEst}, in conjunction with
the supplementary lemma below, shows that if $\{A_\alpha\}$ converges to
$(A_0,\bx)$ in the Uhlenbeck sense, then indeed $v_\alpha =
v(A_\alpha)$ converges to $v_0(A_0)$ in $L^2_{2,A_0,\loc}(X\less\bx)$ and
$w_\alpha = w(A_\alpha)$ converges to $w_0(A_0)$ in
$L^{\sharp,2}_{\loc}(X\less\bx)$.

\begin{lem}
\label{lem:SecondOrderExtASDDegConvergence}
Let the hypotheses be as in Theorem \ref{thm:UhlContExtSecOrderASD}. Then
the following hold:
\begin{itemize}
\item
The sequence of solutions $\{v_\alpha\}$ to equation
\eqref{eq:FirstOrderQuasiLinearA} converges to the solution $v_0$ to
equation \eqref{eq:FirstOrderQuasiLinearA0} in $L^\8\cap L^2_{2,A_0,\loc}$ on
$X\less\bx$,
\item
The sequence of terms $\{w_\alpha\}$ converges to the term $w_0$ in
$L^{\sharp,2}_{\loc}$ on $X\less\bx$,
\item
The sequence of solutions $a_\alpha = d_{A_\alpha}^{+,*}v_\alpha$ to
equation \eqref{eq:FirstOrderQuasiLinearA} converges to the solution
$a_0 = d_{A_0}^{+,*}v_0$ to \eqref{eq:FirstOrderQuasiLinearA0} in
$L^2_{1,A_0,\loc}$ on $X\less\bx$,
\item
The sequence of connections $\{A_\alpha+a_\alpha\}$ converges to $A_0+a_0$ in
$L^2_{1,A_0,\loc}$ on $X\less\bx$.
\end{itemize}
\end{lem}

\begin{proof}
We first consider the terms contributed to the right-hand side of
inequality \eqref{eq:SecondOrderExtASDDegEst}  in Proposition
\ref{prop:SecondOrderExtASDDegEst}
by Proposition \ref{prop:SecondOrderASDDegEst}.
Let $C$ be the constant of Proposition \ref{prop:SecondOrderASDDegEst}.
Suppose $U\Subset X\less\bx$ is a precompact open subset and that
$\eta$ is a positive constant. By Lemma \ref{lem:GoodCutoffSequence} we can
choose a sequence of cutoff functions $\{\psi_n\}$ on $X$ which are
$\bx$-good in the sense that
\begin{itemize}
\item
$\psi_n =1$ on $U$, 
\item
$\supp\psi_n \subset V_n \Subset X\less\bx$, where $\{V_n\}$ is
a collection of open subsets such that $\cup V_n = X\less\bx$,
\item
$\psi_n $ converges pointwise to the constant $1$ on $X\less\bx$, and
\item
$\lim_{n\to\8}
\|\cov\psi_n \|_{L^{2\sharp,4}(X)} 
+ \|\cov^2\psi_n \|_{L^{\sharp,2}(X)} = 0$.
\end{itemize}
We shall use the cutoff functions $\psi_n$ in place of the $\psi_0$
appearing in the statement of Proposition \ref{prop:SecondOrderASDDegEst}.
In particular, we may suppose that for large enough $n_0(\eta,C,M)$
(letting $M$ be the constant in the inequalities
\eqref{eq:UniformBoundvav0}),
\begin{equation}
\label{eq:LowerBoundNSqrtLambda}
\|\cov^2\psi_n \|_{L^{\sharp,2}(X)} 
+
\|\cov\psi_n \|_{L^{2\sharp,4}(X)} \leq \eta/(1000 CM^2),
\quad n\geq n_0.
\end{equation}
For any fixed $n$, the sequence $\{A_\alpha\}$ converges to $A_0$ in
$L^2_{k,A_0}(V_n)$ and so, for large enough
$\alpha_0=\alpha_0(C,\eta,n)$, we have
\begin{equation}
\label{eq:AalphaA0Close}
\|A_\alpha-A_0\|_{L^2_{k,A_0}(V_n)} \leq \eta/(1000 CM^2),
\quad \alpha\geq \alpha_0.
\end{equation}
For the remainder of the section we employ the abbreviation
\begin{equation}
\label{eq:WhatsbalphaAnyway}
b_\alpha := A_\alpha - A_0 \quad\text{on $X\less\bx$.} 
\end{equation}
Using $F(A_\alpha) = F(A_0 + b_\alpha) =
F_{A_0}+d_{A_0}b_\alpha + b_\alpha\wedge b_\alpha$ on $ X\less\bx$
and the Sobolev embeddings of \cite[Lemma 4.1]{FeehanSlice}, we see that
\begin{align*}
\|\psi_n (F^+_{A_\alpha}-F^+_{A_0})\|_{L^{\sharp,2}(X)}
&\leq
\|\cov_{A_0}b_\alpha\|_{L^{\sharp,2}(V_n)} 
+ \|b_\alpha\|_{L^{2\sharp,4}(V_n)}
\leq
c\|b_\alpha\|_{L^2_{2,A_0}(V_n)} 
\\
&\leq
c\|A_\alpha-A_0\|_{L^2_{2,A_0}(V_n)}\leq \eta/(1000 CM^2).
\end{align*}
We note that Lemma \ref{lem:GlobalASDEst} (with right-hand side data
$-\Pi_{A_\alpha,\mu}^\perp F_{A_\alpha}^+$) provides uniform bounds $M$ on
$v_\alpha$ which only depend on $\sup_\alpha \|F_{A_\alpha}^+\|_{L^\8}$,
$\mu$, $g$, and similarly for $v_0$, so
\begin{equation}
\label{eq:UniformBoundvav0}
\|v_\alpha\|_{L^\8\cap L^2_{2,A_\alpha}(X)} \leq M
\quad\text{and}\quad
\|v_0\|_{L^\8\cap L^2_{2,A_0}(X)} \leq M.
\end{equation}
Since $U\subset\{\psi_n =1\}$, by definition of $\psi_n $, we have
$$
\|v_\alpha-v_0\|_{L^\8\cap L^2_{2,A_0}(U)}
\leq
\|\psi_n (v_\alpha-v_0)\|_{L^\8\cap L^2_{2,A_0}(X)}.
$$
Hence, for large enough $n_0$,
the sum of the terms contributed to the right-hand side of inequality
\eqref{eq:SecondOrderExtASDDegEst} in Proposition
\ref{prop:SecondOrderExtASDDegEst} by \eqref{eq:SecondOrderASDDegEst}, is 
bounded by $\eta/3$, for all $\alpha\geq\alpha_0$. 

We now consider the terms on the right-hand side of inequality
\eqref{eq:SecondOrderExtASDDegEst} which are contributed by Lemma 
\ref{lem:2ndOrderExtASDBasicComp}.  {}From the definition of
the sequence of cutoff functions $\{\psi_n\}$, we have
$\lim_{n\to\8}\|1-\psi_n \|_{L^2(X)}=0$. By \S
\ref{sec:Eigenvalue2} we know that
$\|\Pi_{A_\alpha,\mu}\psi_n -\Pi_{A_0,\mu}\psi_n \|$ converges to
zero as $\alpha\to\8$ (for any fixed $n$), where $\|\cdot\|$ denotes
the operator norm in $\End(L^2(\Lambda^+\otimes\fg_E))$. Therefore, for
large enough $n_0$, the sum of the terms
contributed to the right-hand side of inequality
\eqref{eq:SecondOrderExtASDDegEst}  in Proposition
\ref{prop:SecondOrderExtASDDegEst}
by Lemma \ref{lem:2ndOrderExtASDBasicComp}, may be bounded by $\eta/3$ for
all $\alpha\geq\alpha_0$. (The
term $d_{A_0}^+d_{A_0}^{+,*}\psi_n(v_\alpha-v_0)$ is excluded, as it
has already been accounted for by Proposition
\ref{prop:SecondOrderASDDegEst}.)

Lastly, we consider the terms in inequality
\eqref{eq:SecondOrderExtASDDegEst}
contributed by Lemma \ref{lem:2ndOrderExtASDQuadDataComp}. {}From the
definition of $u_\alpha = u(A_\alpha)\in
L^{\sharp,2}(\Lambda^+\otimes\fg_E)$ in
\eqref{eq:ChoiceOfuu0} and our $L^\8\cap L^2_{2,A_\alpha}$ estimates for 
$v_\alpha$ in \eqref{eq:UniformBoundvav0}, we have a uniform $L^2$ bound on
$u_\alpha$.  Hence, we can see that the term
$\|\Pi_{A_\alpha,\mu}(1-\psi_n )u_\alpha\|_{L^2(X)}$ converges to zero
either by (i) approximating $\Pi_{A_\alpha,\mu}$ by $L^2$-orthogonal
projection onto a subspace of $L^2 (\Lambda^+\otimes\fg_E)$ spanned by
approximate eigenvectors of $d_{A_\alpha}^+d_{A_\alpha}^{+,*}$ which are
supported in $\{x\in X:1-\psi_n (x)=0\}$, or (ii) applying the Lebesgue
dominated convergence theorem to the sequence $(1-\psi_n )u_\alpha$,
noting that the support of $1-\psi_n $ shrinks. The same remarks apply
to the term $\|\Pi_{A_0,\mu}(1-\psi_n )u_0\|_{L^2(X)}$. As with the
corresponding term involving $v_\alpha$ in Lemma
\ref{lem:2ndOrderExtASDBasicComp}, the term
$\|(\Pi_{A_\alpha,\mu}-\Pi_{A_0,\mu})\psi_n u_\alpha\|_{L^2(X)}$ also
converges to zero as $\alpha\to\8$, for any fixed $n$.  Thus, the sum
of the terms contributed to the right-hand side of inequality
\eqref{eq:SecondOrderExtASDDegEst}  in Proposition
\ref{prop:SecondOrderASDDegEst}
by Lemma \ref{lem:2ndOrderExtASDQuadDataComp}, may be bounded by $\eta/3$, as
\begin{equation}
\label{eq:LimitCutoffDiffProjuu0}
\lim_{n\to\8}\lim_{\alpha\to\8}
\|\psi_n (\Pi_{A_\alpha,\mu}u_\alpha - \Pi_{A_0,\mu}u_0)\|_{L^2(X)}
=
0.
\end{equation}
Consequently, summing all three contributions to the right-hand side of
inequality \eqref{eq:SecondOrderExtASDDegEst}, and taking first $n$,
then $\alpha$ sufficiently large yields
$$
\|v_\alpha-v_0\|_{L^\8\cap L^2_{2,A_0}(U)} \leq \eta.
$$
Hence, the sequence $\{v_\alpha\}$ converges to $v_0$ in $L^\8\cap
L^2_{2,A_0}(U)$. Lemma \ref{lem:FirstOrderASDDegEst} then implies that
$\{a_\alpha\}$ converges to $a_0$ in $L^2_{1,A_0}(U)$ and thus the sequence of
connections $\{A_\alpha+a_\alpha\}$ converges to $A_0$ in $L^2_{1,A_0}(U)$.
\end{proof}

\subsection{Uhlenbeck convergence of solutions to the anti-self-dual
equation: \boldmath{$L^2_{k,\loc}$} convergence}
\label{subsec:FullUhlenbeckContinuityASD}
The estimates of the preceding subsection imply that
$\psi_n a_\alpha$ converges to $a_0$ in $L^2_{1,A_0}(X)$ if the
sequence of connections $A_\alpha$ converges to $(A_0,\bx)$ in
$L^2_{k,A_0,\loc}(X\less\bx)$. It remains to show, with the aid of our
higher-order interior elliptic estimates, that (i) the sequence of
one-forms $a_\alpha$ converges to $a_0$ in $L^2_{\ell,A_0,\loc}$ over
$X\less\bx$ for all $2\leq\ell\leq k$ and that (ii) if the sequence of
measures $|F_{A_\alpha}|^2$ converges to $|F_{A_0}|^2 +
8\pi^2\sum_{x\in\bx}\delta_x$ in the weak-* topology on measures, then
the sequence of measures $|F_{A_\alpha+a_\alpha}|^2$ also converges to
$|F_{A_0}|^2 + 8\pi^2\sum_{x\in\bx}\delta_x$. Note that given
$A_\alpha$ converges to $(A_0,\bx)$, we are not at liberty to pass to
a subsequence in order to show that $A_\alpha+a_\alpha$ converges to
$(A_0,\bx)$.

\begin{lem}
\label{lem:DiffASDEquations}
Suppose $v$, $v_0$ are solutions to equations
\eqref{eq:SecondOrderQuasiLinearA} and \eqref{eq:SecondOrderQuasiLinearA0},
respectively, defined by connections $A$ on $\fg_E$ and $A_0$ on
$\fg_{E_0}$, with $w$ and $w_0$ defined by
\eqref{eq:DefnOfww_0ForExtASDCompared}. Suppose $U\subset X$ is an
open subset such that $E|_U = E_0|_U$. Then $v':=v-v_0$ solves
\begin{align}
\label{eq:DiffASDEquations}
&d_{A_0}^+d_{A_0}^{+,*}v'
+ 
(d_{A_0}^{+,*}v'\wedge d_{A_0}^{+,*}v')^+ 
+
\{d_{A_0}^{+,*}v_0,d_{A_0}^{+,*}v'\}
\\
\notag
&= w' =: w-w_0 + b\otimes v\otimes d_{A_0}^{+,*}v 
+ b\otimes b\otimes v\otimes v,
\end{align}
where $b := (A-A_0)|_U$ and $\{\alpha,\beta\} := (\alpha\wedge\beta +
\beta\wedge\alpha)^+$. 
\end{lem}

\begin{proof}
We work over the open set $U$ where $E=E_0$.
Using $d_A^{+,*}v = -*d_Av$ and $d_A^{+,*}v = -*d_{A_0}v - *(b\wedge v)$,
we observe that
\begin{equation}
\begin{aligned}
\label{eq:DiffConnQuad}
d_A^{+,*}v\wedge d_A^{+,*}v
&=
d_{A_0}^{+,*}v\wedge d_{A_0}^{+,*}v
- *(b\wedge v)\wedge d_{A_0}^{+,*}v
- d_{A_0}^{+,*}v\wedge *(b\wedge v)
\\
&\quad + *(b\wedge v)\wedge *(b\wedge v).
\end{aligned}
\end{equation}
By equation \eqref{eq:ExpandLaplacianConnection} we see that
\begin{equation}
\label{eq:DiffConnLaplacian}
d_A^+d_A^{+,*}v 
= 
d_{A_0}^+d_{A_0}^{+,*}v - (d_{A_0}^+*(b\wedge v)^+ -
(b\wedge *d_{A_0} v)^+ + (b\wedge *(b\wedge v))^+.
\end{equation}
Noting that $v'=v-v_0$, we obtain the following difference between
quadratic terms,
\begin{equation}
\begin{aligned}
\label{eq:Diff2FormQuad}
&d_{A_0}^{+,*}v\wedge d_{A_0}^{+,*}v
- d_{A_0}^{+,*}v_0\wedge d_{A_0}^{+,*}v_0
\\
&=d_{A_0}^{+,*}v'\wedge d_{A_0}^{+,*}v'
+ d_{A_0}^{+,*}v_0\wedge d_{A_0}^{+,*}v'
+ d_{A_0}^{+,*}v'\wedge d_{A_0}^{+,*}v_0.
\end{aligned}
\end{equation}
Subtracting equation \eqref{eq:SecondOrderQuasiLinearA0} from equation
\eqref{eq:SecondOrderQuasiLinearA} and using the identities
\eqref{eq:DiffConnQuad}, \eqref{eq:DiffConnLaplacian}, and
\eqref{eq:Diff2FormQuad}, we obtain
\begin{align*}
&d_{A_0}^+d_{A_0}^{+,*}v' + (d_{A_0}^{+,*}v'\wedge d_{A_0}^{+,*}v')^+ 
+ (d_{A_0}^{+,*}v_0\wedge d_{A_0}^{+,*}v')^+
+ (d_{A_0}^{+,*}v'\wedge d_{A_0}^{+,*}v_0)^+
\\
&= w-w_0 + d_{A_0}^+*(b\wedge v) + (b\wedge *d_{A_0} v)^+ 
- (b\wedge *(b\wedge v))^+. 
\end{align*}
This gives the desired result.
\end{proof}

We can now complete the proofs of our main gluing-map continuity results.

\begin{proof}[Proof of Theorem \ref{thm:UhlContSecOrderQuasiLin}]
The first step is to show that the sequence of connections
$A_\alpha+a_\alpha$ converges to $A_0+a_0$ in
$L^2_{\ell,A_0,\loc}(X\less\bx)$, when $2\leq \ell\leq k$; Lemma
\ref{lem:SecondOrderExtASDDegConvergence} addresses the case $\ell=1$.

Let $U''\Subset U'\Subset U\Subset X\less\bx$ be precompact open subsets:
we will show that $A_\alpha+a_\alpha$ converges to $A_0+a_0$ in
$L^2_{\ell,A_0}(U'')$.  Observe that $v_\alpha' := v_\alpha-v_0$ is a
solution to equation
\eqref{eq:GeneralSecondOrderASD} with $A$ replaced by $A_0$, $\alpha$
replaced by $d_{A_0}^{+,*}v_0$, and $w$ replaced by the right-hand
side, say $w_\alpha'$, of identity
\eqref{eq:DiffASDEquations}. Then applying Corollary
\ref{cor:L2_2SecondOrderInhomoRegLocal} to equation
\eqref{eq:GeneralSecondOrderASD}, with the preceding replacements,
provides the estimate
\begin{equation}
\label{eq:valphav0diffEst}
\|v_\alpha-v_0\|_{L^2_{\ell+1,A_0}(U'')}
\leq
Q_\ell(\|w_\alpha'\|_{L^2_{\ell-1,A_0}(U')}, \|v_\alpha-v_0\|_{L^2(U')}),
\quad \ell \geq 2.
\end{equation}
The positive, real coefficients of the universal polynomial $Q_\ell$ depend
at most on $A_0$, $v_0$, $\ell$, $U'$, $U''$, and the metric $g$ on $X$;
also, $Q(0,0)=0$.  We have already seen, 
from Lemma \ref{lem:SecondOrderExtASDDegConvergence},
that $\|v_\alpha-v_0\|_{L^2(U)}$ converges to zero. {}From equations
\eqref{eq:DefnOfww_0ForExtASDCompared} and \eqref{eq:DiffASDEquations},
\begin{equation}
\begin{aligned}
\label{eq:walphaEst}
&\|w_\alpha'\|_{L^2_{\ell-1,A_0}(U')}
\\
&\quad\leq
\|b_\alpha\otimes v_\alpha\otimes 
d_{A_0}^{+,*}v_\alpha\|_{L^2_{\ell-1,A_0}(U')} 
+ \|b_\alpha\otimes b_\alpha\otimes v_\alpha\otimes 
v_\alpha\|_{L^2_{\ell-1,A_0}(U')}
\\
&\qquad + \|\Pi_{A_\alpha,\mu}u_\alpha 
- \Pi_{A_0,\mu}u_0\|_{L^2_{\ell-1,A_0}(U')}
+ \|F_{A_\alpha}^+-F_{A_0}^+\|_{L^2_{\ell-1,A_0}(U')},
\end{aligned}
\end{equation}
where $u_\alpha$, $u_0$ are defined by equation
\eqref{eq:ChoiceOfuu0} and $b_\alpha$ by \eqref{eq:WhatsbalphaAnyway}.
The Sobolev multiplication theorems and definition
\eqref{eq:WhatsbalphaAnyway} imply that
\begin{equation}
\begin{aligned}
\label{eq:TwoTermsbavava}
&\|b_\alpha\otimes v_\alpha\otimes 
d_{A_0}^{+,*}v_\alpha\|_{L^2_{\ell-1,A_0}(U')} 
+ \|b_\alpha\otimes b_\alpha\otimes v_\alpha\otimes 
v_\alpha\|_{L^2_{\ell-1,A_0}(U')}
\\
&\quad\leq
\begin{cases}
C\|A_\alpha-A_0\|_{L^2_3(U')}
\|v_\alpha\|_{L^\8\cap L^2_{2,A_0}(U')}^2, &\ell=2,
\\
C\|A_\alpha-A_0\|_{L^2_\ell(U')}
\|v_\alpha\|_{L^2_{\ell,A_0}(U')}^2, &\ell\geq 3,
\end{cases}
\end{aligned}
\end{equation}
for large enough $\alpha$ (to bound quadratic terms in
$\|A_\alpha-A_0\|_{L^2_\ell(U')}$ by $\|A_\alpha-A_0\|_{L^2_\ell(U')}$)
while
\begin{equation}
\label{eq:SequenceCurvDiff}
\|F_{A_\alpha}^+-F_{A_0}^+\|_{L^2_{\ell-1,A_0}(U')}
\leq 
C\|A_\alpha-A_0\|_{L^2_\ell(U')}.
\end{equation}
Using Lemma \ref{lem:GoodCutoffSequence}, choose a sequence of $\bx$-good
cutoff functions $\psi_n $ on $X$, so
\begin{itemize}
\item
$\supp\psi_n \Subset X\less\bx$,
\item
$\psi_n =1$ on $U$, and 
\item
$\psi_n $ converges pointwise to the constant $1$ on $X\less\bx$.  
\end{itemize}
We now exploit the fact that, for all $\alpha$ sufficiently large, the
sequence $A_\alpha$ is uniformly $L^2_{k,A_0}(U)$ close to $A_0$ over the
precompact open subset $U$: with this observation, our standard estimates
in Lemma \ref{lem:LInftyEstDe2AEvec} for eigenvectors of
$d_{A_\alpha}^+d_{A_\alpha}^{+,*}$ and $d_{A_0}^+d_{A_0}^{+,*}$ yield the
first inequality below (the second following by definition of the
$\psi_n $),
\begin{equation}
\begin{aligned}
\label{eq:L2ellDiffOrthogProjuau0}
\|\Pi_{A_\alpha,\mu}u_\alpha - \Pi_{A_0,\mu}u_0\|_{L^2_{\ell-1,A_0}(U')}
&\leq
C\|\Pi_{A_\alpha,\mu}u_\alpha - \Pi_{A_0,\mu}u_0\|_{L^2(U)}
\\
&\leq
C\|\psi_n (\Pi_{A_\alpha,\mu}u_\alpha - \Pi_{A_0,\mu}u_0)\|_{L^2(X)},
\end{aligned}
\end{equation}
where $u_\al=u(A_\alpha)$, $u_0=u(A_0)$ are defined by \eqref{eq:ChoiceOfuu0},
and $C=C([A_0],\mu,U,U',\ell,\alpha_0,g)$.

Combining estimates \eqref{eq:valphav0diffEst}, \eqref{eq:walphaEst},
\eqref{eq:TwoTermsbavava}, \eqref{eq:SequenceCurvDiff},
and \eqref{eq:L2ellDiffOrthogProjuau0} yields the 
bound, for $\ell \geq 2$,
\begin{equation}
\begin{aligned}
\label{eq:Diffvalphav0BigSobolev}
&\|v_\alpha-v_0\|_{L^2_{\ell+1,A_0}(U'')}
\\
&\leq
Q_\ell\left(\|v_\alpha-v_0\|_{L^2_{2,A_0}(U')} 
+ \|A_\alpha-A_0\|_{L^2_{\ell,A_0}(U')}\right.
\\
&\qquad\quad\left. + \|\psi_n (\Pi_{A_\alpha,\mu}u_\alpha 
- \Pi_{A_0,\mu}u_0)\|_{L^2(X)},\|v_\alpha-v_0\|_{L^2(U')}\right).
\end{aligned}
\end{equation}
Lemma \ref{lem:SecondOrderExtASDDegConvergence} --- which we use to get
convergence-to-zero for the terms $\|v_\alpha-v_0\|_{L^2_{2,A_0}(U')}$ ---
and the limit result
\eqref{eq:LimitCutoffDiffProjuu0} (given our choice of cutoff functions 
$\psi_n $) imply that the right-hand side of
\eqref{eq:Diffvalphav0BigSobolev} converges to zero as $\alpha,n\to\8$.
Hence, the sequence $v_\alpha$ converges to $v_0$ in
$L^2_{\ell+1,A_0}(U'')$.  The definition of the $a_\alpha$ and $a_0$
implies that the sequence $a_\alpha$ converges to $a_0$ in
$L^2_{\ell,A_0}(U'')$.  Consequently, the sequence of connections
$A_\alpha+a_\alpha$ converges to $A_0+a_0$ in $L^2_{\ell,A_0}(U'')$. We
conclude that $A_\alpha+a_\alpha$ converges to $A_0+a_0$ in
$L^2_{k,A_0,\loc}(X\less\bx)$, as desired.

The second (and final) step is to show that the sequence of measures
$|F(A_\alpha+a_\alpha)|^2$ converges to
$|F(A_0+a_0)|^2+8\pi^2\sum_{x\in\bx}\delta_x$ on $X$.  Remark that, by
hypothesis, the sequence of measures $|F(A_\alpha)|^2$ converges to
$|F(A_0)|^2+8\pi^2\sum_{x\in\bx}\delta_x$ on $X$. Over $X\less\bx$,
the sequence of functions $|F(A_\alpha+a_\alpha)|^2$ converges in
$C^\8$ to $|F(A_0+a_0)|^2$ so, for any $f\in C^0(X)$ and $r>0$,
\begin{equation}
\label{eq:ConvMeasXlessBalls}
\lim_{\alpha\to\8}\frac{1}{8\pi^2}\int_{X\less\cup_i B(x_i,r)} 
f|F(A_\alpha+a_\alpha)|^2\,dV
=
\frac{1}{8\pi^2}\int_{X\less\cup_i B(x_i,r)} f|F(A_0+a_0)|^2\,dV.
\end{equation}
It remains to consider the corresponding limits over the balls $B(x_i,r)$.

We shall adapt the argument of
\cite[p. 165]{DK} and consider the Chern-Simons invariants of the
restrictions of the extended anti-self-dual connections $A_0+a_0$, and
$A_\alpha+a_\alpha$ to three-spheres $S(x_i,r) = \rd B(x_i,r)$, where $r$
is a small positive constant. Recall that if $B$ is a connection on the
product $\U(2)$ bundle over a three-manifold $Y$, then its Chern-Simons
invariant is given by \cite[p. 164]{DK}
$$
\CS_Y(B)
:=
\frac{1}{8\pi^2}\int_Y\tr\left(dB + \frac{2}{3}B\wedge B\wedge B\right).
$$
If $A$ is any extension of $B$ to a connection over a four-manifold $Z$
with boundary $Y$, then
\begin{equation}
\label{eq:ChernSimonsModZ}
\CS_Y(B) = -\frac{1}{8\pi^2}\int_Z\tr(F_A\wedge F_A)\pmod{\ZZ}.
\end{equation}
Because $(|F_A^+|^2 -  |F_A^+|^2)\,dV = -\tr(F_A\wedge F_A)$ for any
connection $A$, we see that
\begin{equation}
\begin{aligned}
\label{eq:CurvatureMeasIntegralToCS}
&\lim_{\alpha\to\8}
\frac{1}{8\pi^2}\int_{B(x_i,r)}
\left(|F^-(A_\alpha+a_\alpha)|^2 - |F^+(A_\alpha+a_\alpha)|^2\right)\,dV
\\
&\quad =
\lim_{\alpha\to\8}
\CS_{S(x_i,r)}(A_\alpha+a_\alpha)
\pmod{\ZZ}
\\
&\quad =
\CS_{S(x_i,r)}(A_0 + a_0)
\\
&\quad =
\frac{1}{8\pi^2}\int_{B(x_i,r)}
\left(|F^-(A_0+a_0)|^2-|F^+(A_0+a_0)|^2\right)\,dV
\pmod{\ZZ}.
\end{aligned}
\end{equation}
By construction \eqref{eq:QuickExtASDEqnForv} of the solutions
$a_\alpha=d_{A_\alpha}^{+,*}v_\alpha$, we have $F^+(A_\alpha+a_\alpha) =
\Pi_{A_\alpha,\mu}F^+(A_\alpha+a_\alpha)$ and so 
Lemma \ref{lem:LInftyEstDe2AEvec} implies that
$$
\|F^+(A_\alpha+a_\alpha)\|_{C^0(X)}
\leq
C\|F^+(A_\alpha+a_\alpha)\|_{L^2(X)},
$$
where $C$ is a constant independent of $\alpha$ by Lemma
\ref{lem:LInftyEstDe2AEvec} and hypothesis
\eqref{eq:UniformF+BoundOnConnSequence}. Therefore, the Lebesgue bounded
convergence theorem yields
$$
\lim_{\alpha\to\8}
\frac{1}{8\pi^2}\int_{B(x_i,r)}|F^+(A_\alpha+a_\alpha)|^2\,dV
=
\frac{1}{8\pi^2}\int_{B(x_i,r)}|F^+(A_0+a_0)|^2\,dV,
$$
and so the preceding limit identity, the fact that $|F_A|^2 =
|F_A^-|^2+|F_A^+|^2$ for any connection $A$, and
\eqref{eq:CurvatureMeasIntegralToCS} imply that
\begin{equation}
\begin{aligned}
\label{eq:LimitCurvatureMeasureModZ}
&\lim_{\alpha\to\8}
\frac{1}{8\pi^2}\int_{B(x_i,r)}|F(A_\alpha+a_\alpha)|^2\,dV
\\
&\lim_{\alpha\to\8}
\frac{1}{8\pi^2}\int_{B(x_i,r)}
\left((|F^-(A_\alpha+a_\alpha)|^2-|F^+(A_\alpha+a_\alpha)|^2)
+ 2|F^+(A_\alpha+a_\alpha)|^2\right)\,dV
\\
&\quad =
\frac{1}{8\pi^2}\int_{B(x_i,r)}
\left((|F^-(A_0+a_0)|^2-|F^+(A_0+a_0)|^2) + 2|F^+(A_0+a_0)|^2\right)\,dV
\\
&\quad =
\frac{1}{8\pi^2}\int_{B(x_i,r)}|F(A_0+a_0)|^2\,dV
\pmod{\ZZ}.
\end{aligned}
\end{equation}
Note that $F(A_\alpha+a_\alpha) = F(A_\alpha) + d_{A_\alpha}a_\alpha +
a_\alpha\wedge a_\alpha$, so we have the estimate
\begin{align*}
&\int_X\left||F(A_\alpha+a_\alpha)|^2-|F(A_\alpha)|^2\right|\,dV
\\
&\leq 
2\int_X(F(A_\alpha),d_{A_\alpha}a_\alpha + a_\alpha\wedge
a_\alpha)\,dV
+
2\int_X|d_{A_\alpha}a_\alpha + a_\alpha\wedge a_\alpha|^2\,dV.
\end{align*}
Therefore,
\begin{align*}
&\int_X\left||F(A_\alpha+a_\alpha)|^2-|F(A_\alpha)|^2\right|\,dV
\\
&\quad\leq  2(\|F(A_\alpha)\|_{L^2(X)}+\|a_\alpha\|_{L^2_{1,A_\alpha}(X)})
\|a_\alpha\|_{L^2_{1,A_\alpha}(X)}.
\end{align*}
The terms $\|F(A_\alpha)\|_{L^2(X)}^2$ are bounded by a
universal constant times $\|F(A_0)\|_{L^2(X)}^2 + 8\pi^2\ell$.
The terms $\|a_\alpha\|_{L^2_{1,A_\alpha}(X)}$ are bounded by a universal
constant times $\|F^+(A_\alpha)\|_{L^{\sharp,2}(X)}$, by our {\em a
priori\/} estimates of Lemma \ref{lem:GlobalASDEst}
for our solutions to the anti-self-dual equation
\eqref{eq:FirstOrderQuasiLinearA}. Similar remarks apply
to the connections $A_0$ and $A_0+a_0$, so we obtain
\begin{equation}
\begin{aligned}
\label{eq:CurvatureMeasureBounds}
\frac{1}{8\pi^2}
\int_X\left||F(A_\alpha+a_\alpha)|^2-|F(A_\alpha)|^2\right|\,dV
&\leq \quarter,
\\
\frac{1}{8\pi^2}
\int_X\left||F(A_0+a_0)|^2-|F(A_0)|^2\right|\,dV
&\leq \quarter,
\end{aligned}
\end{equation}
as our hypotheses \eqref{eq:UniformF+BoundOnConnSequence} imply that the
factors $\sup_\alpha
\|F^+(A_\alpha)\|_{L^2(X)}$ and $\|F^+(A_0)\|_{L^2(X)}$
are small. Since
\begin{equation}
\label{eq:ApproxASDCurvatureMeasureLimit}
\lim_{\alpha\to\8}
\frac{1}{8\pi^2}\int_{B(x_i,r)}|F(A_\alpha)|^2\,dV
= \frac{1}{8\pi^2}\int_{B(x_i,r)}|F(A_0)|^2\,dV + \kappa_i,
\end{equation}
it follows from the limit \eqref{eq:LimitCurvatureMeasureModZ}, the
estimates \eqref{eq:CurvatureMeasureBounds}, and the limit 
\eqref{eq:ApproxASDCurvatureMeasureLimit} that
\begin{equation}
\begin{aligned}
\label{eq:CurvatureLimitIdentity}
&\lim_{\alpha\to\8}
\frac{1}{8\pi^2}\int_{B(x_i,r)}|F(A_\alpha+a_\alpha)|^2\,dV
\\
&\quad =
\frac{1}{8\pi^2}\int_{B(x_i,r)}
\left(|F(A_0+a_0)|^2\,dV + 8\pi^2\kappa_i\delta_{x_i}\right).
\end{aligned}
\end{equation}
Finally, observe that for all $\alpha$ and $f\in C^0(X)$ we have
\begin{equation}
\label{eq:MeasureConvergenceEstimate}
\frac{1}{8\pi^2}\int_{B(x_i,r)}|f(x_i) - f(x)||F(A_\alpha+a_\alpha)|^2\,dV(x)
\leq
\kappa\sup_{x\in B(x_i,r)}|f(x_i)-f(x)|,
\end{equation}
where $\kappa = -\quarter p_1(\fg_E)$.
By choosing a small enough constant $r>0$, we can make the right-hand side of 
\eqref{eq:MeasureConvergenceEstimate} arbitrarily small, uniformly with respect
to $\alpha$. Hence, combining the limit identities
\eqref{eq:ConvMeasXlessBalls} and
\eqref{eq:CurvatureLimitIdentity} with the estimate
\eqref{eq:MeasureConvergenceEstimate} yields
$$
\lim_{\alpha\to\8}\frac{1}{8\pi^2}\int_X f|F(A_\alpha+a_\alpha)|^2\,dV
=
\frac{1}{8\pi^2}\int_X f|F(A_0+a_0)|^2\,dV
+
\sum_{i=1}^m\kappa_i f(x_i),
$$
as required, where $\bx\in\Sym^\ell(X)$ is represented by the $m$-tuple
$(x_1,\dots,x_m)$. The existence of the limit asserted above follows from
its existence mod $\ZZ$ via the application of the Chern-Simons identity
\eqref{eq:ChernSimonsModZ} in \eqref{eq:CurvatureMeasIntegralToCS} and the
curvature estimates \eqref{eq:CurvatureMeasureBounds}. 
\end{proof}

\ifx\undefined\bysame
\newcommand{\bysame}{\leavevmode\hbox to3em{\hrulefill}\,}
\fi

\end{document}